%% file: anti-Zener-and-Zener.tex
\documentclass[a4paper]{article}
\usepackage{amsmath, amsfonts, amsthm,amssymb}
\usepackage{mathrsfs}
\usepackage[margin=2cm]{geometry}
\usepackage{xcolor}
\usepackage{colortbl}
\usepackage{graphicx}
\usepackage{caption}
\usepackage{subfig}
\usepackage{multirow,makecell,multicol}
\usepackage{hhline}
\usepackage{stackengine}
\usepackage{lscape}
\usepackage[T1]{fontenc}

\newcommand\xrowht[2][0]{\addstackgap[.5\dimexpr#2\relax]{\vphantom{#1}}}
\allowdisplaybreaks
\makeatletter
\renewcommand*{\@fnsymbol}[1]{\ifcase#1 \or \@arabic{\numexpr#1 \relax}  \or \@arabic{\numexpr#1 \relax} \or \@arabic{\numexpr#1 \relax} \else \ * \fi}
\makeatother
\input{tcilatex}
\begin{document}

\title{Energy balance for fractional anti-Zener and Zener models\\
in terms of relaxation modulus and creep compliance}
\author{Sla\dj an Jeli\'{c}\thanks{
Department of Physics, Faculty of Sciences, University of Novi Sad, Trg D.
Obradovi\'{c}a 4, 21000 Novi Sad, Serbia, sladjan.jelic@df.uns.ac.rs} \quad
Du\v{s}an Zorica\thanks{
Department of Physics, Faculty of Sciences, University of Novi Sad, Trg D.
Obradovi\'{c}a 4, 21000 Novi Sad, Serbia, dusan.zorica@df.uns.ac.rs} \ 
\thanks{%
Mathematical Institute, Serbian Academy of Arts and Sciences, Kneza Mihaila
36, 11000 Belgrade, Serbia, dusan\textunderscore zorica@mi.sanu.ac.rs} \ 
\thanks{Corresponding author}}

\maketitle

\begin{abstract}
\noindent Relaxation modulus and creep compliance corresponding to
fractional anti-Zener and Zener models are calculated and restrictions on
model parameters narrowing thermodynamical constraints are posed in order to
ensure relaxation modulus and creep compliance to be completely monotone and Bernstein function respectively,
that a priori guarantee the positivity of stored energy and dissipated
power per unit volume, derived in time domain by considering the power per
unit volume. Both relaxation modulus and creep compliance for model
parameters obeying thermodynamical constraints, proved that can also be
oscillatory functions with decreasing amplitude. Model used in numerical
examples of relaxation modulus and creep compliance is also analyzed for the
asymptotic behavior near the initial time instant and for large time.

\end{abstract}

\section{Introduction}

Fractional anti-Zener and Zener models, that are the subject of the energy
balance analysis in time domain, are formulated in \cite{SD-1} by
considering the rheological schemes corresponding to the classical
anti-Zener and Zener models having the classical spring replaced by the
fractional element whose stress-strain relation is assumed in terms of
fractional integral, as well as having the classical dash-pot replaced by
the Scott-Blair element, i.e., element whose stress-strain relation is
assumed in terms of fractional derivative. Note, the classical anti-Zener
model is also referred to as the Jeffrey model, as discussed in \cite%
{Shitikova} reviewing classical and fractional models of viscoelasticity.
Some of the first to use the Scott-Blair element in posing the constitutive
equations using classical rheological schemes are \cite%
{SchiesselMetzlerBlumenNonnemacher}, while some of the most widely used
fractional order models of viscoelastic body derived using the rheological
schemes containing Scott-Blair element are reviewed in \cite{Mai-10,hilf}.
Rheological scheme of the Burgers model with the classical dash-pot replaced
by the Scott-Blair element is adopted in \cite{OZ-1}, in order to
fractionalize the classical Burgers model. Infinite number of springs and
dash-pots in the rheological scheme is used in \cite%
{ColombaroGiustiMainardi,Giusti,GiustiMainardi} in order to formulate the
Bessel model, while the setting of fractal rheological models is used in 
\cite{HeymansBauwens}.

The analysis of thermodynamical consistency of constitutive models in the
steady state regime, i.e., by requesting non-negativity of the storage and
loss modulus for any frequency, is given in \cite{b-t}. This method is used
in \cite{SD-1} for formulation of thermodynamical restrictions on model
parameters appearing in the fractional anti-Zener and Zener models, as well
as in \cite{OZ-1} and \cite{BazhlekovaTsocheva} for the fractional Burgers
model, where the fractionalization is performed by replacing the classical
dash-pot element with the fractional one in the former, as well as by
replacing the ordinary derivatives with the fractional ones in the classical
Burgers model in the latter one. Constitutive models of viscoelastic body
having the orders of fractional derivatives not exceeding the first order
are examined for thermodynamical consistency in \cite{AKOZ}. Considering
time domain, the dissipation inequality is also used in order to pose
thermodynamical restrictions on model parameters, see for example \cite%
{a-2002,AJP}, while the energy balance properties of the fractional wave
equations are examined in \cite{ZO}. Dissipation of energy in viscoelastic
media of Bessel type, the wave equation thermodynamical consistency and the
requirements guaranteeing the complete monotonicity, as well as the analysis
of the Zener wave equation are presented in \cite%
{ColombaroGiustiMentrelli,EndeLionLammering,HolmHolm,NH}. The extensive
account on thermodynamics of hereditary materials can be found in \cite%
{AmendolaFabrizioGolden}.

A priori energy balance properties of fractional anti-Zener and Zener
models, that are formulated in \cite{SD-1} and analyzed for thermodynamical
consistency in the steady-state regime, are discussed in time domain and it
is proved that the power per unit volume consists of two terms: the time
derivative of the energy per unit volume stored in the viscoelastic body and
dissipated power per unit volume. The stored energy per unit volume,
expressed through strain, consists of the instantaneous term, that is of the
same form as the potential energy of elastic body in which the relaxation
modulus plays the role of the Young modulus, as well as of the hereditary
type term, where the time derivative of the relaxation modulus represents
the memory kernel of the difference of current strain and strain in the
previous time instant. Dissipated power per unit volume, expressed in terms
of strain, also consists of the instantaneous and hereditary type terms,
respectively containing the first and second time derivative of the
relaxation modulus. The positivity of stored energy and dissipated power per
unit volume is guaranteed if the relaxation modulus is a positive
monotonically decreasing and convex function, or moreover a completely
monotone function. Note, the relaxation modulus is a function representing
the time evolution of stress in a stress relaxation experiment, occurring as
a consequence of the sudden and later constant strain, i.e., strain
prescribed as the Heaviside step function.

The stored energy per unit volume can also be expressed through stress, when
there is only the hereditary type term, having the time derivative of the
creep compliance as the memory kernel multiplying the time evolution of
stress, while the dissipated power per unit volume consists of the
instantaneous and hereditary type terms, respectively containing the first
and second time derivative of the creep compliance. Again, the positivity of
stored energy and dissipated power per unit volume is guaranteed if the
creep compliance is a positive monotonically increasing and concave
function, or moreover a Bernstein function. Note, the creep compliance is a
function representing the time evolution of strain in a creep experiment,
occurring as a consequence of the sudden and later constant stress, i.e.,
stress prescribed as the Heaviside step function.

The important property of relaxation modulus being completely monotonic and
creep compliance being a Bernstein function is discussed in \cite{Mai-10}.
On the other hand, the complete monotonicity of relaxation modus in the case
of distributed-order fractional Zener model is studied in \cite%
{BazhlekovaBazhlekov}, while \cite{MainardiSpada} deals with the fractional
Burgers model in creep and stress relaxation tests. The relaxation modulus
and creep compliance corresponding to the eight thermodynamically consistent
fractional Burgers models formulated in \cite{OZ-1} are examined in \cite%
{OZ-2}, where the thermodynamical requirements are narrowed in order to
ensure the complete monotonicity of the relaxation modulus and to ensure
that creep compliance is a Bernstein function. The role of Fox function is
underlined in \cite{GlockleNonnenmacher1,GlockleNonnenmacher2} in seeking
the responses of fractional constitutive models if stress or strain is
assumed to be given. The experimental data, obtained in creep and stress
relaxation experiments on biological tissues, is used for curve fitting in 
\cite{DemirciTonuk,grah}. In \cite{Makris}, the creep compliances in the
steady state regime are reviewed in the case of integer-order models of
viscoelasticity. Considering the inertial effects of the one-dimensional
viscoelastic body, the displacement and stress in stress relaxation and
creep tests are analyzed in \cite{APZ-4,APZ-3} when material is described by
the distributed-order fractional model, while in \cite{SD} the
thermodynamically consistent fractional Burgers models are used.

The qualitative characteristics of relaxation modulus and creep compliance
play a crucial role in energy balance properties of the viscoelastic body,
and therefore their explicit form is calculated and used in order to pose
the conditions guaranteeing that the relaxation modulus is a completely
monotone function and that creep compliance is a Bernstein function.
Previously mentioned conditions prove to narrow the thermodynamical
restrictions on model parameters appearing in the fractional anti-Zener and
Zener models, that are obtained in \cite{SD-1}. The thermodynamical
conditions actually allow for the relaxation modulus and creep compliance to
be non-monotonic and even oscillatory functions with exponentially
decreasing amplitude, which does not necessarily violate the positivity of
stored energy and dissipated power per unit volume.

Fractional anti-Zener and Zener constitutive equations, listed in Appendix %
\ref{FAZ-ZM}, take the following form%
\begin{equation}
\tilde{\sigma}\left( s\right) =s^{\xi }\frac{\phi _{\varepsilon }\left(
s\right) }{\phi _{\sigma }\left( s\right) }\tilde{\varepsilon}\left(
s\right)   \label{konst-jednacina}
\end{equation}%
after the Laplace transform, defined as%
\begin{equation*}
\hat{f}(s)=\mathcal{L}[f(t)](s)=\int_{0}^{\infty }f(t)\,\mathrm{e}^{-st}%
\mathrm{d}t,\quad \text{for}\quad \mathrm{Re}\,s>0,
\end{equation*}%
is applied, where $\phi _{\varepsilon }$ and $\phi _{\sigma }$ are model
dependent constitutive functions in the Laplace domain, that are, along with
the model dependent constitutive parameter $\xi $, listed in Table \ref%
{skupina}. 
\begin{landscape}
\input{tabela-svih-modela.tex} 
\end{landscape}Constitutive model in the Laplace domain (\ref%
{konst-jednacina}) is transformed into the relaxation modulus in the Laplace
domain%
\begin{equation}
\tilde{\sigma}_{sr}\left( s\right) =\frac{1}{s^{1-\xi }}\frac{\phi
_{\varepsilon }\left( s\right) }{\phi _{\sigma }\left( s\right) },
\label{sigma-sr-ld}
\end{equation}%
when the strain is prescribed in the form of Heaviside step function, i.e.,
as $\varepsilon =H,$ while if the stress is prescribed in the form of
Heaviside step function, i.e., as $\sigma =H$, then the constitutive model (%
\ref{konst-jednacina}) becomes the creep compliance in the Laplace domain,
expressed as%
\begin{equation}
\tilde{\varepsilon}_{cr}\left( s\right) =\frac{1}{s^{1+\xi }}\frac{\phi
_{\sigma }\left( s\right) }{\phi _{\varepsilon }\left( s\right) }.
\label{epsilon-cr-ld}
\end{equation}

The use of constitutive models in modeling wave propagation is extensive,
where some of the first studies of viscoelastic materials of fractional
order is performed in \cite{CaputoMainardi-1971b,CaputoMainardi-1971a}.
Damped oscillations and wave propagation in viscoelastic materials modeled
by the Zener, modified Zener, and modified Maxwell models is conducted in 
\cite{R-S1,R-S-2001,R-S,R-S2,R-S-2008}. Waves in fractional Zener and
distributed type media are investigated in \cite{OparnicaBroucke,KOZ10} and
in \cite{OparnicaBroucke1,KOZ19}, while in \cite{OZO} wave propagation in
the fractional Burgers type material is considered. Multidimensional
generalizations of the fractional Zener model and its application in wave
propagation is studied in \cite{AitIchou,CunhaFilho,OparnicaSuli}.
Fractional order constitutive equations of viscoelastic media, along with
the wave propagation, dispersion, and attenuation processes are reviewed in 
\cite{APSZ-1,APSZ-2,Holm-book,Mai-10,R-S-2010}, while \cite{Cai2018} surveys
the propagation of acoustic waves in complex media.

\section{Models' dissipativity properties in the time domain\label{power}}

In order to investigate energy balance properties of the one-dimensional
viscoelastic body, the constitutive equation in the Laplace domain (\ref%
{konst-jednacina}) is rewritten either as 
\begin{equation}
\tilde{\sigma}\left( s\right) =s\tilde{\sigma}_{sr}\left( s\right) \tilde{%
\varepsilon}\left( s\right) ,\quad \text{or as}\quad \tilde{\varepsilon}%
\left( s\right) =s\tilde{\varepsilon}_{cr}\left( s\right) \tilde{\sigma}%
\left( s\right) ,  \label{sigma-epsilon-ld}
\end{equation}%
using the relaxation modulus (\ref{sigma-sr-ld}) and creep compliance (\ref%
{epsilon-cr-ld}) in the Laplace domain, so that, after performing the
inverse Laplace transform in (\ref{sigma-epsilon-ld}), the constitutive
equation in time domain reads either%
\begin{equation}
\sigma \left( t\right) =\frac{\mathrm{d}}{\mathrm{d}t}\left( \sigma
_{sr}\left( t\right) \ast \varepsilon \left( t\right) \right) ,\quad \text{or%
}\quad \varepsilon \left( t\right) =\dot{\varepsilon}_{cr}\left( t\right)
\ast \sigma \left( t\right) ,  \label{sigma-epsilon}
\end{equation}%
since $\left. \left( \sigma _{sr}\left( t\right) \ast \varepsilon \left(
t\right) \right) \right\vert _{t=0}=0$ and $\varepsilon _{cr}^{\left(
g\right) }=\varepsilon _{cr}\left( 0\right) =0$. The former holds true if
the relaxation modulus tends to infinity as a power-type function for small
time, i.e., if $\sigma _{sr}\left( t\right) \sim Kt^{-\alpha },$ as $%
t\rightarrow 0$ with $\alpha \in \left( 0,1\right) $, as well as if strain
is bounded at zero by $\varepsilon _{0}$, since then $\sigma _{sr}\left(
t\right) \ast \varepsilon \left( t\right) =\int_{0}^{t}\sigma _{sr}\left(
t^{\prime }\right) \varepsilon \left( t-t^{\prime }\right) \mathrm{d}%
t^{\prime }\sim K\varepsilon _{0}\int_{0}^{t}\left( t^{\prime }\right)
^{-\alpha }\mathrm{d}t^{\prime }=K\varepsilon _{0}\frac{t^{1-\alpha }}{%
1-\alpha }\rightarrow 0$ as $t\rightarrow 0$, while the latter holds true,
since the glass compliance $\varepsilon _{cr}^{\left( g\right) }$ is zero in 
$\mathcal{L}^{-1}\left[ s\tilde{\varepsilon}_{cr}\left( s\right) \right] =%
\dot{\varepsilon}_{cr}\left( t\right) +\varepsilon _{cr}^{\left( g\right)
}\delta \left( t\right) $ for all considered models.

The power per unit volume, expressed by%
\begin{equation}
P\left( t\right) =\sigma \left( t\right) \dot{\varepsilon}\left( t\right) ,
\label{snaga}
\end{equation}%
can be written either in terms of strain$\ $as%
\begin{align}
P\left( t\right) & =\frac{\mathrm{d}}{\mathrm{d}t}\left( \frac{1}{2}\sigma
_{sr}\left( t\right) \varepsilon ^{2}\left( t\right) +\frac{1}{2}%
\int_{0}^{t}\left( -\dot{\sigma}_{sr}\left( t-t^{\prime }\right) \right)
\left( \varepsilon \left( t\right) -\varepsilon \left( t^{\prime }\right)
\right) ^{2}\mathrm{d}t^{\prime }\right)  \notag \\
&\quad +\frac{1}{2}\left( -\dot{\sigma}_{sr}\left( t\right) \right)
\varepsilon ^{2}\left( t\right) +\frac{1}{2}\int_{0}^{t}\ddot{\sigma}%
_{sr}\left( t-t^{\prime }\right) \left( \varepsilon \left( t\right)
-\varepsilon \left( t^{\prime }\right) \right) ^{2}\mathrm{d}t^{\prime },
\label{snaga-epsilon}
\end{align}%
using (\ref{sigma-epsilon})$_{1},$ or in terms of stress%
\begin{align}
P\left( t\right) & =\frac{\mathrm{d}}{\mathrm{d}t}\left( \frac{1}{2}%
\int_{0}^{t}\dot{\varepsilon}_{cr}\left( t-t^{\prime }\right) \sigma
^{2}\left( t^{\prime }\right) \mathrm{d}t^{\prime }\right)  \notag \\
&\quad +\frac{1}{2}\dot{\varepsilon}_{cr}\left( t\right) \sigma ^{2}\left(
t\right) +\frac{1}{2}\int_{0}^{t}\left( -\ddot{\varepsilon}_{cr}\left(
t\right) \left( t-t^{\prime }\right) \right) \left( \sigma \left( t\right)
-\sigma \left( t^{\prime }\right) \right) ^{2}\mathrm{d}t^{\prime },
\label{snaga-sigma}
\end{align}%
using (\ref{sigma-epsilon})$_{2},$ where the terms%
\begin{align}
W\left( t\right) & =\frac{1}{2}\sigma _{sr}\left( t\right) \varepsilon
^{2}\left( t\right) +\frac{1}{2}\int_{0}^{t}\left( -\dot{\sigma}_{sr}\left(
t-t^{\prime }\right) \right) \left( \varepsilon \left( t\right) -\varepsilon
\left( t^{\prime }\right) \right) ^{2}\mathrm{d}t^{\prime }>0\quad \text{and}
\label{pot-en-epsilon} \\
W\left( t\right) & =\frac{1}{2}\int_{0}^{t}\dot{\varepsilon}_{cr}\left(
t-t^{\prime }\right) \sigma ^{2}\left( t^{\prime }\right) \mathrm{d}%
t^{\prime }>0,  \label{pot-en-sigma}
\end{align}%
occurring in (\ref{snaga-epsilon}) and (\ref{snaga-sigma}) respectively, can
be interpreted as the energy per unit volume stored in the viscoelastic
body, taking into account the memory of strain, respectively stress,
weighted by the derivative of relaxation modulus, respectively creep
compliance, carrying the information about material properties, while the
first term in (\ref{pot-en-epsilon}) resembles to the potential energy of
the elastic body with Young's modulus replaced by the relaxation modulus. On
the other hand, the terms%
\begin{align}
\mathcal{P}\left( t\right) & =\frac{1}{2}\left( -\dot{\sigma}_{sr}\left(
t\right) \right) \varepsilon ^{2}\left( t\right) +\frac{1}{2}\int_{0}^{t}%
\ddot{\sigma}_{sr}\left( t-t^{\prime }\right) \left( \varepsilon \left(
t\right) -\varepsilon \left( t^{\prime }\right) \right) ^{2}\mathrm{d}%
t^{\prime }>0\quad \text{and}  \label{P-epsilon} \\
\mathcal{P}\left( t\right) & =\frac{1}{2}\dot{\varepsilon}_{cr}\left(
t\right) \sigma ^{2}\left( t\right) +\frac{1}{2}\int_{0}^{t}\left( -\ddot{%
\varepsilon}_{cr}\left( t\right) \left( t-t^{\prime }\right) \right) \left(
\sigma \left( t\right) -\sigma \left( t^{\prime }\right) \right) ^{2}\mathrm{%
d}t^{\prime }>0,  \label{P-sigma}
\end{align}%
occurring in (\ref{snaga-epsilon}) and (\ref{snaga-sigma}) respectively, can
be interpreted as the dissipated power per unit volume having two types of
contributions: the first one being instantaneous and depending on a material
dependent positive and decreasing function $-\dot{\sigma}_{sr}$,
respectively $\dot{\varepsilon}_{cr}$, and the second contribution being of
hereditary type with material dependent kernel. The positivity of energy and
dissipativity of power is guaranteed by the properties of the relaxation
modulus and creep compliance, namely by requesting the relaxation modulus to
be completely monotonic function, i.e.,%
\begin{equation*}
\sigma _{sr}\left( t\right) \geqslant 0\quad \text{and}\quad \left(
-1\right) ^{k}\frac{\mathrm{d}^{k}}{\mathrm{d}t^{k}}\dot{\sigma}_{sr}\left(
t\right) \leqslant 0,\quad \text{for}\quad t>0,\quad k\in 
%TCIMACRO{\U{2115} }%
%BeginExpansion
\mathbb{N}
%EndExpansion
_{0},
\end{equation*}%
as well as by requesting the creep compliance to be Bernstein's function,
i.e., a positive function with completely monotonic first derivative.

Introducing the energy and dissipated power per unit volume, the power per
unit volume can be rewritten as%
\begin{equation*}
P\left( t\right) =\frac{\mathrm{d}}{\mathrm{d}t}W\left( t\right) +\mathcal{P}%
\left( t\right) ,
\end{equation*}%
where the first term corresponds to the elastic type properties of
viscoelastic body, while the second term corresponds to its properties of
viscous type.

The expression (\ref{snaga-epsilon}) for the power per unit volume follows
from 
\begin{equation*}
P\left( t\right) =\frac{\mathrm{d}}{\mathrm{d}t}\left( \sigma _{sr}\left(
t\right) \ast \varepsilon \left( t\right) \right) \dot{\varepsilon}\left(
t\right) ,
\end{equation*}%
obtained by (\ref{snaga}) and (\ref{sigma-epsilon})$_{1},$ transforming into%
\begin{align}
P\left( t\right) &=\frac{\mathrm{d}}{\mathrm{d}t}\left( \frac{\mathrm{d}}{%
\mathrm{d}t}\left( \sigma _{sr}\left( t\right) \ast \varepsilon \left(
t\right) \right) \varepsilon \left( t\right) \right) -\frac{\mathrm{d}}{%
\mathrm{d}t}\left( \dot{\sigma}_{sr}\left( t\right) \ast \varepsilon \left(
t\right) +\sigma _{sr}^{\left( g\right) }\varepsilon \left( t\right) \right)
\varepsilon \left( t\right)  \notag \\
&=\frac{\mathrm{d}}{\mathrm{d}t}\left( \frac{\mathrm{d}}{\mathrm{d}t}\left(
\sigma _{sr}\left( t\right) \ast \varepsilon \left( t\right) \right)
\varepsilon \left( t\right) \right) -\frac{\mathrm{d}}{\mathrm{d}t}\left( 
\dot{\sigma}_{sr}\left( t\right) \ast \varepsilon \left( t\right) \right)
\varepsilon \left( t\right) -\frac{\mathrm{d}}{\mathrm{d}t}\left( \frac{1}{2}%
\sigma _{sr}^{\left( g\right) }\varepsilon ^{2}\left( t\right) \right) ,
\label{p-eps}
\end{align}%
using the derivative of a product of two functions and a derivative of a
convolution%
\begin{equation*}
\frac{\mathrm{d}}{\mathrm{d}t}\left( f\left( t\right) \ast g\left( t\right)
\right) =\dot{f}\left( t\right) \ast g\left( t\right) +f\left( 0\right)
g\left( t\right) ,
\end{equation*}%
since the first two terms in (\ref{p-eps}) become%
\begin{eqnarray*}
&&\frac{\mathrm{d}}{\mathrm{d}t}\left( \sigma _{sr}\left( t\right) \ast
\varepsilon \left( t\right) \right) \varepsilon \left( t\right) \\
&&\qquad =\frac{1}{2}\sigma _{sr}\left( t\right) \varepsilon ^{2}\left(
t\right) +\frac{1}{2}\frac{\mathrm{d}}{\mathrm{d}t}\left( \sigma _{sr}\left(
t\right) \ast \varepsilon ^{2}\left( t\right) \right) -\frac{1}{2}%
\int_{0}^{t}\dot{\sigma}_{sr}\left( t-t^{\prime }\right) \left( \varepsilon
\left( t\right) -\varepsilon \left( t^{\prime }\right) \right) ^{2}\mathrm{d}%
t^{\prime } \\
&&\qquad =\frac{1}{2}\sigma _{sr}\left( t\right) \varepsilon ^{2}\left(
t\right) +\frac{1}{2}\dot{\sigma}_{sr}\left( t\right) \ast \varepsilon
^{2}\left( t\right) +\frac{1}{2}\sigma _{sr}^{\left( g\right) }\varepsilon
^{2}\left( t\right) -\frac{1}{2}\int_{0}^{t}\dot{\sigma}_{sr}\left(
t-t^{\prime }\right) \left( \varepsilon \left( t\right) -\varepsilon \left(
t^{\prime }\right) \right) ^{2}\mathrm{d}t^{\prime }
\end{eqnarray*}%
and%
\begin{eqnarray*}
&&\frac{\mathrm{d}}{\mathrm{d}t}\left( \dot{\sigma}_{sr}\left( t\right) \ast
\varepsilon \left( t\right) \right) \varepsilon \left( t\right) \\
&&\qquad =\frac{1}{2}\dot{\sigma}_{sr}\left( t\right) \varepsilon ^{2}\left(
t\right) +\frac{1}{2}\frac{\mathrm{d}}{\mathrm{d}t}\left( \dot{\sigma}%
_{sr}\left( t\right) \ast \varepsilon ^{2}\left( t\right) \right) -\frac{1}{2%
}\int_{0}^{t}\ddot{\sigma}_{sr}\left( t-t^{\prime }\right) \left(
\varepsilon \left( t\right) -\varepsilon \left( t^{\prime }\right) \right)
^{2}\mathrm{d}t^{\prime },
\end{eqnarray*}%
according to 
\begin{equation}
\frac{\mathrm{d}}{\mathrm{d}t}\left( k\left( t\right) \ast u\left( t\right)
\right) u\left( t\right) =\frac{1}{2}k\left( t\right) u^{2}\left( t\right) +%
\frac{1}{2}\frac{\mathrm{d}}{\mathrm{d}t}\left( k\left( t\right) \ast
u^{2}\left( t\right) \right) -\frac{1}{2}\int_{0}^{t}\dot{k}\left(
t-t^{\prime }\right) \left( u\left( t\right) -u\left( t^{\prime }\right)
\right) ^{2}\mathrm{d}t^{\prime }  \label{gomila}
\end{equation}%
and a derivative of a convolution.

On the other hand, the expression (\ref{snaga-sigma}) for the power per unit
volume follows from 
\begin{equation*}
P\left( t\right) =\frac{\mathrm{d}}{\mathrm{d}t}\left( \dot{\varepsilon}%
_{cr}\left( t\right) \ast \sigma \left( t\right) \right) \sigma \left(
t\right) ,
\end{equation*}%
obtained by (\ref{snaga}) and (\ref{sigma-epsilon})$_{2},$ transforming into%
\begin{equation*}
P\left( t\right) =\frac{1}{2}\dot{\varepsilon}_{cr}\left( t\right) \sigma
^{2}\left( t\right) +\frac{1}{2}\frac{\mathrm{d}}{\mathrm{d}t}\left( \dot{%
\varepsilon}_{cr}\left( t\right) \ast \sigma ^{2}\left( t\right) \right) -%
\frac{1}{2}\int_{0}^{t}\ddot{\varepsilon}_{cr}\left( t-t^{\prime }\right)
\left( \sigma \left( t\right) -\sigma \left( t^{\prime }\right) \right) ^{2}%
\mathrm{d}t^{\prime },
\end{equation*}%
according to (\ref{gomila}).

\section{Stress relaxation and creep\label{sr-and-cr}}

Relaxation modulus and creep compliance play a crucial role in energy
balance properties of viscoelastic body, since being a completely monotonic
function and a Bernstein function respectively, they guarantee the
positivity of the energy stored in viscoelastic body, see (\ref%
{pot-en-epsilon}) and (\ref{pot-en-sigma}), as well as the positivity of the
dissipated power, see (\ref{P-epsilon}) and (\ref{P-sigma}). Therefore, the
relaxation modulus and creep compliance, given by (\ref{sigma-sr-ld}) and (%
\ref{epsilon-cr-ld}) in the Laplace domain, are calculated in the time
domain and their properties are investigated, yielding that the relaxation
modulus and creep compliance are completely monotonic and Bernstein
functions respectively, if the model parameters satisfy narrowed
thermodynamical restrictions in addition to request that the relaxation
modulus and creep compliance in the Laplace domain (\ref{sigma-sr-ld}) and (%
\ref{epsilon-cr-ld}) do not have poles, while the character of relaxation
modulus and creep compliance in time domain changes if their counterparts in
the Laplace domain have poles. Namely, in the case of negative real pole
there exists an additional term decaying exponentially in time, while if
there exists a pair of complex conjugated poles with negative real part,
then there is an additional term displaying damped oscillatory character.

\subsection{Relaxation modulus and creep compliance}

\subsubsection{Relaxation modulus}

The relaxation modulus is obtained as%
\begin{equation}
\sigma _{sr}(t)=\sigma _{sr}^{\left( \mathrm{NP}\right) }\left( t\right)
+\left\{ \!\!\!%
\begin{tabular}{ll}
$0$, & if $\tilde{\sigma}_{sr}$ has no poles, \smallskip \\ 
$\sigma _{sr}^{\left( \mathrm{RP}\right) }\left( t\right) $, & if $\tilde{%
\sigma}_{sr}$ has a negative real pole, \smallskip \\ 
$\sigma _{sr}^{\left( \mathrm{CCP}\right) }\left( t\right) $, & if $\tilde{%
\sigma}_{sr}$ has a pair of complex conjugated poles,%
\end{tabular}%
\ \right.  \label{sr-opste}
\end{equation}%
by inverting the Laplace transform in the relaxation modulus in the Laplace
domain $\tilde{\sigma}_{sr},$ given by (\ref{sigma-sr-ld}), using the
definition and integration in the complex plane along the Bromwich contour,
with the functions $\sigma _{sr}^{\left( \mathrm{NP}\right) },$ $\sigma
_{sr}^{\left( \mathrm{RP}\right) },$ and $\sigma _{sr}^{\left( \mathrm{CCP}%
\right) }$ respectively taking the forms%
\begin{align}
\sigma _{sr}^{\left( \mathrm{NP}\right) }\left( t\right) & =\frac{1}{\pi }%
\int_{0}^{\infty }\frac{1}{\rho ^{1-\xi }}\frac{\left\vert \phi
_{\varepsilon }\left( \rho \mathrm{e}^{\mathrm{i}\pi }\right) \right\vert }{%
\left\vert \phi _{\sigma }\left( \rho \mathrm{e}^{\mathrm{i}\pi }\right)
\right\vert }\sin \left( \arg \phi _{\varepsilon }\left( \rho \mathrm{e}^{%
\mathrm{i}\pi }\right) -\arg \phi _{\sigma }\left( \rho \mathrm{e}^{\mathrm{i%
}\pi }\right) +\xi \pi \right) \mathrm{e}^{-\rho t}\mathrm{d}\rho ,
\label{sigma-NP} \\
\sigma _{sr}^{\left( \mathrm{RP}\right) }\left( t\right) & =-\frac{1}{\rho _{%
\scriptscriptstyle{\mathrm{RP}}}^{1-\xi }}\frac{\left\vert \phi
_{\varepsilon }\left( s_{\scriptscriptstyle{\mathrm{RP}}}\right) \right\vert 
}{\left\vert \phi _{\sigma }^{\prime }\left( s_{\scriptscriptstyle{\mathrm{RP%
}}}\right) \right\vert }\cos \left( \arg \phi _{\varepsilon }\left( s_{%
\scriptscriptstyle{\mathrm{RP}}}\right) -\arg \phi _{\sigma }^{\prime
}\left( s_{\scriptscriptstyle{\mathrm{RP}}}\right) +\xi \pi \right) \mathrm{e%
}^{-\rho _{\scriptscriptstyle{\mathrm{RP}}}t},  \label{sigma-RP} \\
\sigma _{sr}^{\left( \mathrm{CCP}\right) }\left( t\right) & =2\frac{1}{\rho
_{\scriptscriptstyle{\mathrm{CCP}}}^{1-\xi }}\frac{\left\vert \phi
_{\varepsilon }\left( s_{\scriptscriptstyle{\mathrm{CCP}}}\right)
\right\vert }{\left\vert \phi _{\sigma }^{\prime }\left( s_{%
\scriptscriptstyle{\mathrm{CCP}}}\right) \right\vert }\mathrm{e}%
^{-\left\vert \func{Re}s_{\scriptscriptstyle{\mathrm{CCP}}}\right\vert
t}\cos \left( \func{Im}s_{\scriptscriptstyle{\mathrm{CCP}}}t+\arg \phi
_{\varepsilon }\left( s_{\scriptscriptstyle{\mathrm{CCP}}}\right) -\arg \phi
_{\sigma }^{\prime }\left( s_{\scriptscriptstyle{\mathrm{CCP}}}\right)
-\left( 1-\xi \right) \varphi _{\scriptscriptstyle{\mathrm{CCP}}}\right) ,
\label{sigma-CCP}
\end{align}%
where $\phi _{\sigma }^{\prime }\left( s\right) =\frac{\mathrm{d}}{\mathrm{d}%
s}\phi _{\sigma }\left( s\right) $ and where poles of function $\tilde{\sigma%
}_{sr}$: $s_{\scriptscriptstyle{\mathrm{RP}}}=\rho _{\scriptscriptstyle{%
\mathrm{RP}}}\,\mathrm{e}^{\mathrm{i}\pi }$ and $s_{\scriptscriptstyle{%
\mathrm{CCP}}}=\rho _{\scriptscriptstyle{\mathrm{CCP}}}\,\mathrm{e}^{\mathrm{%
i}\varphi _{\scriptscriptstyle{\mathrm{CCP}}}}$ are respectively a negative
real zero of function $\phi _{\sigma }$ and its complex zero having negative
real part.

The equivalent form of function $\sigma _{sr}^{\left( \mathrm{NP}\right) },$
given by (\ref{sigma-NP}), is%
\begin{equation}
\sigma _{sr}^{\left( \mathrm{NP}\right) }\left( t\right) =\frac{1}{\pi }%
\int_{0}^{\infty }\frac{1}{\rho ^{1-\xi }}\frac{K\left( \rho \right) }{%
\left\vert \phi _{\sigma }\left( \rho \mathrm{e}^{\mathrm{i}\pi }\right)
\right\vert ^{2}}\mathrm{e}^{-\rho t}\mathrm{d}\rho ,  \label{sigma-NP-1}
\end{equation}%
with 
\begin{align}
K\left( \rho \right) & =\frac{1}{2\mathrm{i}}\left( \mathrm{e}^{\mathrm{i}%
\xi \pi }\phi _{\varepsilon }\left( \rho \mathrm{e}^{\mathrm{i}\pi }\right) 
\bar{\phi}_{\sigma }\left( \rho \mathrm{e}^{\mathrm{i}\pi }\right) -\mathrm{e%
}^{-\mathrm{i}\xi \pi }\bar{\phi}_{\varepsilon }\left( \rho \mathrm{e}^{%
\mathrm{i}\pi }\right) \phi _{\sigma }\left( \rho \mathrm{e}^{\mathrm{i}\pi
}\right) \right)  \notag \\
& =\left\vert \phi _{\varepsilon }\left( \rho \mathrm{e}^{\mathrm{i}\pi
}\right) \right\vert \left\vert \phi _{\sigma }\left( \rho \mathrm{e}^{%
\mathrm{i}\pi }\right) \right\vert \sin \left( \arg \phi _{\varepsilon
}\left( \rho \mathrm{e}^{\mathrm{i}\pi }\right) -\arg \phi _{\sigma }\left(
\rho \mathrm{e}^{\mathrm{i}\pi }\right) +\xi \pi \right)  \notag \\
& =\cos \left( \xi \pi \right) \left( \mathrm{\func{Im}}\phi _{\varepsilon
}\left( \rho \mathrm{e}^{\mathrm{i}\pi }\right) \func{Re}\phi _{\sigma
}\left( \rho \mathrm{e}^{\mathrm{i}\pi }\right) -\func{Re}\phi _{\varepsilon
}\left( \rho \mathrm{e}^{\mathrm{i}\pi }\right) \func{Im}\phi _{\sigma
}\left( \rho \mathrm{e}^{\mathrm{i}\pi }\right) \right)  \notag \\
& \quad +\sin \left( \xi \pi \right) \left( \mathrm{\func{Re}}\phi
_{\varepsilon }\left( \rho \mathrm{e}^{\mathrm{i}\pi }\right) \func{Re}\phi
_{\sigma }\left( \rho \mathrm{e}^{\mathrm{i}\pi }\right) +\func{Im}\phi
_{\varepsilon }\left( \rho \mathrm{e}^{\mathrm{i}\pi }\right) \func{Im}\phi
_{\sigma }\left( \rho \mathrm{e}^{\mathrm{i}\pi }\right) \right) ,  \label{K}
\end{align}%
being convenient for examining whether the function $\sigma _{sr}^{\left( 
\mathrm{NP}\right) }$ is completely monotonic, that is guaranteed by
requiring%
\begin{equation}
K\left( \rho \right) \geqslant 0,\quad \text{i.e.,}\quad \sin \left( \arg
\phi _{\varepsilon }\left( \rho \mathrm{e}^{\mathrm{i}\pi }\right) -\arg
\phi _{\sigma }\left( \rho \mathrm{e}^{\mathrm{i}\pi }\right) +\xi \pi
\right) \geqslant 0\quad \text{for}\quad \rho \geqslant 0,  \label{kondisns}
\end{equation}%
while if the condition (\ref{kondisns}) is not satisfied, then function $%
\sigma _{sr}^{\left( \mathrm{NP}\right) }$ can be non-monotonic. Function $%
\sigma _{sr}^{\left( \mathrm{RP}\right) }$ is either a positive
exponentially decreasing function tending to zero or a negative
exponentially increasing function also tending to zero, with time constant $%
\rho _{\scriptscriptstyle{\mathrm{RP}}},$ obtained as a negative real pole
of function $\tilde{\sigma}_{sr},$ see (\ref{sigma-RP}), while function $%
\sigma _{sr}^{\left( \mathrm{CCP}\right) }$ is an oscillatory function of an
exponentially decreasing amplitude, having angular frequency defined by the
imaginary part of a complex pole of function $\tilde{\sigma}_{sr},$ i.e., by 
$\func{Im}s_{\scriptscriptstyle{\mathrm{CCP}}}$, and damping parameter
defined by the real part of a complex pole of function $\tilde{\sigma}_{sr}$%
, i.e., by $\left\vert \func{Re}s_{\scriptscriptstyle{\mathrm{CCP}}%
}\right\vert ,$ see (\ref{sigma-CCP}).

If the function $\phi _{\sigma },$ appearing in the relaxation modulus in
the Laplace domain (\ref{sigma-sr-ld}), consists of two terms, as it is the
case for models ID.ID, ID.DD$^{{}^{+}}$, ID.IDD, ID.DDD$^{{}^{+}}$, and
ID.IDD$^{{}^{+}}$, see Table \ref{skupina}, then the relaxation modulus can
be written in terms of the two-parameter Mittag-Leffler function $e_{\xi
,\zeta ,\lambda }$ as%
\begin{equation}
\sigma _{sr}\left( t\right) =\frac{b_{1}}{a_{2}}e_{\alpha +\beta ,1-\xi
+\alpha +\beta ,\frac{a_{1}}{a_{2}}}\left( t\right) +\frac{b_{2}}{a_{2}}%
e_{\alpha +\beta ,1-\xi +\alpha +\beta -\lambda ,\frac{a_{1}}{a_{2}}}\left(
t\right) +\frac{b_{3}}{a_{2}}e_{\alpha +\beta ,1-\xi +\alpha +\beta -\kappa ,%
\frac{a_{1}}{a_{2}}}\left( t\right) ,  \label{sr-ML}
\end{equation}%
where one takes $b_{3}=0$ for models ID.ID and ID.DD$^{{}^{+}}$, with the
notation%
\begin{equation}
e_{\xi ,\zeta ,\lambda }\left( t\right) =t^{\zeta -1}E_{\xi ,\zeta }\left(
-\lambda t^{\xi }\right) ,\quad \text{where}\quad E_{\xi ,\zeta }\left(
z\right) =\sum_{n=0}^{\infty }\frac{z^{n}}{\Gamma \left( \xi n+\zeta \right) 
},  \label{ML2}
\end{equation}%
see \cite{GoreMai}. Namely, the relaxation modulus in the Laplace domain (%
\ref{sigma-sr-ld}) for the mentioned models can be written as%
\begin{equation*}
\tilde{\sigma}_{sr}\left( s\right) =\frac{1}{s^{1-\xi }}\frac{%
b_{1}+b_{2}s^{\lambda }+b_{3}s^{\kappa }}{a_{1}+a_{2}s^{\alpha +\beta }},
\end{equation*}%
see Table \ref{skupina}, transforming into%
\begin{equation*}
\tilde{\sigma}_{sr}\left( s\right) =\frac{b_{1}}{a_{2}}\frac{s^{-\left(
1-\xi \right) }}{s^{\alpha +\beta }+\frac{a_{1}}{a_{2}}}+\frac{b_{2}}{a_{2}}%
\frac{s^{-\left( 1-\xi -\lambda \right) }}{s^{\alpha +\beta }+\frac{a_{1}}{%
a_{2}}}+\frac{b_{3}}{a_{2}}\frac{s^{-\left( 1-\xi -\kappa \right) }}{%
s^{\alpha +\beta }+\frac{a_{1}}{a_{2}}},
\end{equation*}%
that yields the relaxation modulus in the form (\ref{sr-ML}) by the Laplace
transform of the two-parameter Mittag-Leffler function%
\begin{equation}
e_{\xi ,\zeta ,\lambda }\left( t\right) =\mathcal{L}^{-1}\left[ \frac{s^{\xi
-\zeta }}{s^{\xi }+\lambda }\right] \left( t\right) .  \label{ML2-ld}
\end{equation}

The integral representation of the relaxation modulus, already given by (\ref%
{sigma-NP-1}) in the case when $\tilde{\sigma}_{sr}$ has no poles, using the
relaxation modulus (\ref{sr-ML}) expressed in terms of the two-parameter
Mittag-Leffler function (\ref{ML2}), according to its integral
representation 
\begin{equation}
e_{\xi ,\zeta ,\lambda }\left( t\right) =\frac{1}{\pi }\int_{0}^{\infty }%
\frac{\lambda \sin \left( \left( \zeta -\xi \right) \pi \right) +\rho ^{\xi
}\sin \left( \zeta \pi \right) }{\left\vert \rho ^{\xi }\mathrm{e}^{\mathrm{i%
}\xi \pi }+\lambda \right\vert ^{2}}\rho ^{\xi -\zeta }\mathrm{e}^{-\rho t}%
\mathrm{d}\rho ,  \label{ML-IR}
\end{equation}%
see \cite{GoreMai}, provided that $\xi \in \left( 0,1\right) \ $and $\zeta
<1+\xi $, is reobtained in the form%
\begin{align}
\sigma _{sr}\left( t\right) & =\frac{1}{a_{2}\pi }\int_{0}^{\infty }\frac{1}{%
\rho ^{1-\xi }}\left( \frac{a_{1}}{a_{2}}\frac{b_{1}\sin \left( \xi \pi
\right) +b_{2}\rho ^{\lambda }\sin \left( \left( \xi +\lambda \right) \pi
\right) +b_{3}\rho ^{\kappa }\sin \left( \left( \xi +\kappa \right) \pi
\right) }{\left\vert \rho ^{\alpha +\beta }\mathrm{e}^{\mathrm{i}\left(
\alpha +\beta \right) \pi }+\frac{a_{1}}{a_{2}}\right\vert ^{2}}\right.  
\notag \\
& \quad +\left. \rho ^{\alpha +\beta }\frac{b_{1}\sin \left( \left( \xi
-\alpha -\beta \right) \pi \right) +b_{2}\rho ^{\lambda }\sin \left( \left(
\xi +\lambda -\alpha -\beta \right) \pi \right) +b_{3}\rho ^{\kappa }\sin
\left( \left( \xi +\kappa -\alpha -\beta \right) \pi \right) }{\left\vert
\rho ^{\alpha +\beta }\mathrm{e}^{\mathrm{i}\left( \alpha +\beta \right) \pi
}+\frac{a_{1}}{a_{2}}\right\vert ^{2}}\right) \mathrm{e}^{-\rho t}\mathrm{d}%
\rho ,  \label{sigma-NP-2}
\end{align}%
so that its complete monotonicity, provided that $\alpha +\beta \in \left(
0,1\right) ,$ is ensured by requiring the function%
\begin{align*}
K\left( \rho \right) & =a_{1}b_{1}\sin \left( \xi \pi \right)
+a_{1}b_{2}\rho ^{\lambda }\sin \left( \left( \xi +\lambda \right) \pi
\right) +a_{1}b_{3}\rho ^{\kappa }\sin \left( \left( \xi +\kappa \right) \pi
\right) +a_{2}b_{1}\rho ^{\alpha +\beta }\sin \left( \left( \xi -\alpha
-\beta \right) \pi \right)  \\
& \quad +a_{2}b_{2}\rho ^{\alpha +\beta +\lambda }\sin \left( \left( \xi
+\lambda -\alpha -\beta \right) \pi \right) +a_{2}b_{3}\rho ^{\alpha +\beta
+\kappa }\sin \left( \left( \xi +\kappa -\alpha -\beta \right) \pi \right) 
\end{align*}%
to be non-negative for $\rho \geqslant 0$. Note, the forms of the relaxation
modulus (\ref{sigma-NP-1}) and (\ref{sigma-NP-2}) coincide.

\subsubsection{Creep compliance}

The creep compliance is obtained as%
\begin{equation}
\varepsilon _{cr}(t)=\varepsilon _{cr}^{\left( \mathrm{NP}\right) }\left(
t\right) +\left\{ \!\!\!%
\begin{tabular}{ll}
$0$, & if $\tilde{\varepsilon}_{cr}$ has no poles, \smallskip \\ 
$\varepsilon _{cr}^{\left( \mathrm{RP}\right) }\left( t\right) $, & if $%
\tilde{\varepsilon}_{cr}$ has a negative real pole, \smallskip \\ 
$\varepsilon _{cr}^{\left( \mathrm{CCP}\right) }\left( t\right) $, & if $%
\tilde{\varepsilon}_{cr}$ has a pair of complex conjugated poles,%
\end{tabular}%
\ \right.  \label{cr-opste}
\end{equation}%
with the functions $\varepsilon _{cr}^{\left( \mathrm{NP}\right) },$ $%
\varepsilon _{cr}^{\left( \mathrm{RP}\right) },$ and $\varepsilon
_{cr}^{\left( \mathrm{CCP}\right) }$ respectively taking the forms%
\begin{align}
\varepsilon _{cr}^{\left( \mathrm{NP}\right) }\left( t\right) & =\frac{1}{%
\pi }\int_{0}^{\infty }\frac{1}{\rho ^{1+\xi }}\frac{\left\vert \phi
_{\sigma }\left( \rho \mathrm{e}^{\mathrm{i}\pi }\right) \right\vert }{%
\left\vert \phi _{\varepsilon }\left( \rho \mathrm{e}^{\mathrm{i}\pi
}\right) \right\vert }\sin \left( \arg \phi _{\varepsilon }\left( \rho 
\mathrm{e}^{\mathrm{i}\pi }\right) -\arg \phi _{\sigma }\left( \rho \mathrm{e%
}^{\mathrm{i}\pi }\right) +\xi \pi \right) \left( 1-\mathrm{e}^{-\rho
t}\right) \mathrm{d}\rho ,  \label{epsilon-NP} \\
\varepsilon _{cr}^{\left( \mathrm{RP}\right) }\left( t\right) & =-\frac{1}{%
\rho _{\scriptscriptstyle{\mathrm{RP}}}^{1+\xi }}\frac{\left\vert \phi
_{\sigma }\left( s_{\scriptscriptstyle{\mathrm{RP}}}\right) \right\vert }{%
\left\vert \phi _{\varepsilon }^{\prime }\left( s_{\scriptscriptstyle{%
\mathrm{RP}}}\right) \right\vert }\cos \left( \arg \phi _{\varepsilon
}^{\prime }\left( s_{\scriptscriptstyle{\mathrm{RP}}}\right) -\arg \phi
_{\sigma }\left( s_{\scriptscriptstyle{\mathrm{RP}}}\right) +\xi \pi \right)
\left( 1-\mathrm{e}^{-\rho _{\scriptscriptstyle{\mathrm{RP}}}t}\right) ,
\label{epsilon-RP} \\
\varepsilon _{cr}^{\left( \mathrm{CCP}\right) }\left( t\right) & =2\frac{1}{%
\rho _{\scriptscriptstyle{\mathrm{CCP}}}^{1+\xi }}\frac{\left\vert \phi
_{\sigma }\left( s_{\scriptscriptstyle{\mathrm{CCP}}}\right) \right\vert }{%
\left\vert \phi _{\varepsilon }^{\prime }\left( s_{\scriptscriptstyle{%
\mathrm{CCP}}}\right) \right\vert }  \notag \\
& \qquad \times \bigg(\mathrm{e}^{-\left\vert \func{Re}s_{\scriptscriptstyle{%
\mathrm{CCP}}}\right\vert t}\cos \left( \func{Im}s_{\scriptscriptstyle{%
\mathrm{CCP}}}t-\arg \phi _{\varepsilon }^{\prime }\left( s_{%
\scriptscriptstyle{\mathrm{CCP}}}\right) +\arg \phi _{\sigma }\left( s_{%
\scriptscriptstyle{\mathrm{CCP}}}\right) -\left( 1+\xi \right) \varphi _{%
\scriptscriptstyle{\mathrm{CCP}}}\right)  \notag \\
& \qquad \qquad -\cos \left( \arg \phi _{\varepsilon }^{\prime }\left( s_{%
\scriptscriptstyle{\mathrm{CCP}}}\right) -\arg \phi _{\sigma }\left( s_{%
\scriptscriptstyle{\mathrm{CCP}}}\right) +\left( 1+\xi \right) \varphi _{%
\scriptscriptstyle{\mathrm{CCP}}}\right) \bigg),  \label{epsilon-CCP}
\end{align}%
where $\phi _{\varepsilon }^{\prime }\left( s\right) =\frac{\mathrm{d}}{%
\mathrm{d}s}\phi _{\varepsilon }\left( s\right) $ and where $s_{%
\scriptscriptstyle{\mathrm{RP}}}=\rho _{\scriptscriptstyle{\mathrm{RP}}}\,%
\mathrm{e}^{\mathrm{i}\pi }$ and $s_{\scriptscriptstyle{\mathrm{CCP}}}=\rho
_{\scriptscriptstyle{\mathrm{CCP}}}\,\mathrm{e}^{\mathrm{i}\varphi _{%
\scriptscriptstyle{\mathrm{CCP}}}}$ are respectively a negative real zero of
function $\phi _{\varepsilon }$ and its complex zero having negative real
part.

The equivalent form of function $\varepsilon _{cr}^{\left( \mathrm{NP}%
\right) },$ given by (\ref{epsilon-NP}), is%
\begin{equation}
\varepsilon _{cr}^{\left( \mathrm{NP}\right) }\left( t\right) =\frac{1}{\pi }%
\int_{0}^{\infty }\frac{1}{\rho ^{1+\xi }}\frac{K\left( \rho \right) }{%
\left\vert \phi _{\varepsilon }\left( \rho \mathrm{e}^{\mathrm{i}\pi
}\right) \right\vert ^{2}}\left( 1-\mathrm{e}^{-\rho t}\right) \mathrm{d}%
\rho ,  \label{epsilon-NP-1}
\end{equation}%
with function $K$ defined by (\ref{K}), so that the conditions (\ref%
{kondisns}), which guaranteed that the relaxation modulus $\sigma
_{sr}^{\left( \mathrm{NP}\right) },$ see (\ref{sigma-NP-1}), is completely
monotonic, also guarantee that the creep compliance $\varepsilon
_{cr}^{\left( \mathrm{NP}\right) }$ is a Bernstein function. On the other
hand, if the condition (\ref{kondisns}) is not satisfied, then function $%
\varepsilon _{cr}^{\left( \mathrm{NP}\right) }$ can be non-monotonic.
Function $\varepsilon _{cr}^{\left( \mathrm{RP}\right) }$ is a function
starting from zero, that either exponentially increases to a positive
horizontal asymptote, or exponentially decreases to a negative horizontal
asymptote, having the time constant $\rho _{\scriptscriptstyle{\mathrm{RP}}%
}, $ obtained as a negative real pole of function $\tilde{\varepsilon}_{cr},$
see (\ref{epsilon-RP}), while function $\varepsilon _{cr}^{\left( \mathrm{CCP%
}\right) }$ is a vertically shifted oscillatory function of an exponentially
decreasing amplitude, having angular frequency defined by the imaginary part
of a complex pole of function $\tilde{\varepsilon}_{cr},$ i.e., by $\func{Im}%
s_{\scriptscriptstyle{\mathrm{CCP}}}$, and damping parameter defined by the
real part of a complex pole of function $\tilde{\varepsilon}_{cr}$, i.e., by 
$\left\vert \func{Re}s_{\scriptscriptstyle{\mathrm{CCP}}}\right\vert ,$ see (%
\ref{epsilon-CCP}).

If the function $\phi _{\varepsilon },$ appearing in the creep compliance in
the Laplace domain (\ref{epsilon-cr-ld}), consists of two terms, as it is
the case for models ID.ID, ID.DD$^{{}^{+}}$, IID.ID, IDD.DD$^{{}^{+}}$, I$%
^{{}^{+}}$ID.ID, and IDD$^{{}^{+}}$.DD$^{{}^{+}}$, see Table \ref{skupina},
then the creep compliance can be written in terms of the two-parameter
Mittag-Leffler function $e_{\xi ,\zeta ,\lambda },$ given by (\ref{ML2}), as%
\begin{equation}
\varepsilon _{cr}\left( t\right) =\frac{a_{1}}{b_{2}}e_{\alpha +\beta ,1+\xi
+\alpha +\beta ,\frac{b_{1}}{b_{2}}}\left( t\right) +\frac{a_{2}}{b_{2}}%
e_{\alpha +\beta ,1+\xi +\alpha +\beta -\lambda ,\frac{b_{1}}{b_{2}}}\left(
t\right) +\frac{a_{3}}{b_{2}}e_{\alpha +\beta ,1+\xi +\alpha +\beta -\kappa ,%
\frac{b_{1}}{b_{2}}}\left( t\right) ,  \label{cr-ML}
\end{equation}%
where one takes $a_{3}=0$ for models ID.ID and ID.DD$^{{}^{+}}$. Namely, the
creep compliance in the Laplace domain (\ref{epsilon-cr-ld}) for the
mentioned models can be written as%
\begin{equation*}
\tilde{\varepsilon}_{cr}\left( s\right) =\frac{1}{s^{1+\xi }}\frac{%
a_{1}+a_{2}s^{\lambda }+a_{3}s^{\kappa }}{b_{1}+b_{2}s^{\alpha +\beta }},
\end{equation*}%
see Table \ref{skupina}, transforming into%
\begin{equation*}
\tilde{\varepsilon}_{cr}\left( s\right) =\frac{a_{1}}{b_{2}}\frac{s^{-\left(
1+\xi \right) }}{s^{\alpha +\beta }+\frac{b_{1}}{b_{2}}}+\frac{a_{2}}{b_{2}}%
\frac{s^{-\left( 1+\xi -\lambda \right) }}{s^{\alpha +\beta }+\frac{b_{1}}{%
b_{2}}}+\frac{a_{3}}{b_{2}}\frac{s^{-\left( 1+\xi -\kappa \right) }}{%
s^{\alpha +\beta }+\frac{b_{1}}{b_{2}}},
\end{equation*}%
that yields the creep compliance in the form (\ref{cr-ML}) by the Laplace
transform of the two-parameter Mittag-Leffler function, given by (\ref%
{ML2-ld}).

Note, the creep compliance written in terms of the two-parameter
Mittag-Leffler function (\ref{cr-ML}), does not admit the integral
representation in the form (\ref{ML-IR}), since the condition $\zeta <1+\xi $
is not satisfied for the first term (\ref{cr-ML}). However, the derivative
of creep compliance (\ref{cr-ML}), obtained as%
\begin{equation*}
\dot{\varepsilon}_{cr}\left( t\right) =\frac{a_{1}}{b_{2}}e_{\alpha +\beta
,\xi +\alpha +\beta ,\frac{b_{1}}{b_{2}}}\left( t\right) +\frac{a_{2}}{b_{2}}%
e_{\alpha +\beta ,\xi +\alpha +\beta -\lambda ,\frac{b_{1}}{b_{2}}}\left(
t\right) +\frac{a_{3}}{b_{2}}e_{\alpha +\beta ,\xi +\alpha +\beta -\kappa ,%
\frac{b_{1}}{b_{2}}}\left( t\right) ,
\end{equation*}%
according to the property $\frac{\mathrm{d}}{\mathrm{d}t}e_{\xi ,\zeta
,\lambda }\left( t\right) =e_{\xi ,\zeta -1,\lambda }\left( t\right) $ of
the two-parameter Mittag-Leffler function, admits the integral
representation in the form%
\begin{align}
\dot{\varepsilon}_{cr}\left( t\right) & =\frac{1}{b_{2}\pi }\int_{0}^{\infty
}\frac{1}{\rho ^{\xi }}\left( \frac{b_{1}}{b_{2}}\frac{a_{1}\sin \left( \xi
\pi \right) +a_{2}\rho ^{\lambda }\sin \left( \left( \xi -\lambda \right)
\pi \right) +a_{3}\rho ^{\kappa }\sin \left( \left( \xi -\kappa \right) \pi
\right) }{\left\vert \rho ^{\alpha +\beta }\mathrm{e}^{\mathrm{i}\left(
\alpha +\beta \right) \pi }+\frac{b_{1}}{b_{2}}\right\vert ^{2}}\right. 
\notag \\
& \quad +\left. \rho ^{\alpha +\beta }\frac{a_{1}\sin \left( \left( \xi
+\alpha +\beta \right) \pi \right) +a_{2}\rho ^{\lambda }\sin \left( \left(
\xi -\lambda +\alpha +\beta \right) \pi \right) +a_{3}\rho ^{\kappa }\sin
\left( \left( \xi -\kappa +\alpha +\beta \right) \pi \right) }{\left\vert
\rho ^{\alpha +\beta }\mathrm{e}^{\mathrm{i}\left( \alpha +\beta \right) \pi
}+\frac{b_{1}}{b_{2}}\right\vert ^{2}}\right) \mathrm{e}^{-\rho t}\mathrm{d}%
\rho ,  \label{epsilon-NP-2}
\end{align}%
according to the expression (\ref{ML-IR}), provided that $\alpha +\beta \in
\left( 0,1\right) ,$ so that by requiring the function%
\begin{align*}
K\left( \rho \right) & =a_{1}b_{1}\sin \left( \xi \pi \right)
+a_{2}b_{1}\rho ^{\lambda }\sin \left( \left( \xi -\lambda \right) \pi
\right) +a_{3}b_{1}\rho ^{\kappa }\sin \left( \left( \xi -\kappa \right) \pi
\right) \\
& \quad +a_{1}b_{2}\sin \left( \left( \xi +\alpha +\beta \right) \pi \right)
+a_{2}b_{2}\rho ^{\lambda }\sin \left( \left( \xi -\lambda +\alpha +\beta
\right) \pi \right) +a_{3}b_{2}\rho ^{\kappa }\sin \left( \left( \xi -\kappa
+\alpha +\beta \right) \pi \right)
\end{align*}%
to be non-negative, it is ensured that function $\dot{\varepsilon}_{cr},$
see (\ref{ML-IR}), is completely monotonic and therefore the creep
compliance is a Bernstein function. Note, the creep compliance in the form (%
\ref{epsilon-NP-1}) coincides with the integral of (\ref{epsilon-NP-2}),
since $\varepsilon _{cr}\left( 0\right) =0$.

\subsection{Calculation of relaxation modulus and creep compliance \label%
{CRM copy(1)}}

\subsubsection{Relaxation modulus calculation \label{CRM}}

The relaxation modulus is obtained in the form (\ref{sr-opste}), containing
functions given by (\ref{sigma-NP}), (\ref{sigma-RP}), and (\ref{sigma-CCP}%
), by inverting the relaxation modulus in the Laplace domain (\ref%
{sigma-sr-ld}) according to the definition of inverse Laplace transform%
\begin{equation*}
\sigma _{sr}\left( t\right) =\mathcal{L}^{-1}\left[ \tilde{\sigma}%
_{sr}\left( s\right) \right] \left( t\right) =\frac{1}{2\pi \mathrm{i}}%
\int_{p_{0}-\mathrm{i}\infty }^{p_{0}+\mathrm{i}\infty }\tilde{\sigma}%
_{sr}\left( s\right) \mathrm{e}^{st}\mathrm{d}s.  \label{inv-laplas}
\end{equation*}%
More precisely, the relaxation modulus $\sigma _{sr}$ in the forms $\sigma
_{sr}=\sigma _{sr}^{\left( \mathrm{NP}\right) }$ and $\sigma _{sr}=\sigma
_{sr}^{\left( \mathrm{NP}\right) }+\sigma _{sr}^{\left( \mathrm{RP}\right) }$
is obtained by the Cauchy integral theorem%
\begin{equation}
\oint_{\Gamma ^{(\mathrm{I,II})}}\tilde{\sigma}_{sr}\left( s\right) \mathrm{e%
}^{st}\mathrm{d}s=0,  \label{Kosijeva-teorema-sigma}
\end{equation}%
since function $\tilde{\sigma}_{sr}$ either does not have poles and then
integration is performed along the contour $\Gamma ^{(\mathrm{I})}$,
depicted in Figure \ref{nemaTG}, or has a negative real pole lying outside
of the contour $\Gamma ^{(\mathrm{II})}$, depicted in Figure \ref%
{negativnaTG}, while the relaxation modulus $\sigma _{sr}$ in the form $%
\sigma _{sr}=\sigma _{sr}^{\left( \mathrm{NP}\right) }+\sigma _{sr}^{\left( 
\mathrm{CCP}\right) }$ is obtained by the Cauchy residue theorem%
\begin{equation}
\oint_{\Gamma ^{(\mathrm{I})}}\tilde{\sigma}_{sr}\left( s\right) \mathrm{e}%
^{st}\mathrm{d}s=2\pi \mathrm{i}\left( \func{Res}\left( \tilde{\sigma}%
_{sr}\left( s\right) \mathrm{e}^{st},s_{\scriptscriptstyle{\mathrm{CCP}}%
}\right) +\func{Res}\left( \tilde{\sigma}_{sr}\left( s\right) \mathrm{e}%
^{st},\bar{s}_{\scriptscriptstyle{\mathrm{CCP}}}\right) \right) ,
\label{Kosijeva-teorema-rezduumi-sigma}
\end{equation}%
since function $\tilde{\sigma}_{sr}$ has $s_{\scriptscriptstyle{\mathrm{CCP}}%
}$ and its complex conjugate $\bar{s}_{\scriptscriptstyle{\mathrm{CCP}}}$ as
poles lying within the contour $\Gamma ^{(\mathrm{I})},$ depicted in Figure %
\ref{nemaTG}.

\noindent 
\begin{minipage}{\columnwidth}
\begin{minipage}[c]{0.4\columnwidth}
\centering
\includegraphics[width=0.7\columnwidth]{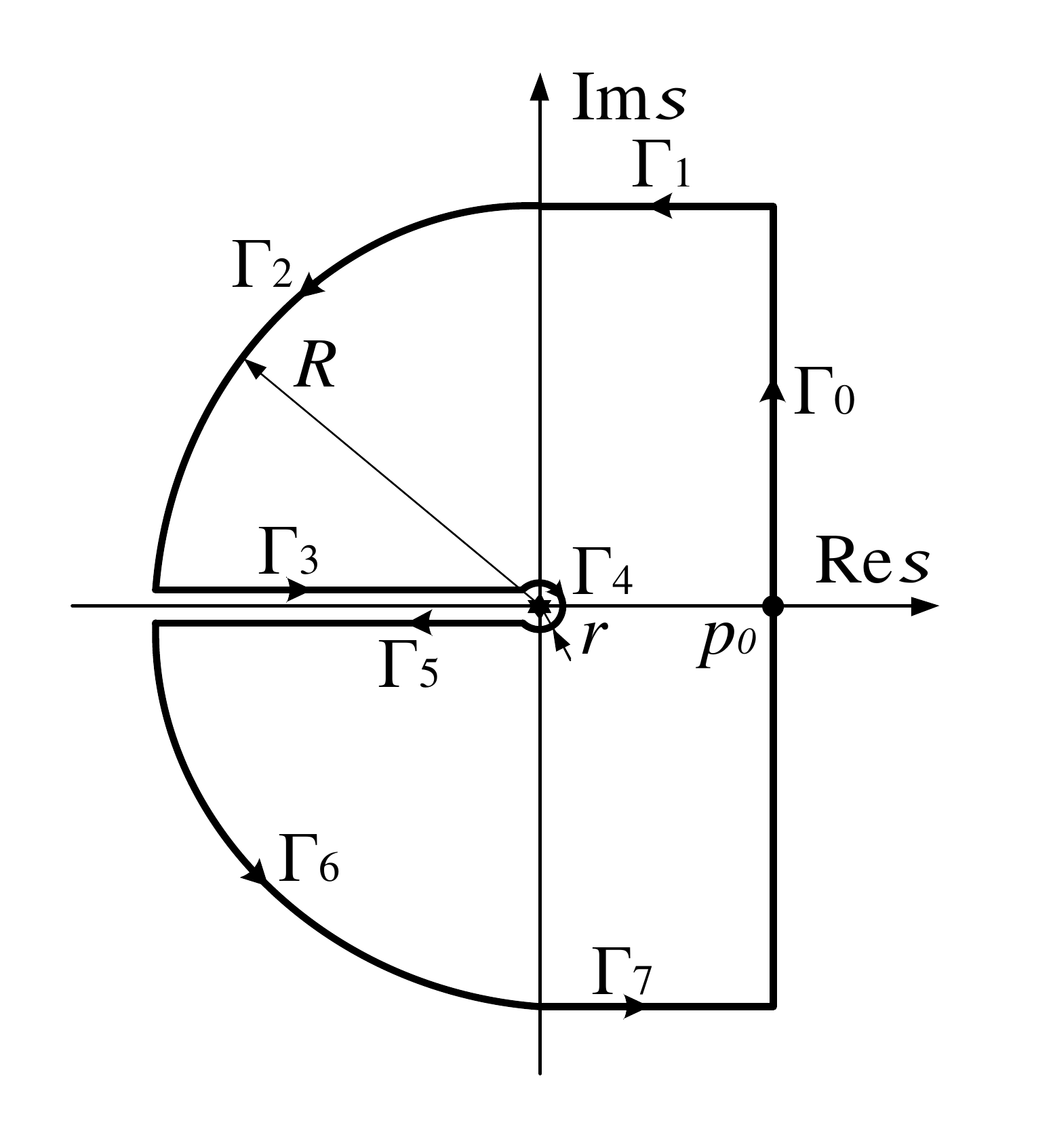}
\captionof{figure}{Integration contour $\Gamma^{(\mathrm{I})}$.}
\label{nemaTG}
\end{minipage}
\hfil
\begin{minipage}[c]{0.55\columnwidth}
\centering
\begin{tabular}{rll}
$\Gamma _{0}:$ & Bromwich path, &  \\ 
$\Gamma _{1}:$ & $s=p+\mathrm{i}R,$ & $p\in \left[ 0,p_{0}\right],\, p_0\geq 0$ arbitrary, \\ 
$\Gamma _{2}:$ & $s=R\mathrm{e}^{\mathrm{i}\varphi },$ & $\varphi \in \left[ 
\frac{\pi }{2},\pi \right] ,$ \\ 
$\Gamma _{3}:$ & $s=\rho \mathrm{e}^{\mathrm{i}\pi },$ & $\rho \in \left[ r,R%
\right] ,$ \\ 
$\Gamma _{4}:$ & $s=r\mathrm{e}^{\mathrm{i}\varphi },$ & $\varphi \in \left[ -\pi
,\pi \right] ,$ \\ 
$\Gamma _{5}:$ & $s=\rho \mathrm{e}^{-\mathrm{i}\pi },$ & $\rho \in \left[ r,R%
\right] ,$ \\
$\Gamma _{6}:$  & $s=R\mathrm{e}^{\mathrm{i}\varphi },$ & $\varphi \in \left[ 
-\pi, -\frac{\pi }{2} \right] ,$ \\
$\Gamma _{7}:$ & $s=p-\mathrm{i}R,$ & $p\in \left[ 0,p_{0}\right],\, p_0\geq 0$ arbitrary.  
\end{tabular}
\captionof{table}{Parametrization of integration contour $\Gamma^{(\mathrm{I})}$.}
\label{nemaTG-param}
\end{minipage}
\end{minipage}\smallskip

\noindent 
\begin{minipage}{\columnwidth}
\begin{minipage}[c]{0.4\columnwidth}
\centering
\includegraphics[width=0.7\columnwidth]{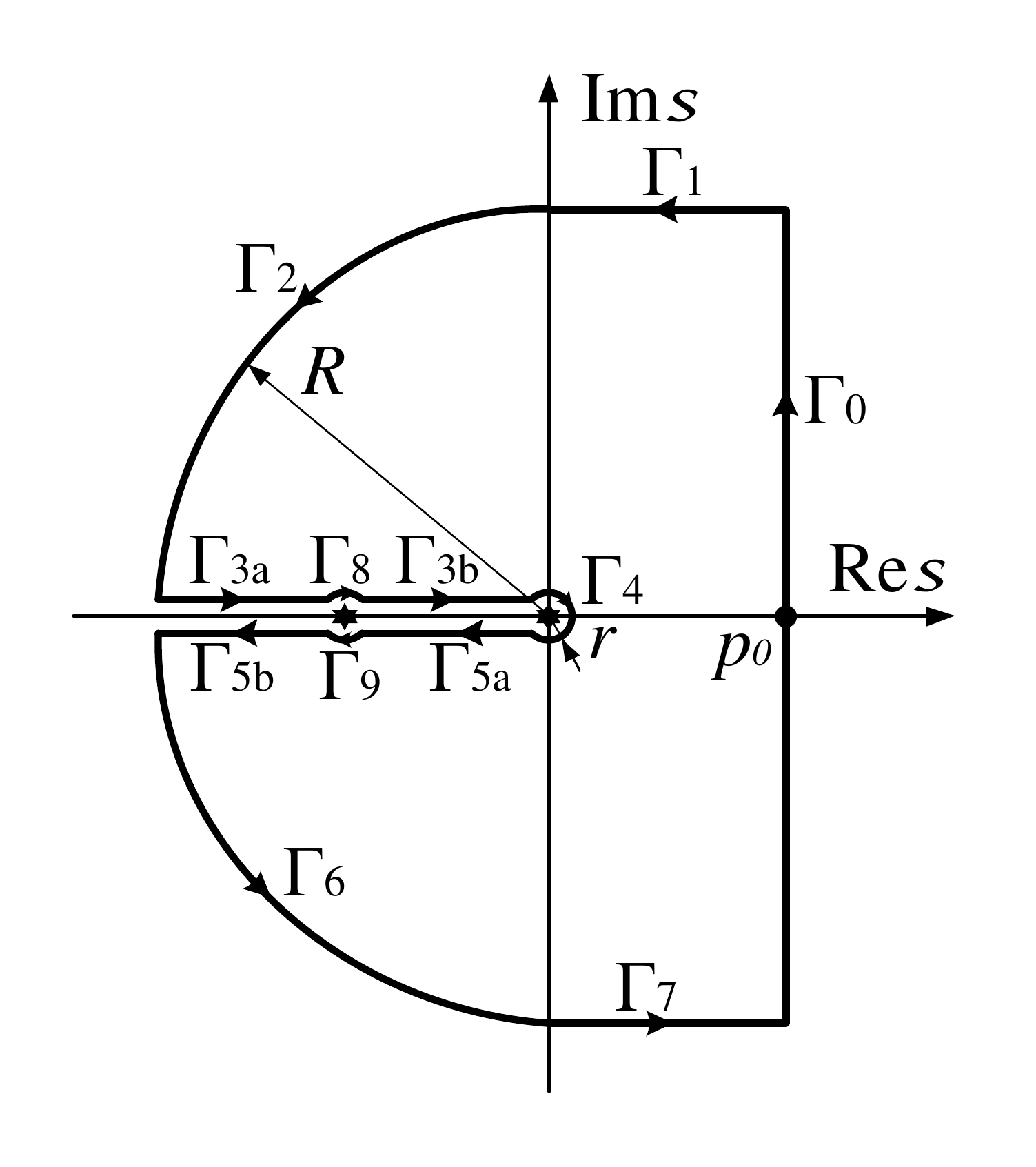}
\captionof{figure}{Integration contour $\Gamma^{(\mathrm{II})}$.}
\label{negativnaTG}
\end{minipage}
\hfil
\begin{minipage}[c]{0.55\columnwidth}
\centering
\begin{tabular}{rll}
$\Gamma _{0}:$ & Bromwich path, &  \\ 
$\Gamma _{1}:$ & $s=p+\mathrm{i}R,$ & $p\in \left[ 0,p_{0}\right],\, p_0\geq 0$ arbitrary, \\ 
$\Gamma _{2}:$ & $s=R\mathrm{e}^{\mathrm{i}\varphi },$ & $\varphi \in \left[ 
\frac{\pi }{2},\pi \right] ,$ \\ 
$\Gamma _{3a}\cup\Gamma _{3b}:$ & $s=\rho \mathrm{e}^{\mathrm{i}\pi },$ & $\rho \in \left[ r,R%
\right] ,$ \\ 
$\Gamma _{4}:$ & $s=r\mathrm{e}^{\mathrm{i}\varphi },$ & $\varphi \in \left[ -\pi
,\pi \right] ,$ \\ 
$\Gamma _{5a}\cup\Gamma _{5b}:$ & $s=\rho \mathrm{e}^{-\mathrm{i}\pi },$ & $\rho \in \left[ r,R%
\right] ,$ \\
$\Gamma _{6}:$  & $s=R\mathrm{e}^{\mathrm{i}\varphi },$ & $\varphi \in \left[ 
-\pi, -\frac{\pi }{2} \right] ,$ \\
$\Gamma _{7}:$ & $s=p-\mathrm{i}R,$ & $p\in \left[ 0,p_{0}\right],\, p_0\geq 0$ arbitrary,\\
$\Gamma _{8}:$  & $s=-\rho^*+r\mathrm{e}^{\mathrm{i}\varphi },$ & $\varphi \in \left[0,
\pi\right] ,$ \\
$\Gamma _{9}:$  & $s=-\rho^*+r\mathrm{e}^{\mathrm{i}\varphi },$ & $\varphi \in \left[ -\pi,
0 \right]$.  
\end{tabular}
\captionof{table}{Parametrization of integration contour $\Gamma^{(\mathrm{II})}$.}
\label{negativnaTG-param}
\end{minipage}
\end{minipage}\smallskip

In the case when relaxation modulus in the Laplace domain (\ref{sigma-sr-ld}%
) does not have poles, the Cauchy integral theorem (\ref%
{Kosijeva-teorema-sigma}), with the contour $\Gamma ^{(\mathrm{I})}$
depicted in Figure \ref{nemaTG} and by taking into account integrals having
non-zero contributions, yields%
\begin{equation}
\int_{\Gamma _{0}}\tilde{\sigma}_{sr}\left( s\right) \mathrm{e}^{st}\mathrm{d%
}s+\int_{\Gamma _{3}}\tilde{\sigma}_{sr}\left( s\right) \mathrm{e}^{st}%
\mathrm{d}s+\int_{\Gamma _{5}}\tilde{\sigma}_{sr}\left( s\right) \mathrm{e}%
^{st}\mathrm{d}s=0,  \label{equation}
\end{equation}%
that, according to the parameterization of contours $\Gamma _{3}$ and $%
\Gamma _{5},$ given in Table \ref{nemaTG-param}, in the limit when $%
r\rightarrow 0$ and $R\rightarrow \infty $ becomes%
\begin{equation}
2\pi \mathrm{i\,}\sigma _{sr}\left( t\right) +\int_{\infty }^{0}\frac{1}{%
\rho ^{1-\xi }\mathrm{e}^{\mathrm{i}\left( 1-\xi \right) \pi }}\frac{\phi
_{\varepsilon }\left( \rho \mathrm{e}^{\mathrm{i}\pi }\right) }{\phi
_{\sigma }\left( \rho \mathrm{e}^{\mathrm{i}\pi }\right) }\mathrm{e}^{\rho t%
\mathrm{e}^{^{\mathrm{i}\pi }}}\mathrm{e}^{\mathrm{i}\pi }\mathrm{d}\rho
+\int_{0}^{\infty }\frac{1}{\rho ^{1-\xi }\mathrm{e}^{-\mathrm{i}\left(
1-\xi \right) \pi }}\frac{\phi _{\varepsilon }\left( \rho \mathrm{e}^{-%
\mathrm{i}\pi }\right) }{\phi _{\sigma }\left( \rho \mathrm{e}^{-\mathrm{i}%
\pi }\right) }\mathrm{e}^{\rho t\mathrm{e}^{-^{\mathrm{i}\pi }}}\mathrm{e}^{-%
\mathrm{i}\pi }\mathrm{d}\rho =0,  \label{KTS-NP}
\end{equation}%
transforming into%
\begin{align*}
\sigma _{sr}\left( t\right) & =\sigma _{sr}^{\left( \mathrm{NP}\right)
}\left( t\right) \\
& =\frac{1}{2\pi \mathrm{i}}\int_{0}^{\infty }\frac{1}{\rho ^{1-\xi }}\frac{%
\mathrm{e}^{\mathrm{i}\left( 1-\xi \right) \pi }\bar{\phi}_{\varepsilon
}\left( \rho \mathrm{e}^{\mathrm{i}\pi }\right) \phi _{\sigma }\left( \rho 
\mathrm{e}^{\mathrm{i}\pi }\right) -\mathrm{e}^{-\mathrm{i}\left( 1-\xi
\right) \pi }\phi _{\varepsilon }\left( \rho \mathrm{e}^{\mathrm{i}\pi
}\right) \bar{\phi}_{\sigma }\left( \rho \mathrm{e}^{\mathrm{i}\pi }\right) 
}{\left\vert \phi _{\sigma }\left( \rho \mathrm{e}^{\mathrm{i}\pi }\right)
\right\vert ^{2}}\mathrm{e}^{-\rho t}\mathrm{d}\rho ,\quad \text{i.e.,} \\
\sigma _{sr}\left( t\right) & =\sigma _{sr}^{\left( \mathrm{NP}\right)
}\left( t\right) \\
& =\frac{1}{2\pi \mathrm{i}}\int_{0}^{\infty }\frac{1}{\rho ^{1-\xi }}\frac{%
\left\vert \phi _{\varepsilon }\left( \rho \mathrm{e}^{\mathrm{i}\pi
}\right) \right\vert }{\left\vert \phi _{\sigma }\left( \rho \mathrm{e}^{%
\mathrm{i}\pi }\right) \right\vert }\left( \mathrm{e}^{\mathrm{i}\left(
1-\xi \right) \pi -\arg \phi _{\varepsilon }\left( \rho \mathrm{e}^{\mathrm{i%
}\pi }\right) +\arg \phi _{\sigma }\left( \rho \mathrm{e}^{\mathrm{i}\pi
}\right) }-\mathrm{e}^{-\mathrm{i}\left( 1-\xi \right) \pi +\arg \phi
_{\varepsilon }\left( \rho \mathrm{e}^{\mathrm{i}\pi }\right) -\arg \phi
_{\sigma }\left( \rho \mathrm{e}^{\mathrm{i}\pi }\right) }\right) \mathrm{e}%
^{-\rho t}\mathrm{d}\rho ,
\end{align*}%
and becoming of the form (\ref{sigma-NP}) or (\ref{sigma-NP-1}) having the
function $K,$ given by (\ref{K}), taken into account.

In addition to the integrals along contours $\Gamma _{3a}\cup \Gamma _{3b}$
and $\Gamma _{5a}\cup \Gamma _{5b},$ that are parts of contour $\Gamma ^{(%
\mathrm{II})}$ from Figure \ref{negativnaTG}, there are additional integrals
along contours $\Gamma _{8}$ and $\Gamma _{9}$ having non-zero contributions
in the case when relaxation modulus in the Laplace domain (\ref{sigma-sr-ld}%
) has a negative real pole $s_{\scriptscriptstyle{\mathrm{RP}}}=\rho _{%
\scriptscriptstyle{\mathrm{RP}}}\,\mathrm{e}^{\mathrm{i}\pi },$ so that the
Cauchy integral theorem (\ref{Kosijeva-teorema-sigma}) yields%
\begin{equation*}
\int_{\Gamma _{0}}\tilde{\sigma}_{sr}\left( s\right) \mathrm{e}^{st}\mathrm{d%
}s+\int_{\Gamma _{3a}\cup \Gamma _{3b}}\tilde{\sigma}_{sr}\left( s\right) 
\mathrm{e}^{st}\mathrm{d}s+\int_{\Gamma _{5a}\cup \Gamma _{5b}}\tilde{\sigma}%
_{sr}\left( s\right) \mathrm{e}^{st}\mathrm{d}s+\int_{\Gamma _{8}}\tilde{%
\sigma}_{sr}\left( s\right) \mathrm{e}^{st}\mathrm{d}s+\int_{\Gamma _{9}}%
\tilde{\sigma}_{sr}\left( s\right) \mathrm{e}^{st}\mathrm{d}s=0,
\end{equation*}%
with the first three terms being already defined by (\ref{KTS-NP}) and with
the remaining terms, containing integrals along contours $\Gamma _{8}$ and $%
\Gamma _{9}$ parameterized as in Table \ref{negativnaTG-param}, transforming
the previous expression into%
\begin{align*}
2\pi \mathrm{i\,}\sigma _{sr}\left( t\right) +2\pi \mathrm{i\,}\sigma
_{sr}^{\left( \mathrm{NP}\right) }\left( t\right) +& \int_{\pi }^{0}\frac{1}{%
\left( s_{\scriptscriptstyle{\mathrm{RP}}}+r\mathrm{e}^{\mathrm{i}\varphi
}\right) ^{1-\xi }}\frac{\phi _{\varepsilon }\left( s_{\scriptscriptstyle{%
\mathrm{RP}}}+r\mathrm{e}^{\mathrm{i}\varphi }\right) }{\phi _{\sigma
}\left( s_{\scriptscriptstyle{\mathrm{RP}}}+r\mathrm{e}^{\mathrm{i}\varphi
}\right) }\mathrm{e}^{\left( s_{\scriptscriptstyle{\mathrm{RP}}}+r\mathrm{e}%
^{\mathrm{i}\varphi }\right) t}\mathrm{i}r\mathrm{e}^{\mathrm{i}\varphi }%
\mathrm{d}\varphi \\
+& \int_{0}^{-\pi }\frac{1}{\left( \bar{s}_{\scriptscriptstyle{\mathrm{RP}}%
}+r\mathrm{e}^{\mathrm{i}\varphi }\right) ^{1-\xi }}\frac{\phi _{\varepsilon
}\left( \bar{s}_{\scriptscriptstyle{\mathrm{RP}}}+r\mathrm{e}^{\mathrm{i}%
\varphi }\right) }{\phi _{\sigma }\left( \bar{s}_{\scriptscriptstyle{\mathrm{%
RP}}}+r\mathrm{e}^{\mathrm{i}\varphi }\right) }\mathrm{e}^{\left( \bar{s}_{%
\scriptscriptstyle{\mathrm{RP}}}+r\mathrm{e}^{\mathrm{i}\varphi }\right) t}%
\mathrm{i}r\mathrm{e}^{\mathrm{i}\varphi }\mathrm{d}\varphi =0,
\end{align*}%
in the limit when $r\rightarrow 0$ and $R\rightarrow \infty $, so that%
\begin{equation*}
\sigma _{sr}\left( t\right) =\sigma _{sr}^{\left( \mathrm{NP}\right) }\left(
t\right) +\sigma _{sr}^{\left( \mathrm{RP}\right) }\left( t\right) ,
\end{equation*}%
where 
\begin{align*}
\sigma _{sr}^{\left( \mathrm{RP}\right) }\left( t\right) & =\frac{1}{2\pi 
\mathrm{i}}\left( -\mathrm{i}\pi \frac{1}{s_{\scriptscriptstyle{\mathrm{RP}}%
}^{1-\xi }}\frac{\phi _{\varepsilon }\left( s_{\scriptscriptstyle{\mathrm{RP}%
}}\right) }{\phi _{\sigma }^{\prime }\left( s_{\scriptscriptstyle{\mathrm{RP}%
}}\right) }\mathrm{e}^{s_{\scriptscriptstyle{\mathrm{RP}}}t}-\mathrm{i}\pi 
\frac{1}{\bar{s}_{\scriptscriptstyle{\mathrm{RP}}}^{1-\xi }}\frac{\phi
_{\varepsilon }\left( \bar{s}_{\scriptscriptstyle{\mathrm{RP}}}\right) }{%
\phi _{\sigma }^{\prime }\left( \bar{s}_{\scriptscriptstyle{\mathrm{RP}}%
}\right) }\mathrm{e}^{\bar{s}_{\scriptscriptstyle{\mathrm{RP}}}t}\right) \\
& =-\frac{1}{2}\frac{1}{\rho _{\scriptscriptstyle{\mathrm{RP}}}^{1-\xi }}%
\frac{\mathrm{e}^{-\mathrm{i}\left( 1-\xi \right) \pi }\phi _{\varepsilon
}\left( s_{\scriptscriptstyle{\mathrm{RP}}}\right) \bar{\phi}_{\sigma
}^{\prime }\left( s_{\scriptscriptstyle{\mathrm{RP}}}\right) +\mathrm{e}^{%
\mathrm{i}\left( 1-\xi \right) \pi }\bar{\phi}_{\varepsilon }\left( s_{%
\scriptscriptstyle{\mathrm{RP}}}\right) \phi _{\sigma }^{\prime }\left( s_{%
\scriptscriptstyle{\mathrm{RP}}}\right) }{\left\vert \phi _{\sigma }^{\prime
}\left( s_{\scriptscriptstyle{\mathrm{RP}}}\right) \right\vert ^{2}}\mathrm{e%
}^{-\rho _{\scriptscriptstyle{\mathrm{RP}}}t}
\end{align*}%
becoming of the form given by (\ref{sigma-RP}). The integral along contour $%
\Gamma _{8}$ is calculated as 
\begin{align*}
\lim_{r\rightarrow 0}\int_{\Gamma _{8}}\tilde{\sigma}_{sr}\left( s\right) 
\mathrm{e}^{st}\mathrm{d}s& =\lim_{r\rightarrow 0}\int_{\pi }^{0}\frac{1}{%
\left( s_{\scriptscriptstyle{\mathrm{RP}}}+r\mathrm{e}^{\mathrm{i}\varphi
}\right) ^{1-\xi }}\frac{\phi _{\varepsilon }\left( s_{\scriptscriptstyle{%
\mathrm{RP}}}+r\mathrm{e}^{\mathrm{i}\varphi }\right) }{\phi _{\sigma
}\left( s_{\scriptscriptstyle{\mathrm{RP}}}+r\mathrm{e}^{\mathrm{i}\varphi
}\right) }\mathrm{e}^{\left( s_{\scriptscriptstyle{\mathrm{RP}}}+r\mathrm{e}%
^{\mathrm{i}\varphi }\right) t}\mathrm{i}r\mathrm{e}^{\mathrm{i}\varphi }%
\mathrm{d}\varphi \\
& =\lim_{r\rightarrow 0}\int_{\pi }^{0}\frac{1}{\left( s_{\scriptscriptstyle{%
\mathrm{RP}}}+r\mathrm{e}^{\mathrm{i}\varphi }\right) ^{1-\xi }}\frac{\phi
_{\varepsilon }\left( s_{\scriptscriptstyle{\mathrm{RP}}}+r\mathrm{e}^{%
\mathrm{i}\varphi }\right) }{\phi _{\sigma }\left( s_{\scriptscriptstyle{%
\mathrm{RP}}}\right) +\left. \phi _{\sigma }^{\prime }\left( s\right) \left(
s-s_{\scriptscriptstyle{\mathrm{RP}}}\right) \right\vert _{s=s_{%
\scriptscriptstyle{\mathrm{RP}}}+r\mathrm{e}^{\mathrm{i}\varphi }}+\ldots }%
\mathrm{e}^{\left( s_{\scriptscriptstyle{\mathrm{RP}}}+r\mathrm{e}^{\mathrm{i%
}\varphi }\right) t}\mathrm{i}r\mathrm{e}^{\mathrm{i}\varphi }\mathrm{d}%
\varphi \\
& =-\mathrm{i}\pi \frac{1}{s_{\scriptscriptstyle{\mathrm{RP}}}^{1-\xi }}%
\frac{\phi _{\varepsilon }\left( s_{\scriptscriptstyle{\mathrm{RP}}}\right) 
}{\phi _{\sigma }^{\prime }\left( s_{\scriptscriptstyle{\mathrm{RP}}}\right) 
}\mathrm{e}^{s_{\scriptscriptstyle{\mathrm{RP}}}t},
\end{align*}%
by expanding function $\phi _{\sigma }$ into the series and by taking into
account $\phi _{\sigma }\left( s_{\scriptscriptstyle{\mathrm{RP}}}\right) =0$%
, while the similar calculation yields the integral along contour $\Gamma
_{9}$ in the form%
\begin{equation*}
\lim_{r\rightarrow 0}\int_{\Gamma _{9}}\tilde{\sigma}_{sr}\left( s\right) 
\mathrm{e}^{st}\mathrm{d}s=-\mathrm{i}\pi \frac{1}{\bar{s}_{%
\scriptscriptstyle{\mathrm{RP}}}^{1-\xi }}\frac{\phi _{\varepsilon }\left( 
\bar{s}_{\scriptscriptstyle{\mathrm{RP}}}\right) }{\phi _{\sigma }^{\prime
}\left( \bar{s}_{\scriptscriptstyle{\mathrm{RP}}}\right) }\mathrm{e}^{\bar{s}%
_{\scriptscriptstyle{\mathrm{RP}}}t}.
\end{equation*}

On the right-hand side of equation (\ref{equation}), the term originating
from the residues appears when the function $\tilde{\sigma}_{sr},$ given by (%
\ref{sigma-sr-ld}), has a pair of complex conjugated poles $s_{%
\scriptscriptstyle{\mathrm{RP}}}$ and $\bar{s}_{\scriptscriptstyle{\mathrm{RP%
}}}$, since the Cauchy residues theorem (\ref%
{Kosijeva-teorema-rezduumi-sigma}), with the contour $\Gamma ^{(\mathrm{I})}$
depicted in Figure \ref{nemaTG}, where the poles are located within the area
bounded by the contour, by taking into account integrals having non-zero
contributions, yields%
\begin{equation*}
\int_{\Gamma _{0}}\tilde{\sigma}_{sr}\left( s\right) \mathrm{e}^{st}\mathrm{d%
}s+\int_{\Gamma _{3}}\tilde{\sigma}_{sr}\left( s\right) \mathrm{e}^{st}%
\mathrm{d}s+\int_{\Gamma _{5}}\tilde{\sigma}_{sr}\left( s\right) \mathrm{e}%
^{st}\mathrm{d}s=2\pi \mathrm{i}\left( \func{Res}\left( \tilde{\sigma}%
_{sr}\left( s\right) \mathrm{e}^{st},s_{\scriptscriptstyle{\mathrm{CCP}}%
}\right) +\func{Res}\left( \tilde{\sigma}_{sr}\left( s\right) \mathrm{e}%
^{st},\bar{s}_{\scriptscriptstyle{\mathrm{CCP}}}\right) \right) ,
\end{equation*}%
that becomes%
\begin{equation*}
\sigma _{sr}\left( t\right) -\sigma _{sr}^{\left( \mathrm{NP}\right) }\left(
t\right) =\sigma _{sr}^{\left( \mathrm{CCP}\right) }\left( t\right) ,
\end{equation*}%
where 
\begin{align*}
\sigma _{sr}^{\left( \mathrm{CCP}\right) }\left( t\right) & =\frac{1}{\rho _{%
\scriptscriptstyle{\mathrm{CCP}}}^{1-\xi }\mathrm{e}^{\mathrm{i}\left( 1-\xi
\right) \varphi _{\scriptscriptstyle{\mathrm{CCP}}}}}\frac{\phi
_{\varepsilon }\left( s_{\scriptscriptstyle{\mathrm{CCP}}}\right) }{\phi
_{\sigma }^{\prime }\left( s_{\scriptscriptstyle{\mathrm{CCP}}}\right) }%
\mathrm{e}^{\func{Re}s_{\scriptscriptstyle{\mathrm{CCP}}}t+\mathrm{i}\func{Im%
}s_{\scriptscriptstyle{\mathrm{CCP}}}t} \\
& \quad \quad +\frac{1}{\rho _{\scriptscriptstyle{\mathrm{CCP}}}^{1-\xi }%
\mathrm{e}^{-\mathrm{i}\left( 1-\xi \right) \varphi _{\scriptscriptstyle{%
\mathrm{CCP}}}}}\frac{\phi _{\varepsilon }\left( \bar{s}_{\scriptscriptstyle{%
\mathrm{CCP}}}\right) }{\phi _{\sigma }^{\prime }\left( \bar{s}_{%
\scriptscriptstyle{\mathrm{CCP}}}\right) }\mathrm{e}^{\func{Re}\bar{s}_{%
\scriptscriptstyle{\mathrm{CCP}}}t+\mathrm{i}\func{Im}\bar{s}_{%
\scriptscriptstyle{\mathrm{CCP}}}t} \\
& =\frac{1}{\rho _{\scriptscriptstyle{\mathrm{CCP}}}^{1-\xi }}\mathrm{e}^{%
\func{Re}s_{\scriptscriptstyle{\mathrm{CCP}}}t}\left( \frac{\phi
_{\varepsilon }\left( s_{\scriptscriptstyle{\mathrm{CCP}}}\right) }{\phi
_{\sigma }^{\prime }\left( s_{\scriptscriptstyle{\mathrm{CCP}}}\right) }%
\frac{\mathrm{e}^{\mathrm{i}\func{Im}s_{\scriptscriptstyle{\mathrm{CCP}}}t}}{%
\mathrm{e}^{\mathrm{i}\left( 1-\xi \right) \varphi _{\scriptscriptstyle{%
\mathrm{CCP}}}}}+\frac{\phi _{\varepsilon }\left( \bar{s}_{\scriptscriptstyle%
{\mathrm{CCP}}}\right) }{\phi _{\sigma }^{\prime }\left( \bar{s}_{%
\scriptscriptstyle{\mathrm{CCP}}}\right) }\frac{\mathrm{e}^{\mathrm{i}\func{%
Im}\bar{s}_{\scriptscriptstyle{\mathrm{CCP}}}t}}{\mathrm{e}^{-\mathrm{i}%
\left( 1-\xi \right) \varphi _{\scriptscriptstyle{\mathrm{CCP}}}}}\right) \\
& =\frac{1}{\rho _{\scriptscriptstyle{\mathrm{CCP}}}^{1-\xi }}\mathrm{e}%
^{-\left\vert \func{Re}s_{\scriptscriptstyle{\mathrm{CCP}}}\right\vert t} \\
& \quad \quad \times \frac{\mathrm{e}^{\mathrm{i}\left( \func{Im}s_{%
\scriptscriptstyle{\mathrm{CCP}}}t-\left( 1-\xi \right) \varphi _{%
\scriptscriptstyle{\mathrm{CCP}}}\right) }\phi _{\varepsilon }\left( s_{%
\scriptscriptstyle{\mathrm{CCP}}}\right) \bar{\phi}_{\sigma }^{\prime
}\left( s_{\scriptscriptstyle{\mathrm{CCP}}}\right) +\mathrm{e}^{-\mathrm{i}%
\left( \func{Im}s_{\scriptscriptstyle{\mathrm{CCP}}}t-\left( 1-\xi \right)
\varphi _{\scriptscriptstyle{\mathrm{CCP}}}\right) }\bar{\phi}_{\varepsilon
}\left( s_{\scriptscriptstyle{\mathrm{CCP}}}\right) \phi _{\sigma }^{\prime
}\left( s_{\scriptscriptstyle{\mathrm{CCP}}}\right) }{\left\vert \phi
_{\sigma }^{\prime }\left( s_{\scriptscriptstyle{\mathrm{CCP}}}\right)
\right\vert ^{2}},
\end{align*}%
becoming of the form (\ref{sigma-CCP}).

It is left to be shown that the integrals along contours $\Gamma _{1},$ $%
\Gamma _{2},$ $\Gamma _{4},$ $\Gamma _{6},$ and $\Gamma _{7},$ as parts of
contours $\Gamma ^{(\mathrm{I})}$ and $\Gamma ^{(\mathrm{II})}$ from Figures %
\ref{nemaTG} and \ref{negativnaTG}, have zero contributions in both Cauchy
integral and residues theorems, see (\ref{Kosijeva-teorema-sigma}) and (\ref%
{Kosijeva-teorema-rezduumi-sigma}), in the limit when $r\rightarrow 0$ and $%
R\rightarrow \infty .$

The integral along the contour $\Gamma _{1}$ reads%
\begin{equation*}
I_{\Gamma _{1}}=\int_{p_{0}}^{0}\frac{1}{\left( p+\mathrm{i}R\right) ^{1-\xi
}}\frac{\phi _{\varepsilon }\left( p+\mathrm{i}R\right) }{\phi _{\sigma
}\left( p+\mathrm{i}R\right) }\mathrm{e}^{\left( p+\mathrm{i}R\right) t}%
\mathrm{d}p,  \label{I-gamma-1}
\end{equation*}%
according to the contour parametrization given in Tables \ref{nemaTG-param}
and \ref{negativnaTG-param}, so that it has zero contribution, since its
absolute value yields%
\begin{align*}
\left\vert I_{\Gamma _{1}}\right\vert & \leqslant \int_{0}^{p_{0}}\frac{1}{%
\left\vert p+\mathrm{i}R\right\vert ^{1-\xi }}\frac{\left\vert \phi
_{\varepsilon }\left( p+\mathrm{i}R\right) \right\vert }{\left\vert \phi
_{\sigma }\left( p+\mathrm{i}R\right) \right\vert }\mathrm{e}^{pt}\mathrm{d}p
\\
& \leqslant \int_{0}^{p_{0}}\frac{1}{R^{1-\xi }\left\vert 1-\mathrm{i}\frac{p%
}{R}\right\vert ^{1-\xi }}\frac{\left\vert \phi _{\varepsilon }\left( p+%
\mathrm{i}R\right) \right\vert }{\left\vert \phi _{\sigma }\left( p+\mathrm{i%
}R\right) \right\vert }\mathrm{e}^{pt}\mathrm{d}p \\
& \leqslant \int_{0}^{p_{0}}\frac{1}{R^{1-\xi }}\frac{\left\vert \phi
_{\varepsilon }\left( p+\mathrm{i}R\right) \right\vert }{\left\vert \phi
_{\sigma }\left( p+\mathrm{i}R\right) \right\vert }\mathrm{e}^{pt}\mathrm{d}p
\\
& \leqslant \int_{0}^{p_{0}}\frac{1}{R^{1-\xi -\zeta _{R}}}\mathrm{e}^{pt}%
\mathrm{d}p\rightarrow 0,\quad \text{for}\quad R\rightarrow \infty ,
\end{align*}%
if the functions $\phi _{\varepsilon }$ and $\phi _{\sigma },$ that are
power-type functions, see Table \ref{skupina}, satisfy $\frac{\left\vert
\phi _{\varepsilon }\left( p+\mathrm{i}R\right) \right\vert }{\left\vert
\phi _{\sigma }\left( p+\mathrm{i}R\right) \right\vert }\sim R^{\zeta _{R}},$
where $\zeta _{R}<1-\xi $ when $R\rightarrow \infty $ for $p\in \left[
0,p_{0}\right] .$ Using the similar argumentation, one can prove that the
integral along contour $\Gamma _{7}$ has zero contribution as well.

The integral along contour $\Gamma _{2},$ parametrized as in Tables \ref%
{nemaTG-param} and \ref{negativnaTG-param}, takes the form%
\begin{equation*}
I_{\Gamma _{2}}=\int_{\frac{\pi }{2}}^{\pi }\frac{1}{R^{1-\xi }\mathrm{e}^{%
\mathrm{i}\left( 1-\xi \right) \varphi }}\frac{\phi _{\varepsilon }\left( R%
\mathrm{e}^{\mathrm{i}\varphi }\right) }{\phi _{\sigma }\left( R\mathrm{e}^{%
\mathrm{i}\varphi }\right) }\mathrm{e}^{Rt\mathrm{e}^{\mathrm{i}\varphi }}%
\mathrm{i}R\mathrm{e}^{\mathrm{i}\varphi }\mathrm{d}\varphi
\end{equation*}%
and has zero contribution, that can be proved by considering its absolute
value%
\begin{equation*}
\left\vert I_{\Gamma _{2}}\right\vert \leqslant \int_{\frac{\pi }{2}}^{\pi
}R^{\xi }\frac{\left\vert \phi _{\varepsilon }\left( R\mathrm{e}^{\mathrm{i}%
\varphi }\right) \right\vert }{\left\vert \phi _{\sigma }\left( R\mathrm{e}^{%
\mathrm{i}\varphi }\right) \right\vert }\mathrm{e}^{Rt\mathrm{\cos }\varphi }%
\mathrm{d}\varphi \rightarrow 0,\quad \text{for}\quad R\rightarrow \infty ,
\end{equation*}%
since the integrand term $\mathrm{e}^{Rt\mathrm{\cos }\varphi }\rightarrow 0$
for $\varphi \in \left[ \frac{\pi }{2},\pi \right] $ when $R\rightarrow
\infty $ and since functions $\phi _{\varepsilon }$ and $\phi _{\sigma }$
are of power-type, see Table \ref{skupina}. The integral along contour $%
\Gamma _{6}$ also has zero contribution, that is proved by the similar
argumentation.

If the power-type functions $\phi _{\varepsilon }$ and $\phi _{\sigma },$
see Table \ref{skupina}, are such that $\frac{\left\vert \phi _{\varepsilon
}\left( r\mathrm{e}^{\mathrm{i}\varphi }\right) \right\vert }{\left\vert
\phi _{\sigma }\left( r\mathrm{e}^{\mathrm{i}\varphi }\right) \right\vert }%
\sim \frac{1}{r^{\zeta _{r}}},$ where $\zeta _{r}<\xi $ when $r\rightarrow 0$
for $\varphi \in \left[ -\pi ,\pi \right] ,$ it can be proved that the
integral along contour $\Gamma _{4}$ has zero contribution, since, according
to the parametrization given in Tables \ref{nemaTG-param} and \ref%
{negativnaTG-param}, the integral is%
\begin{equation*}
I_{\Gamma _{4}}=\int_{\pi }^{-\pi }\frac{1}{r^{1-\xi }\mathrm{e}^{\mathrm{i}%
\left( 1-\xi \right) \varphi }}\frac{\phi _{\varepsilon }\left( r\mathrm{e}^{%
\mathrm{i}\varphi }\right) }{\phi _{\sigma }\left( r\mathrm{e}^{\mathrm{i}%
\varphi }\right) }\mathrm{e}^{rt\mathrm{e}^{\mathrm{i}\varphi }}\mathrm{i}r%
\mathrm{e}^{\mathrm{i}\varphi }\mathrm{d}\varphi ,
\end{equation*}%
so that%
\begin{align*}
\left\vert I_{\Gamma _{4}}\right\vert & \leqslant \int_{-\pi }^{\pi }r^{\xi }%
\frac{\left\vert \phi _{\varepsilon }\left( r\mathrm{e}^{\mathrm{i}\varphi
}\right) \right\vert }{\left\vert \phi _{\sigma }\left( r\mathrm{e}^{\mathrm{%
i}\varphi }\right) \right\vert }\mathrm{e}^{rt\mathrm{\cos }\varphi }\mathrm{%
d}\varphi \\
& \leqslant \int_{-\pi }^{\pi }r^{\xi -\zeta _{r}}\mathrm{d}\varphi
\rightarrow 0,\quad \text{for}\quad r\rightarrow 0.
\end{align*}

\subsubsection{Creep compliance calculation \label{CCC}}

The creep compliance in the form (\ref{cr-opste}), containing functions
given by (\ref{epsilon-NP}), (\ref{epsilon-RP}), and (\ref{epsilon-CCP}), is
obtained according to 
\begin{equation}
\varepsilon _{cr}\left( t\right) =\int_{0}^{t}\epsilon _{cr}\left( t^{\prime
}\right) \mathrm{d}t^{\prime },  \label{epsilon-creep-integral}
\end{equation}%
that is a consequence of considering the function%
\begin{equation}
\tilde{\epsilon}_{cr}\left( s\right) =s\tilde{\varepsilon}_{cr}\left(
s\right) =\frac{1}{s^{\xi }}\frac{\phi _{\sigma }\left( s\right) }{\phi
_{\varepsilon }\left( s\right) },  \label{s-epsilon}
\end{equation}%
having the inverse Laplace transform obtained as%
\begin{equation*}
\epsilon _{cr}\left( t\right) =\frac{\mathrm{d}}{\mathrm{d}t}\varepsilon
_{cr}\left( t\right) +\varepsilon _{cr}^{\left( g\right) }\delta \left(
t\right) =\frac{\mathrm{d}}{\mathrm{d}t}\varepsilon _{cr}\left( t\right) ,
\label{izvod-epsilon-krip}
\end{equation*}%
with the glass compliance $\varepsilon _{cr}^{\left( g\right)
}=\lim_{t\rightarrow 0}\varepsilon _{cr}\left( t\right) =\lim_{s\rightarrow
\infty }s\tilde{\varepsilon}_{cr}\left( s\right) =\lim_{s\rightarrow \infty }%
\frac{1}{s^{\xi }}\frac{\phi _{\sigma }\left( s\right) }{\phi _{\varepsilon
}\left( s\right) }=0$ for all considered models, rather than the creep
compliance in Laplace domain (\ref{epsilon-cr-ld}). The function $\epsilon
_{cr}$ is calculated using the definition of the inverse Laplace transform
applied to function $\tilde{\epsilon}_{cr},$ see (\ref{s-epsilon}), and
integration in the complex plane along the Bromwich contour, i.e., by%
\begin{equation}
\epsilon _{cr}\left( t\right) =\frac{1}{2\pi \mathrm{i}}\int_{p_{0}-\mathrm{i%
}\infty }^{p_{0}+\mathrm{i}\infty }\tilde{\epsilon}_{cr}\left( s\right) 
\mathrm{e}^{st}\mathrm{d}s.  \label{inv-laplas-epsilon}
\end{equation}

More precisely, the creep compliance $\varepsilon _{cr}$ in the forms $%
\varepsilon _{cr}=\varepsilon _{cr}^{\left( \mathrm{NP}\right) }$ and $%
\varepsilon _{cr}=\varepsilon _{cr}^{\left( \mathrm{NP}\right) }+\varepsilon
_{cr}^{\left( \mathrm{RP}\right) }$ is obtained according to (\ref%
{epsilon-creep-integral}), with the function $\epsilon _{cr}$ calculated by
the Cauchy integral theorem 
\begin{equation}
\oint_{\Gamma ^{(\mathrm{I,II})}}\tilde{\epsilon}_{cr}\left( s\right) 
\mathrm{e}^{st}\mathrm{d}s=0,  \label{Kosijeva-teorema-epsilon}
\end{equation}%
where the integration is performed either along the contour $\Gamma ^{(%
\mathrm{I})}$, depicted in Figure \ref{nemaTG}, if the function $\tilde{%
\epsilon}_{cr}$ does not have poles, or along the contour $\Gamma ^{(\mathrm{%
II})}$, depicted in Figure \ref{negativnaTG}, if the function $\tilde{%
\epsilon}_{cr}$ has a negative real pole that lies outside of contour $%
\Gamma ^{(\mathrm{II})},$ while the creep compliance $\varepsilon _{cr}$ in
the form $\varepsilon _{cr}=\varepsilon _{cr}^{\left( \mathrm{NP}\right)
}+\varepsilon _{cr}^{\left( \mathrm{CCP}\right) }$ is obtained according to (%
\ref{epsilon-creep-integral}), with the function $\epsilon _{cr}$ calculated
by the Cauchy residue theorem 
\begin{equation}
\oint_{\Gamma ^{(\mathrm{I})}}\tilde{\epsilon}_{cr}\left( s\right) \mathrm{e}%
^{st}\mathrm{d}s=2\pi \mathrm{i}\left( \func{Res}\left( \tilde{\epsilon}%
_{cr}\left( s\right) \mathrm{e}^{st},s_{\scriptscriptstyle{\mathrm{CCP}}%
}\right) +\func{Res}\left( \tilde{\epsilon}_{cr}\left( s\right) \mathrm{e}%
^{st},\bar{s}_{\scriptscriptstyle{\mathrm{CCP}}}\right) \right) ,
\label{Kosijeva-teorema-rezduumi-epsilon}
\end{equation}%
since function $\tilde{\epsilon}_{cr}$ has $s_{\scriptscriptstyle{\mathrm{CCP%
}}}$ and its complex conjugate $\bar{s}_{\scriptscriptstyle{\mathrm{CCP}}}$
as poles lying within the contour $\Gamma ^{(\mathrm{I})},$ depicted in
Figure \ref{nemaTG}.

Namely, in the case when function $\tilde{\epsilon}_{cr},$ given by (\ref%
{s-epsilon}), does not have poles, the Cauchy integral theorem (\ref%
{Kosijeva-teorema-epsilon}), with the contour $\Gamma ^{(\mathrm{I})}$
depicted in Figure \ref{nemaTG}, taking into account integrals that have
non-zero contributions, yields%
\begin{equation*}
\int_{\Gamma _{0}}\tilde{\epsilon}_{cr}\left( s\right) \mathrm{e}^{st}%
\mathrm{d}s+\int_{\Gamma _{3}}\tilde{\epsilon}_{cr}\left( s\right) \mathrm{e}%
^{st}\mathrm{d}s+\int_{\Gamma _{5}}\tilde{\epsilon}_{cr}\left( s\right) 
\mathrm{e}^{st}\mathrm{d}s=0,
\end{equation*}%
where the integrals along contours $\Gamma _{0},$ $\Gamma _{3},$ and $\Gamma
_{5},$ parameterized as in Table \ref{nemaTG-param}, in the limit when $%
r\rightarrow 0$ and $R\rightarrow \infty $ become%
\begin{equation}
2\pi \mathrm{i\,}\mathcal{\epsilon }_{cr}\left( t\right) +\int_{\infty }^{0}%
\frac{1}{\rho ^{\xi }\mathrm{e}^{\mathrm{i}\xi \pi }}\frac{\phi _{\sigma
}\left( \rho \mathrm{e}^{\mathrm{i}\pi }\right) }{\phi _{\varepsilon }\left(
\rho \mathrm{e}^{\mathrm{i}\pi }\right) }\mathrm{e}^{\rho t\mathrm{e}^{^{%
\mathrm{i}\pi }}}\mathrm{e}^{\mathrm{i}\pi }\mathrm{d}\rho +\int_{0}^{\infty
}\frac{1}{\rho ^{\xi }\mathrm{e}^{-\mathrm{i}\xi \pi }}\frac{\phi _{\sigma
}\left( \rho \mathrm{e}^{-\mathrm{i}\pi }\right) }{\phi _{\varepsilon
}\left( \rho \mathrm{e}^{-\mathrm{i}\pi }\right) }\mathrm{e}^{\rho t\mathrm{e%
}^{-^{\mathrm{i}\pi }}}\mathrm{e}^{-\mathrm{i}\pi }\mathrm{d}\rho =0
\label{KTE-NP}
\end{equation}%
transforming into%
\begin{align*}
\mathcal{\epsilon }_{cr}\left( t\right) & =\mathcal{\epsilon }_{cr}^{\left( 
\mathrm{NP}\right) }\left( t\right) \\
& =\frac{1}{2\pi \mathrm{i}}\int_{0}^{\infty }\frac{1}{\rho ^{\xi }}\frac{%
\mathrm{e}^{\mathrm{i}\xi \pi }\phi _{\varepsilon }\left( \rho \mathrm{e}^{%
\mathrm{i}\pi }\right) \bar{\phi}_{\sigma }\left( \rho \mathrm{e}^{\mathrm{i}%
\pi }\right) -\mathrm{e}^{-\mathrm{i}\xi \pi }\bar{\phi}_{\varepsilon
}\left( \rho \mathrm{e}^{\mathrm{i}\pi }\right) \phi _{\sigma }\left( \rho 
\mathrm{e}^{\mathrm{i}\pi }\right) }{\left\vert \phi _{\varepsilon }\left(
\rho \mathrm{e}^{\mathrm{i}\pi }\right) \right\vert ^{2}}\mathrm{e}^{-\rho t}%
\mathrm{d}\rho \\
& =\frac{1}{\pi }\int_{0}^{\infty }\frac{1}{\rho ^{\xi }}\frac{K\left( \rho
\right) }{\left\vert \phi _{\varepsilon }\left( \rho \mathrm{e}^{\mathrm{i}%
\pi }\right) \right\vert ^{2}}\mathrm{e}^{-\rho t}\mathrm{d}\rho ,\quad 
\text{i.e.,} \\
\mathcal{\epsilon }_{cr}\left( t\right) & =\mathcal{\epsilon }_{cr}^{\left( 
\mathrm{NP}\right) }\left( t\right) \\
& =\frac{1}{2\pi \mathrm{i}}\int_{0}^{\infty }\frac{1}{\rho ^{\xi }}\frac{%
\left\vert \phi _{\sigma }\left( \rho \mathrm{e}^{\mathrm{i}\pi }\right)
\right\vert }{\left\vert \phi _{\varepsilon }\left( \rho \mathrm{e}^{\mathrm{%
i}\pi }\right) \right\vert }\left( \mathrm{e}^{\mathrm{i}\left( \xi \pi
+\arg \phi _{\varepsilon }\left( \rho \mathrm{e}^{\mathrm{i}\pi }\right)
-\arg \phi _{\sigma }\left( \rho \mathrm{e}^{\mathrm{i}\pi }\right) \right)
}-\mathrm{e}^{-\mathrm{i}\left( \xi \pi +\arg \phi _{\varepsilon }\left(
\rho \mathrm{e}^{\mathrm{i}\pi }\right) -\arg \phi _{\sigma }\left( \rho 
\mathrm{e}^{\mathrm{i}\pi }\right) \right) }\right) \mathrm{e}^{-\rho t}%
\mathrm{d}\rho \\
& =\frac{1}{\pi }\int_{0}^{\infty }\frac{1}{\rho ^{\xi }}\frac{\left\vert
\phi _{\sigma }\left( \rho \mathrm{e}^{\mathrm{i}\pi }\right) \right\vert }{%
\left\vert \phi _{\varepsilon }\left( \rho \mathrm{e}^{\mathrm{i}\pi
}\right) \right\vert }\sin \left( \arg \phi _{\varepsilon }\left( \rho 
\mathrm{e}^{\mathrm{i}\pi }\right) -\arg \phi _{\sigma }\left( \rho \mathrm{e%
}^{\mathrm{i}\pi }\right) +\xi \pi \right) \mathrm{e}^{-\rho t}\mathrm{d}%
\rho ,
\end{align*}%
that, according to (\ref{epsilon-creep-integral}) and according to the form (%
\ref{K}) of function $K,$ become creep compliance $\varepsilon _{cr}^{\left( 
\mathrm{NP}\right) }$ in the equivalent forms (\ref{epsilon-NP}) and (\ref%
{epsilon-NP-1}).

On the other hand, if the function $\tilde{\epsilon}_{cr},$ given by (\ref%
{s-epsilon}), has a negative real pole $s_{\scriptscriptstyle{\mathrm{RP}}%
}=\rho _{\scriptscriptstyle{\mathrm{RP}}}\,\mathrm{e}^{\mathrm{i}\pi },$
then the Cauchy integral theorem (\ref{Kosijeva-teorema-epsilon}) is
considered for contour $\Gamma ^{(\mathrm{II})}$ from Figure \ref%
{negativnaTG}, so that, in addition to the integrals along contours $\Gamma
_{3a}\cup \Gamma _{3b}$ and $\Gamma _{5a}\cup \Gamma _{5b},$ there are also
integrals along contours $\Gamma _{8}$ and $\Gamma _{9},$ that have non-zero
contributions, implying%
\begin{equation*}
\int_{\Gamma _{0}}\tilde{\epsilon}_{cr}\left( s\right) \mathrm{e}^{st}%
\mathrm{d}s+\int_{\Gamma _{3a}\cup \Gamma _{3b}}\tilde{\epsilon}_{cr}\left(
s\right) \mathrm{e}^{st}\mathrm{d}s+\int_{\Gamma _{5a}\cup \Gamma _{5b}}%
\tilde{\epsilon}_{cr}\left( s\right) \mathrm{e}^{st}\mathrm{d}s+\int_{\Gamma
_{8}}\tilde{\epsilon}_{cr}\left( s\right) \mathrm{e}^{st}\mathrm{d}%
s+\int_{\Gamma _{9}}\tilde{\epsilon}_{cr}\left( s\right) \mathrm{e}^{st}%
\mathrm{d}s=0,
\end{equation*}%
with the first three terms being already defined by (\ref{KTE-NP}) and with
the remaining terms calculated along the contours $\Gamma _{8}$ and $\Gamma
_{9},$ parameterized as in Table \ref{negativnaTG-param}, transforming the
previous expression in the limit when $r\rightarrow 0$ and $R\rightarrow
\infty $ into%
\begin{align*}
2\pi \mathrm{i\,}\mathcal{\epsilon }_{cr}\left( t\right) +2\pi \mathrm{i\,}%
\mathcal{\epsilon }_{cr}^{\left( \mathrm{NP}\right) }\left( t\right) &
+\int_{\pi }^{0}\frac{1}{\left( s_{\scriptscriptstyle{\mathrm{RP}}}+r\mathrm{%
e}^{\mathrm{i}\varphi }\right) ^{\xi }}\frac{\phi _{\sigma }\left( s_{%
\scriptscriptstyle{\mathrm{RP}}}+r\mathrm{e}^{\mathrm{i}\varphi }\right) }{%
\phi _{\varepsilon }\left( s_{\scriptscriptstyle{\mathrm{RP}}}+r\mathrm{e}^{%
\mathrm{i}\varphi }\right) }\mathrm{e}^{\left( s_{\scriptscriptstyle{\mathrm{%
RP}}}+r\mathrm{e}^{\mathrm{i}\varphi }\right) t}\mathrm{i}r\mathrm{e}^{%
\mathrm{i}\varphi }\mathrm{d}\varphi \\
& +\int_{0}^{-\pi }\frac{1}{\left( \bar{s}_{\scriptscriptstyle{\mathrm{RP}}%
}+r\mathrm{e}^{\mathrm{i}\varphi }\right) ^{\xi }}\frac{\phi _{\sigma
}\left( \bar{s}_{\scriptscriptstyle{\mathrm{RP}}}+r\mathrm{e}^{\mathrm{i}%
\varphi }\right) }{\phi _{\varepsilon }\left( \bar{s}_{\scriptscriptstyle{%
\mathrm{RP}}}+r\mathrm{e}^{\mathrm{i}\varphi }\right) }\mathrm{e}^{\left( 
\bar{s}_{\scriptscriptstyle{\mathrm{RP}}}+r\mathrm{e}^{\mathrm{i}\varphi
}\right) t}\mathrm{i}r\mathrm{e}^{\mathrm{i}\varphi }\mathrm{d}\varphi =0,
\end{align*}%
that becomes%
\begin{align*}
\mathcal{\epsilon }_{cr}\left( t\right) & =\mathcal{\epsilon }_{cr}^{\left( 
\mathrm{NP}\right) }\left( t\right) +\frac{1}{2\pi \mathrm{i}}\left( -%
\mathrm{i}\pi \frac{1}{s_{\scriptscriptstyle{\mathrm{RP}}}^{\xi }}\frac{\phi
_{\sigma }\left( s_{\scriptscriptstyle{\mathrm{RP}}}\right) }{\phi
_{\varepsilon }^{\prime }\left( s_{\scriptscriptstyle{\mathrm{RP}}}\right) }%
\mathrm{e}^{s_{\scriptscriptstyle{\mathrm{RP}}}t}-\mathrm{i}\pi \frac{1}{%
\bar{s}_{\scriptscriptstyle{\mathrm{RP}}}^{\xi }}\frac{\phi _{\sigma }\left( 
\bar{s}_{\scriptscriptstyle{\mathrm{RP}}}\right) }{\phi _{\varepsilon
}^{\prime }\left( \bar{s}_{\scriptscriptstyle{\mathrm{RP}}}\right) }\mathrm{e%
}^{\bar{s}_{\scriptscriptstyle{\mathrm{RP}}}t}\right) \\
& =\mathcal{\epsilon }_{cr}^{\left( \mathrm{NP}\right) }\left( t\right) -%
\frac{1}{2}\frac{1}{\rho _{\scriptscriptstyle{\mathrm{RP}}}^{\xi }}\frac{%
\mathrm{e}^{-\mathrm{i}\xi \pi }\bar{\phi}_{\varepsilon }^{\prime }\left( s_{%
\scriptscriptstyle{\mathrm{RP}}}\right) \phi _{\sigma }\left( s_{%
\scriptscriptstyle{\mathrm{RP}}}\right) +\mathrm{e}^{\mathrm{i}\xi \pi }\phi
_{\varepsilon }^{\prime }\left( s_{\scriptscriptstyle{\mathrm{RP}}}\right) 
\bar{\phi}_{\sigma }\left( s_{\scriptscriptstyle{\mathrm{RP}}}\right) }{%
\left\vert \phi _{\varepsilon }^{\prime }\left( s_{\scriptscriptstyle{%
\mathrm{RP}}}\right) \right\vert ^{2}}\mathrm{e}^{-\rho _{\scriptscriptstyle{%
\mathrm{RP}}}t},
\end{align*}%
where the last term transforms into%
\begin{equation*}
\mathcal{\epsilon }_{cr}^{\left( \mathrm{RP}\right) }\left( t\right) =-\frac{%
1}{\rho _{\scriptscriptstyle{\mathrm{RP}}}^{\xi }}\frac{\left\vert \phi
_{\sigma }\left( s_{\scriptscriptstyle{\mathrm{RP}}}\right) \right\vert }{%
\left\vert \phi _{\varepsilon }^{\prime }\left( s_{\scriptscriptstyle{%
\mathrm{RP}}}\right) \right\vert }\cos \left( \arg \phi _{\varepsilon
}^{\prime }\left( s_{\scriptscriptstyle{\mathrm{RP}}}\right) -\arg \phi
_{\sigma }\left( s_{\scriptscriptstyle{\mathrm{RP}}}\right) +\xi \pi \right) 
\mathrm{e}^{-\rho _{\scriptscriptstyle{\mathrm{RP}}}t},
\end{equation*}%
yielding the function $\varepsilon _{cr}^{\left( \mathrm{RP}\right) }$ in
the form (\ref{epsilon-RP}), according to (\ref{epsilon-creep-integral}).
The integral along contour $\Gamma _{8}$ is calculated as 
\begin{align*}
\lim_{r\rightarrow 0}\int_{\Gamma _{8}}\tilde{\epsilon}_{cr}\left( s\right) 
\mathrm{e}^{st}\mathrm{d}s& =\lim_{r\rightarrow 0}\int_{\pi }^{0}\frac{1}{%
\left( s_{\scriptscriptstyle{\mathrm{RP}}}+r\mathrm{e}^{\mathrm{i}\varphi
}\right) ^{\xi }}\frac{\phi _{\sigma }\left( s_{\scriptscriptstyle{\mathrm{RP%
}}}+r\mathrm{e}^{\mathrm{i}\varphi }\right) }{\phi _{\varepsilon }\left( s_{%
\scriptscriptstyle{\mathrm{RP}}}+r\mathrm{e}^{\mathrm{i}\varphi }\right) }%
\mathrm{e}^{\left( s_{\scriptscriptstyle{\mathrm{RP}}}+r\mathrm{e}^{\mathrm{i%
}\varphi }\right) t}\mathrm{i}r\mathrm{e}^{\mathrm{i}\varphi }\mathrm{d}%
\varphi \\
& =\lim_{r\rightarrow 0}\int_{\pi }^{0}\frac{1}{\left( s_{\scriptscriptstyle{%
\mathrm{RP}}}+r\mathrm{e}^{\mathrm{i}\varphi }\right) ^{\xi }}\frac{\phi
_{\sigma }\left( s_{\scriptscriptstyle{\mathrm{RP}}}+r\mathrm{e}^{\mathrm{i}%
\varphi }\right) }{\phi _{\varepsilon }\left( s_{\scriptscriptstyle{\mathrm{%
RP}}}\right) +\left. \phi _{\varepsilon }^{\prime }\left( s\right) \left(
s-s_{\scriptscriptstyle{\mathrm{RP}}}\right) \right\vert _{s=s_{%
\scriptscriptstyle{\mathrm{RP}}}+r\mathrm{e}^{\mathrm{i}\varphi }}+\ldots }%
\mathrm{e}^{\left( s_{\scriptscriptstyle{\mathrm{RP}}}+r\mathrm{e}^{\mathrm{i%
}\varphi }\right) t}\mathrm{i}r\mathrm{e}^{\mathrm{i}\varphi }\mathrm{d}%
\varphi \\
& =-\mathrm{i}\pi \frac{1}{s_{\scriptscriptstyle{\mathrm{RP}}}^{\xi }}\frac{%
\phi _{\sigma }\left( s_{\scriptscriptstyle{\mathrm{RP}}}\right) }{\phi
_{\varepsilon }^{\prime }\left( s_{\scriptscriptstyle{\mathrm{RP}}}\right) }%
\mathrm{e}^{s_{\scriptscriptstyle{\mathrm{RP}}}t},
\end{align*}%
by expanding function $\phi _{\varepsilon }$ into the series and by taking
into account $\phi _{\varepsilon }\left( s_{\scriptscriptstyle{\mathrm{RP}}%
}\right) =0$, while the similar calculation yields the integral along
contour $\Gamma _{9}$ in the form%
\begin{equation*}
\lim_{r\rightarrow 0}\int_{\Gamma _{9}}\tilde{\epsilon}_{cr}\left( s\right) 
\mathrm{e}^{st}\mathrm{d}s=-\mathrm{i}\pi \frac{1}{\bar{s}_{%
\scriptscriptstyle{\mathrm{RP}}}^{\xi }}\frac{\phi _{\sigma }\left( \bar{s}_{%
\scriptscriptstyle{\mathrm{RP}}}\right) }{\phi _{\varepsilon }^{\prime
}\left( \bar{s}_{\scriptscriptstyle{\mathrm{RP}}}\right) }\mathrm{e}^{\bar{s}%
_{\scriptscriptstyle{\mathrm{RP}}}t}.
\end{equation*}

In addition, when the function $\tilde{\epsilon}_{cr}$, see (\ref{s-epsilon}%
), has a pair of complex conjugated poles $s_{\scriptscriptstyle{\mathrm{RP}}%
}$ and $\bar{s}_{\scriptscriptstyle{\mathrm{RP}}}$, located within the area
bounded by the contour $\Gamma ^{(\mathrm{I})}$ depicted in Figure \ref%
{nemaTG}, the Cauchy residues theorem (\ref%
{Kosijeva-teorema-rezduumi-epsilon}), by taking into account integrals
having non-zero contributions, yields%
\begin{equation*}
\int_{\Gamma _{0}}\tilde{\epsilon}_{cr}\left( s\right) \mathrm{e}^{st}%
\mathrm{d}s+\int_{\Gamma _{3}}\tilde{\epsilon}_{cr}\left( s\right) \mathrm{e}%
^{st}\mathrm{d}s+\int_{\Gamma _{5}}\tilde{\epsilon}_{cr}\left( s\right) 
\mathrm{e}^{st}\mathrm{d}s=2\pi \mathrm{i}\left( \func{Res}\left( \tilde{%
\epsilon}_{cr}\left( s\right) \mathrm{e}^{st},s_{\scriptscriptstyle{\mathrm{%
CCP}}}\right) +\func{Res}\left( \tilde{\epsilon}_{cr}\left( s\right) \mathrm{%
e}^{st},\bar{s}_{\scriptscriptstyle{\mathrm{CCP}}}\right) \right) ,
\end{equation*}%
that, in the limit when $r\rightarrow 0$ and $R\rightarrow \infty $ and with
the terms on the left-hand-side of the previous expression already defined
by (\ref{KTE-NP}), becomes%
\begin{equation*}
\mathcal{\epsilon }_{cr}\left( t\right) -\mathcal{\epsilon }_{cr}^{\left( 
\mathrm{NP}\right) }\left( t\right) =\mathcal{\epsilon }_{cr}^{\left( 
\mathrm{CCP}\right) }\left( t\right) ,
\end{equation*}%
with%
\begin{align*}
\mathcal{\epsilon }_{cr}^{\left( \mathrm{CCP}\right) }\left( t\right) & =%
\frac{1}{\rho _{\scriptscriptstyle{\mathrm{CCP}}}^{\xi }\mathrm{e}^{\mathrm{i%
}\xi \varphi _{\scriptscriptstyle{\mathrm{CCP}}}}}\frac{\phi _{\sigma
}\left( s_{\scriptscriptstyle{\mathrm{CCP}}}\right) }{\phi _{\varepsilon
}^{\prime }\left( s_{\scriptscriptstyle{\mathrm{CCP}}}\right) }\mathrm{e}^{%
\func{Re}s_{\scriptscriptstyle{\mathrm{CCP}}}t+\mathrm{i}\func{Im}s_{%
\scriptscriptstyle{\mathrm{CCP}}}t} \\
& \quad \quad +\frac{1}{\rho _{\scriptscriptstyle{\mathrm{CCP}}}^{\xi }%
\mathrm{e}^{-\mathrm{i}\xi \varphi _{\scriptscriptstyle{\mathrm{CCP}}}}}%
\frac{\phi _{\sigma }\left( \bar{s}_{\scriptscriptstyle{\mathrm{CCP}}%
}\right) }{\phi _{\varepsilon }^{\prime }\left( \bar{s}_{\scriptscriptstyle{%
\mathrm{CCP}}}\right) }\mathrm{e}^{\func{Re}\bar{s}_{\scriptscriptstyle{%
\mathrm{CCP}}}t+\mathrm{i}\func{Im}\bar{s}_{\scriptscriptstyle{\mathrm{CCP}}%
}t} \\
& =\frac{1}{\rho _{\scriptscriptstyle{\mathrm{CCP}}}^{\xi }}\mathrm{e}^{%
\func{Re}s_{\scriptscriptstyle{\mathrm{CCP}}}t}\left( \frac{1}{\mathrm{e}^{%
\mathrm{i}\xi \varphi _{\scriptscriptstyle{\mathrm{CCP}}}}}\frac{\phi
_{\sigma }\left( s_{\scriptscriptstyle{\mathrm{CCP}}}\right) }{\phi
_{\varepsilon }^{\prime }\left( s_{\scriptscriptstyle{\mathrm{CCP}}}\right) }%
\mathrm{e}^{\mathrm{i}\func{Im}s_{\scriptscriptstyle{\mathrm{CCP}}}t}+\frac{1%
}{\mathrm{e}^{-\mathrm{i}\xi \varphi _{\scriptscriptstyle{\mathrm{CCP}}}}}%
\frac{\phi _{\sigma }\left( \bar{s}_{\scriptscriptstyle{\mathrm{CCP}}%
}\right) }{\phi _{\varepsilon }^{\prime }\left( \bar{s}_{\scriptscriptstyle{%
\mathrm{CCP}}}\right) }\mathrm{e}^{\mathrm{i}\func{Im}\bar{s}_{%
\scriptscriptstyle{\mathrm{CCP}}}t}\right) \\
& =\frac{1}{\rho _{\scriptscriptstyle{\mathrm{CCP}}}^{\xi }}\mathrm{e}%
^{-\left\vert \func{Re}s_{\scriptscriptstyle{\mathrm{CCP}}}\right\vert t} \\
& \quad \quad \times \frac{\mathrm{e}^{\mathrm{i}\left( \func{Im}s_{%
\scriptscriptstyle{\mathrm{CCP}}}t-\xi \varphi _{\scriptscriptstyle{\mathrm{%
CCP}}}\right) }\bar{\phi}_{\varepsilon }^{\prime }\left( s_{%
\scriptscriptstyle{\mathrm{CCP}}}\right) \phi _{\sigma }\left( s_{%
\scriptscriptstyle{\mathrm{CCP}}}\right) +\mathrm{e}^{-\mathrm{i}\left( 
\func{Im}s_{\scriptscriptstyle{\mathrm{CCP}}}t-\xi \varphi _{%
\scriptscriptstyle{\mathrm{CCP}}}\right) }\phi _{\varepsilon }^{\prime
}\left( s_{\scriptscriptstyle{\mathrm{CCP}}}\right) \bar{\phi}_{\sigma
}\left( s_{\scriptscriptstyle{\mathrm{CCP}}}\right) }{\left\vert \phi
_{\varepsilon }^{\prime }\left( s_{\scriptscriptstyle{\mathrm{CCP}}}\right)
\right\vert ^{2}} \\
& =2\frac{1}{\rho _{\scriptscriptstyle{\mathrm{CCP}}}^{\xi }}\frac{%
\left\vert \phi _{\sigma }\left( s_{\scriptscriptstyle{\mathrm{CCP}}}\right)
\right\vert }{\left\vert \phi _{\varepsilon }^{\prime }\left( s_{%
\scriptscriptstyle{\mathrm{CCP}}}\right) \right\vert }\mathrm{e}%
^{-\left\vert \func{Re}s_{\scriptscriptstyle{\mathrm{CCP}}}\right\vert
t}\cos \left( \func{Im}s_{\scriptscriptstyle{\mathrm{CCP}}}t-\arg \phi
_{\varepsilon }^{\prime }\left( s_{\scriptscriptstyle{\mathrm{CCP}}}\right)
+\arg \phi _{\sigma }\left( s_{\scriptscriptstyle{\mathrm{CCP}}}\right) -\xi
\varphi _{\scriptscriptstyle{\mathrm{CCP}}}\right) ,
\end{align*}%
implying, according to (\ref{epsilon-creep-integral}), 
\begin{equation*}
\varepsilon _{cr}\left( t\right) =\varepsilon _{cr}^{\left( \mathrm{NP}%
\right) }\left( t\right) +\varepsilon _{cr}^{\left( \mathrm{CCP}\right)
}\left( t\right) ,
\end{equation*}%
with 
\begin{align*}
\varepsilon _{cr}^{\left( \mathrm{CCP}\right) }\left( t\right) & =2\frac{1}{%
\rho _{\scriptscriptstyle{\mathrm{CCP}}}^{1+\xi }}\frac{\left\vert \phi
_{\sigma }\left( s_{\scriptscriptstyle{\mathrm{CCP}}}\right) \right\vert }{%
\left\vert \phi _{\varepsilon }^{\prime }\left( s_{\scriptscriptstyle{%
\mathrm{CCP}}}\right) \right\vert } \\
& \qquad \times \bigg(\mathrm{e}^{-\left\vert \func{Re}s_{\scriptscriptstyle{%
\mathrm{CCP}}}\right\vert t}\cos \left( \func{Im}s_{\scriptscriptstyle{%
\mathrm{CCP}}}t-\arg \phi _{\varepsilon }^{\prime }\left( s_{%
\scriptscriptstyle{\mathrm{CCP}}}\right) +\arg \phi _{\sigma }\left( s_{%
\scriptscriptstyle{\mathrm{CCP}}}\right) -\left( 1+\xi \right) \varphi _{%
\scriptscriptstyle{\mathrm{CCP}}}\right) \\
& \qquad \qquad -\cos \left( \arg \phi _{\varepsilon }^{\prime }\left( s_{%
\scriptscriptstyle{\mathrm{CCP}}}\right) -\arg \phi _{\sigma }\left( s_{%
\scriptscriptstyle{\mathrm{CCP}}}\right) +\left( 1+\xi \right) \varphi _{%
\scriptscriptstyle{\mathrm{CCP}}}\right) \bigg),
\end{align*}%
where the integral%
\begin{align*}
\int_{0}^{t}\mathrm{e}^{\mu t^{\prime }}\cos \left( \omega t^{\prime }+\phi
\right) \mathrm{d}t^{\prime }& =\frac{1}{\mu ^{2}+\omega ^{2}}\left( \mathrm{%
e}^{\mu t}\left( \omega \sin \left( \omega t+\phi \right) +\mu \cos \left(
\omega t+\phi \right) \right) -\left( \omega \sin \phi +\mu \cos \phi
\right) \right) \\
& =\frac{1}{\sqrt{\mu ^{2}+\omega ^{2}}}\left( \mathrm{e}^{\mu t}\cos \left(
\omega t+\phi -\varphi \right) -\cos \left( \phi -\varphi \right) \right)
\end{align*}%
is used with substitutions $\mu =\func{Re}s_{\scriptscriptstyle{\mathrm{CCP}}%
},$ $\omega =\func{Im}s_{\scriptscriptstyle{\mathrm{CCP}}},$ $\phi =-\arg
\phi _{\varepsilon }^{\prime }\left( s_{\scriptscriptstyle{\mathrm{CCP}}%
}\right) +\arg \phi _{\sigma }\left( s_{\scriptscriptstyle{\mathrm{CCP}}%
}\right) -\xi \varphi _{\scriptscriptstyle{\mathrm{CCP}}},$ and $\varphi
=\arctan \frac{\omega }{\mu }=\varphi _{\scriptscriptstyle{\mathrm{CCP}}}$.

The integrals along the contours $\Gamma _{1}$ ($\Gamma _{7}$), $\Gamma _{2}$
($\Gamma _{6}$), and $\Gamma _{4}$ tend to zero in the Cauchy integral and
residue theorems (\ref{Kosijeva-teorema-epsilon}) and (\ref%
{Kosijeva-teorema-rezduumi-epsilon}) in the limit when $r\rightarrow 0$ and $%
R\rightarrow \infty $. Namely, if the power-type functions $\phi
_{\varepsilon }$ and $\phi _{\sigma },$ see Table \ref{skupina}, satisfy
conditions 
\begin{gather*}
\frac{\left\vert \phi _{\varepsilon }\left( p+\mathrm{i}R\right) \right\vert 
}{\left\vert \phi _{\sigma }\left( p+\mathrm{i}R\right) \right\vert }\sim
R^{\zeta _{R}},\quad \text{when}\quad R\rightarrow \infty \quad \text{for}%
\quad p\in \left[ 0,p_{0}\right] , \\
\frac{\left\vert \phi _{\varepsilon }\left( r\mathrm{e}^{\mathrm{i}\varphi
}\right) \right\vert }{\left\vert \phi _{\sigma }\left( r\mathrm{e}^{\mathrm{%
i}\varphi }\right) \right\vert }\sim \frac{1}{r^{\zeta _{r}}},\quad \text{%
when}\quad r\rightarrow 0\quad \text{for}\quad \varphi \in \left[ -\pi ,\pi %
\right] ,
\end{gather*}%
that are already posed in Section \ref{CRM} in order to ensure zero
contributions of the integrals of function $\tilde{\sigma}_{sr}\left(
s\right) \mathrm{e}^{st}$ along the contours $\Gamma _{1}$ ($\Gamma _{7}$), $%
\Gamma _{2}$ ($\Gamma _{6}$), and $\Gamma _{4}$, then one has%
\begin{gather*}
\frac{\left\vert \phi _{\sigma }\left( p+\mathrm{i}R\right) \right\vert }{%
\left\vert \phi _{\varepsilon }\left( p+\mathrm{i}R\right) \right\vert }\sim 
\frac{1}{R^{\zeta _{R}}},\quad \text{when}\quad R\rightarrow \infty \quad 
\text{for}\quad p\in \left[ 0,p_{0}\right] , \\
\frac{\left\vert \phi _{\sigma }\left( r\mathrm{e}^{\mathrm{i}\varphi
}\right) \right\vert }{\left\vert \phi _{\varepsilon }\left( r\mathrm{e}^{%
\mathrm{i}\varphi }\right) \right\vert }\sim r^{\zeta _{r}},\quad \text{when}%
\quad r\rightarrow 0\quad \text{for}\quad \varphi \in \left[ -\pi ,\pi %
\right] ,
\end{gather*}%
and due to the similar forms of the relaxation modulus in the Laplace domain
and function $\tilde{\epsilon}_{cr},$ see (\ref{sigma-sr-ld}) and (\ref%
{s-epsilon}), one has the zero contributions of integrals along the contours 
$\Gamma _{1}$ ($\Gamma _{7}$) and $\Gamma _{4}$ if $\zeta _{R}>-\xi $ and $%
\zeta _{r}>-\left( 1-\xi \right) $, yielding%
\begin{equation*}
-\xi <\zeta _{R}<1-\xi \quad \text{and}\quad -\left( 1-\xi \right) <\zeta
_{r}<\xi
\end{equation*}%
with the previously posed conditions $\zeta _{R}<1-\xi $ and $\zeta _{r}<\xi
,$ while the zero contributions of integrals along the contours $\Gamma _{2}$
($\Gamma _{6}$) is simply guaranteed by the fact that $\phi _{\varepsilon }$
and $\phi _{\sigma }$ are power-type functions.

\section{Narrowing the thermodynamical requirements\label{narrovlje}}

\subsection{Symmetric models}

\subsubsection{Model ID.ID}

Thermodynamical, see (\ref{TD-AZ-ID-ID-less-1}) and (\ref{TD-AZ-ID-ID-less-2}%
), along with the narrowed thermodynamical restrictions on the parameters
for model ID.ID, given by the expression%
\begin{equation*}
\left( a_{1}\,_{0}\mathrm{I}_{t}^{\alpha }+a_{2}\,_{0}\mathrm{D}_{t}^{\beta
}\right) \sigma \left( t\right) =\left( b_{1}\,_{0}\mathrm{I}_{t}^{\mu
}+a_{2}\,_{0}\mathrm{D}_{t}^{\alpha +\beta -\mu }\right) \varepsilon \left(
t\right) ,
\end{equation*}%
are of the following form%
\begin{gather}
0\leqslant \alpha +\beta -\mu \leqslant 1,\quad \mu \leqslant \alpha ,\quad
\beta +\mu \leqslant 1,  \label{STD-ID.DD-1} \\
-\frac{a_{1}}{a_{2}}\frac{\cos \frac{\left( 2\alpha +\beta -\mu \right) \pi 
}{2}}{\cos \frac{\left( \beta +\mu \right) \pi }{2}}\leqslant \frac{b_{1}}{%
b_{2}}\leqslant \frac{a_{1}}{a_{2}}\frac{\sin \frac{\left( 2\alpha +\beta
-\mu \right) \pi }{2}}{\sin \frac{\left( \beta +\mu \right) \pi }{2}}\frac{%
\cos \frac{\left( 2\alpha +\beta -\mu \right) \pi }{2}}{\cos \frac{\left(
\beta +\mu \right) \pi }{2}}\leqslant \frac{a_{1}}{a_{2}}\frac{\sin \frac{%
\left( 2\alpha +\beta -\mu \right) \pi }{2}}{\sin \frac{\left( \beta +\mu
\right) \pi }{2}},  \label{STD-ID.DD-2}
\end{gather}%
where the inequality (\ref{STD-ID.DD-2})$_{3}$ narrows the thermodynamical
restriction (\ref{TD-AZ-ID-ID-less-2}) if $\alpha \leqslant 2\alpha +\beta
-\mu <1,$ and it is obtained by requesting the function $K$, given in the
general form by (\ref{K}) and reducing to the expression%
\begin{align}
K\left( \rho \right) & =a_{1}b_{1}\sin \left( \left( \alpha -\mu \right) \pi
\right) +a_{1}b_{2}\rho ^{\alpha +\beta }\sin \left( \left( 2\alpha +\beta
-\mu \right) \pi \right)  \notag \\
& \quad -a_{2}b_{1}\rho ^{\alpha +\beta }\sin \left( \left( \beta +\mu
\right) \pi \right) +a_{2}b_{2}\rho ^{2\left( \alpha +\beta \right) }\sin
\left( \left( \alpha -\mu \right) \pi \right)  \label{K-ID.DD}
\end{align}%
for the model ID.ID, to have a non-negative value. The second term in (\ref%
{K-ID.DD}), that can be either positive for $2\alpha +\beta -\mu \in \left[
\alpha ,1\right) ,$ or non-positive for $2\alpha +\beta -\mu \in \left[
1,1+\alpha \right) ,$ see (\ref{STD-ID.DD-1})$_{1}$, combined with the third
non-positive term in (\ref{K-ID.DD}) yields narrowed thermodynamical
restriction (\ref{STD-ID.DD-2})$_{3}$ if $2\alpha +\beta -\mu \in \left[
\alpha ,1\right) ,$ since $\beta +\mu \leqslant 2\alpha +\beta -\mu <1$
reducing to $\mu \leqslant \alpha ,$ see (\ref{STD-ID.DD-1})$_{2}$, implies $%
\frac{\cos \frac{\left( 2\alpha +\beta -\mu \right) \pi }{2}}{\cos \frac{%
\left( \beta +\mu \right) \pi }{2}}\leqslant 1$, while if $2\alpha +\beta
-\mu \in \left[ 1,1+\alpha \right) ,$ then one has the second and third term
in (\ref{K-ID.DD}) non-positive and therefore non-negativity of function $K\ 
$cannot be guaranteed, see (\ref{K-ID.DD}).

In conclusion, the non-negativity of function $K,$ corresponding to the
model ID.ID and given by (\ref{K-ID.DD}), is ensured if $2\alpha +\beta -\mu
\in \left[ \alpha ,1\right) $ by requesting (\ref{STD-ID.DD-2}) in addition
to (\ref{STD-ID.DD-1}) implying the corresponding relaxation modulus to be a
completely monotonic function and creep compliance a Bernstein function,
while if $2\alpha +\beta -\mu \in \left[ 1,1+\alpha \right) ,$ then the
mentioned properties cannot be ensured.

\subsubsection{Model ID.DD$^{{}^{+}}$}

Thermodynamical, see (\ref{TD-AZ-ID-DD-over-1}) and (\ref{TD-AZ-ID-DD-over-2}%
), along with the narrowed thermodynamical restrictions on the parameters
for model ID.DD$^{{}^{+}}$, given by the expression%
\begin{equation*}
\left( a_{1}\,_{0}\mathrm{I}_{t}^{\alpha }+a_{2}\,_{0}\mathrm{D}_{t}^{\beta
}\right) \sigma \left( t\right) =\left( b_{1}\,_{0}\mathrm{D}_{t}^{\mu
}+b_{2}\,_{0}\mathrm{D}_{t}^{\alpha +\beta +\mu }\right) \varepsilon \left(
t\right) ,  \label{ID.DD+}
\end{equation*}%
are of the following form%
\begin{gather}
1\leqslant \alpha +\beta +\mu \leqslant 2,\quad \beta \leqslant \mu
\leqslant 1-\alpha ,  \label{STD-ID.DD+-1} \\
\frac{a_{1}}{a_{2}}\frac{\left\vert \cos \frac{\left( 2\alpha +\beta +\mu
\right) \pi }{2}\right\vert }{\cos \frac{\left( \mu -\beta \right) \pi }{2}}%
\leqslant \frac{a_{1}}{a_{2}}\frac{\left\vert \cos \frac{\left( 2\alpha
+\beta +\mu \right) \pi }{2}\right\vert }{\cos \frac{\left( \mu -\beta
\right) \pi }{2}}\frac{\sin \frac{\left( 2\alpha +\beta +\mu \right) \pi }{2}%
}{\sin \frac{\left( \mu -\beta \right) \pi }{2}}\leqslant \frac{b_{1}}{b_{2}}%
,  \label{STD-ID.DD+-2}
\end{gather}%
where the inequality (\ref{STD-ID.DD+-2})$_{1}$ narrows the thermodynamical
restriction (\ref{TD-AZ-ID-DD-over-2}) and it is obtained by requesting the
function $K$, given in the general form by (\ref{K}) and reducing to the
expression%
\begin{align}
K\left( \rho \right) & =a_{1}b_{1}\sin \left( \left( \alpha +\mu \right) \pi
\right) +a_{1}b_{2}\rho ^{\alpha +\beta }\sin \left( \left( 2\alpha +\beta
+\mu \right) \pi \right)  \notag \\
& \quad +a_{2}b_{1}\rho ^{\alpha +\beta }\sin \left( \left( \mu -\beta
\right) \pi \right) +a_{2}b_{2}\rho ^{2\left( \alpha +\beta \right) }\sin
\left( \left( \alpha +\mu \right) \pi \right) ,  \label{K-ID.DD+}
\end{align}%
for the model ID.DD$^{{}^{+}}$, to have a non-negative value. The second
term in (\ref{K-ID.DD+}), that is non-positive since (\ref{STD-ID.DD+-1})$%
_{1}$ implies $2\alpha +\beta +\mu \geqslant 1+\alpha $ and (\ref%
{STD-ID.DD+-1})$_{2}$ implies $2\alpha +\beta +\mu =\left( \alpha +\beta
\right) +\left( \alpha +\mu \right) \leqslant 2$ and therefore $2\alpha
+\beta +\mu \in \left[ 1+\alpha ,2\right] $, combined with the third
non-negative term in (\ref{K-ID.DD+}) yields narrowed thermodynamical
restriction (\ref{STD-ID.DD+-2})$_{1},$ since $\frac{\left( \mu -\beta
\right) \pi }{2}\leqslant \pi -\frac{\left( 2\alpha +\beta +\mu \right) \pi 
}{2}<\frac{\pi }{2}$ reducing to $\alpha +\mu \leqslant 1$, see (\ref%
{STD-ID.DD+-1})$_{2}$, implies $\frac{\sin \frac{\left( 2\alpha +\beta +\mu
\right) \pi }{2}}{\sin \frac{\left( \mu -\beta \right) \pi }{2}}\geqslant 1$.

In conclusion, the non-negativity of function $K,$ corresponding to the
model ID.DD$^{{}^{+}}$ and given by (\ref{K-ID.DD+}), is ensured by
requesting (\ref{STD-ID.DD+-2}) in addition to (\ref{STD-ID.DD+-1}),
implying the corresponding relaxation modulus to be a completely monotonic
function and creep compliance a Bernstein function.

\subsubsection{Model IID.IID}

Function $K$, given by (\ref{K}), in the case of model IID.IID, having the
constitutive equation given in the form%
\begin{equation*}
\left( a_{1}\,_{0}\mathrm{I}_{t}^{\alpha }+a_{2}\,_{0}\mathrm{I}_{t}^{\beta
}+a_{3}\,_{0}\mathrm{D}_{t}^{\gamma }\right) \sigma \left( t\right) =\left(
b_{1}\,_{0}\mathrm{I}_{t}^{\alpha +\gamma -\eta }+b_{2}\,_{0}\mathrm{I}%
_{t}^{\beta +\gamma -\eta }+b_{3}\,_{0}\mathrm{D}_{t}^{\eta }\right)
\varepsilon \left( t\right) ,
\end{equation*}%
and with corresponding thermodynamical, see (\ref{franjo-less-less-TD-1}) - (%
\ref{franjo-less-less-TD-3}), as well as with the narrowed thermodynamical
requirements%
\begin{gather}
\beta <\alpha ,\quad \gamma \leqslant \eta ,\quad 0\leqslant \beta +\gamma
-\eta \leqslant \alpha +2\gamma -\eta \leqslant 1,\quad \alpha +\gamma
\leqslant \beta +\eta ,  \label{STD-IID.IID-1} \\
-\frac{b_{3}}{b_{1}}\frac{\cos \frac{\left( \alpha +\eta \right) \pi }{2}}{%
\cos \frac{\left( \alpha +2\gamma -\eta \right) \pi }{2}}\leqslant \frac{%
a_{3}}{a_{1}}\leqslant \frac{b_{3}}{b_{1}}\frac{\sin \frac{\left( \alpha
+\eta \right) \pi }{2}}{\sin \frac{\left( \alpha +2\gamma -\eta \right) \pi 
}{2}}\frac{\cos \frac{\left( \alpha +\eta \right) \pi }{2}}{\cos \frac{%
\left( \alpha +2\gamma -\eta \right) \pi }{2}}\leqslant \frac{b_{3}}{b_{1}}%
\frac{\sin \frac{\left( \alpha +\eta \right) \pi }{2}}{\sin \frac{\left(
\alpha +2\gamma -\eta \right) \pi }{2}},  \label{STD-IID.IID-2} \\
-\frac{b_{3}}{b_{2}}\frac{\cos \frac{\left( \beta +\eta \right) \pi }{2}}{%
\cos \frac{\left( \beta +2\gamma -\eta \right) \pi }{2}}\leqslant \frac{a_{3}%
}{a_{2}}\leqslant \frac{b_{3}}{b_{2}}\frac{\sin \frac{\left( \beta +\eta
\right) \pi }{2}}{\sin \frac{\left( \beta +2\gamma -\eta \right) \pi }{2}}%
\frac{\cos \frac{\left( \beta +\eta \right) \pi }{2}}{\cos \frac{\left(
\beta +2\gamma -\eta \right) \pi }{2}}\leqslant \frac{b_{3}}{b_{2}}\frac{%
\sin \frac{\left( \beta +\eta \right) \pi }{2}}{\sin \frac{\left( \beta
+2\gamma -\eta \right) \pi }{2}},  \label{STD-IID.IID-3}
\end{gather}%
valid if $\beta +\eta <\alpha +\eta <1,$ becomes%
\begin{align}
K\left( \rho \right) &=a_{1}b_{1}\sin \left( \left( \eta -\gamma \right) \pi
\right) +a_{1}b_{2}\rho ^{\alpha -\beta }\sin \left( \left( \alpha +\eta
-\beta -\gamma \right) \pi \right) +a_{1}b_{3}\rho ^{\alpha +\gamma }\sin
\left( \left( \alpha +\eta \right) \pi \right)  \notag \\
&\quad+a_{2}b_{1}\rho ^{\alpha -\beta }\sin \left( \left( \beta +\eta
-\alpha -\gamma \right) \pi \right) +a_{2}b_{2}\rho ^{2\left( \alpha -\beta
\right) }\sin \left( \left( \eta -\gamma \right) \pi \right) +a_{2}b_{3}\rho
^{2\alpha -\beta +\gamma }\sin \left( \left( \beta +\eta \right) \pi \right)
\notag \\
&\quad-a_{3}b_{1}\rho ^{\alpha +\gamma }\sin \left( \left( \alpha +2\gamma
-\eta \right) \pi \right) -a_{3}b_{2}\rho ^{2\alpha -\beta +\gamma }\sin
\left( \left( \beta +2\gamma -\eta \right) \pi \right) +a_{3}b_{3}\rho
^{2\left( \alpha +\gamma \right) }\sin \left( \left( \eta -\gamma \right)
\pi \right) ,  \label{K-IID.IID}
\end{align}%
so that the non-negativity requirement of function $K$ implies%
\begin{gather}
a_{1}b_{2}\sin \left( \left( \alpha +\eta -\beta -\gamma \right) \pi \right)
+a_{2}b_{1}\sin \left( \left( \beta +\eta -\alpha -\gamma \right) \pi
\right) \geqslant 0,  \label{NN-of-K-IID.IID-1} \\
a_{1}b_{3}\sin \left( \left( \alpha +\eta \right) \pi \right)
-a_{3}b_{1}\sin \left( \left( \alpha +2\gamma -\eta \right) \pi \right)
\geqslant 0,  \label{NN-of-K-IID.IID-2} \\
a_{2}b_{3}\sin \left( \left( \beta +\eta \right) \pi \right) -a_{3}b_{2}\sin
\left( \left( \beta +2\gamma -\eta \right) \pi \right) \geqslant 0.
\label{NN-of-K-IID.IID-3}
\end{gather}

According to the thermodynamical requirements (\ref{STD-IID.IID-1})$_{1,2},$
one finds that $\alpha +\eta -\beta -\gamma =\left( \alpha -\beta \right)
+\left( \eta -\gamma \right) \geqslant 0$, while the thermodynamical
requirement $\beta +\gamma -\eta \geqslant 0$, given by (\ref{STD-IID.IID-1})%
$_{3}$, can be modified into $\alpha +\eta -\beta -\gamma \leqslant \alpha
\leqslant 1,$ and therefore the first term in (\ref{NN-of-K-IID.IID-1}) has
a non-negative value. The argument $\left( \beta +\eta -\alpha -\gamma
\right) \pi =\left( \left( \eta -\gamma \right) -\left( \alpha -\beta
\right) \right) \pi $ of the second term in (\ref{NN-of-K-IID.IID-1}),
firstly has a value that is less than a value of the argument $\left( \alpha
+\eta -\beta -\gamma \right) \pi =\left( \left( \eta -\gamma \right) +\left(
\alpha -\beta \right) \right) \pi $ of the first term in (\ref%
{NN-of-K-IID.IID-1}), more precisely $\beta +\eta -\alpha -\gamma =\left(
\eta -\gamma \right) -\left( \alpha -\beta \right) \leqslant \alpha +\eta
-\beta -\gamma =\left( \eta -\gamma \right) +\left( \alpha -\beta \right)
\leqslant \alpha \leqslant 1,$ and, secondly it has a non-negative value,
i.e., $\beta +\eta -\alpha -\gamma =\left( \beta +\eta \right) -\left(
\alpha +\gamma \right) \geqslant 0,$\ according to the thermodynamical
requirement (\ref{STD-IID.IID-1})$_{4}$. Therefore, the inequality (\ref%
{NN-of-K-IID.IID-1}) is trivially satisfied.

According to (\ref{STD-IID.IID-1})$_{3},$ the sine in the second term in (%
\ref{NN-of-K-IID.IID-2}) has a non-negative value, while the first term has
either a positive value if $\alpha +\eta <1,$ or a non-positive value if $%
\alpha +\eta \geqslant 1,$ so that if $\alpha +\eta <1,$ then the
requirement (\ref{NN-of-K-IID.IID-2}) is transformed into%
\begin{equation*}
\frac{a_{3}}{a_{1}}\leqslant \frac{b_{3}}{b_{1}}\frac{\sin \frac{\left(
\alpha +\eta \right) \pi }{2}}{\sin \frac{\left( \alpha +2\gamma -\eta
\right) \pi }{2}}\frac{\cos \frac{\left( \alpha +\eta \right) \pi }{2}}{\cos 
\frac{\left( \alpha +2\gamma -\eta \right) \pi }{2}},
\end{equation*}%
representing the narrowed thermodynamical restriction (\ref{STD-IID.IID-2})$%
_{3}$, since $\alpha +2\gamma -\eta \leqslant \alpha +\eta <1,$ reducing to $%
\gamma \leqslant \eta ,$ see (\ref{STD-IID.IID-1})$_{2}$, implies $\frac{%
\cos \frac{\left( \alpha +\eta \right) \pi }{2}}{\cos \frac{\left( \alpha
+2\gamma -\eta \right) \pi }{2}}\leqslant 1$. Similarly, according to (\ref%
{STD-IID.IID-1})$_{1,3},$ the sine in the second term in (\ref%
{NN-of-K-IID.IID-3}) has a non-negative value, since $0\leqslant \beta
+2\gamma -\eta \leqslant \alpha +2\gamma -\eta \leqslant 1$, while the first
term has either a positive value if $\beta +\eta <1,$ or a non-positive
value if $\beta +\eta \geqslant 1,$ so that if $\beta +\eta <1,$ then the
requirement (\ref{NN-of-K-IID.IID-3}) is transformed into%
\begin{equation*}
\frac{a_{3}}{a_{2}}\leqslant \frac{b_{3}}{b_{2}}\frac{\sin \frac{\left(
\beta +\eta \right) \pi }{2}}{\sin \frac{\left( \beta +2\gamma -\eta \right)
\pi }{2}}\frac{\cos \frac{\left( \beta +\eta \right) \pi }{2}}{\cos \frac{%
\left( \beta +2\gamma -\eta \right) \pi }{2}},
\end{equation*}%
representing the narrowed thermodynamical restriction (\ref{STD-IID.IID-3})$%
_{3}$, since $\beta +2\gamma -\eta \leqslant \beta +\eta <1$ reducing to $%
\gamma \leqslant \eta ,$ see (\ref{STD-IID.IID-1})$_{2}$, implies $\frac{%
\cos \frac{\left( \beta +\eta \right) \pi }{2}}{\cos \frac{\left( \beta
+2\gamma -\eta \right) \pi }{2}}\leqslant 1.$

In conclusion, the non-negativity of function $K,$ corresponding to the
model IID.IID and given by (\ref{K-IID.IID}), is ensured if $\beta +\eta
<\alpha +\eta <1,$ see (\ref{STD-IID.IID-1})$_{1}$, by requesting (\ref%
{STD-IID.IID-2}) and (\ref{STD-IID.IID-3}) in addition to (\ref%
{STD-IID.IID-1}) implying the corresponding relaxation modulus to be a
completely monotonic function and creep compliance a Bernstein function,
while if either $\alpha +\eta \geqslant 1$ or $\alpha +\eta >\beta +\eta
\geqslant 1,$ then the mentioned properties cannot be ensured, since either (%
\ref{NN-of-K-IID.IID-2}), or both (\ref{NN-of-K-IID.IID-2}) and (\ref%
{NN-of-K-IID.IID-3}) cannot be satisfied.

\subsubsection{Model IDD.IDD}

Function $K$, given by (\ref{K}), in the case of the model IDD.IDD, having
the constitutive equation given in the form%
\begin{equation*}
\left( a_{1}\,_{0}\mathrm{I}_{t}^{\alpha }+a_{2}\,_{0}\mathrm{D}_{t}^{\beta
}+a_{3}\,_{0}\mathrm{D}_{t}^{\gamma }\right) \sigma \left( t\right) =\left(
b_{1}\,_{0}\mathrm{I}_{t}^{\mu }+b_{2}\,_{0}\mathrm{D}_{t}^{\alpha +\beta
-\mu }+b_{3}\,_{0}\mathrm{D}_{t}^{\alpha +\gamma -\mu }\right) \varepsilon
\left( t\right) ,
\end{equation*}%
and with corresponding thermodynamical, see (\ref{IDD-IDD-less-less-TD-1}) -
(\ref{IDD-IDD-less-less-TD-3}), as well as with the narrowed thermodynamical
requirements%
\begin{gather}
0\leqslant \alpha +\gamma -\mu \leqslant 1,\quad \beta <\gamma ,\quad \mu
\leqslant \alpha ,\quad \gamma +\mu \leqslant \alpha +\beta ,\quad \gamma
+\mu \leqslant 1,  \label{STD-IDD.IDD-1} \\
-\frac{b_{2}}{b_{1}}\frac{\cos \frac{\left( 2\alpha +\beta -\mu \right) \pi 
}{2}}{\cos \frac{\left( \beta +\mu \right) \pi }{2}}\leqslant \frac{a_{2}}{%
a_{1}}\leqslant \frac{b_{2}}{b_{1}}\frac{\sin \frac{\left( 2\alpha +\beta
-\mu \right) \pi }{2}}{\sin \frac{\left( \beta +\mu \right) \pi }{2}}\frac{%
\cos \frac{\left( 2\alpha +\beta -\mu \right) \pi }{2}}{\cos \frac{\left(
\beta +\mu \right) \pi }{2}}\leqslant \frac{b_{2}}{b_{1}}\frac{\sin \frac{%
\left( 2\alpha +\beta -\mu \right) \pi }{2}}{\sin \frac{\left( \beta +\mu
\right) \pi }{2}},  \label{STD-IDD.IDD-2} \\
-\frac{b_{3}}{b_{1}}\frac{\cos \frac{\left( 2\alpha +\gamma -\mu \right) \pi 
}{2}}{\cos \frac{\left( \gamma +\mu \right) \pi }{2}}\leqslant \frac{a_{3}}{%
a_{1}}\leqslant \frac{b_{3}}{b_{1}}\frac{\sin \frac{\left( 2\alpha +\gamma
-\mu \right) \pi }{2}}{\sin \frac{\left( \gamma +\mu \right) \pi }{2}}\frac{%
\cos \frac{\left( 2\alpha +\gamma -\mu \right) \pi }{2}}{\cos \frac{\left(
\gamma +\mu \right) \pi }{2}}\leqslant \frac{b_{3}}{b_{1}}\frac{\sin \frac{%
\left( 2\alpha +\gamma -\mu \right) \pi }{2}}{\sin \frac{\left( \gamma +\mu
\right) \pi }{2}},  \label{STD-IDD.IDD-3}
\end{gather}%
valid if $2\alpha +\beta -\mu <2\alpha +\gamma -\mu <1$, becomes%
\begin{align}
K\left( \rho \right) &=a_{1}b_{1}\sin \left( \left( \alpha -\mu \right) \pi
\right) +a_{1}b_{2}\rho ^{\alpha +\beta }\sin \left( \left( 2\alpha +\beta
-\mu \right) \pi \right) +a_{1}b_{3}\rho ^{\alpha +\gamma }\sin \left(
\left( 2\alpha +\gamma -\mu \right) \pi \right)  \notag \\
&\quad-a_{2}b_{1}\rho ^{\alpha +\beta }\sin \left( \left( \beta +\mu \right)
\pi \right) +a_{2}b_{2}\rho ^{2\left( \alpha +\beta \right) }\sin \left(
\left( \alpha -\mu \right) \pi \right) +a_{2}b_{3}\rho ^{2\alpha +\beta
+\gamma }\sin \left( \left( \alpha +\gamma -\beta -\mu \right) \pi \right) 
\notag \\
&\quad-a_{3}b_{1}\rho ^{\alpha +\gamma }\sin \left( \left( \gamma +\mu
\right) \pi \right) +a_{3}b_{2}\rho ^{2\alpha +\beta +\gamma }\sin \left(
\left( \alpha +\beta -\gamma -\mu \right) \pi \right) +a_{3}b_{3}\rho
^{2\left( \alpha +\gamma \right) }\sin \left( \left( \alpha -\mu \right) \pi
\right) ,  \label{K-IDD.IDD}
\end{align}%
so that the non-negativity requirement of function $K$ implies%
\begin{gather}
a_{2}b_{3}\sin \left( \left( \alpha +\gamma -\beta -\mu \right) \pi \right)
+a_{3}b_{2}\sin \left( \left( \alpha +\beta -\gamma -\mu \right) \pi \right)
\geqslant 0,  \label{NN-of-K-IDD.IDD-1} \\
a_{1}b_{2}\sin \left( \left( 2\alpha +\beta -\mu \right) \pi \right)
-a_{2}b_{1}\sin \left( \left( \beta +\mu \right) \pi \right) \geqslant 0,
\label{NN-of-K-IDD.IDD-2} \\
a_{1}b_{3}\sin \left( \left( 2\alpha +\gamma -\mu \right) \pi \right)
-a_{3}b_{1}\sin \left( \left( \gamma +\mu \right) \pi \right) \geqslant 0.
\label{NN-of-K-IDD.IDD-3}
\end{gather}

The first term in (\ref{NN-of-K-IDD.IDD-1}) has a non-negative value since $%
\alpha +\gamma -\beta -\mu =\left( \alpha -\mu \right) +\left( \gamma -\beta
\right) \geqslant 0$, according to the thermodynamical requirements (\ref%
{STD-IDD.IDD-1})$_{2,3}$, as well as $\alpha +\gamma -\beta -\mu =\left(
\alpha +\gamma -\mu \right) -\beta \leqslant 1-\beta ,$ according to the
thermodynamical requirement (\ref{STD-IDD.IDD-1})$_{1}$. The second term in (%
\ref{NN-of-K-IDD.IDD-1}) also has a non-negative value, since $\alpha +\beta
-\gamma -\mu =\left( \alpha +\beta \right) -\left( \gamma +\mu \right)
\geqslant 0$, according to the thermodynamical requirement (\ref%
{STD-IDD.IDD-1})$_{4}$, as well as $\alpha +\beta -\gamma -\mu =\alpha
-\left( \gamma +\mu -\beta \right) \leqslant \alpha \leqslant 1,$ since $%
\left( \gamma -\beta \right) +\mu \geqslant 0$, according to (\ref%
{STD-IDD.IDD-1})$_{2}.$ Therefore, the inequality (\ref{NN-of-K-IDD.IDD-1})
is trivially satisfied.

According to (\ref{STD-IDD.IDD-1})$_{2,5}$, yielding $0\leqslant \beta +\mu
<\gamma +\mu \leqslant 1,$ the sine in the second term in (\ref%
{NN-of-K-IDD.IDD-2}) as well as in (\ref{NN-of-K-IDD.IDD-3}) has a
non-negative value, while the first term has either a positive value if $%
2\alpha +\beta -\mu <1,$ or a non-positive value if $2\alpha +\beta -\mu
\geqslant 1,$ so that if $2\alpha +\beta -\mu <1,$ then the requirement (\ref%
{NN-of-K-IDD.IDD-2}) is transformed into%
\begin{equation*}
\frac{a_{2}}{a_{1}}\leqslant \frac{b_{2}}{b_{1}}\frac{\sin \frac{\left(
2\alpha +\beta -\mu \right) \pi }{2}}{\sin \frac{\left( \beta +\mu \right)
\pi }{2}}\frac{\cos \frac{\left( 2\alpha +\beta -\mu \right) \pi }{2}}{\cos 
\frac{\left( \beta +\mu \right) \pi }{2}},
\end{equation*}%
representing the narrowed thermodynamical restriction (\ref{STD-IDD.IDD-2})$%
_{3}$, since $\beta +\mu \leqslant 2\alpha +\beta -\mu <1,$ reducing to $\mu
\leqslant \alpha ,$ see (\ref{STD-IDD.IDD-1})$_{3},$ implies $\frac{\cos 
\frac{\left( 2\alpha +\beta -\mu \right) \pi }{2}}{\cos \frac{\left( \beta
+\mu \right) \pi }{2}}\leqslant 1$. Similarly, the second term in (\ref%
{NN-of-K-IDD.IDD-3}), as already mentioned, has a non-negative value, while
the first term has either a positive value if $2\alpha +\gamma -\mu <1$ or a
non-positive value if $2\alpha +\gamma -\mu \geqslant 1,$ so that if $%
2\alpha +\gamma -\mu <1,$ then the requirement (\ref{NN-of-K-IDD.IDD-3}) is
transformed into%
\begin{equation*}
\frac{a_{3}}{a_{1}}\leqslant \frac{b_{3}}{b_{1}}\frac{\sin \frac{\left(
2\alpha +\gamma -\mu \right) \pi }{2}}{\sin \frac{\left( \gamma +\mu \right)
\pi }{2}}\frac{\cos \frac{\left( 2\alpha +\gamma -\mu \right) \pi }{2}}{\cos 
\frac{\left( \gamma +\mu \right) \pi }{2}},
\end{equation*}%
representing the narrowed thermodynamical restriction (\ref{STD-IDD.IDD-3})$%
_{3}$, since $\gamma +\mu \leqslant 2\alpha +\gamma -\mu <1$ reducing to $%
\mu \leqslant \alpha ,$ see (\ref{STD-IDD.IDD-1})$_{3}$, implies $\frac{\cos 
\frac{\left( 2\alpha +\gamma -\mu \right) \pi }{2}}{\cos \frac{\left( \gamma
+\mu \right) \pi }{2}}\leqslant 1.$

In conclusion, the non-negativity of function $K,$ corresponding to the
model IDD.IDD and given by (\ref{K-IDD.IDD}), is ensured if $2\alpha +\beta
-\mu <2\alpha +\gamma -\mu <1$, see (\ref{STD-IDD.IDD-1})$_{2}$, by
requesting (\ref{STD-IDD.IDD-2}) and (\ref{STD-IDD.IDD-3}) in addition to (%
\ref{STD-IDD.IDD-1}) implying the corresponding relaxation modulus to be a
completely monotonic function and creep compliance a Bernstein function,
while if either $2\alpha +\gamma -\mu \geqslant 1$ or $2\alpha +\gamma -\mu
>2\alpha +\beta -\mu \geqslant 1,$ then the mentioned properties cannot be
ensured, since either (\ref{NN-of-K-IDD.IDD-3}), or both (\ref%
{NN-of-K-IDD.IDD-2}) and (\ref{NN-of-K-IDD.IDD-3}) cannot be satisfied.

\subsubsection{Model IID.IDD}

Function $K$, given by (\ref{K}), in the case of the model IID.IDD, having
the constitutive equation given in the form%
\begin{equation*}
\left( a_{1}\,_{0}\mathrm{I}_{t}^{\alpha }+a_{2}\,_{0}\mathrm{I}_{t}^{\beta
}+a_{3}\,_{0}\mathrm{D}_{t}^{\gamma }\right) \sigma \left( t\right) =\left(
b_{1}\,_{0}\mathrm{I}_{t}^{\mu }+b_{2}\,_{0}\mathrm{D}_{t}^{\nu }+b_{3}\,_{0}%
\mathrm{D}_{t}^{\alpha +\gamma -\mu }\right) \varepsilon \left( t\right) ,
\end{equation*}%
and with corresponding thermodynamical, see (\ref%
{TD-frale-less-less-less-less-1}) and (\ref{TD-frale-less-less-less-less-2}%
), as well as with the narrowed thermodynamical requirements%
\begin{gather}
\mu \leqslant \beta <\alpha ,\quad \gamma \leqslant \nu ,\quad \alpha +\beta
+\gamma \leqslant 1+\mu ,\quad \mu +\nu -\gamma <\alpha \leqslant 1-\nu ,
\label{STD-IID.IDD-1} \\
0\leqslant \left\{ 
\begin{tabular}{l}
$\alpha -\beta -\gamma -\mu $ \smallskip \\ 
$\alpha -2\mu -\nu $ \smallskip%
\end{tabular}%
\right\} \leqslant 2\alpha -\beta -2\mu -\nu \leqslant \left\{ 
\begin{tabular}{l}
$2\alpha -\beta -\mu $ \smallskip \\ 
$2\alpha +\gamma -2\mu -\nu $ \smallskip%
\end{tabular}%
\right\} \leqslant 2\alpha +\gamma -\mu <1,  \label{STD-IID.IDD-2} \\
-\frac{b_{3}}{b_{1}}\frac{\cos \frac{\left( 2\alpha +\gamma -\mu \right) \pi 
}{2}}{\cos \frac{\left( \gamma +\mu \right) \pi }{2}}\leqslant \frac{a_{3}}{%
a_{1}}\leqslant \frac{b_{3}}{b_{1}}\frac{\sin \frac{\left( 2\alpha +\gamma
-\mu \right) \pi }{2}}{\sin \frac{\left( \gamma +\mu \right) \pi }{2}}\frac{%
\cos \frac{\left( 2\alpha +\gamma -\mu \right) \pi }{2}}{\cos \frac{\left(
\gamma +\mu \right) \pi }{2}}\leqslant \frac{b_{3}}{b_{1}}\frac{\sin \frac{%
\left( 2\alpha +\gamma -\mu \right) \pi }{2}}{\sin \frac{\left( \gamma +\mu
\right) \pi }{2}},  \label{STD-IID.IDD-3}
\end{gather}%
becomes%
\begin{align}
K\left( \rho \right) & =a_{1}b_{1}\sin \left( \left( \alpha -\mu \right) \pi
\right) +a_{1}b_{2}\rho ^{\alpha -\beta }\sin \left( \left( 2\alpha -\beta
-\mu \right) \pi \right)  \notag \\
& \quad +a_{1}b_{3}\rho ^{\alpha +\gamma }\sin \left( \left( 2\alpha +\gamma
-\mu \right) \pi \right) +a_{2}b_{1}\rho ^{\mu +\nu }\sin \left( \left(
\alpha -2\mu -\nu \right) \pi \right)  \notag \\
& \quad +a_{2}b_{2}\rho ^{\alpha -\beta +\mu +\nu }\sin \left( \left(
2\alpha -\beta -2\mu -\nu \right) \pi \right) +a_{2}b_{3}\rho ^{\alpha
+\gamma +\mu +\nu }\sin \left( \left( 2\alpha +\gamma -2\mu -\nu \right) \pi
\right)  \notag \\
& \quad -a_{3}b_{1}\rho ^{\alpha +\gamma }\sin \left( \left( \mu +\gamma
\right) \pi \right) +a_{3}b_{2}\rho ^{2\alpha -\beta +\gamma }\sin \left(
\left( \alpha -\beta -\gamma -\mu \right) \pi \right) +a_{3}b_{3}\rho
^{2\left( \alpha +\gamma \right) }\sin \left( \left( \alpha -\mu \right) \pi
\right) ,  \label{K-IID.IDD}
\end{align}%
so that the non-negativity requirement of function $K$ implies%
\begin{gather}
a_{1}b_{3}\sin \left( \left( 2\alpha +\gamma -\mu \right) \pi \right)
-a_{3}b_{1}\sin \left( \left( \gamma +\mu \right) \pi \right) \geqslant 0,
\label{NN-of-K-IID.IDD-1} \\
\sin \left( \left( 2\alpha -\beta -\mu \right) \pi \right) \geqslant 0,\quad
\sin \left( \left( \alpha -2\mu -\nu \right) \pi \right) \geqslant 0,
\label{NN-of-K-IID.IDD-2} \\
\sin \left( \left( 2\alpha -\beta -2\mu -\nu \right) \pi \right) \geqslant
0,\quad \sin \left( \left( 2\alpha +\gamma -2\mu -\nu \right) \pi \right)
\geqslant 0,\quad \sin \left( \left( \alpha -\beta -\gamma -\mu \right) \pi
\right) \geqslant 0.  \label{NN-of-K-IID.IDD-3}
\end{gather}%
The notation used in the narrowed thermodynamical requirement (\ref%
{STD-IID.IDD-2})$_{1,2}$ means that both $\alpha -\beta -\gamma -\mu $ and $%
\alpha -2\mu -\nu $ are non-negative and not greater than $2\alpha -\beta
-2\mu -\nu ,$ while the relation between $\alpha -\beta -\gamma -\mu $ and $%
\alpha -2\mu -\nu $ cannot be established and the same interpretation holds
for the notation used in (\ref{STD-IID.IDD-2})$_{3,4}$.

The first term in (\ref{NN-of-K-IID.IDD-1}) has either a positive value if $%
2\alpha +\gamma -\mu <1,$ or a non-positive value if $2\alpha +\gamma -\mu
\geqslant 1,$ since $2\alpha +\gamma -\mu =\left( \alpha +\gamma \right)
+\left( \alpha -\mu \right) \in \left( 0,2\right) ,$ due to $\alpha +\gamma
\leqslant \alpha +\nu \leqslant 1,$ according to the thermodynamical
requirements (\ref{STD-IID.IDD-1})$_{1,2,4}$, while the second term in (\ref%
{NN-of-K-IID.IDD-1}) has a non-positive value, since $\gamma +\mu \leqslant
\alpha +\nu \leqslant 1$, according to the thermodynamical requirements (\ref%
{STD-IID.IDD-1})$_{1,2,4},$ so that if $2\alpha +\gamma -\mu <1,$ then the
requirement (\ref{NN-of-K-IID.IDD-1}) is transformed into%
\begin{equation*}
\frac{a_{3}}{a_{1}}\leqslant \frac{b_{3}}{b_{1}}\frac{\sin \frac{\left(
2\alpha +\gamma -\mu \right) \pi }{2}}{\sin \frac{\left( \gamma +\mu \right)
\pi }{2}}\frac{\cos \frac{\left( 2\alpha +\gamma -\mu \right) \pi }{2}}{\cos 
\frac{\left( \gamma +\mu \right) \pi }{2}},
\end{equation*}%
representing the narrowed thermodynamical restriction (\ref{STD-IID.IDD-3})$%
_{3}$, since $\gamma +\mu \leqslant 2\alpha +\gamma -\mu <1$ reducing to $%
\mu \leqslant \alpha ,$ see (\ref{STD-IID.IDD-1})$_{1}$, implies $\frac{\cos 
\frac{\left( 2\alpha +\gamma -\mu \right) \pi }{2}}{\cos \frac{\left( \gamma
+\mu \right) \pi }{2}}\leqslant 1.$

The chain of inequalities 
\begin{equation*}
\left\{ 
\begin{tabular}{l}
$\alpha -\beta -\gamma -\mu $ \smallskip \\ 
$\alpha -2\mu -\nu $ \smallskip%
\end{tabular}%
\right\} \leqslant 2\alpha -\beta -2\mu -\nu \leqslant \left\{ 
\begin{tabular}{l}
$2\alpha -\beta -\mu $ \smallskip \\ 
$2\alpha +\gamma -2\mu -\nu $ \smallskip%
\end{tabular}%
\right\} \leqslant 2\alpha +\gamma -\mu ,
\end{equation*}%
holds since inequality $\alpha -\beta -\gamma -\mu \leqslant 2\alpha -\beta
-2\mu -\nu $ reduces to the thermodynamical requirement (\ref{STD-IID.IDD-1})%
$_{4},$ inequality $\alpha -2\mu -\nu \leqslant 2\alpha -\beta -2\mu -\nu $
reduces to the requirement (\ref{STD-IID.IDD-1})$_{1},$ while all other
inequalities are trivially satisfied. On the other hand, neither can be
stated that $\alpha -2\mu -\nu \leqslant \alpha -\beta -\gamma -\mu ,$ nor
that $2\alpha +\gamma -2\mu -\nu \leqslant 2\alpha -\beta -\mu ,$ since both
inequalities reduce to $\beta +\gamma \leqslant \mu +\nu ,$ while the
thermodynamical requirements (\ref{STD-IID.IDD-1})$_{1,2}$ imply $\beta
+\gamma \leqslant \alpha +\nu $.

In addition to the request $2\alpha +\gamma -\mu <1,$ that guarantees
non-negativity of (\ref{NN-of-K-IID.IDD-1}), in order to satisfy
inequalities (\ref{NN-of-K-IID.IDD-2}) and (\ref{NN-of-K-IID.IDD-3}), one
requests that $\alpha -\beta -\gamma -\mu \geqslant 0$ as well as that $%
\alpha -2\mu -\nu \geqslant 0,$ since both $\alpha -\beta -\gamma -\mu
=\alpha -\beta -\left( \gamma +\mu \right) \geqslant -1+\alpha -\beta
\geqslant -1$ and $\alpha -2\mu -\nu =\alpha -\mu -\left( \mu +\nu \right)
\geqslant -1+\alpha -\mu \geqslant -1,$ according to the thermodynamical
requirement (\ref{STD-IID.IDD-1})$_{1,2,4}.$

In conclusion, the non-negativity of function $K,$ corresponding to the
model IID.IDD and given by (\ref{K-IID.IDD}), is ensured if the request (\ref%
{STD-IID.IDD-2}) on the orders of fractional integrals and derivatives,
along with the thermodynamical restriction (\ref{STD-IID.IDD-1}), is
satisfied in addition to the request (\ref{STD-IID.IDD-3}) on the model
parameters, implying the corresponding relaxation modulus to be a completely
monotonic function and creep compliance a Bernstein function, while if the
request (\ref{STD-IID.IDD-2}) is violated$,$ then the mentioned properties
cannot be ensured, since some or all of the conditions in (\ref%
{NN-of-K-IID.IDD-1}), (\ref{NN-of-K-IID.IDD-2}), and (\ref{NN-of-K-IID.IDD-3}%
) may not satisfied.

\subsubsection{Model I$^{{}^{+}}$ID.I$^{{}^{+}}$ID}

Function $K$, given by (\ref{K}), in the case of model I$^{{}^{+}}$ID.I$%
^{{}^{+}}$ID, having the constitutive equation given in the form%
\begin{equation*}
\left( a_{1}\,_{0}\mathrm{I}_{t}^{1+\alpha }+a_{2}\,_{0}\mathrm{I}_{t}^{%
\frac{1+\alpha -\gamma }{2}}+a_{3}\,_{0}\mathrm{D}_{t}^{\gamma }\right)
\sigma \left( t\right) =\left( b_{1}\,_{0}\mathrm{I}_{t}^{1+\mu }+b_{2}\,_{0}%
\mathrm{I}_{t}^{\frac{1+\mu -\left( \alpha +\gamma -\mu \right) }{2}%
}+b_{3}\,_{0}\mathrm{D}_{t}^{\alpha +\gamma -\mu }\right) \varepsilon \left(
t\right) ,
\end{equation*}%
and with corresponding thermodynamical, see (\ref{franjo-over-over-TD-1}) - (%
\ref{franjo-over-over-TD-4}), as well as with the narrowed thermodynamical
requirements%
\begin{gather}
\mu \leqslant \alpha ,\quad \alpha +\gamma +2\left( \alpha -\mu \right)
=3\alpha +\gamma -2\mu \leqslant 1,  \label{STD-I+ID.I+ID-1} \\
\frac{a_{2}}{a_{1}}\leqslant \frac{b_{2}}{b_{1}}\frac{\cos \frac{\left(
1-3\alpha -\gamma +2\mu \right) \pi }{4}}{\cos \frac{\left( 1+\alpha -\gamma
-2\mu \right) \pi }{4}}\frac{\sin \frac{\left( 1-3\alpha -\gamma +2\mu
\right) \pi }{4}}{\sin \frac{\left( 1+\alpha -\gamma -2\mu \right) \pi }{4}}%
\leqslant \frac{b_{2}}{b_{1}}\frac{\cos \frac{\left( 1-3\alpha -\gamma +2\mu
\right) \pi }{4}}{\cos \frac{\left( 1+\alpha -\gamma -2\mu \right) \pi }{4}},
\label{STD-I+ID.I+ID-2.1} \\
\frac{a_{3}}{a_{2}}\leqslant \frac{b_{3}}{b_{2}}\frac{\sin \frac{\left(
1+3\alpha +\gamma -2\mu \right) \pi }{4}}{\sin \frac{\left( 1-\alpha +\gamma
+2\mu \right) \pi }{4}}\frac{\cos \frac{\left( 1+3\alpha +\gamma -2\mu
\right) \pi }{4}}{\cos \frac{\left( 1-\alpha +\gamma +2\mu \right) \pi }{4}}%
\leqslant \frac{b_{3}}{b_{2}}\frac{\sin \frac{\left( 1+3\alpha +\gamma -2\mu
\right) \pi }{4}}{\sin \frac{\left( 1-\alpha +\gamma +2\mu \right) \pi }{4}},
\label{STD-I+ID.I+ID-2.2} \\
a_{3}b_{1}\cos \frac{\left( \gamma +\mu \right) \pi }{2}\leqslant
a_{2}b_{2}\sin \frac{\left( \alpha -\mu \right) \pi }{2}+a_{1}b_{3}\cos 
\frac{\left( 2\alpha +\gamma -\mu \right) \pi }{2},  \label{STD-I+ID.I+ID-3}
\\
a_{1}b_{3}\sin \frac{\left( 2\alpha +\gamma -\mu \right) \pi }{2}\leqslant
a_{2}b_{2}\cos \frac{\left( \alpha -\mu \right) \pi }{2}-a_{3}b_{1}\sin 
\frac{\left( \gamma +\mu \right) \pi }{2},  \label{STD-I+ID.I+ID-7} \\
a_{1}b_{3}\sin \frac{\left( 2\alpha +\gamma -\mu \right) \pi }{2}\leqslant
a_{2}b_{2}\cos \frac{\left( \alpha -\mu \right) \pi }{2}\frac{\sin \frac{%
\left( \alpha -\mu \right) \pi }{2}}{\cos \frac{\left( 2\alpha +\gamma -\mu
\right) \pi }{2}}+a_{3}b_{1}\sin \frac{\left( \gamma +\mu \right) \pi }{2}%
\frac{\cos \frac{\left( \gamma +\mu \right) \pi }{2}}{\cos \frac{\left(
2\alpha +\gamma -\mu \right) \pi }{2}},  \label{STD-I+ID.I+ID-8}
\end{gather}%
becomes%
\begin{align}
K\left( \rho \right) &=a_{1}b_{1}\sin \left( \left( \alpha -\mu \right) \pi
\right) +a_{1}b_{2}\rho ^{\frac{1+\alpha +\gamma }{2}}\sin \frac{\left(
1+3\alpha +\gamma -2\mu \right) \pi }{2}+a_{1}b_{3}\rho ^{1+\alpha +\gamma
}\sin \left( \left( 1+2\alpha +\gamma -\mu \right) \pi \right)  \notag \\
&\quad-a_{2}b_{1}\rho ^{\frac{1+\alpha +\gamma }{2}}\sin \frac{\left(
1-\alpha +\gamma +2\mu \right) \pi }{2}+a_{2}b_{2}\rho ^{1+\alpha +\gamma
}\sin \left( \left( \alpha -\mu \right) \pi \right)  \notag \\
&\quad+a_{2}b_{3}\rho ^{3\frac{1+\alpha +\gamma }{2}}\sin \frac{\left(
1+3\alpha +\gamma -2\mu \right) \pi }{2}-a_{3}b_{1}\rho ^{1+\alpha +\gamma
}\sin \left( \left( 1+\gamma +\mu \right) \pi \right)  \notag \\
&\quad-a_{3}b_{2}\rho ^{3\frac{1+\alpha +\gamma }{2}}\sin \frac{\left(
1-\alpha +\gamma +2\mu \right) \pi }{2}+a_{3}b_{3}\rho ^{2\left( 1+\alpha
+\gamma \right) }\sin \left( \left( \alpha -\mu \right) \pi \right) ,
\label{K-I+ID.I+ID}
\end{align}%
so that the non-negativity requirement of function $K$ implies%
\begin{gather}
a_{1}b_{2}\sin \frac{\left( 1+3\alpha +\gamma -2\mu \right) \pi }{2}%
-a_{2}b_{1}\sin \frac{\left( 1-\alpha +\gamma +2\mu \right) \pi }{2}%
\geqslant 0,  \label{NN-of-K-I+ID.I+ID-1} \\
a_{2}b_{3}\sin \frac{\left( 1+3\alpha +\gamma -2\mu \right) \pi }{2}%
-a_{3}b_{2}\sin \frac{\left( 1-\alpha +\gamma +2\mu \right) \pi }{2}%
\geqslant 0,  \label{NN-of-K-I+ID.I+ID-2} \\
-a_{1}b_{3}\sin \left( \left( 2\alpha +\gamma -\mu \right) \pi \right)
+a_{2}b_{2}\sin \left( \left( \alpha -\mu \right) \pi \right)
+a_{3}b_{1}\sin \left( \left( \gamma +\mu \right) \pi \right) \geqslant 0.
\label{NN-of-K-I+ID.I+ID-3}
\end{gather}

The sines appearing in the first term in (\ref{NN-of-K-I+ID.I+ID-1}) and (%
\ref{NN-of-K-I+ID.I+ID-2}) have a non-negative value, since $3\alpha +\gamma
-2\mu \in \left( 0,1\right) ,$ i.e., $1+3\alpha +\gamma -2\mu \in \left(
1,2\right) ,$ according to the thermodynamical requirement (\ref%
{STD-I+ID.I+ID-1})$_{2},$ as well as the sines appearing in the second term
in (\ref{NN-of-K-I+ID.I+ID-1}) and in (\ref{NN-of-K-I+ID.I+ID-2}), since $%
1-\alpha +\gamma +2\mu \in \left( 0,2\right) $, because of $-\alpha +\gamma
+2\mu \in \left( -1,1\right) $, that is valid since $-\alpha +\gamma +2\mu
=\left( \alpha +\gamma \right) -2\left( \alpha -\mu \right) \leqslant
1-4\left( \alpha -\mu \right) \leqslant 1$ according to the thermodynamical
restrictions (\ref{STD-I+ID.I+ID-1}), as well as since the thermodynamical
requirement (\ref{STD-I+ID.I+ID-1})$_{2}$ transforms into $-\alpha +\gamma
+2\mu -2\left( \alpha +\gamma \right) =-\left( 3\alpha +\gamma -2\mu \right)
\geqslant -1\,$implying $-\alpha +\gamma +2\mu \geqslant -1+2\left( \alpha
+\gamma \right) \geqslant -1,$ so that $1-\alpha +\gamma +2\mu \geqslant
2\left( \alpha +\gamma \right) \geqslant 0$. Therefore, the requirements (%
\ref{NN-of-K-I+ID.I+ID-1}) and (\ref{NN-of-K-I+ID.I+ID-2}) respectively
transform into%
\begin{equation}
\frac{a_{2}}{a_{1}}\leqslant \frac{b_{2}}{b_{1}}\frac{\cos \frac{\left(
1-3\alpha -\gamma +2\mu \right) \pi }{4}}{\cos \frac{\left( 1+\alpha -\gamma
-2\mu \right) \pi }{4}}\frac{\sin \frac{\left( 1-3\alpha -\gamma +2\mu
\right) \pi }{4}}{\sin \frac{\left( 1+\alpha -\gamma -2\mu \right) \pi }{4}}%
\quad \text{and}\quad \frac{a_{3}}{a_{2}}\leqslant \frac{b_{3}}{b_{2}}\frac{%
\sin \frac{\left( 1+3\alpha +\gamma -2\mu \right) \pi }{4}}{\sin \frac{%
\left( 1-\alpha +\gamma +2\mu \right) \pi }{4}}\frac{\cos \frac{\left(
1+3\alpha +\gamma -2\mu \right) \pi }{4}}{\cos \frac{\left( 1-\alpha +\gamma
+2\mu \right) \pi }{4}},  \label{komplikovano-NN-of-K-I+ID.I+ID-1}
\end{equation}%
where the requirement (\ref{NN-of-K-I+ID.I+ID-1}) is rewritten as%
\begin{equation*}
a_{1}b_{2}\sin \frac{\left( 1-3\alpha -\gamma +2\mu \right) \pi }{2}%
-a_{2}b_{1}\sin \frac{\left( 1+\alpha -\gamma -2\mu \right) \pi }{2}%
\geqslant 0,
\end{equation*}%
since $\sin \left( \frac{\pi }{2}+\phi \right) =\sin \left( \frac{\pi }{2}%
-\phi \right) $ and then transformed into (\ref%
{komplikovano-NN-of-K-I+ID.I+ID-1})$_{1}$ in order to be combined with the
thermodynamical restriction (\ref{franjo-over-over-TD-2})$_{1}$. Expressions
given by (\ref{komplikovano-NN-of-K-I+ID.I+ID-1}) represent the narrowed
thermodynamical restrictions (\ref{STD-I+ID.I+ID-2.1})$_{2}$ and (\ref%
{STD-I+ID.I+ID-2.2})$_{2},$ since $\frac{\sin \frac{\left( 1-3\alpha -\gamma
+2\mu \right) \pi }{4}}{\sin \frac{\left( 1+\alpha -\gamma -2\mu \right) \pi 
}{4}}\leqslant 1$ and $\frac{\cos \frac{\left( 1+3\alpha +\gamma -2\mu
\right) \pi }{4}}{\cos \frac{\left( 1-\alpha +\gamma +2\mu \right) \pi }{4}}%
\leqslant 1$, due to $1-3\alpha -\gamma +2\mu \leqslant 1+\alpha -\gamma
-2\mu $ and $1-\alpha +\gamma +2\mu \leqslant 1+3\alpha +\gamma -2\mu ,$
both reducing to $\mu \leqslant \alpha $, see (\ref{STD-I+ID.I+ID-1})$_{1}$.

The requirement (\ref{STD-I+ID.I+ID-1})$_{2},$ reading $3\alpha +\gamma
-2\mu =2\alpha +\gamma -\mu +\left( \alpha -\mu \right) \leqslant 1$, i.e., $%
2\alpha +\gamma -\mu \leqslant 1-\left( \alpha -\mu \right) \leqslant 1$,
combined with the requirement (\ref{STD-I+ID.I+ID-1})$_{1}$ guarantees a
non-negative value of the sine in the first term in (\ref%
{NN-of-K-I+ID.I+ID-3}),\ while the requirement (\ref{STD-I+ID.I+ID-1})$_{1}$
ensures non-negativity of the sine in the second term in (\ref%
{NN-of-K-I+ID.I+ID-3}), and finally non-negativity of the sine in the third
term in (\ref{NN-of-K-I+ID.I+ID-3}) is guaranteed by the requirement (\ref%
{STD-I+ID.I+ID-1})$_{2}$, so that the requirement (\ref{NN-of-K-I+ID.I+ID-3}%
) is transformed into the narrowed thermodynamical restriction (\ref%
{STD-I+ID.I+ID-8}) if the right-hand-side of (\ref{STD-I+ID.I+ID-8}) is less
than the right-hand-side of (\ref{STD-I+ID.I+ID-7}), since $\frac{\sin \frac{%
\left( \alpha -\mu \right) \pi }{2}}{\cos \frac{\left( 2\alpha +\gamma -\mu
\right) \pi }{2}}=\frac{\cos \frac{\left( 1-\left( \alpha -\mu \right)
\right) \pi }{2}}{\cos \frac{\left( 2\alpha +\gamma -\mu \right) \pi }{2}}%
\leqslant 1$ and $\frac{\cos \frac{\left( \gamma +\mu \right) \pi }{2}}{\cos 
\frac{\left( 2\alpha +\gamma -\mu \right) \pi }{2}}\geqslant 1,$ due to $%
2\alpha +\gamma -\mu \leqslant 1-\left( \alpha -\mu \right) <1$ and $\gamma
+\mu \leqslant 2\alpha +\gamma -\mu <1$, respectively reducing to $2\alpha
+\gamma -\mu \leqslant 1$ and $\mu \leqslant \alpha ,$ see (\ref%
{STD-I+ID.I+ID-1}).

In conclusion, the non-negativity of function $K,$ corresponding to the
model I$^{{}^{+}}$ID.I$^{{}^{+}}$ID and given by (\ref{K-I+ID.I+ID}), is
ensured by requesting (\ref{STD-I+ID.I+ID-2.1}), (\ref{STD-I+ID.I+ID-2.2})$,$
and possibly (\ref{STD-I+ID.I+ID-8}) in addition to (\ref{STD-I+ID.I+ID-1}),
(\ref{STD-I+ID.I+ID-3}), and (\ref{STD-I+ID.I+ID-7}) implying the
corresponding relaxation modulus to be a completely monotonic function and
creep compliance a Bernstein function.

\subsubsection{Model IDD$^{{}^{+}}$.IDD$^{{}^{+}}$}

Function $K$, given by (\ref{K}), in the case of the model IDD$^{{}^{+}}$.IDD%
$^{{}^{+}}$, having the constitutive equation given in the form%
\begin{equation*}
\left( a_{1}\,_{0}\mathrm{I}_{t}^{\alpha }+a_{2}\,_{0}\mathrm{D}_{t}^{\frac{%
1+\gamma -\alpha }{2}}+a_{3}\,_{0}\mathrm{D}_{t}^{1+\gamma }\right) \sigma
\left( t\right) =\left( b_{1}\,_{0}\mathrm{I}_{t}^{\alpha +\gamma -\eta
}+b_{2}\,_{0}\mathrm{D}_{t}^{\frac{1+\eta -\left( \alpha +\gamma -\eta
\right) }{2}}+b_{3}\,_{0}\mathrm{D}_{t}^{1+\eta }\right) \varepsilon \left(
t\right) ,
\end{equation*}%
and with corresponding thermodynamical, see (\ref{TD-IDD-IDD-over-over-1}) -
(\ref{TD-IDD-IDD-over-over-5}), as well as with the narrowed thermodynamical
requirements%
\begin{gather}
0\leqslant \alpha +\gamma -\eta \leqslant 1,\quad \gamma \leqslant \eta
,\quad \alpha +\eta +\left( \eta -\gamma \right) =\alpha -\gamma +2\eta
\leqslant 1,  \label{STD-IDD+.IDD+-1} \\
\frac{a_{2}}{a_{1}}\leqslant \frac{b_{2}}{b_{1}}\frac{\sin \frac{\left(
1+\alpha -\gamma +2\eta \right) \pi }{4}}{\sin \frac{\left( 1+\alpha
+3\gamma -2\eta \right) \pi }{4}}\frac{\cos \frac{\left( 1+\alpha -\gamma
+2\eta \right) \pi }{4}}{\cos \frac{\left( 1+\alpha +3\gamma -2\eta \right)
\pi }{4}}\leqslant \frac{b_{2}}{b_{1}}\frac{\sin \frac{\left( 1+\alpha
-\gamma +2\eta \right) \pi }{4}}{\sin \frac{\left( 1+\alpha +3\gamma -2\eta
\right) \pi }{4}},  \label{STD-IDD+.IDD+-2} \\
\frac{a_{3}}{a_{2}}\leqslant \frac{b_{3}}{b_{2}}\frac{\cos \frac{\left(
1-\alpha +\gamma -2\eta \right) \pi }{4}}{\cos \frac{\left( 1-\alpha
-3\gamma +2\eta \right) \pi }{4}}\frac{\sin \frac{\left( 1-\alpha +\gamma
-2\eta \right) \pi }{4}}{\sin \frac{\left( 1-\alpha -3\gamma +2\eta \right)
\pi }{4}}\leqslant \frac{b_{3}}{b_{2}}\frac{\cos \frac{\left( 1-\alpha
+\gamma -2\eta \right) \pi }{4}}{\cos \frac{\left( 1-\alpha -3\gamma +2\eta
\right) \pi }{4}},  \label{STD-IDD+.IDD+-2a} \\
a_{3}b_{1}\cos \frac{\left( \alpha +2\gamma -\eta \right) \pi }{2}%
-a_{2}b_{2}\sin \frac{\left( \eta -\gamma \right) \pi }{2}\leqslant
a_{1}b_{3}\cos \frac{\left( \alpha +\eta \right) \pi }{2},
\label{STD-IDD+.IDD+-3} \\
a_{1}b_{3}\sin \frac{\left( \alpha +\eta \right) \pi }{2}\leqslant
a_{2}b_{2}\cos \frac{\left( \eta -\gamma \right) \pi }{2}-a_{3}b_{1}\sin 
\frac{\left( \alpha +2\gamma -\eta \right) \pi }{2},  \label{STD-IDD+.IDD+-4}
\\
a_{1}b_{3}\sin \frac{\left( \alpha +\eta \right) \pi }{2}\leqslant
a_{2}b_{2}\cos \frac{\left( \eta -\gamma \right) \pi }{2}\frac{\sin \frac{%
\left( \eta -\gamma \right) \pi }{2}}{\cos \frac{\left( \alpha +\eta \right)
\pi }{2}}+a_{3}b_{1}\sin \frac{\left( \alpha +2\gamma -\eta \right) \pi }{2}%
\frac{\cos \frac{\left( \alpha +2\gamma -\eta \right) \pi }{2}}{\cos \frac{%
\left( \alpha +\eta \right) \pi }{2}},  \label{STD-IDD+.IDD+-5}
\end{gather}%
becomes%
\begin{align}
K\left( \rho \right) & =a_{1}b_{1}\sin \left( \left( \eta -\gamma \right)
\pi \right) +a_{1}b_{2}\rho ^{\frac{1+\alpha +\gamma }{2}}\sin \frac{\left(
1+\alpha -\gamma +2\eta \right) \pi }{2}-a_{1}b_{3}\rho ^{1+\alpha +\gamma
}\sin \left( \left( \alpha +\eta \right) \pi \right)  \notag \\
& \quad -a_{2}b_{1}\rho ^{\frac{1+\alpha +\gamma }{2}}\sin \frac{\left(
1+\alpha +3\gamma -2\eta \right) \pi }{2}+a_{2}b_{2}\rho ^{1+\alpha +\gamma
}\sin \left( \left( \eta -\gamma \right) \pi \right)  \notag \\
& \quad +a_{2}b_{3}\rho ^{3\frac{1+\alpha +\gamma }{2}}\sin \frac{\left(
1+\alpha +2\eta -\gamma \right) \pi }{2}+a_{3}b_{1}\rho ^{1+\alpha +\gamma
}\sin \left( \left( \alpha +2\gamma -\eta \right) \pi \right)  \notag \\
& \quad -a_{3}b_{2}\rho ^{3\frac{1+\alpha +\gamma }{2}}\sin \frac{\left(
1+\alpha +3\gamma -2\eta \right) \pi }{2}+a_{3}b_{3}\rho ^{2\left( 1+\alpha
+\gamma \right) }\sin \left( \left( \eta -\gamma \right) \pi \right) ,
\label{K-IDD+.IDD+}
\end{align}%
so that the non-negativity requirement of function $K$ implies%
\begin{gather}
a_{1}b_{2}\sin \frac{\left( 1+\alpha -\gamma +2\eta \right) \pi }{2}%
-a_{2}b_{1}\sin \frac{\left( 1+\alpha +3\gamma -2\eta \right) \pi }{2}%
\geqslant 0  \label{NN-of-K-IDD+.IDD+-1} \\
a_{2}b_{3}\sin \frac{\left( 1+\alpha -\gamma +2\eta \right) \pi }{2}%
-a_{3}b_{2}\sin \frac{\left( 1+\alpha +3\gamma -2\eta \right) \pi }{2}%
\geqslant 0  \label{NN-of-K-IDD+.IDD+-2} \\
-a_{1}b_{3}\sin \left( \left( \alpha +\eta \right) \pi \right)
+a_{2}b_{2}\sin \left( \left( \eta -\gamma \right) \pi \right)
+a_{3}b_{1}\sin \left( \left( \alpha +2\gamma -\eta \right) \pi \right)
\geqslant 0.  \label{NN-of-K-IDD+.IDD+-3}
\end{gather}

According to the thermodynamical requirement (\ref{STD-IDD+.IDD+-1})$_{3},$
the first term in (\ref{NN-of-K-IDD+.IDD+-1}), as well as in (\ref%
{NN-of-K-IDD+.IDD+-2}) has a non-negative value, since $1+\alpha -\gamma
+2\eta \in \left( 1,2\right) $, while the second term in (\ref%
{NN-of-K-IDD+.IDD+-1}), as well as in (\ref{NN-of-K-IDD+.IDD+-2}) is non-
positive, since one has $1+\alpha +3\gamma -2\eta \in \left( 0,2\right) $,
due to $1+\alpha +3\gamma -2\eta =1-\left( \alpha -\gamma +2\eta \right)
+2\left( \alpha +\gamma \right) \geqslant 2\left( \alpha +\gamma \right)
\geqslant 0,$\ where the restriction (\ref{STD-IDD+.IDD+-1})$_{3}$ is used
in the form $-\left( \alpha -\gamma +2\eta \right) \geqslant -1$, as well as 
$1+\alpha +3\gamma -2\eta =1+\left( \alpha -\gamma +2\eta \right) -4\left(
\eta -\gamma \right) \leqslant 2-4\left( \eta -\gamma \right) \leqslant 2$,
with the restriction (\ref{STD-IDD+.IDD+-1})$_{3}$ used again, so that the
requirements (\ref{NN-of-K-IDD+.IDD+-1}) and (\ref{NN-of-K-IDD+.IDD+-2}) are
transformed into%
\begin{equation}
\frac{a_{2}}{a_{1}}\leqslant \frac{b_{2}}{b_{1}}\frac{\sin \frac{\left(
1+\alpha -\gamma +2\eta \right) \pi }{4}}{\sin \frac{\left( 1+\alpha
+3\gamma -2\eta \right) \pi }{4}}\frac{\cos \frac{\left( 1+\alpha -\gamma
+2\eta \right) \pi }{4}}{\cos \frac{\left( 1+\alpha +3\gamma -2\eta \right)
\pi }{4}}\quad \text{and}\quad \frac{a_{3}}{a_{2}}\leqslant \frac{b_{3}}{%
b_{2}}\frac{\cos \frac{\left( 1-\alpha +\gamma -2\eta \right) \pi }{4}}{\cos 
\frac{\left( 1-\alpha -3\gamma +2\eta \right) \pi }{4}}\frac{\sin \frac{%
\left( 1-\alpha +\gamma -2\eta \right) \pi }{4}}{\sin \frac{\left( 1-\alpha
-3\gamma +2\eta \right) \pi }{4}},  \label{komplikovano-NN-of-K-IDD+.IDD+-1}
\end{equation}%
and represent narrowed thermodynamical requirements (\ref{STD-IDD+.IDD+-2})
and (\ref{STD-IDD+.IDD+-2a}),\ where the requirement (\ref%
{NN-of-K-IDD+.IDD+-2}) becomes 
\begin{equation*}
a_{2}b_{3}\sin \frac{\left( 1-\alpha +\gamma -2\eta \right) \pi }{2}%
-a_{3}b_{2}\sin \frac{\left( 1-\alpha -3\gamma +2\eta \right) \pi }{2}%
\geqslant 0,
\end{equation*}%
since $\sin \left( \frac{\pi }{2}+\phi \right) =\sin \left( \frac{\pi }{2}%
-\phi \right) $ and then transformed into (\ref%
{komplikovano-NN-of-K-IDD+.IDD+-1})$_{2}$ in order to be combined with the
thermodynamical restriction (\ref{TD-IDD-IDD-over-over-3})$_{2}$.
Expressions given by (\ref{komplikovano-NN-of-K-IDD+.IDD+-1}) represent the
narrowed thermodynamical restrictions (\ref{STD-IDD+.IDD+-2})$_{2}$ and (\ref%
{STD-IDD+.IDD+-2a})$_{2}$, since $\frac{\cos \frac{\left( 1+\alpha -\gamma
+2\eta \right) \pi }{4}}{\cos \frac{\left( 1+\alpha +3\gamma -2\eta \right)
\pi }{4}}\leqslant 1$ and $\frac{\sin \frac{\left( 1-\alpha +\gamma -2\eta
\right) \pi }{4}}{\sin \frac{\left( 1-\alpha -3\gamma +2\eta \right) \pi }{4}%
}\leqslant 1$, due to $1+\alpha +3\gamma -2\eta \leqslant 1+\alpha -\gamma
+2\eta $ and $1-\alpha +\gamma -2\eta \leqslant 1-\alpha -3\gamma +2\eta ,$
both reducing to $\gamma \leqslant \eta $, see (\ref{STD-IDD+.IDD+-1})$_{2}$.

The requirements (\ref{STD-IDD+.IDD+-1}) ensure non-negative values of sines
in the first and third term in (\ref{NN-of-K-IDD+.IDD+-3}), since $\alpha
+\eta \leqslant 1-\left( \eta -\gamma \right) \leqslant 1$ and $\alpha
+2\gamma -\eta =\left( \alpha -\gamma +2\eta \right) -3\left( \eta -\gamma
\right) \leqslant 1$ as well as $\alpha +2\gamma -\eta =\left( \alpha
+\gamma -\eta \right) +\gamma \geqslant \gamma ,$ while the second term in (%
\ref{NN-of-K-IDD+.IDD+-3}) has a non-negative value according to (\ref%
{STD-IDD+.IDD+-1})$_{2}$, so that the requirement (\ref{NN-of-K-IDD+.IDD+-3}%
) is transformed into the narrowed thermodynamical restriction (\ref%
{STD-IDD+.IDD+-5}) if the right-hand-side of (\ref{STD-IDD+.IDD+-5}) is less
than the right-hand-side of (\ref{STD-IDD+.IDD+-4}), since $\frac{\sin \frac{%
\left( \eta -\gamma \right) \pi }{2}}{\cos \frac{\left( \alpha +\eta \right)
\pi }{2}}=\frac{\cos \frac{\left( 1-\left( \eta -\gamma \right) \right) \pi 
}{2}}{\cos \frac{\left( \alpha +\eta \right) \pi }{2}}\leqslant 1$ and $%
\frac{\cos \frac{\left( \alpha +2\gamma -\eta \right) \pi }{2}}{\cos \frac{%
\left( \alpha +\eta \right) \pi }{2}}\geqslant 1,$ due to $\alpha +\eta
\leqslant 1-\left( \eta -\gamma \right) <1$ and $\alpha +2\gamma -\eta
\leqslant \alpha +\eta <1$, respectively reducing to $\alpha -\gamma +2\eta
\leqslant 1$ and $\gamma \leqslant \eta ,$ see (\ref{STD-IDD+.IDD+-1})$%
_{3,2} $.

In conclusion, the non-negativity of function $K,$ corresponding to the
model IDD$^{{}^{+}}$.IDD$^{{}^{+}}$ and given by (\ref{K-IDD+.IDD+}), is
ensured by requesting (\ref{STD-IDD+.IDD+-2}), (\ref{STD-IDD+.IDD+-2a}), and
possibly (\ref{STD-IDD+.IDD+-5}) in addition to (\ref{STD-IDD+.IDD+-1}), (%
\ref{STD-IDD+.IDD+-3}), and (\ref{STD-IDD+.IDD+-4}) implying the
corresponding relaxation modulus to be a completely monotonic function and
creep compliance a Bernstein function.

\subsubsection{Model I$^{{}^{+}}$ID.IDD$^{{}^{+}}$}

Function $K$, given by (\ref{K}), in the case of the model I$^{{}^{+}}$ID.IDD%
$^{{}^{+}}$, having the constitutive equation given in the form%
\begin{equation*}
\left( a_{1}\,_{0}\mathrm{I}_{t}^{1+\alpha }+a_{2}\,_{0}\mathrm{I}_{t}^{%
\frac{1+\alpha -\gamma }{2}}+a_{3}\,_{0}\mathrm{D}_{t}^{\gamma }\right)
\sigma \left( t\right) =\left( b_{1}\,_{0}\mathrm{I}_{t}^{\alpha +\gamma
-\eta }+b_{2}\,_{0}\mathrm{D}_{t}^{\frac{1+\eta -\left( \alpha +\gamma -\eta
\right) }{2}}+b_{3}\,_{0}\mathrm{D}_{t}^{1+\eta }\right) \varepsilon \left(
t\right) ,
\end{equation*}%
and with corresponding thermodynamical, see (\ref%
{TD-frale-over-less-less-over-1}) - (\ref{TD-frale-over-less-less-over-4}),
as well as with the narrowed thermodynamical requirements%
\begin{gather}
\eta \leqslant \gamma ,\quad \alpha +\gamma +2\left( \gamma -\eta \right)
=\alpha +3\gamma -2\eta \leqslant 1,  \label{STD-frale-over-less-less-over-1}
\\
\frac{a_{1}}{b_{1}}\frac{\sin \frac{\left( 1+\alpha -\gamma +2\eta \right)
\pi }{4}}{\cos \frac{\left( 1-\alpha -3\gamma +2\eta \right) \pi }{4}}%
\leqslant \frac{a_{1}}{b_{1}}\frac{\sin \frac{\left( 1+\alpha -\gamma +2\eta
\right) \pi }{4}}{\cos \frac{\left( 1-\alpha -3\gamma +2\eta \right) \pi }{4}%
}\frac{\cos \frac{\left( 1+\alpha -\gamma +2\eta \right) \pi }{4}}{\sin 
\frac{\left( 1-\alpha -3\gamma +2\eta \right) \pi }{4}}\leqslant \frac{a_{2}%
}{b_{2}},  \label{STD-frale-over-less-less-over-2} \\
\frac{a_{2}}{b_{2}}\leqslant \frac{a_{3}}{b_{3}}\frac{\cos \frac{\left(
1-\alpha -3\gamma +2\eta \right) \pi }{4}}{\sin \frac{\left( 1+\alpha
-\gamma +2\eta \right) \pi }{4}}\frac{\sin \frac{\left( 1-\alpha -3\gamma
+2\eta \right) \pi }{4}}{\cos \frac{\left( 1+\alpha -\gamma +2\eta \right)
\pi }{4}}\leqslant \frac{a_{3}}{b_{3}}\frac{\cos \frac{\left( 1-\alpha
-3\gamma +2\eta \right) \pi }{4}}{\sin \frac{\left( 1+\alpha -\gamma +2\eta
\right) \pi }{4}},  \label{STD-frale-over-less-less-over-2a} \\
a_{1}b_{3}\cos \frac{\left( \alpha +\eta \right) \pi }{2}-a_{2}b_{2}\sin 
\frac{\left( \gamma -\eta \right) \pi }{2}\leqslant a_{3}b_{1}\cos \frac{%
\left( \alpha +2\gamma -\eta \right) \pi }{2},
\label{STD-frale-over-less-less-over-3} \\
a_{3}b_{1}\sin \frac{\left( \alpha +2\gamma -\eta \right) \pi }{2}\leqslant
a_{2}b_{2}\cos \frac{\left( \gamma -\eta \right) \pi }{2}-a_{1}b_{3}\sin 
\frac{\left( \alpha +\eta \right) \pi }{2},
\label{STD-frale-over-less-less-over-4} \\
a_{3}b_{1}\sin \frac{\left( \alpha +2\gamma -\eta \right) \pi }{2}\leqslant
a_{2}b_{2}\cos \frac{\left( \gamma -\eta \right) \pi }{2}\frac{\sin \frac{%
\left( \gamma -\eta \right) \pi }{2}}{\cos \frac{\left( \alpha +2\gamma
-\eta \right) \pi }{2}}+a_{1}b_{3}\sin \frac{\left( \alpha +\eta \right) \pi 
}{2}\frac{\cos \frac{\left( \alpha +\eta \right) \pi }{2}}{\cos \frac{\left(
\alpha +2\gamma -\eta \right) \pi }{2}},
\label{STD-frale-over-less-less-over-5}
\end{gather}%
becomes%
\begin{align}
K\left( \rho \right) &=a_{1}b_{1}\sin \left( \left( \gamma -\eta \right) \pi
\right) -a_{1}b_{2}\rho ^{\frac{1+\alpha +\gamma }{2}}\sin \frac{\left(
1+\alpha -\gamma +2\eta \right) \pi }{2}+a_{1}b_{3}\rho ^{1+\alpha +\gamma
}\sin \left( \left( \alpha +\eta \right) \pi \right)  \notag \\
&\quad+a_{2}b_{1}\rho ^{\frac{1+\alpha +\gamma }{2}}\sin \frac{\left(
1-\alpha -3\gamma +2\eta \right) \pi }{2}+a_{2}b_{2}\rho ^{1+\alpha +\gamma
}\sin \left( \left( \gamma -\eta \right) \pi \right)  \notag \\
&\quad -a_{2}b_{3}\rho ^{3\frac{1+\alpha +\gamma }{2}}\sin \frac{\left(
1+\alpha -\gamma +2\eta \right) \pi }{2} -a_{3}b_{1}\rho ^{1+\alpha +\gamma
}\sin \left( \left( \alpha +2\gamma -\eta \right) \pi \right)  \notag \\
&\quad+a_{3}b_{2}\rho ^{3\frac{1+\alpha +\gamma }{2}}\sin \frac{\left(
1-\alpha -3\gamma +2\eta \right) \pi }{2}+a_{3}b_{3}\rho ^{2\left( 1+\alpha
+\gamma \right) }\sin \left( \left( \gamma -\eta \right) \pi \right) ,
\label{K-I+ID.IDD+}
\end{align}%
so that the non-negativity requirement of function $K$ implies%
\begin{gather}
-a_{1}b_{2}\sin \frac{\left( 1+\alpha -\gamma +2\eta \right) \pi }{2}%
+a_{2}b_{1}\sin \frac{\left( 1-\alpha -3\gamma +2\eta \right) \pi }{2}%
\geqslant 0,  \label{NN-of-K-I+ID.IDD+-1} \\
-a_{2}b_{3}\sin \frac{\left( 1+\alpha -\gamma +2\eta \right) \pi }{2}%
+a_{3}b_{2}\sin \frac{\left( 1-\alpha -3\gamma +2\eta \right) \pi }{2}%
\geqslant 0,  \label{NN-of-K-I+ID.IDD+-2} \\
a_{1}b_{3}\sin \left( \left( \alpha +\eta \right) \pi \right)
+a_{2}b_{2}\sin \left( \left( \gamma -\eta \right) \pi \right)
-a_{3}b_{1}\sin \left( \left( \alpha +2\gamma -\eta \right) \pi \right)
\geqslant 0.  \label{NN-of-K-I+ID.IDD+-3}
\end{gather}

According to the thermodynamical requirements (\ref%
{STD-frale-over-less-less-over-1}), the sine in the first term in (\ref%
{NN-of-K-I+ID.IDD+-1}), as well as in (\ref{NN-of-K-I+ID.IDD+-2}), has a
non-negative value, since $1+\alpha -\gamma +2\eta \geqslant \left( 1-\left(
\gamma -\eta \right) \right) +\alpha +\eta \geqslant 0$ and$\ 1+\alpha
-\gamma +2\eta =1+\left( \alpha +3\gamma -2\eta \right) -4\left( \gamma
-\eta \right) \leqslant 2-4\left( \gamma -\eta \right) \leqslant 2$, while
the argument $\frac{\left( 1-\alpha -3\gamma +2\eta \right) \pi }{2}$ in
second term in (\ref{NN-of-K-I+ID.IDD+-1}) and (\ref{NN-of-K-I+ID.IDD+-2})
is in the interval $\left( 0,\frac{\pi }{2}\right) $, since $\alpha +3\gamma
-2\eta \in \left( 0,1\right) $, according to (\ref%
{STD-frale-over-less-less-over-1})$_{2}$, so that the requirements (\ref%
{NN-of-K-I+ID.IDD+-1}) and (\ref{NN-of-K-I+ID.IDD+-2}) are transformed into%
\begin{equation*}
\frac{a_{1}}{b_{1}}\frac{\sin \frac{\left( 1+\alpha -\gamma +2\eta \right)
\pi }{4}}{\cos \frac{\left( 1-\alpha -3\gamma +2\eta \right) \pi }{4}}\frac{%
\cos \frac{\left( 1+\alpha -\gamma +2\eta \right) \pi }{4}}{\sin \frac{%
\left( 1-\alpha -3\gamma +2\eta \right) \pi }{4}}\leqslant \frac{a_{2}}{b_{2}%
}\leqslant \frac{a_{3}}{b_{3}}\frac{\cos \frac{\left( 1-\alpha -3\gamma
+2\eta \right) \pi }{4}}{\sin \frac{\left( 1+\alpha -\gamma +2\eta \right)
\pi }{4}}\frac{\sin \frac{\left( 1-\alpha -3\gamma +2\eta \right) \pi }{4}}{%
\cos \frac{\left( 1+\alpha -\gamma +2\eta \right) \pi }{4}},
\end{equation*}%
and represent narrowed thermodynamical requirements (\ref%
{STD-frale-over-less-less-over-2})$_{1}$ and (\ref%
{STD-frale-over-less-less-over-2a})$_{2}$, when combined with the
thermodynamical restriction (\ref{TD-frale-over-less-less-over-2}), since $%
\frac{\cos \frac{\left( 1+\alpha -\gamma +2\eta \right) \pi }{4}}{\sin \frac{%
\left( 1-\alpha -3\gamma +2\eta \right) \pi }{4}}=\frac{\cos \frac{\left(
1+\alpha -\gamma +2\eta \right) \pi }{4}}{\cos \frac{\left( 1+\alpha
+3\gamma -2\eta \right) \pi }{4}}\geqslant 1$ and $\frac{\sin \frac{\left(
1-\alpha -3\gamma +2\eta \right) \pi }{4}}{\cos \frac{\left( 1+\alpha
-\gamma +2\eta \right) \pi }{4}}=\frac{\cos \frac{\left( 1+\alpha +3\gamma
-2\eta \right) \pi }{4}}{\cos \frac{\left( 1+\alpha -\gamma +2\eta \right)
\pi }{4}}\leqslant 1$, due to $1+\alpha -\gamma +2\eta \leqslant 1+\alpha
+3\gamma -2\eta $ reducing to $\eta \leqslant \gamma ,$ see (\ref%
{STD-frale-over-less-less-over-2})$_{1}.$

The requirements (\ref{STD-frale-over-less-less-over-1}) ensure non-negative
values of sines in the first and third term in (\ref{NN-of-K-I+ID.IDD+-3}),
since $0\leqslant \alpha +\eta \leqslant \alpha +\gamma +2\left( \gamma
-\eta \right) \leqslant 1$ and $0\leqslant \alpha +2\gamma -\eta =\left(
\alpha +3\gamma -2\eta \right) -\left( \gamma -\eta \right) \leqslant 1,$
while sine in the second term in (\ref{NN-of-K-I+ID.IDD+-3}) has a
non-negative value according to the requirement (\ref%
{STD-frale-over-less-less-over-1})$_{1}$, so that the requirement (\ref%
{NN-of-K-I+ID.IDD+-3}) is transformed into the narrowed thermodynamical
restriction (\ref{STD-frale-over-less-less-over-5}) if the right-hand-side
of (\ref{STD-frale-over-less-less-over-5}) is less than the right-hand-side
of (\ref{STD-frale-over-less-less-over-4}), since $\frac{\sin \frac{\left(
\gamma -\eta \right) \pi }{2}}{\cos \frac{\left( \alpha +2\gamma -\eta
\right) \pi }{2}}=\frac{\cos \frac{\left( 1-\left( \gamma -\eta \right)
\right) \pi }{2}}{\cos \frac{\left( \alpha +2\gamma -\eta \right) \pi }{2}}%
\leqslant 1$ and $\frac{\cos \frac{\left( \alpha +\eta \right) \pi }{2}}{%
\cos \frac{\left( \alpha +2\gamma -\eta \right) \pi }{2}}\geqslant 1,$ due
to $\alpha +2\gamma -\eta \leqslant 1-\left( \gamma -\eta \right) <1$ and $%
\alpha +\eta \leqslant \alpha +2\gamma -\eta <1$, respectively reducing to $%
\alpha +3\gamma -2\eta \leqslant 1$ and $\eta \leqslant \gamma ,$ see (\ref%
{STD-frale-over-less-less-over-1}).

In conclusion, the non-negativity of function $K,$ corresponding to the
model I$^{{}^{+}}$ID.IDD$^{{}^{+}}$ and given by (\ref{K-I+ID.IDD+}), is
ensured by requesting (\ref{STD-frale-over-less-less-over-2}), (\ref%
{STD-frale-over-less-less-over-2a}) and possibly (\ref%
{STD-frale-over-less-less-over-5}) in addition to (\ref%
{STD-frale-over-less-less-over-1}), (\ref{STD-frale-over-less-less-over-3})
and (\ref{STD-frale-over-less-less-over-4}) implying the corresponding
relaxation modulus to be a completely monotonic function and creep
compliance a Bernstein function.

\subsection{Asymmetric models}

\subsubsection{Model IID.ID}

Function $K$, given by (\ref{K}), in the case of model IID.ID, having the
constitutive equation given in the form%
\begin{equation*}
\left( a_{1}\,_{0}\mathrm{I}_{t}^{\alpha +\beta -\gamma }+a_{2}\,_{0}\mathrm{%
I}_{t}^{\nu }+a_{3}\,_{0}\mathrm{D}_{t}^{\gamma }\right) \sigma \left(
t\right) =\left( b_{1}\,_{0}\mathrm{I}_{t}^{\alpha }+b_{2}\,_{0}\mathrm{D}%
_{t}^{\beta }\right) \varepsilon \left( t\right) ,
\end{equation*}%
and with corresponding thermodynamical, see (\ref{TD-(A)Z-IID-ID-less-1})
and (\ref{TD-(A)Z-IID-ID-less-2}), as well as with the narrowed
thermodynamical requirements%
\begin{gather}
0\leqslant \alpha \leqslant \nu <\alpha +\beta -\gamma \leqslant 1,\quad
\beta +\nu \leqslant 1,  \label{STD-IID.ID-1} \\
\frac{b_{1}}{b_{2}}\leqslant \frac{a_{1}}{a_{3}}\frac{\sin \frac{\left(
\alpha +2\beta -\gamma \right) \pi }{2}}{\sin \frac{\left( \alpha +\gamma
\right) \pi }{2}}\frac{\cos \frac{\left( \alpha +2\beta -\gamma \right) \pi 
}{2}}{\cos \frac{\left( \alpha +\gamma \right) \pi }{2}}\leqslant \frac{a_{1}%
}{a_{3}}\frac{\sin \frac{\left( \alpha +2\beta -\gamma \right) \pi }{2}}{%
\sin \frac{\left( \alpha +\gamma \right) \pi }{2}},  \label{STD-IID.ID-2}
\end{gather}%
valid if $\alpha +2\beta -\gamma <1,$ becomes%
\begin{align}
K\left( \rho \right) & =a_{1}b_{1}\sin \left( \left( \beta -\gamma \right)
\pi \right) +a_{1}b_{2}\rho ^{\alpha +\beta }\sin \left( \left( \alpha
+2\beta -\gamma \right) \pi \right)  \notag \\
& \quad +a_{2}b_{1}\rho ^{\alpha +\beta -\gamma -\nu }\sin \left( \left( \nu
-\alpha \right) \pi \right) +a_{2}b_{2}\rho ^{2\left( \alpha +\beta \right)
-\gamma -\nu }\sin \left( \left( \beta +\nu \right) \pi \right)  \notag \\
& \quad -a_{3}b_{1}\rho ^{\alpha +\beta }\sin \left( \left( \alpha +\gamma
\right) \pi \right) +a_{3}b_{2}\rho ^{2\left( \alpha +\beta \right) }\sin
\left( \left( \beta -\gamma \right) \pi \right) ,  \label{K-IID.ID}
\end{align}%
so that the non-negativity requirement of function $K$ implies%
\begin{equation}
a_{1}b_{2}\sin \left( \left( \alpha +2\beta -\gamma \right) \pi \right)
-a_{3}b_{1}\sin \left( \left( \alpha +\gamma \right) \pi \right) \geqslant 0.
\label{NN-of-K-IID.ID}
\end{equation}

The first term in (\ref{NN-of-K-IID.ID}) has either a positive value if $%
\alpha +2\beta -\gamma <1$, or a non-positive value if $\alpha +2\beta
-\gamma \geqslant 1$, since $\alpha +2\beta -\gamma =\left( \alpha +\beta
-\gamma \right) +\beta \in \left( 0,2\right) $, according to the
thermodynamical requirement (\ref{STD-IID.ID-1})$_{1},$ while the second
term in (\ref{NN-of-K-IID.ID}) has a non-positive value, since the
thermodynamical requirement (\ref{STD-IID.ID-1})$_{1}$ gives $\alpha +\gamma
\leqslant \nu +\gamma <\alpha +\beta $ implying $\gamma <\beta ,$ that
combined with $\alpha \leqslant \nu ,$ see (\ref{STD-IID.ID-1})$_{1},$
yields $\alpha +\gamma \leqslant \beta +\nu \leqslant 1$ according to the
thermodynamical requirement (\ref{STD-IID.ID-1})$_{2},$ so that the
non-negativity of requirement (\ref{NN-of-K-IID.ID}), as well as of the
function $K$, given by (\ref{K-IID.ID}), can be guaranteed only if $\alpha
+2\beta -\gamma <1,$ and than the requirement (\ref{NN-of-K-IID.ID}) is
transformed into%
\begin{equation*}
\frac{b_{1}}{b_{2}}\leqslant \frac{a_{1}}{a_{3}}\frac{\sin \frac{\left(
\alpha +2\beta -\gamma \right) \pi }{2}}{\sin \frac{\left( \alpha +\gamma
\right) \pi }{2}}\frac{\cos \frac{\left( \alpha +2\beta -\gamma \right) \pi 
}{2}}{\cos \frac{\left( \alpha +\gamma \right) \pi }{2}},
\end{equation*}%
representing the narrowed thermodynamical restriction (\ref{STD-IID.ID-2})$%
_{2}$, since $\alpha +\gamma \leqslant \alpha +2\beta -\gamma <1$ reducing
to $\gamma <\beta ,$ see (\ref{STD-IID.ID-1})$_{1}$, implies $\frac{\cos 
\frac{\left( \alpha +2\beta -\gamma \right) \pi }{2}}{\cos \frac{\left(
\alpha +\gamma \right) \pi }{2}}\leqslant 1$.

In conclusion, the non-negativity of function $K,$ corresponding to the
model IID.ID and given by (\ref{K-IID.ID}), is ensured if $\alpha +2\beta
-\gamma <1$ by requesting (\ref{STD-IID.ID-2})$,$ in addition to (\ref%
{STD-IID.ID-1}) implying the corresponding relaxation modulus to be a
completely monotonic function and creep compliance a Bernstein function,
while if $\alpha +2\beta -\gamma \geqslant 1,$ then the mentioned properties
cannot be ensured, since (\ref{NN-of-K-IID.ID}) cannot be satisfied.

\subsubsection{Model IDD.DD$^{{}^{+}}$}

Function $K$, given by (\ref{K}), in the case of model IDD.DD$^{{}^{+}}$,
having the constitutive equation given in the form%
\begin{equation*}
\left( a_{1}\,_{0}\mathrm{I}_{t}^{\alpha }+a_{2}\,_{0}\mathrm{D}_{t}^{\beta
}+a_{3}\,_{0}\mathrm{D}_{t}^{\gamma }\right) \sigma \left( t\right) =\left(
b_{1}\,_{0}\mathrm{D}_{t}^{\mu }+b_{2}\,_{0}\mathrm{D}_{t}^{\alpha +\beta
+\mu }\right) \varepsilon \left( t\right) ,
\end{equation*}%
and with corresponding thermodynamical, see (\ref{TD-AZ-IDD-DD-less-over-1})
and (\ref{TD-AZ-IDD-DD-less-over-2}), as well as with the narrowed
thermodynamical requirements%
\begin{gather}
1\leqslant \alpha +\beta +\mu \leqslant 2,\quad \beta <\gamma \leqslant \mu
\leqslant 1-\alpha ,  \label{STD-IDD.DD+-1} \\
\frac{a_{1}}{a_{2}}\frac{\left\vert \cos \frac{\left( 2\alpha +\beta +\mu
\right) \pi }{2}\right\vert }{\cos \frac{\left( \mu -\beta \right) \pi }{2}}%
\leqslant \frac{a_{1}}{a_{2}}\frac{\left\vert \cos \frac{\left( 2\alpha
+\beta +\mu \right) \pi }{2}\right\vert }{\cos \frac{\left( \mu -\beta
\right) \pi }{2}}\frac{\sin \frac{\left( 2\alpha +\beta +\mu \right) \pi }{2}%
}{\sin \frac{\left( \mu -\beta \right) \pi }{2}}\leqslant \frac{b_{1}}{b_{2}}%
,  \label{STD-IDD.DD+-2}
\end{gather}%
becomes%
\begin{align}
K\left( \rho \right) & =a_{1}b_{1}\sin \left( \left( \alpha +\mu \right) \pi
\right) +a_{1}b_{2}\rho ^{\alpha +\beta }\sin \left( \left( 2\alpha +\beta
+\mu \right) \pi \right)  \notag \\
& \quad +a_{2}b_{1}\rho ^{\alpha +\beta }\sin \left( \left( \mu -\beta
\right) \pi \right) +a_{2}b_{2}\rho ^{2\left( \alpha +\beta \right) }\sin
\left( \left( \alpha +\mu \right) \pi \right) +a_{3}b_{2}\rho ^{2\alpha
+\beta +\mu }\sin \left( \left( \alpha +\beta \right) \pi \right) ,
\label{K-IDD.DD+}
\end{align}%
so that the non-negativity requirement of function $K$ implies%
\begin{equation}
a_{1}b_{2}\sin \left( \left( 2\alpha +\beta +\mu \right) \pi \right)
+a_{2}b_{1}\sin \left( \left( \mu -\beta \right) \pi \right) >0.
\label{NN-of-K-IDD.DD+}
\end{equation}

According to the thermodynamical requirement (\ref{STD-IDD.DD+-1})$_{2},$
one finds $2\alpha +\beta +\mu =\left( \alpha +\beta \right) +\left( \alpha
+\mu \right) \leqslant 2$, as well as according to the thermodynamical
requirement (\ref{STD-IDD.DD+-1})$_{1},$ one has $2\alpha +\beta +\mu
=\left( \alpha +\beta +\mu \right) +\alpha \geqslant 1$, so that the first
term in (\ref{NN-of-K-IDD.DD+}) has a non-positive value, while the second
term in (\ref{NN-of-K-IDD.DD+}) has a non-negative value, according the
thermodynamical requirement (\ref{STD-IDD.DD+-1})$_{2}$. Therefore, the
requirement (\ref{NN-of-K-IDD.DD+}) is transformed into%
\begin{equation*}
\frac{a_{1}}{a_{2}}\frac{\left\vert \cos \frac{\left( 2\alpha +\beta +\mu
\right) \pi }{2}\right\vert }{\cos \frac{\left( \mu -\beta \right) \pi }{2}}%
\frac{\sin \frac{\left( 2\alpha +\beta +\mu \right) \pi }{2}}{\sin \frac{%
\left( \mu -\beta \right) \pi }{2}}\leqslant \frac{b_{1}}{b_{2}},
\end{equation*}%
representing the narrowed thermodynamical restriction (\ref{STD-IDD.DD+-2})$%
_{1}$, since $\pi -\frac{\left( 2\alpha +\beta +\mu \right) \pi }{2}%
\geqslant \frac{\left( \mu -\beta \right) \pi }{2}$, reducing to $\alpha
+\mu \leqslant 1,$ see (\ref{STD-IDD.DD+-1})$_{2},$ implies $\frac{\sin 
\frac{\left( 2\alpha +\beta +\mu \right) \pi }{2}}{\sin \frac{\left( \mu
-\beta \right) \pi }{2}}\geqslant 1$.

In conclusion, the non-negativity of function $K,$ corresponding to the
model IDD.DD$^{{}^{+}}$ and given by (\ref{K-IDD.DD+}), is ensured by
requesting (\ref{STD-IDD.DD+-2}) in addition to (\ref{STD-IDD.DD+-1})
implying the corresponding relaxation modulus to be a completely monotonic
function and creep compliance a Bernstein function.

\subsubsection{Model I$^{{}^{+}}$ID.ID\label{model-I+ID.ID-narrowed}}

Function $K$, given by (\ref{K}), in the case of the model I$^{{}^{+}}$%
ID.ID, having the constitutive equation given in the form%
\begin{equation}
\left( a_{1}\,{}_{0}\mathrm{I}_{t}^{\alpha +\beta +\nu }+a_{2}\,{}_{0}%
\mathrm{I}_{t}^{\nu }+a_{3}\,{}_{0}\mathrm{D}_{t}^{\alpha +\beta -\nu
}\right) \sigma \left( t\right) =\left( b_{1}\,{}_{0}\mathrm{I}_{t}^{\alpha
}+b_{2}\,{}_{0}\mathrm{D}_{t}^{\beta }\right) \varepsilon \left( t\right) ,
\label{I+ID.ID}
\end{equation}%
and with corresponding thermodynamical, see (\ref{TD-(A)Z-IID-ID-over-1})
and (\ref{TD-(A)Z-IID-ID-over-2}), as well as with the narrowed
thermodynamical requirements%
\begin{gather}
0\leqslant \alpha +\beta -\nu \leqslant 1,\quad 1\leqslant \alpha +\beta
+\nu \leqslant 2,\quad \alpha \leqslant \nu \leqslant 1-\beta ,
\label{STD-I+ID.ID-1} \\
\frac{a_{1}}{a_{2}}\frac{\left\vert \cos \frac{\left( \alpha +2\beta +\nu
\right) \pi }{2}\right\vert }{\cos \frac{\left( \nu -\alpha \right) \pi }{2}}%
\leqslant \frac{a_{1}}{a_{2}}\frac{\left\vert \cos \frac{\left( \alpha
+2\beta +\nu \right) \pi }{2}\right\vert }{\cos \frac{\left( \nu -\alpha
\right) \pi }{2}}\frac{\sin \frac{\left( \alpha +2\beta +\nu \right) \pi }{2}%
}{\sin \frac{\left( \nu -\alpha \right) \pi }{2}}\leqslant \frac{b_{1}}{b_{2}%
},  \label{STD-I+ID.ID-2} \\
\frac{b_{1}}{b_{2}}\leqslant \frac{a_{2}}{a_{3}}\frac{\sin \frac{\left(
\beta +\nu \right) \pi }{2}}{\sin \frac{\left( 2\alpha +\beta -\nu \right)
\pi }{2}}\frac{\cos \frac{\left( \beta +\nu \right) \pi }{2}}{\cos \frac{%
\left( 2\alpha +\beta -\nu \right) \pi }{2}}\leqslant \frac{a_{2}}{a_{3}}%
\frac{\sin \frac{\left( \beta +\nu \right) \pi }{2}}{\sin \frac{\left(
2\alpha +\beta -\nu \right) \pi }{2}},  \label{STD-I+ID.ID-3}
\end{gather}%
becomes%
\begin{align}
K\left( \rho \right) & =a_{1}b_{1}\sin \left( \left( \beta +\nu \right) \pi
\right) +a_{1}b_{2}\rho ^{\alpha +\beta }\sin \left( \left( \alpha +2\beta
+\nu \right) \pi \right)  \notag \\
& \quad +a_{2}b_{1}\rho ^{\alpha +\beta }\sin \left( \left( \nu -\alpha
\right) \pi \right) +a_{2}b_{2}\rho ^{2\left( \alpha +\beta \right) }\sin
\left( \left( \beta +\nu \right) \pi \right)  \notag \\
& \quad -a_{3}b_{1}\rho ^{2\left( \alpha +\beta \right) }\sin \left( \left(
2\alpha +\beta -\nu \right) \pi \right) +a_{3}b_{2}\rho ^{3\left( \alpha
+\beta \right) }\sin \left( \left( \nu -\alpha \right) \pi \right) ,
\label{K-I+ID.ID}
\end{align}%
so that the non-negativity requirement of function $K$ implies%
\begin{gather}
a_{1}b_{2}\sin \left( \left( \alpha +2\beta +\nu \right) \pi \right)
+a_{2}b_{1}\sin \left( \left( \nu -\alpha \right) \pi \right) \geqslant 0,
\label{NN-of-K-I+ID.ID-1} \\
a_{2}b_{2}\sin \left( \left( \beta +\nu \right) \pi \right) -a_{3}b_{1}\sin
\left( \left( 2\alpha +\beta -\nu \right) \pi \right) \geqslant 0.
\label{NN-of-K-I+ID.ID-2}
\end{gather}

According to the thermodynamical requirement (\ref{STD-I+ID.ID-1})$_{3},$
one finds $\alpha +2\beta +\nu =\left( \alpha +\beta \right) +\left( \beta
+\nu \right) \leqslant 2$, while according to the thermodynamical
requirement (\ref{STD-I+ID.ID-1})$_{2}$ one finds $\alpha +2\beta +\nu
=\left( \alpha +\beta +\nu \right) +\beta \geqslant 1$, so that the first
term in (\ref{NN-of-K-I+ID.ID-1}) has a non-positive value, while the second
term in (\ref{NN-of-K-I+ID.ID-1}) has a non-negative value, according to the
thermodynamical requirement (\ref{STD-I+ID.ID-1})$_{3}$. Therefore, the
inequality (\ref{NN-of-K-I+ID.ID-1}) is transformed into%
\begin{equation*}
\frac{a_{1}}{a_{2}}\frac{\left\vert \cos \frac{\left( \alpha +2\beta +\nu
\right) \pi }{2}\right\vert }{\cos \frac{\left( \nu -\alpha \right) \pi }{2}}%
\frac{\sin \frac{\left( \alpha +2\beta +\nu \right) \pi }{2}}{\sin \frac{%
\left( \nu -\alpha \right) \pi }{2}}\leqslant \frac{b_{1}}{b_{2}},
\end{equation*}%
representing the narrowed thermodynamical restriction (\ref{STD-I+ID.ID-2})$%
_{1}$, since $\pi -\frac{\left( \alpha +2\beta +\nu \right) \pi }{2}%
\geqslant \frac{\left( \nu -\alpha \right) \pi }{2}$, reducing to $\beta
+\nu \leqslant 1,$ see (\ref{STD-I+ID.ID-1})$_{3},$ implies $\frac{\sin 
\frac{\left( \alpha +2\beta +\nu \right) \pi }{2}}{\sin \frac{\left( \nu
-\alpha \right) \pi }{2}}\geqslant 1$.

According to the thermodynamical requirement (\ref{STD-I+ID.ID-1})$_{3},$
the first term in (\ref{NN-of-K-I+ID.ID-2}) has a non-negative value, while
according to the thermodynamical requirement (\ref{STD-I+ID.ID-1})$_{1},$
one finds $2\alpha +\beta -\nu =\left( \alpha +\beta -\nu \right) +\alpha
\geqslant 0,$ as well as according to the requirement (\ref{STD-I+ID.ID-1})$%
_{3},$ one has $2\alpha +\beta -\nu =\left( \alpha +\beta \right) -\left(
\nu -\alpha \right) \leqslant 1,$ and therefore the second term in (\ref%
{NN-of-K-I+ID.ID-2}) has a non-positive value, so that the inequality (\ref%
{NN-of-K-I+ID.ID-2}) is transformed into%
\begin{equation*}
\frac{b_{1}}{b_{2}}\leqslant \frac{a_{2}}{a_{3}}\frac{\sin \frac{\left(
\beta +\nu \right) \pi }{2}}{\sin \frac{\left( 2\alpha +\beta -\nu \right)
\pi }{2}}\frac{\cos \frac{\left( \beta +\nu \right) \pi }{2}}{\cos \frac{%
\left( 2\alpha +\beta -\nu \right) \pi }{2}},
\end{equation*}%
representing the narrowed thermodynamical restriction (\ref{STD-I+ID.ID-3})$%
_{2}$, since $2\alpha +\beta -\nu \leqslant \beta +\nu <1$, reducing to $%
\alpha \leqslant \nu ,$ see (\ref{STD-I+ID.ID-1})$_{3},$ implies $\frac{\cos 
\frac{\left( \beta +\nu \right) \pi }{2}}{\cos \frac{\left( 2\alpha +\beta
-\nu \right) \pi }{2}}\leqslant 1$.

In conclusion, the non-negativity of function $K,$ corresponding to the
model I$^{{}^{+}}$ID.ID and given by (\ref{K-I+ID.ID}), is ensured by
requesting (\ref{STD-I+ID.ID-2}) and (\ref{STD-I+ID.ID-3}) in addition to (%
\ref{STD-I+ID.ID-1}) implying the corresponding relaxation modulus to be a
completely monotonic function and creep compliance a Bernstein function.

\subsubsection{Model IDD$^{{}^{+}}$.DD$^{{}^{+}}$}

Function $K$, given by (\ref{K}), in the case of model IDD$^{{}^{+}}$.DD$%
^{{}^{+}}$, having the constitutive equation given in the form%
\begin{equation*}
\left( a_{1}\,_{0}\mathrm{I}_{t}^{\alpha }+a_{2}\,_{0}\mathrm{D}_{t}^{\beta
}+a_{3}\,_{0}\mathrm{D}_{t}^{\alpha +2\beta }\right) \sigma \left( t\right)
=\left( b_{1}\,_{0}\mathrm{D}_{t}^{\mu }+b_{2}\,_{0}\mathrm{D}_{t}^{\alpha
+\beta +\mu }\right) \varepsilon \left( t\right) ,
\end{equation*}%
and with corresponding thermodynamical, see (\ref{TD-AZ-IDD-DD-over-over-1})
and (\ref{TD-AZ-IDD-DD-over-over-2}), as well as with the narrowed
thermodynamical requirements%
\begin{gather}
1\leqslant \alpha +2\beta \leqslant 2,\quad 1\leqslant \alpha +\beta +\mu
\leqslant 2,\quad \beta \leqslant \mu \leqslant 1-\alpha ,
\label{STD-IDD+.DD+-1} \\
\frac{a_{1}}{a_{2}}\frac{\left\vert \cos \frac{\left( 2\alpha +\beta +\mu
\right) \pi }{2}\right\vert }{\cos \frac{\left( \mu -\beta \right) \pi }{2}}%
\leqslant \frac{a_{1}}{a_{2}}\frac{\left\vert \cos \frac{\left( 2\alpha
+\beta +\mu \right) \pi }{2}\right\vert }{\cos \frac{\left( \mu -\beta
\right) \pi }{2}}\frac{\sin \frac{\left( 2\alpha +\beta +\mu \right) \pi }{2}%
}{\sin \frac{\left( \mu -\beta \right) \pi }{2}}\leqslant \frac{b_{1}}{b_{2}}%
,  \label{STD-IDD+.DD+-2} \\
\frac{b_{1}}{b_{2}}\leqslant \frac{a_{2}}{a_{3}}\frac{\sin \frac{\left(
\alpha +\mu \right) \pi }{2}}{\sin \frac{\left( \alpha +2\beta -\mu \right)
\pi }{2}}\frac{\cos \frac{\left( \alpha +\mu \right) \pi }{2}}{\cos \frac{%
\left( \alpha +2\beta -\mu \right) \pi }{2}}\leqslant \frac{a_{2}}{a_{3}}%
\frac{\sin \frac{\left( \alpha +\mu \right) \pi }{2}}{\sin \frac{\left(
\alpha +2\beta -\mu \right) \pi }{2}},  \label{STD-IDD+.DD+-3}
\end{gather}%
valid if $2\alpha +\beta -\mu <1,$ becomes%
\begin{align}
K\left( \rho \right) & =a_{1}b_{1}\sin \left( \left( \alpha +\mu \right) \pi
\right) +a_{1}b_{2}\rho ^{\alpha +\beta }\sin \left( \left( 2\alpha +\beta
+\mu \right) \pi \right)  \notag \\
& \quad +a_{2}b_{1}\rho ^{\alpha +\beta }\sin \left( \left( \mu -\beta
\right) \pi \right) +a_{2}b_{2}\rho ^{2\left( \alpha +\beta \right) }\sin
\left( \left( \alpha +\mu \right) \pi \right)  \notag \\
& \quad -a_{3}b_{1}\rho ^{2\left( \alpha +\beta \right) }\sin \left( \left(
\alpha +2\beta -\mu \right) \pi \right) +a_{3}b_{2}\rho ^{3\left( \alpha
+\beta \right) }\sin \left( \left( \mu -\beta \right) \pi \right) ,
\label{K-IDD+.DD+}
\end{align}%
so that the non-negativity requirement of function $K$ implies%
\begin{gather}
a_{1}b_{2}\sin \left( \left( 2\alpha +\beta +\mu \right) \pi \right)
+a_{2}b_{1}\sin \left( \left( \mu -\beta \right) \pi \right) \geqslant 0,
\label{NN-of-K-IDD+.DD+-1} \\
a_{2}b_{2}\sin \left( \left( \alpha +\mu \right) \pi \right) -a_{3}b_{1}\sin
\left( \left( \alpha +2\beta -\mu \right) \pi \right) \geqslant 0.
\label{NN-of-K-IDD+.DD+-2}
\end{gather}

According to the thermodynamical requirement (\ref{STD-IDD+.DD+-1})$_{3},$
one finds $\alpha +2\beta +\mu =\left( \alpha +\beta \right) +\left( \beta
+\mu \right) \leqslant 2$, while according to the thermodynamical
requirement (\ref{STD-IDD+.DD+-1})$_{2},$ one finds $\alpha +2\beta +\mu
=\left( \alpha +\beta +\mu \right) +\beta \geqslant 1$, so that the first
term in (\ref{NN-of-K-IDD+.DD+-1}) has a non-positive value, while the
second term in (\ref{NN-of-K-IDD+.DD+-1}) has a non-negative value according
to the thermodynamical requirement (\ref{STD-IDD+.DD+-1})$_{3}$. Therefore,
the inequality (\ref{NN-of-K-IDD+.DD+-1}) is transformed into%
\begin{equation*}
\frac{a_{1}}{a_{2}}\frac{\left\vert \cos \frac{\left( 2\alpha +\beta +\mu
\right) \pi }{2}\right\vert }{\cos \frac{\left( \mu -\beta \right) \pi }{2}}%
\frac{\sin \frac{\left( 2\alpha +\beta +\mu \right) \pi }{2}}{\sin \frac{%
\left( \mu -\beta \right) \pi }{2}}\leqslant \frac{b_{1}}{b_{2}},
\end{equation*}%
representing the narrowed thermodynamical restriction (\ref{STD-IDD+.DD+-2})$%
_{1}$, since $\pi -\frac{\left( 2\alpha +\beta +\mu \right) \pi }{2}%
\geqslant \frac{\left( \mu -\beta \right) \pi }{2}$, reducing to $\alpha
+\mu \leqslant 1,$ see (\ref{STD-IDD+.DD+-1})$_{3},$ implies $\frac{\sin 
\frac{\left( 2\alpha +\beta +\mu \right) \pi }{2}}{\sin \frac{\left( \mu
-\beta \right) \pi }{2}}\geqslant 1$.

According to the thermodynamical requirement (\ref{STD-IDD+.DD+-1})$_{3},$
the first term in (\ref{NN-of-K-IDD+.DD+-2}) has a non-negative value, while
according to the thermodynamical requirement (\ref{STD-IDD+.DD+-1})$_{1},$
one finds $\alpha +2\beta -\mu =\left( \alpha +2\beta \right) -\mu \geqslant
0,$ as well as according to the requirement (\ref{STD-IDD+.DD+-1})$_{3},$
one finds $\alpha +2\beta -\mu =\left( \alpha +\beta \right) -\left( \mu
-\beta \right) \leqslant 1,$ and therefore the second term in (\ref%
{NN-of-K-IDD+.DD+-2}) has a non-positive value, so that\ the inequality (\ref%
{NN-of-K-IDD+.DD+-2}) is transformed into%
\begin{equation*}
\frac{b_{1}}{b_{2}}\leqslant \frac{a_{2}}{a_{3}}\frac{\sin \frac{\left(
\alpha +\mu \right) \pi }{2}}{\sin \frac{\left( \alpha +2\beta -\mu \right)
\pi }{2}}\frac{\cos \frac{\left( \alpha +\mu \right) \pi }{2}}{\cos \frac{%
\left( \alpha +2\beta -\mu \right) \pi }{2}},
\end{equation*}%
representing the narrowed thermodynamical restriction (\ref{STD-IDD+.DD+-3})$%
_{2}$, since $\alpha +2\beta -\mu \leqslant \alpha +\mu <1$, reducing to $%
\beta \leqslant \mu ,$ see (\ref{STD-IDD+.DD+-1})$_{3},$ implies $\frac{\cos 
\frac{\left( \alpha +\mu \right) \pi }{2}}{\cos \frac{\left( \alpha +2\beta
-\mu \right) \pi }{2}}\leqslant 1$.

In conclusion, the non-negativity of function $K,$ corresponding to the
model IDD$^{{}^{+}}$.DD$^{{}^{+}}$ and given by (\ref{K-IDD+.DD+}), is
ensured by requesting (\ref{STD-IDD+.DD+-2}) and (\ref{STD-IDD+.DD+-3}) in
addition to (\ref{STD-IDD+.DD+-1}), implying the corresponding relaxation
modulus to be a completely monotonic function and creep compliance a
Bernstein function.

\subsubsection{Model ID.IDD}

Function $K$, given by (\ref{K}), in the case of the model ID.IDD, having
the constitutive equation given in the form%
\begin{equation*}
\left( a_{1}\,_{0}\mathrm{I}_{t}^{\mu }+a_{2}\,_{0}\mathrm{D}_{t}^{\nu
}\right) \sigma \left( t\right) =\left( b_{1}\,_{0}\mathrm{I}_{t}^{\alpha
}+b_{2}\,_{0}\mathrm{D}_{t}^{\beta }+b_{3}\,_{0}\mathrm{D}_{t}^{\mu +\nu
-\alpha }\right) \varepsilon \left( t\right) ,
\end{equation*}%
and with corresponding thermodynamical, see (\ref{TD-(A)Z-ID-IDD-less-1})
and (\ref{TD-(A)Z-ID-IDD-less-2}), as well as with the narrowed
thermodynamical requirements%
\begin{gather}
0\leqslant \nu \leqslant \beta <\mu +\nu -\alpha \leqslant 1,\quad \alpha
\leqslant \mu \leqslant 1-\beta ,  \label{STD-ID.IDD-1} \\
-\frac{a_{1}}{a_{2}}\frac{\cos \frac{\left( 2\alpha +\beta -\mu \right) \pi 
}{2}}{\cos \frac{\left( \beta +\mu \right) \pi }{2}}\leqslant \frac{b_{1}}{%
b_{3}}\leqslant \frac{a_{1}}{a_{2}}\frac{\sin \frac{\left( 2\mu +\nu -\alpha
\right) \pi }{2}}{\sin \frac{\left( \alpha +\nu \right) \pi }{2}}\frac{\cos 
\frac{\left( 2\mu +\nu -\alpha \right) \pi }{2}}{\cos \frac{\left( \alpha
+\nu \right) \pi }{2}}\leqslant \frac{a_{1}}{a_{2}}\frac{\sin \frac{\left(
2\alpha +\beta -\mu \right) \pi }{2}}{\sin \frac{\left( \beta +\mu \right)
\pi }{2}},  \label{STD-ID.IDD-2}
\end{gather}%
valid if $2\mu +\nu -\alpha <1$, becomes%
\begin{align}
K\left( \rho \right) & =a_{1}b_{1}\sin \left( \left( \mu -\alpha \right) \pi
\right) +a_{1}b_{2}\rho ^{\alpha +\beta }\sin \left( \left( \beta +\mu
\right) \pi \right) +a_{1}b_{3}\rho ^{\mu +\nu }\sin \left( \left( 2\mu +\nu
-\alpha \right) \pi \right)  \notag \\
& \quad -a_{2}b_{1}\rho ^{\mu +\nu }\sin \left( \left( \alpha +\nu \right)
\pi \right) +a_{2}b_{2}\rho ^{\alpha +\beta +\mu +\nu }\sin \left( \left(
\beta -\nu \right) \pi \right) +a_{2}b_{3}\rho ^{2\left( \mu +\nu \right)
}\sin \left( \left( \mu -\alpha \right) \pi \right) ,  \label{K-ID.IDD}
\end{align}%
so that the non-negativity requirement of function $K$ implies%
\begin{equation}
a_{1}b_{3}\sin \left( \left( 2\mu +\nu -\alpha \right) \pi \right)
-a_{2}b_{1}\sin \left( \left( \alpha +\nu \right) \pi \right) \geqslant 0.
\label{NN-of-K-ID.IDD}
\end{equation}

The first term in (\ref{NN-of-K-ID.IDD}) has either a positive value if $%
2\mu +\nu -\alpha <1$, or non-positive value if $2\mu +\nu -\alpha \geqslant
1$, since $2\mu +\nu -\alpha =\left( \mu +\nu -\alpha \right) +\mu \in
\left( 0,2\right) $ according to the thermodynamical requirement (\ref%
{STD-ID.IDD-1})$_{1}$, while the second term in (\ref{NN-of-K-ID.IDD}) has a
non-positive value, since $\alpha +\nu \leqslant \alpha +\beta \leqslant 1,$
according to the thermodynamical requirements (\ref{STD-ID.IDD-1}).
Therefore, the inequality (\ref{NN-of-K-ID.IDD}) can be satisfied only if $%
2\mu +\nu -\alpha <1$ and it is transformed into%
\begin{equation*}
\frac{b_{1}}{b_{3}}\leqslant \frac{a_{1}}{a_{2}}\frac{\sin \frac{\left( 2\mu
+\nu -\alpha \right) \pi }{2}}{\sin \frac{\left( \alpha +\nu \right) \pi }{2}%
}\frac{\cos \frac{\left( 2\mu +\nu -\alpha \right) \pi }{2}}{\cos \frac{%
\left( \alpha +\nu \right) \pi }{2}},
\end{equation*}%
representing the narrowed thermodynamical restriction (\ref{STD-ID.IDD-2})$%
_{3}$, since $\alpha +\nu \leqslant 2\mu +\nu -\alpha <1,$ reducing to $%
\alpha \leqslant \mu ,$ see (\ref{STD-ID.IDD-1})$_{2},$ implies $\frac{\cos 
\frac{\left( 2\mu +\nu -\alpha \right) \pi }{2}}{\cos \frac{\left( \alpha
+\nu \right) \pi }{2}}\leqslant 1$.

In conclusion, the non-negativity of function $K,$ corresponding to the
model ID.IDD and given by (\ref{K-ID.IDD}), is ensured if $2\mu +\nu -\alpha
<1$ by requesting (\ref{STD-ID.IDD-2}) in addition to (\ref{STD-ID.IDD-1})
implying the corresponding relaxation modulus to be a completely monotonic
function and creep compliance a Bernstein function, while if $2\mu +\nu
-\alpha \geqslant 1$, then the mentioned properties cannot be ensured, since
(\ref{NN-of-K-ID.IDD}) cannot be satisfied.

\subsubsection{Model ID.DDD$^{{}^{+}}$}

Function $K$, given by (\ref{K}), in the case of the model ID.DDD$^{{}^{+}}$%
, having the constitutive equation given in the form%
\begin{equation*}
\left( a_{1}\,_{0}\mathrm{I}_{t}^{\alpha }+a_{2}\,_{0}\mathrm{D}_{t}^{\beta
}\right) \sigma \left( t\right) =\left( b_{1}\,_{0}\mathrm{D}_{t}^{\mu
}+b_{2}\,_{0}\mathrm{D}_{t}^{\nu }+b_{3}\,_{0}\mathrm{D}_{t}^{\alpha +\beta
+\nu }\right) \varepsilon \left( t\right) ,
\end{equation*}%
and with corresponding thermodynamical, see (\ref{TD-(A)Z-ID-DDD-over-1})
and (\ref{TD-(A)Z-ID-DDD-over-2}), as well as with the narrowed
thermodynamical requirements%
\begin{gather}
1\leqslant \alpha +\beta +\nu \leqslant 2,\quad \beta \leqslant \mu <\nu
\leqslant 1-\alpha ,  \label{STD-ID.DDD+-1} \\
\frac{a_{1}}{a_{2}}\frac{\left\vert \cos \frac{\left( 2\alpha +\beta +\nu
\right) \pi }{2}\right\vert }{\cos \frac{\left( \nu -\beta \right) \pi }{2}}%
\leqslant \frac{a_{1}}{a_{2}}\frac{\left\vert \cos \frac{\left( 2\alpha
+\beta +\nu \right) \pi }{2}\right\vert }{\cos \frac{\left( \nu -\beta
\right) \pi }{2}}\frac{\sin \frac{\left( 2\alpha +\beta +\nu \right) \pi }{2}%
}{\sin \frac{\left( \nu -\beta \right) \pi }{2}}\leqslant \frac{b_{2}}{b_{3}}%
,  \label{STD-ID.DDD+-2}
\end{gather}%
becomes%
\begin{align}
K\left( \rho \right) & =a_{1}b_{1}\sin \left( \left( \alpha +\mu \right) \pi
\right) +a_{1}b_{2}\rho ^{\nu -\mu }\sin \left( \left( \alpha +\nu \right)
\pi \right) +a_{1}b_{3}\rho ^{\alpha +\beta +\nu -\mu }\sin \left( \left(
2\alpha +\beta +\nu \right) \pi \right)  \notag \\
& \quad +a_{2}b_{1}\rho ^{\alpha +\beta }\sin \left( \left( \mu -\beta
\right) \pi \right) +a_{2}b_{2}\rho ^{\alpha +\beta +\nu -\mu }\sin \left(
\left( \nu -\beta \right) \pi \right) +a_{2}b_{3}\rho ^{2\alpha +2\beta +\nu
-\mu }\sin \left( \left( \alpha +\nu \right) \pi \right) ,  \label{K-ID.DDD+}
\end{align}%
so that the non-negativity requirement of function $K$ implies%
\begin{equation}
a_{1}b_{3}\sin \left( \left( 2\alpha +\beta +\nu \right) \pi \right)
+a_{2}b_{2}\sin \left( \left( \nu -\beta \right) \pi \right) \geqslant 0.
\label{NN-of-K-ID.DDD+}
\end{equation}

The first term in (\ref{NN-of-K-ID.DDD+}) has a non-positive value, since\ $%
2\alpha +\beta +\nu =\left( \alpha +\beta +\nu \right) +\alpha \geqslant 1$,
according to the thermodynamical requirements (\ref{STD-ID.DDD+-1})$_{1}$,
as well as $2\alpha +\beta +\nu =\left( \alpha +\beta \right) +\left( \alpha
+\nu \right) \leqslant 2$ according to the thermodynamical requirements (\ref%
{STD-ID.DDD+-1})$_{2}$, while the second term in (\ref{NN-of-K-ID.DDD+}) has
a non-negative value, according to the thermodynamical requirement (\ref%
{STD-ID.DDD+-1})$_{2}$. Therefore, the inequality (\ref{NN-of-K-ID.DDD+}) is
transformed into%
\begin{equation*}
\frac{a_{1}}{a_{2}}\frac{\left\vert \cos \frac{\left( 2\alpha +\beta +\nu
\right) \pi }{2}\right\vert }{\cos \frac{\left( \nu -\beta \right) \pi }{2}}%
\frac{\sin \frac{\left( 2\alpha +\beta +\nu \right) \pi }{2}}{\sin \frac{%
\left( \nu -\beta \right) \pi }{2}}\leqslant \frac{b_{2}}{b_{3}},
\end{equation*}%
representing the narrowed thermodynamical restriction (\ref{STD-ID.DDD+-2})$%
_{1}$, since $\frac{\left( \nu -\beta \right) \pi }{2}\leqslant \pi -\frac{%
\left( 2\alpha +\beta +\nu \right) \pi }{2}$ reducing to $\alpha +\nu
\leqslant 1,$ see (\ref{STD-ID.DDD+-1})$_{2},$ implies $\frac{\sin \frac{%
\left( 2\alpha +\beta +\nu \right) \pi }{2}}{\sin \frac{\left( \nu -\beta
\right) \pi }{2}}\geqslant 1$.

In conclusion, the non-negativity of function $K,$ corresponding to the
model ID.DDD$^{{}^{+}}$ and given by (\ref{K-ID.DDD+}), is ensured by
requesting (\ref{STD-ID.DDD+-2}) in addition to (\ref{STD-ID.DDD+-1})
implying the corresponding relaxation modulus to be a completely monotonic
function and creep compliance a Bernstein function.

\subsubsection{Model ID.IDD$^{{}^{+}}$}

Function $K$, given by (\ref{K}), in the case of the model ID.IDD$^{{}^{+}}$%
, having the constitutive equation given in the form%
\begin{equation*}
\left( a_{1}\,_{0}\mathrm{I}_{t}^{\alpha }+a_{2}\,_{0}\mathrm{D}_{t}^{\beta
}\right) \sigma \left( t\right) =\left( b_{1}\,_{0}\mathrm{I}_{t}^{\alpha
+\beta -\nu }+b_{2}\,_{0}\mathrm{D}_{t}^{\nu }+b_{3}\,_{0}\mathrm{D}%
_{t}^{\alpha +\beta +\nu }\right) \varepsilon \left( t\right) ,
\end{equation*}%
and with corresponding thermodynamical, see (\ref{TD-(A)Z-ID-IDD-over-1})
and (\ref{TD-(A)Z-ID-IDD-over-2}), as well as with the narrowed
thermodynamical requirements%
\begin{gather}
0\leqslant \alpha +\beta -\nu \leqslant 1,\quad 1\leqslant \alpha +\beta
+\nu \leqslant 2,\quad \beta \leqslant \nu \leqslant 1-\alpha ,
\label{STD-ID.IDD+-1} \\
\frac{b_{1}}{b_{2}}\frac{\sin \frac{\left( \alpha +2\beta -\nu \right) \pi }{%
2}}{\sin \frac{\left( \alpha +\nu \right) \pi }{2}}\leqslant \frac{b_{1}}{%
b_{2}}\frac{\sin \frac{\left( \alpha +2\beta -\nu \right) \pi }{2}}{\sin 
\frac{\left( \alpha +\nu \right) \pi }{2}}\frac{\cos \frac{\left( \alpha
+2\beta -\nu \right) \pi }{2}}{\cos \frac{\left( \alpha +\nu \right) \pi }{2}%
}\leqslant \frac{a_{1}}{a_{2}},  \label{STD-ID.IDD+-2} \\
\frac{a_{1}}{a_{2}}\leqslant \frac{b_{2}}{b_{3}}\frac{\cos \frac{\left( \nu
-\beta \right) \pi }{2}}{\left\vert \cos \frac{\left( 2\alpha +\beta +\nu
\right) \pi }{2}\right\vert }\frac{\sin \frac{\left( \nu -\beta \right) \pi 
}{2}}{\sin \frac{\left( 2\alpha +\beta +\nu \right) \pi }{2}}\leqslant \frac{%
b_{2}}{b_{3}}\frac{\cos \frac{\left( \nu -\beta \right) \pi }{2}}{\left\vert
\cos \frac{\left( 2\alpha +\beta +\nu \right) \pi }{2}\right\vert },
\label{STD-ID.IDD+-3}
\end{gather}%
becomes%
\begin{align}
K\left( \rho \right) & =a_{1}b_{1}\sin \left( \left( \nu -\beta \right) \pi
\right) +a_{1}b_{2}\rho ^{\alpha +\beta }\sin \left( \left( \alpha +\nu
\right) \pi \right) +a_{1}b_{3}\rho ^{2\left( \alpha +\beta \right) }\sin
\left( \left( 2\alpha +\beta +\nu \right) \pi \right)  \notag \\
& \quad -a_{2}b_{1}\rho ^{\alpha +\beta }\sin \left( \left( \alpha +2\beta
-\nu \right) \pi \right) +a_{2}b_{2}\rho ^{2\left( \alpha +\beta \right)
}\sin \left( \left( \nu -\beta \right) \pi \right) +a_{2}b_{3}\rho ^{3\left(
\alpha +\beta \right) }\sin \left( \left( \alpha +\nu \right) \pi \right) ,
\label{K-ID.IDD+}
\end{align}%
so that the non-negativity requirement of function $K$ implies%
\begin{gather}
a_{1}b_{2}\sin \left( \left( \alpha +\nu \right) \pi \right) -a_{2}b_{1}\sin
\left( \left( \alpha +2\beta -\nu \right) \pi \right) \geqslant 0,
\label{NN-of-K-ID.IDD+-1} \\
a_{1}b_{3}\sin \left( \left( 2\alpha +\beta +\nu \right) \pi \right)
+a_{2}b_{2}\sin \left( \left( \nu -\beta \right) \pi \right) \geqslant 0.
\label{NN-of-K-ID.IDD+-2}
\end{gather}

According to the thermodynamical requirement (\ref{STD-ID.IDD+-1})$_{3},$
the first term in (\ref{NN-of-K-ID.IDD+-1}) has a non-negative value, while
according to the thermodynamical requirement (\ref{STD-ID.IDD+-1})$_{1},$
one finds $\alpha +2\beta -\nu =\left( \alpha +\beta -\nu \right) +\beta
\geqslant 0$, as well as according to the requirement (\ref{STD-ID.IDD+-1})$%
_{3},$ one has $\alpha +2\beta -\nu =\left( \alpha +\beta \right) -\left(
\nu -\beta \right) \leqslant 1,$ and therefore the second term in (\ref%
{NN-of-K-ID.IDD+-1}) has a non-positive value, so that the requirement (\ref%
{NN-of-K-ID.IDD+-1}) is transformed into%
\begin{equation*}
\frac{b_{1}}{b_{2}}\frac{\sin \frac{\left( \alpha +2\beta -\nu \right) \pi }{%
2}}{\sin \frac{\left( \alpha +\nu \right) \pi }{2}}\frac{\cos \frac{\left(
\alpha +2\beta -\nu \right) \pi }{2}}{\cos \frac{\left( \alpha +\nu \right)
\pi }{2}}\leqslant \frac{a_{1}}{a_{2}},
\end{equation*}%
representing the narrowed thermodynamical restriction (\ref{STD-ID.IDD+-2})$%
_{1}$, since $\alpha +2\beta -\nu \leqslant \alpha +\nu <1$ reducing to $%
\beta \leqslant \nu ,$ see (\ref{STD-ID.IDD+-1})$_{3}$, implies $\frac{\cos 
\frac{\left( \alpha +2\beta -\nu \right) \pi }{2}}{\cos \frac{\left( \alpha
+\nu \right) \pi }{2}}\geqslant 1$.

The first term in (\ref{NN-of-K-ID.IDD+-2}) has a non-positive value, since\ 
$2\alpha +\beta +\nu =\left( \alpha +\beta +\nu \right) +\alpha \geqslant 1$%
, according to the thermodynamical requirement (\ref{STD-ID.IDD+-1})$_{2}$,
as well as $2\alpha +\beta +\nu =\left( \alpha +\beta \right) +\left( \alpha
+\nu \right) \leqslant 2$ according to the thermodynamical requirement (\ref%
{STD-ID.IDD+-1})$_{3}$, while the second term in (\ref{NN-of-K-ID.IDD+-2})
has a non-negative value, according to the thermodynamical requirement (\ref%
{STD-ID.IDD+-1})$_{3}$. Therefore, the inequality (\ref{NN-of-K-ID.IDD+-2})
is transformed into%
\begin{equation*}
\frac{a_{1}}{a_{2}}\leqslant \frac{b_{2}}{b_{3}}\frac{\cos \frac{\left( \nu
-\beta \right) \pi }{2}}{\left\vert \cos \frac{\left( 2\alpha +\beta +\nu
\right) \pi }{2}\right\vert }\frac{\sin \frac{\left( \nu -\beta \right) \pi 
}{2}}{\sin \frac{\left( 2\alpha +\beta +\nu \right) \pi }{2}},
\end{equation*}%
representing the narrowed thermodynamical restriction (\ref{STD-ID.IDD+-3})$%
_{2}$, since $\frac{\left( \nu -\beta \right) \pi }{2}\leqslant \pi -\frac{%
\left( 2\alpha +\beta +\nu \right) \pi }{2}$ reducing to $\alpha +\nu
\leqslant 1,$ see (\ref{STD-ID.IDD+-1})$_{3},$ implies $\frac{\sin \frac{%
\left( \nu -\beta \right) \pi }{2}}{\sin \frac{\left( 2\alpha +\beta +\nu
\right) \pi }{2}}\leqslant 1$.

In conclusion, the non-negativity of function $K,$ corresponding to the
model ID.IDD$^{{}^{+}}$ and given by (\ref{K-ID.IDD+}), is ensured by
requesting (\ref{STD-ID.IDD+-2}) and (\ref{STD-ID.IDD+-3}) in addition to (%
\ref{STD-ID.IDD+-1}) implying the corresponding relaxation modulus to be a
completely monotonic function and creep compliance a Bernstein function.

\section{Numerical examples\label{nmericali}}

Thermodynamical and narrowed thermodynamical restrictions on the parameters
for the model I$^{+}$ID.ID, taking the form%
\begin{equation}
\left( a_{1}\,{}_{0}\mathrm{I}_{t}^{\alpha +\beta +\nu }+a_{2}\,{}_{0}%
\mathrm{I}_{t}^{\nu }+a_{3}\,{}_{0}\mathrm{D}_{t}^{\alpha +\beta -\nu
}\right) \sigma \left( t\right) =\left( b_{1}\,{}_{0}\mathrm{I}_{t}^{\alpha
}+b_{2}\,{}_{0}\mathrm{D}_{t}^{\beta }\right) \varepsilon \left( t\right) ,
\label{I+ID.ID-ponovljeni}
\end{equation}%
are given as%
\begin{gather}
0\leqslant \alpha +\beta -\nu \leqslant 1,\quad 1\leqslant \alpha +\beta
+\nu \leqslant 2,\quad \alpha \leqslant \nu \leqslant 1-\beta ,
\label{nerovlje-1} \\
\frac{a_{1}}{a_{2}}\frac{\left\vert \cos \frac{\left( \alpha +2\beta +\nu
\right) \pi }{2}\right\vert }{\cos \frac{\left( \nu -\alpha \right) \pi }{2}}%
\leqslant \frac{a_{1}}{a_{2}}\frac{\sin \frac{\left( \alpha +2\beta +\nu
\right) \pi }{2}}{\sin \frac{\left( \nu -\alpha \right) \pi }{2}}\frac{%
\left\vert \cos \frac{\left( \alpha +2\beta +\nu \right) \pi }{2}\right\vert 
}{\cos \frac{\left( \nu -\alpha \right) \pi }{2}}\leqslant \frac{b_{1}}{b_{2}%
},  \label{nerovlje-2} \\
\frac{b_{1}}{b_{2}}\leqslant \frac{a_{2}}{a_{3}}\frac{\sin \frac{\left(
\beta +\nu \right) \pi }{2}}{\sin \frac{\left( 2\alpha +\beta -\nu \right)
\pi }{2}}\frac{\cos \frac{\left( \beta +\nu \right) \pi }{2}}{\cos \frac{%
\left( 2\alpha +\beta -\nu \right) \pi }{2}}\leqslant \frac{a_{2}}{a_{3}}%
\frac{\sin \frac{\left( \beta +\nu \right) \pi }{2}}{\sin \frac{\left(
2\alpha +\beta -\nu \right) \pi }{2}},  \label{nerovlje-3}
\end{gather}%
see also Section \ref{model-I+ID.ID-narrowed} and Appendix \ref{asimetricni
modeli}.

\subsection{Stress relaxation for the model I$^{+}$ID.ID}

The relaxation modulus for the model I$^{+}$ID.ID, given by (\ref%
{I+ID.ID-ponovljeni}), takes the form%
\begin{equation}
\tilde{\sigma}_{sr}\left( s\right) =\frac{1}{s^{1-\left( \beta +\nu \right) }%
}\frac{\phi _{\varepsilon }\left( s\right) }{\phi _{\sigma }\left( s\right) }%
=\frac{1}{s^{1-\left( \beta +\nu \right) }}\frac{b_{1}+b_{2}s^{\alpha +\beta
}}{a_{1}+a_{2}s^{\alpha +\beta }+a_{3}s^{2\left( \alpha +\beta \right) }}
\label{I+ID.ID-LD}
\end{equation}%
in the Laplace domain, after applying the Laplace transform to the
constitutive equation (\ref{I+ID.ID-ponovljeni}), where displacement $%
\varepsilon $ is prescribed as a Heaviside step function, see also (\ref%
{sigma-sr-ld}) and Table \ref{skupina}, so that, after applying the inverse
Laplace transform to (\ref{I+ID.ID-LD}), according to (\ref{sr-opste}), one
obtains the relaxation modulus for the model I$^{+}$ID.ID in the time domain
as%
\begin{equation}
\sigma _{sr}(t)=\sigma _{sr}^{\left( \mathrm{NP}\right) }\left( t\right)
+\left\{ \!\!\!%
\begin{tabular}{ll}
$0$, & if $\tilde{\sigma}_{sr}$ has no poles, \smallskip \\ 
$\sigma _{sr}^{\left( \mathrm{RP}\right) }\left( t\right) $, & if $\tilde{%
\sigma}_{sr}$ has a negative real pole, \smallskip \\ 
$\sigma _{sr}^{\left( \mathrm{CCP}\right) }\left( t\right) $, & if $\tilde{%
\sigma}_{sr}$ has a pair of complex conjugated poles,%
\end{tabular}%
\ \right.  \label{SRM-I+ID.ID}
\end{equation}%
where functions $\sigma _{sr}^{\left( \mathrm{NP}\right) },$ $\sigma
_{sr}^{\left( \mathrm{RP}\right) },$ and $\sigma _{sr}^{\left( \mathrm{CCP}%
\right) }$ are given by%
\begin{align}
\sigma _{sr}^{\left( \mathrm{NP}\right) }\left( t\right) & =\frac{1}{\pi }%
\int_{0}^{\infty }\frac{1}{\rho ^{1-\left( \beta +\nu \right) }}\frac{%
\left\vert b_{1}+b_{2}\rho ^{\alpha +\beta }\mathrm{e}^{\mathrm{i}\left(
\alpha +\beta \right) \pi }\right\vert }{\left\vert a_{1}+a_{2}\rho ^{\alpha
+\beta }\mathrm{e}^{\mathrm{i}\left( \alpha +\beta \right) \pi }+a_{3}\rho
^{2\left( \alpha +\beta \right) }\mathrm{e}^{2\mathrm{i}\left( \alpha +\beta
\right) \pi }\right\vert }  \notag \\
& \qquad \times \sin \left( \left( \beta +\nu \right) \pi +\arctan \frac{%
b_{2}\rho ^{\alpha +\beta }\sin \left( \left( \alpha +\beta \right) \pi
\right) }{b_{1}+b_{2}\rho ^{\alpha +\beta }\cos \left( \left( \alpha +\beta
\right) \pi \right) }\right.  \notag \\
& \qquad \qquad -\left. \arctan \frac{a_{2}\rho ^{\alpha +\beta }\sin \left(
\left( \alpha +\beta \right) \pi \right) +a_{3}\rho ^{2\left( \alpha +\beta
\right) }\sin \left( 2\left( \alpha +\beta \right) \pi \right) }{%
a_{1}+a_{2}\rho ^{\alpha +\beta }\cos \left( \left( \alpha +\beta \right)
\pi \right) +a_{3}\rho ^{2\left( \alpha +\beta \right) }\cos \left( 2\left(
\alpha +\beta \right) \pi \right) }\right) \mathrm{e}^{-\rho t}\mathrm{d}%
\rho ,  \label{sr-np-numer} \\
\sigma _{sr}^{\left( \mathrm{NP}\right) }\left( t\right) & =\frac{1}{\pi }%
\int_{0}^{\infty }\frac{1}{\rho ^{1-\left( \beta +\nu \right) }}\frac{%
K\left( \rho \right) }{\left\vert a_{1}+a_{2}\rho ^{\alpha +\beta }\mathrm{e}%
^{\mathrm{i}\left( \alpha +\beta \right) \pi }+a_{3}\rho ^{2\left( \alpha
+\beta \right) }\mathrm{e}^{2\mathrm{i}\left( \alpha +\beta \right) \pi
}\right\vert ^{2}}\mathrm{e}^{-\rho t}\mathrm{d}\rho ,  \notag \\
\sigma _{sr}^{\left( \mathrm{RP}\right) }\left( t\right) & =-\rho _{%
\scriptscriptstyle{\mathrm{RP}}}^{\beta -\mu }\frac{\left\vert
b_{1}+b_{2}\rho _{\scriptscriptstyle{\mathrm{RP}}}^{\alpha +\beta }\,\mathrm{%
e}^{\mathrm{i}\left( \alpha +\beta \right) \pi }\right\vert }{\left( \alpha
+\beta \right) \left\vert a_{2}+2a_{3}\rho _{\scriptscriptstyle{\mathrm{RP}}%
}^{\alpha +\beta }\,\mathrm{e}^{\mathrm{i}\left( \alpha +\beta \right) \pi
}\right\vert }  \notag \\
& \qquad \times \cos \left( \left( \beta +\nu \right) \pi +\arctan \frac{%
b_{2}\rho _{\scriptscriptstyle{\mathrm{RP}}}^{\alpha +\beta }\sin \left(
\left( \alpha +\beta \right) \pi \right) }{b_{1}+b_{2}\rho _{%
\scriptscriptstyle{\mathrm{RP}}}^{\alpha +\beta }\cos \left( \left( \alpha
+\beta \right) \pi \right) }\right.  \notag \\
& \qquad \qquad -\left. \arctan \frac{a_{2}\sin \left( \left( \alpha +\beta
\right) \pi \right) +2a_{3}\rho _{\scriptscriptstyle{\mathrm{RP}}}^{\alpha
+\beta }\sin \left( 2\left( \alpha +\beta \right) \pi \right) }{a_{2}\cos
\left( \left( \alpha +\beta \right) \pi \right) +2a_{3}\rho _{%
\scriptscriptstyle{\mathrm{RP}}}^{\alpha +\beta }\cos \left( 2\left( \alpha
+\beta \right) \pi \right) }\right) \mathrm{e}^{-\rho _{\scriptscriptstyle{%
\mathrm{RP}}}t},  \label{sr-rp-numer} \\
\sigma _{sr}^{\left( \mathrm{CCP}\right) }\left( t\right) & =2\rho _{%
\scriptscriptstyle{\mathrm{CCP}}}^{\beta -\mu }\frac{\left\vert
b_{1}+b_{2}\rho _{\scriptscriptstyle{\mathrm{CCP}}}^{\alpha +\beta }\,%
\mathrm{e}^{\mathrm{i}\left( \alpha +\beta \right) \varphi _{%
\scriptscriptstyle{\mathrm{CCP}}}}\right\vert }{\left( \alpha +\beta \right)
\left\vert a_{2}+2a_{3}\rho _{\scriptscriptstyle{\mathrm{CCP}}}^{\alpha
+\beta }\,\mathrm{e}^{\mathrm{i}\left( \alpha +\beta \right) \varphi _{%
\scriptscriptstyle{\mathrm{CCP}}}}\right\vert }  \notag \\
& \qquad \times \mathrm{e}^{-\left\vert \func{Re}s_{\scriptscriptstyle{%
\mathrm{CCP}}}\right\vert t}\cos \left( \func{Im}s_{\scriptscriptstyle{%
\mathrm{CCP}}}t-\left( 1-\left( \beta +\nu \right) \right) \varphi _{%
\scriptscriptstyle{\mathrm{CCP}}}+\arctan \frac{b_{2}\rho _{%
\scriptscriptstyle{\mathrm{CCP}}}^{\alpha +\beta }\sin \left( \left( \alpha
+\beta \right) \varphi _{\scriptscriptstyle{\mathrm{CCP}}}\right) }{%
b_{1}+b_{2}\rho _{\scriptscriptstyle{\mathrm{CCP}}}^{\alpha +\beta }\cos
\left( \left( \alpha +\beta \right) \varphi _{\scriptscriptstyle{\mathrm{CCP}%
}}\right) }\right.  \notag \\
& \qquad \qquad +\left. \arctan \frac{a_{2}\sin \left( \left( 1-\left(
\alpha +\beta \right) \right) \varphi _{\scriptscriptstyle{\mathrm{CCP}}%
}\right) +2a_{3}\rho _{\scriptscriptstyle{\mathrm{CCP}}}^{\alpha +\beta
}\sin \left( \left( 1-2\left( \alpha +\beta \right) \right) \varphi _{%
\scriptscriptstyle{\mathrm{CCP}}}\right) }{a_{2}\cos \left( \left( 1-\left(
\alpha +\beta \right) \right) \varphi _{\scriptscriptstyle{\mathrm{CCP}}%
}\right) +2a_{3}\rho _{\scriptscriptstyle{\mathrm{CCP}}}^{\alpha +\beta
}\cos \left( \left( 1-2\left( \alpha +\beta \right) \right) \varphi _{%
\scriptscriptstyle{\mathrm{CCP}}}\right) }\right)  \label{sr-ccp-numer}
\end{align}%
according to (\ref{sigma-NP}), (\ref{sigma-NP-1}), (\ref{sigma-RP}), and (%
\ref{sigma-CCP}) respectively, with%
\begin{equation}
\phi _{\varepsilon }\left( s\right) =b_{1}+b_{2}s^{\alpha +\beta }\quad 
\text{and}\quad \phi _{\sigma }\left( s\right) =a_{1}+a_{2}s^{\alpha +\beta
}+a_{3}s^{2\left( \alpha +\beta \right) },  \label{fi-sigma-I+ID.ID}
\end{equation}%
where the function $K,$ defined by (\ref{K}), is represented by%
\begin{align}
K\left( \rho \right) & =a_{1}b_{1}\sin \left( \left( \beta +\nu \right) \pi
\right) +\rho ^{\alpha +\beta }\left( a_{1}b_{2}\sin \left( \left( \alpha
+2\beta +\nu \right) \pi \right) +a_{2}b_{1}\sin \left( \left( \nu -\alpha
\right) \pi \right) \right)  \notag \\
& \quad +\rho ^{2\left( \alpha +\beta \right) }\left( a_{2}b_{2}\sin \left(
\left( \beta +\nu \right) \pi \right) -a_{3}b_{1}\sin \left( \left( 2\alpha
+\beta -\nu \right) \pi \right) \right) +\rho ^{3\left( \alpha +\beta
\right) }a_{3}b_{2}\sin \left( \left( \nu -\alpha \right) \pi \right) ,
\label{K-I+ID.ID-num}
\end{align}%
and with poles of the relaxation modulus in Laplace domain $\rho _{%
\scriptscriptstyle{\mathrm{RP}}}=-\sqrt[2\left( \alpha +\beta \right) ]{%
\frac{a_{1}}{a_{3}}},$ if $\tan \left( \left( \alpha +\beta \right) \pi
\right) =-\sqrt{\frac{4a_{1}a_{3}}{a_{2}^{2}}-1}$ and $s_{\scriptscriptstyle{%
\mathrm{CCP}}}=\rho _{\scriptscriptstyle{\mathrm{CCP}}}\,\mathrm{e}^{\mathrm{%
i}\varphi _{\scriptscriptstyle{\mathrm{CCP}}}},$ with $\rho _{%
\scriptscriptstyle{\mathrm{CCP}}}=\sqrt[2\left( \alpha +\beta \right) ]{%
\frac{a_{1}}{a_{3}}}$ and $\varphi _{\scriptscriptstyle{\mathrm{CCP}}%
}=\left( 1-\frac{1}{\pi }\arctan \sqrt{\frac{4a_{1}a_{3}}{a_{2}^{2}}-1}%
\right) \frac{\pi }{\alpha +\beta },$ if $\tan \left( \left( \alpha +\beta
\right) \pi \right) >-\sqrt{\frac{4a_{1}a_{3}}{a_{2}^{2}}-1},$ obtained in
Section \ref{nule-trinom} as zeros of function $\phi _{\sigma },$ given by (%
\ref{fi-sigma-I+ID.ID})$_{2},$ by (\ref{nule-rp}) and (\ref{nule-ccp}).

Considering the model I$^{+}$ID.ID, the asymptotic behavior of the
corresponding relaxation modulus in the Laplace domain, according to (\ref%
{I+ID.ID-LD}), takes the form%
\begin{align*}
\tilde{\sigma}_{sr}\left( s\right) & =\frac{b_{2}}{a_{3}}\frac{1}{%
s^{1-\left( \nu -\alpha \right) }}\frac{1+\frac{b_{1}}{b_{2}}\frac{1}{%
s^{\alpha +\beta }}}{1+\frac{a_{2}}{a_{3}}\frac{1}{s^{\alpha +\beta }}+\frac{%
a_{1}}{a_{3}}\frac{1}{s^{2\left( \alpha +\beta \right) }}} \\
& =\frac{b_{2}}{a_{3}}\frac{1}{s^{1-\left( \nu -\alpha \right) }}\left( 1+%
\frac{b_{1}}{b_{2}}\frac{1}{s^{\alpha +\beta }}\right) \left( 1+\frac{a_{2}}{%
a_{3}}\frac{1}{s^{\alpha +\beta }}+\frac{a_{1}}{a_{3}}\frac{1}{s^{2\left(
\alpha +\beta \right) }}\right) ^{-1} \\
& =\frac{b_{2}}{a_{3}}\frac{1}{s^{1-\left( \nu -\alpha \right) }}\left( 1+%
\frac{b_{1}}{b_{2}}\frac{1}{s^{\alpha +\beta }}\right) \left( 1-\frac{a_{2}}{%
a_{3}}\frac{1}{s^{\alpha +\beta }}+\left( -\frac{a_{1}}{a_{3}}+\left( \frac{%
a_{2}}{a_{3}}\right) ^{2}\right) \frac{1}{s^{2\left( \alpha +\beta \right) }}%
+O\left( s^{-3\left( \alpha +\beta \right) }\right) \right) \\
& =\frac{b_{2}}{a_{3}}\frac{1}{s^{1-\left( \nu -\alpha \right) }}\left(
1+\left( \frac{b_{1}}{b_{2}}-\frac{a_{2}}{a_{3}}\right) \frac{1}{s^{\alpha
+\beta }}+\left( \left( \frac{a_{2}}{a_{3}}\right) ^{2}-\frac{a_{1}}{a_{3}}-%
\frac{a_{2}}{a_{3}}\frac{b_{1}}{b_{2}}\right) \frac{1}{s^{2\left( \alpha
+\beta \right) }}+O\left( s^{-3\left( \alpha +\beta \right) }\right) \right)
\\
& =\frac{b_{2}}{a_{3}}\frac{1}{s^{1-\left( \nu -\alpha \right) }}+\frac{b_{2}%
}{a_{3}}\left( \frac{b_{1}}{b_{2}}-\frac{a_{2}}{a_{3}}\right) \frac{1}{%
s^{1+2\alpha +\beta -\nu }} \\
& \quad +\frac{b_{2}}{a_{3}}\left( \left( \frac{a_{2}}{a_{3}}\right) ^{2}-%
\frac{a_{1}}{a_{3}}-\frac{a_{2}}{a_{3}}\frac{b_{1}}{b_{2}}\right) \frac{1}{%
s^{1+3\alpha +2\beta -\nu }}+O\left( s^{-\left( 1+4\alpha +3\beta -\nu
\right) }\right) ,\quad \text{when}\quad s\rightarrow \infty \text{,}
\end{align*}%
that implies the short time asymptotics of the relaxation modulus as%
\begin{align}
\sigma _{sr}\left( t\right) & =\frac{b_{2}}{a_{3}}\frac{t^{-\left( \nu
-\alpha \right) }}{\Gamma \left( 1-\left( \nu -\alpha \right) \right) }+%
\frac{b_{2}}{a_{3}}\left( \frac{b_{1}}{b_{2}}-\frac{a_{2}}{a_{3}}\right) 
\frac{t^{2\alpha +\beta -\nu }}{\Gamma \left( 1+2\alpha +\beta -\nu \right) }
\notag \\
& \quad +\frac{b_{2}}{a_{3}}\left( \left( \frac{a_{2}}{a_{3}}\right) ^{2}-%
\frac{a_{1}}{a_{3}}-\frac{a_{2}}{a_{3}}\frac{b_{1}}{b_{2}}\right) \frac{%
t^{3\alpha +2\beta -\nu }}{\Gamma \left( 1+3\alpha +2\beta -\nu \right) }%
+O\left( t^{-\left( 4\alpha +3\beta -\nu \right) }\right) ,\quad \text{when}%
\quad t\rightarrow 0,  \label{sr-sha}
\end{align}%
according to the theorem that if $\tilde{f}\left( s\right) \sim \tilde{g}%
\left( s\right) $ as $s\rightarrow \infty ,$ then $f\left( t\right) \sim
g\left( t\right) $ as $t\rightarrow 0,$ while the large time asymptotics of
the relaxation modulus%
\begin{equation}
\sigma _{sr}\left( t\right) =\frac{b_{1}}{a_{1}}\frac{t^{-\left( \beta +\nu
\right) }}{\Gamma \left( 1-\left( \beta +\nu \right) \right) }+O\left(
t^{-\left( 1-\delta \right) }\right) ,\quad \text{when}\quad t\rightarrow
\infty ,  \label{sr-la}
\end{equation}%
follows from the asymptotics of relaxation modulus in the Laplace domain (%
\ref{I+ID.ID-LD}), obtained as 
\begin{align*}
\tilde{\sigma}_{sr}\left( s\right) & =\frac{b_{1}}{a_{1}}\frac{1}{%
s^{1-\left( \beta +\nu \right) }}\frac{1+\frac{b_{2}}{b_{1}}s^{\alpha +\beta
}}{1+\frac{a_{2}}{a_{1}}s^{\alpha +\beta }+\frac{a_{3}}{a_{1}}s^{2\left(
\alpha +\beta \right) }} \\
& =\frac{b_{1}}{a_{1}}\frac{1}{s^{1-\left( \beta +\nu \right) }}\left( 1+%
\frac{b_{2}}{b_{1}}s^{\alpha +\beta }\right) \left( 1+\frac{a_{2}}{a_{1}}%
s^{\alpha +\beta }+\frac{a_{3}}{a_{1}}s^{2\left( \alpha +\beta \right)
}\right) ^{-1} \\
& =\frac{b_{1}}{a_{1}}\frac{1}{s^{1-\left( \beta +\nu \right) }}\left( 1+%
\frac{b_{2}}{b_{1}}s^{\alpha +\beta }\right) \left( 1-\frac{a_{2}}{a_{1}}%
s^{\alpha +\beta }+O\left( s^{2\left( \alpha +\beta \right) }\right) \right)
\\
& =\frac{b_{1}}{a_{1}}\frac{1}{s^{1-\left( \beta +\nu \right) }}\left(
1+\left( \frac{b_{2}}{b_{1}}-\frac{a_{2}}{a_{1}}\right) s^{\alpha +\beta
}+O\left( s^{2\left( \alpha +\beta \right) }\right) \right) \\
& =\frac{b_{1}}{a_{1}}\frac{1}{s^{1-\left( \beta +\nu \right) }}+O\left(
s^{-\delta }\right) ,\quad \text{when}\quad s\rightarrow 0,
\end{align*}%
with $0<\delta <1-\left( \beta +\nu \right) $, since $\alpha +2\beta +\nu
>1, $ according to the theorem that if $\tilde{f}\left( s\right) \sim \tilde{%
g}\left( s\right) $ as $s\rightarrow 0,$ then $f\left( t\right) \sim g\left(
t\right) $ as $t\rightarrow \infty .$

\begin{figure}[h]
\begin{center}
\begin{minipage}{0.46\columnwidth}
  \subfloat[Case when $\tilde{\sigma}_{sr}$ has no poles. Relaxation modulus is completely monotonic.]{
   \includegraphics[width=\columnwidth]{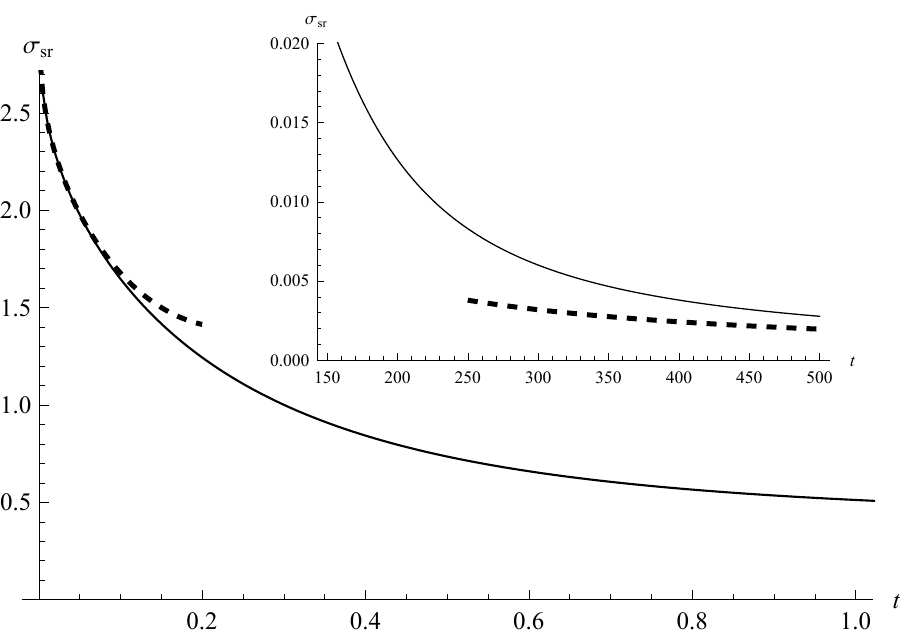}
   \label{sr-np}}
  \end{minipage}
\hfil
\begin{minipage}{0.46\columnwidth}
  \subfloat[Case when $\tilde{\sigma}_{sr}$ has a negative real pole. Complete monotonicity of relaxation modulus is not guaranteed.]{
   \includegraphics[width=\columnwidth]{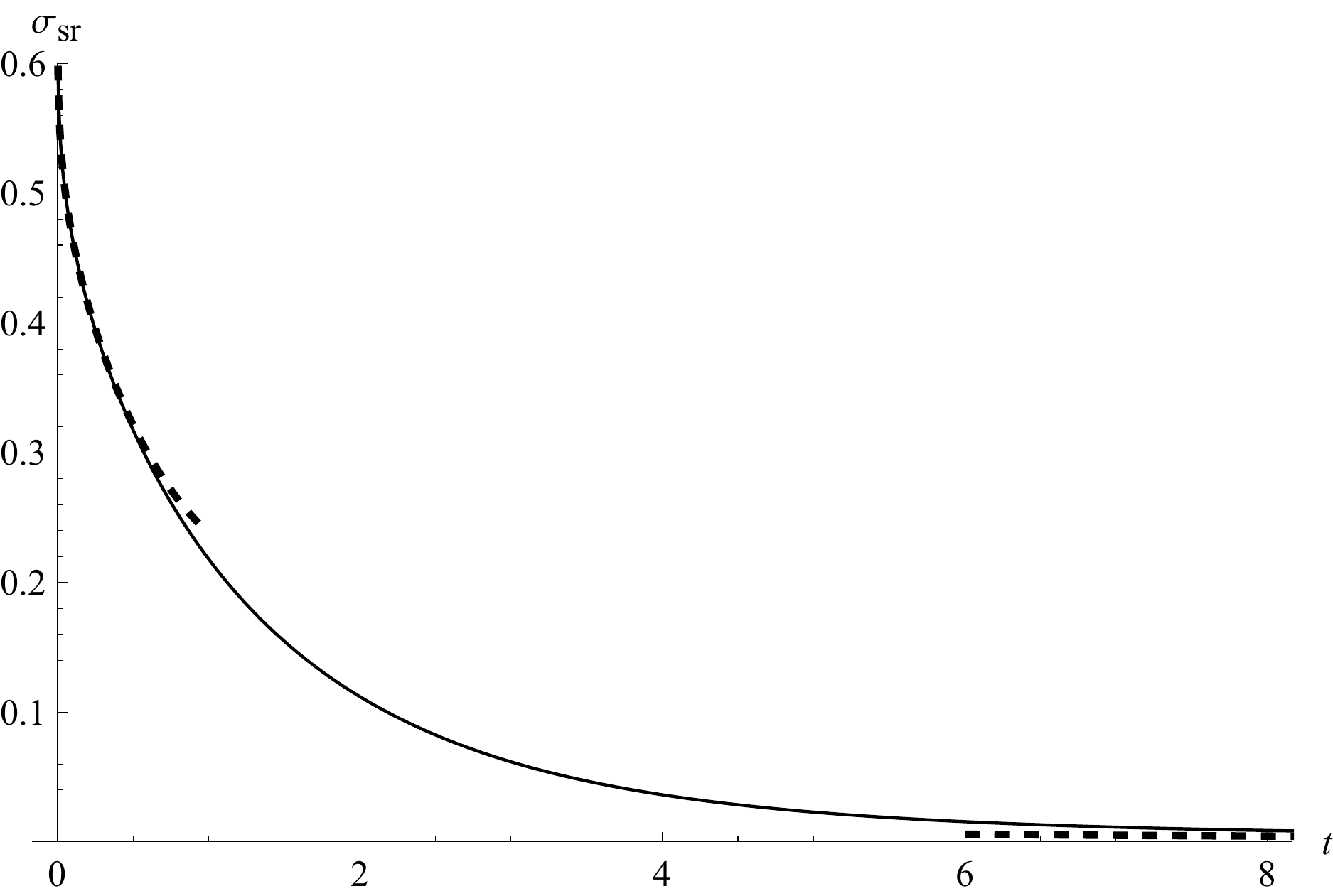}
   \label{sr-rp}}
  \end{minipage}
\vfil  
\begin{minipage}{0.46\columnwidth}
  \subfloat[Case when $\tilde{\sigma}_{sr}$ has a pair of complex conjugated poles.]{
   \includegraphics[width=\columnwidth]{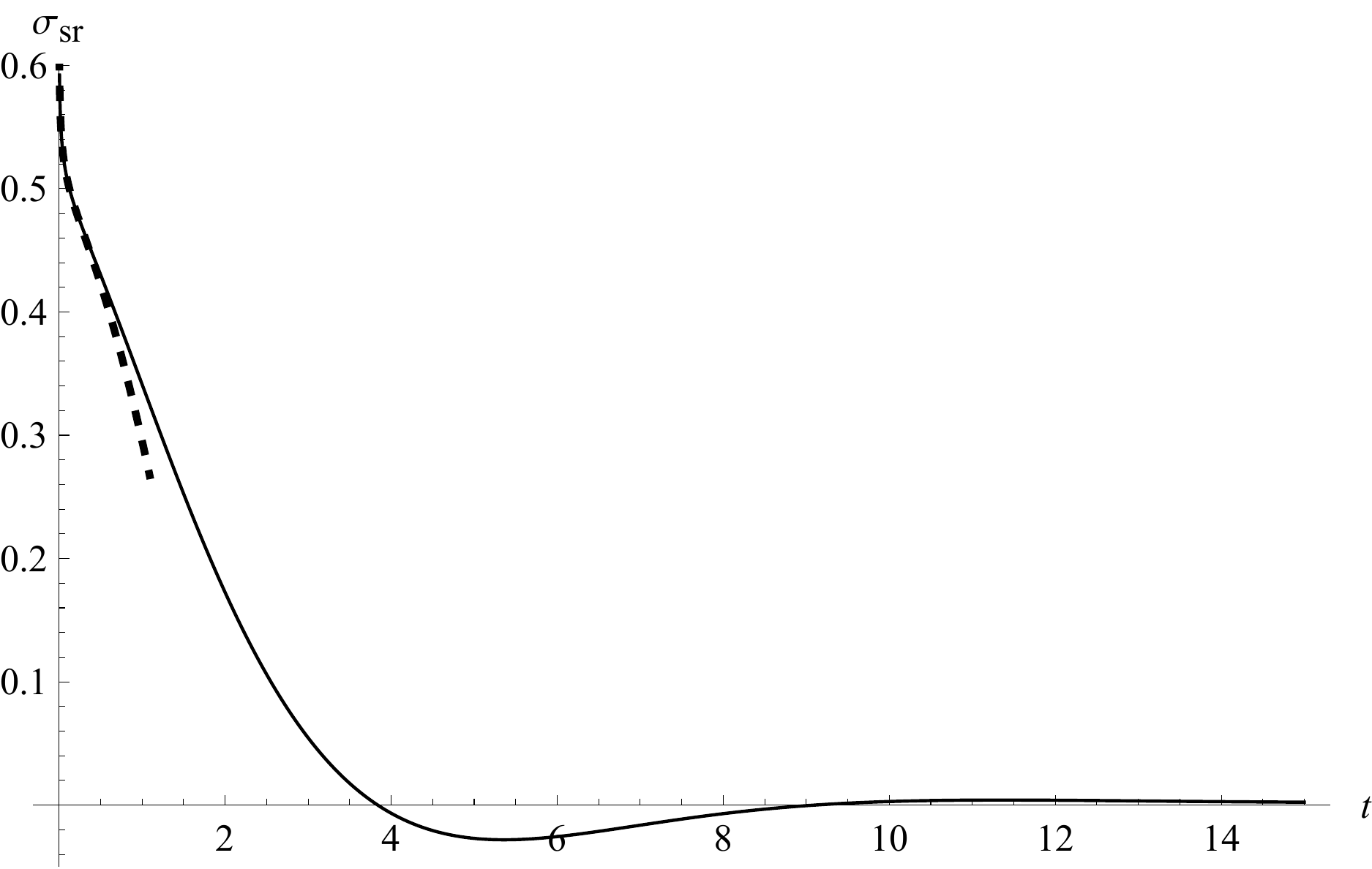}
   \label{sr-ccp-srsha}}
  \end{minipage}
\hfil
\begin{minipage}{0.46\columnwidth}
  \subfloat[Case when $\tilde{\sigma}_{sr}$ has a pair of complex conjugated poles.]{
   \includegraphics[width=\columnwidth]{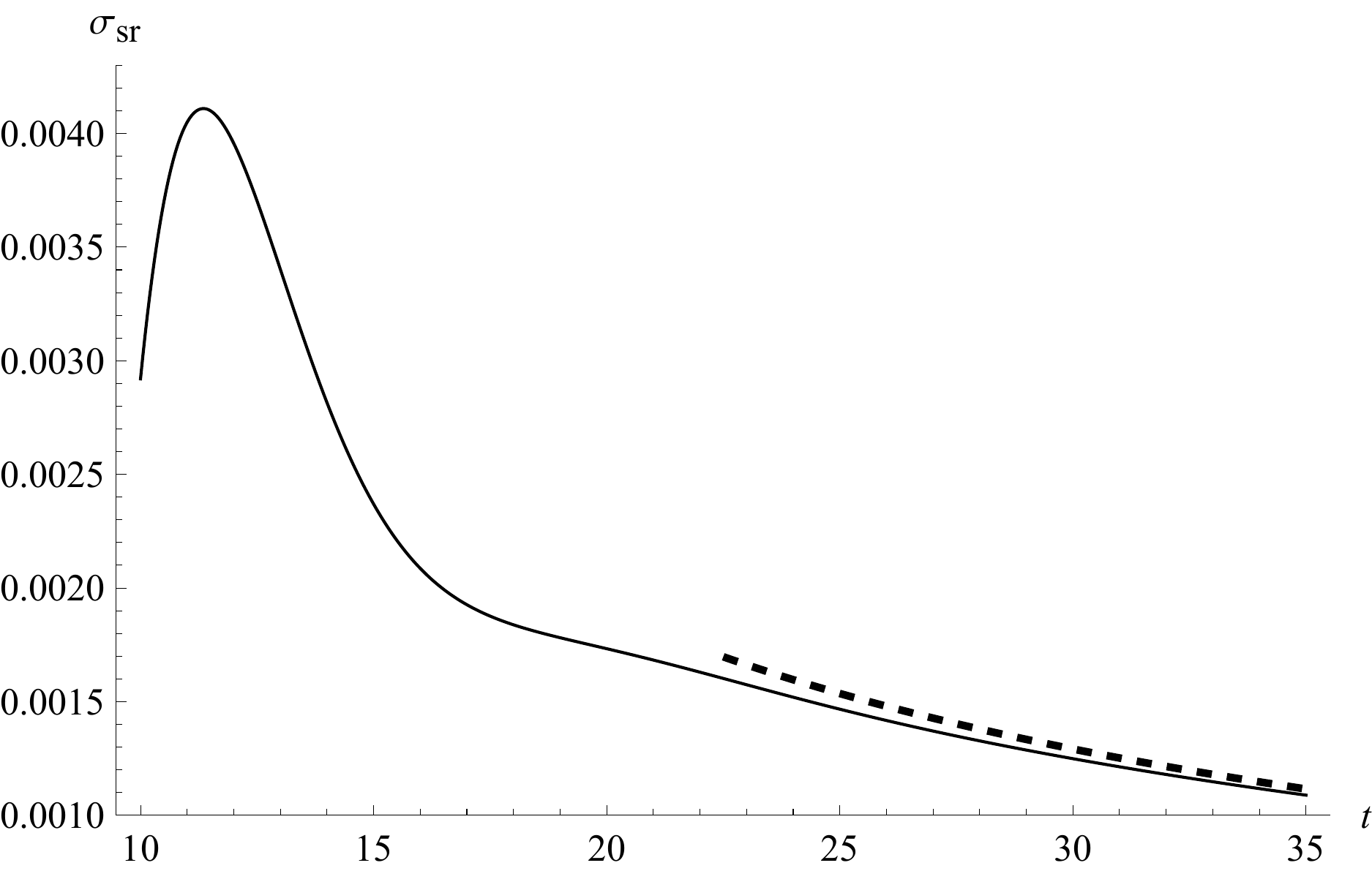}
   \label{sr-ccp-srla}}
  \end{minipage}
\end{center}
\caption{Time evolution of the relaxation modulus. }
\label{sr}
\end{figure}

Figure \ref{sr} presents time evolution of the relaxation modulus for
different values of model parameters given in Table \ref{parametri}, so that
in the case when model parameters are such that the relaxation modulus in
Laplace domain does not have poles and additionally satisfy narrowed
thermodynamical requirements (\ref{nerovlje-1}) - (\ref{nerovlje-3}), the
relaxation modulus from Figure \ref{sr-np}, obtained according to (\ref%
{SRM-I+ID.ID}) with function $\sigma _{sr}^{\left( \mathrm{NP}\right) }$
given by (\ref{sr-np-numer}), is completely monotonic, while the complete
monotonicity of the relaxation modulus from Figure \ref{sr-rp}, obtained
according to (\ref{SRM-I+ID.ID}) with function $\sigma _{sr}^{\left( \mathrm{%
RP}\right) }$ given by (\ref{sr-rp-numer}), is not guaranteed, although its
form is seemingly such, since model parameters do not satisfy narrowed
thermodynamical requirements (\ref{nerovlje-1}) - (\ref{nerovlje-3}).
Figures \ref{sr-ccp-srsha} and \ref{sr-ccp-srla}, obtained according to (\ref%
{SRM-I+ID.ID}) with function $\sigma _{sr}^{\left( \mathrm{CCP}\right) }$
given by (\ref{sr-ccp-numer}), present oscillatory relaxation modulus having
amplitude exponentially decreasing in time, displaying the significant
damping and the negative stress for a certain time interval, however stress
being positive and decreasing to zero for large time. Short and long time
behavior of the relaxation modulus, as obvious from all graphs from Figure %
\ref{sr}, is governed by the asymptotic expressions (\ref{sr-sha}) and (\ref%
{sr-la}). \input{parametri.tex}

\subsection{Creep compliance for the model I$^{+}$ID.ID}

The creep compliance for the model I$^{+}$ID.ID, given by (\ref{I+ID.ID}),
takes the form%
\begin{equation}
\tilde{\varepsilon}_{cr}\left( s\right) =\frac{1}{s^{1+\beta +\nu }}\frac{%
\phi _{\sigma }\left( s\right) }{\phi _{\varepsilon }\left( s\right) }=\frac{%
1}{s^{1+\beta +\nu }}\frac{a_{1}+a_{2}s^{\alpha +\beta }+a_{3}s^{2\left(
\alpha +\beta \right) }}{b_{1}+b_{2}s^{\alpha +\beta }}
\label{I+ID.ID-LD-cr}
\end{equation}%
in the Laplace domain, after applying the Laplace transform to the
constitutive equation (\ref{I+ID.ID-ponovljeni}), where stress $\sigma $ is
prescribed as the Heaviside step function, see also (\ref{epsilon-cr-ld})
and Table \ref{skupina}, so that, after applying the inverse Laplace
transform to (\ref{I+ID.ID-LD-cr}), according to (\ref{cr-opste}), one
obtains the creep compliance for the model I$^{+}$ID.ID in the time domain as%
\begin{align}
\varepsilon _{cr}(t)& =\varepsilon _{cr}^{\left( \mathrm{NP}\right) }\left(
t\right) =\frac{1}{\pi }\int_{0}^{\infty }\frac{1}{\rho ^{1+\beta +\nu }}%
\frac{\left\vert a_{1}+a_{2}\rho ^{\alpha +\beta }\mathrm{e}^{\mathrm{i}%
\left( \alpha +\beta \right) \pi }+a_{3}\rho ^{2\left( \alpha +\beta \right)
}\mathrm{e}^{2\mathrm{i}\left( \alpha +\beta \right) \pi }\right\vert }{%
\left\vert b_{1}+b_{2}\rho ^{\alpha +\beta }\mathrm{e}^{\mathrm{i}\left(
\alpha +\beta \right) \pi }\right\vert }  \notag \\
& \qquad \times \sin \left( \left( \beta +\nu \right) \pi +\arctan \frac{%
b_{2}\rho ^{\alpha +\beta }\sin \left( \left( \alpha +\beta \right) \pi
\right) }{b_{1}+b_{2}\rho ^{\alpha +\beta }\cos \left( \left( \alpha +\beta
\right) \pi \right) }\right.  \notag \\
& \qquad \qquad -\left. \arctan \frac{a_{2}\rho ^{\alpha +\beta }\sin \left(
\left( \alpha +\beta \right) \pi \right) +a_{3}\rho ^{2\left( \alpha +\beta
\right) }\sin \left( 2\left( \alpha +\beta \right) \pi \right) }{%
a_{1}+a_{2}\rho ^{\alpha +\beta }\cos \left( \left( \alpha +\beta \right)
\pi \right) +a_{3}\rho ^{2\left( \alpha +\beta \right) }\cos \left( 2\left(
\alpha +\beta \right) \pi \right) }\right) \left( 1-\mathrm{e}^{-\rho
t}\right) \mathrm{d}\rho ,  \label{cr-preko-fiova}
\end{align}%
or equivalently as 
\begin{equation}
\varepsilon _{cr}\left( t\right) =\varepsilon _{cr}^{\left( \mathrm{NP}%
\right) }\left( t\right) =\frac{1}{\pi }\int_{0}^{\infty }\frac{1}{\rho
^{1+\beta +\nu }}\frac{K\left( \rho \right) }{\left\vert b_{1}+b_{2}\rho
^{\alpha +\beta }\mathrm{e}^{\mathrm{i}\left( \alpha +\beta \right) \pi
}\right\vert ^{2}}\left( 1-\mathrm{e}^{-\rho t}\right) \mathrm{d}\rho ,
\label{cr-preko-K}
\end{equation}%
where the equivalent forms of function $\varepsilon _{cr}^{\left( \mathrm{NP}%
\right) }$ are given by (\ref{epsilon-NP}) and (\ref{epsilon-NP-1}) and
where function $K$ is given by (\ref{K-I+ID.ID-num}), since function $\phi
_{\varepsilon },$ given by (\ref{fi-sigma-I+ID.ID})$_{1},$ has no zeros in
the complex plane, implying that the creep compliance in Laplace domain
admits no poles.

Rather than using any of the expressions (\ref{cr-preko-fiova}) and (\ref%
{cr-preko-K}), for the purpose of better numerical convergence, one rewrites
the creep compliance in Laplace domain (\ref{I+ID.ID-LD-cr}) as 
\begin{align*}
\tilde{\epsilon}_{cr}\left( s\right) &=s\tilde{\varepsilon}_{cr}\left(
s\right) =\frac{1}{s^{\beta +\nu }}\left( \frac{a_{1}}{b_{2}}\frac{1}{%
s^{\alpha +\beta }+\frac{b_{1}}{b_{2}}}+\frac{a_{3}}{b_{2}}s^{\alpha +\beta }%
\frac{s^{\alpha +\beta }+\frac{a_{2}}{a_{3}}}{s^{\alpha +\beta }+\frac{b_{1}%
}{b_{2}}}\right) \\
&=\frac{a_{3}}{b_{2}}\frac{1}{s^{\nu -\alpha }}+\frac{a_{3}}{b_{2}}\left( 
\frac{a_{2}}{a_{3}}-\frac{b_{1}}{b_{2}}\right) \frac{1}{s^{\nu -\alpha }}%
\frac{1}{s^{\alpha +\beta }+\frac{b_{1}}{b_{2}}}+\frac{a_{1}}{b_{2}}\frac{1}{%
s^{\beta +\nu }}\frac{1}{s^{\alpha +\beta }+\frac{b_{1}}{b_{2}}},
\end{align*}%
yielding the creep compliance in time domain in the form%
\begin{align}
\varepsilon _{cr}\left( t\right) &=\int_{0}^{t}\epsilon _{cr}\left(
t^{\prime }\right) \mathrm{d}t^{\prime }  \notag \\
&=\frac{a_{3}}{b_{2}}\frac{t^{-\left( 1-\left( \nu -\alpha \right) \right) }%
}{\Gamma \left( \nu -\alpha \right) }+\frac{1}{\pi }\frac{a_{3}}{b_{2}}%
\left( \frac{a_{2}}{a_{3}}-\frac{b_{1}}{b_{2}}\right) \int_{0}^{\infty }%
\frac{1}{\rho ^{1+\nu -\alpha }}\frac{\rho ^{\alpha +\beta }\sin \left(
\left( \beta +\nu \right) \pi \right) +\frac{b_{1}}{b_{2}}\sin \left( \left(
\nu -\alpha \right) \pi \right) }{\left\vert \rho ^{\alpha +\beta }\mathrm{e}%
^{\mathrm{i}\left( \alpha +\beta \right) \pi }+\frac{b_{1}}{b_{2}}%
\right\vert ^{2}}\left( 1-\mathrm{e}^{-\rho t}\right) \mathrm{d}\rho  \notag
\\
&+\frac{1}{\pi }\frac{a_{1}}{b_{2}}\int_{0}^{\infty }\frac{1}{\rho ^{1+\beta
+\nu }}\frac{\rho ^{\alpha +\beta }\sin \left( \left( \alpha +2\beta +\nu
\right) \pi \right) +\frac{b_{1}}{b_{2}}\sin \left( \left( \beta +\nu
\right) \pi \right) }{\left\vert \rho ^{\alpha +\beta }\mathrm{e}^{\mathrm{i}%
\left( \alpha +\beta \right) \pi }+\frac{b_{1}}{b_{2}}\right\vert ^{2}}%
\left( 1-\mathrm{e}^{-\rho t}\right) \mathrm{d}\rho ,
\label{cr-preko-jedinice}
\end{align}%
since the definition of inverse Laplace transform (\ref{inv-laplas-epsilon}%
), used along with the Cauchy integral theorem (\ref%
{Kosijeva-teorema-epsilon}) where the contour is depicted in Figure \ref%
{nemaTG}, implies%
\begin{equation*}
\mathcal{L}^{-1}\left[ \frac{1}{s^{\xi }}\frac{1}{s^{\alpha +\beta }+\frac{%
b_{1}}{b_{2}}}\right] \left( t\right) =\frac{1}{\pi }\int_{0}^{\infty }\frac{%
1}{\rho ^{\xi }}\frac{\rho ^{\alpha +\beta }\sin \left( \left( \alpha +\beta
+\xi \right) \pi \right) +\frac{b_{1}}{b_{2}}\sin \left( \xi \pi \right) }{%
\left\vert \rho ^{\alpha +\beta }\mathrm{e}^{\mathrm{i}\left( \alpha +\beta
\right) \pi }+\frac{b_{1}}{b_{2}}\right\vert ^{2}}\mathrm{e}^{-\rho t}%
\mathrm{d}\rho .
\end{equation*}%
Equivalently, the creep compliance is also derived in terms of two-parameter
Mittag-Leffler function, see (\ref{cr-ML}), that in the case of model I$^{+}$%
ID.ID yields%
\begin{equation}
\varepsilon _{cr}\left( t\right) =\frac{a_{1}}{b_{2}}e_{\alpha +\beta
,1+\alpha +2\beta +\nu ,\frac{b_{1}}{b_{2}}}\left( t\right) +\frac{a_{2}}{%
b_{2}}e_{\alpha +\beta ,1+\beta +\nu ,\frac{b_{1}}{b_{2}}}\left( t\right) +%
\frac{a_{3}}{b_{2}}e_{\alpha +\beta ,1+\nu -\alpha ,\frac{b_{1}}{b_{2}}%
}\left( t\right) .  \label{cr-preko-ML}
\end{equation}

Considering the model I$^{+}$ID.ID, the asymptotic behavior of the
corresponding creep compliance in Laplace domain, according to (\ref%
{I+ID.ID-LD-cr}), takes the form%
\begin{align*}
\tilde{\varepsilon}_{cr}\left( s\right) & =\frac{a_{3}}{b_{2}}\frac{1}{%
s^{1+\nu -\alpha }}\frac{1+\frac{a_{2}}{a_{3}}\frac{1}{s^{\alpha +\beta }}+%
\frac{a_{1}}{a_{3}}\frac{1}{s^{2\left( \alpha +\beta \right) }}}{1+\frac{%
b_{1}}{b_{2}}\frac{1}{s^{\alpha +\beta }}} \\
& =\frac{a_{3}}{b_{2}}\frac{1}{s^{1+\nu -\alpha }}\left( 1+\frac{a_{2}}{a_{3}%
}\frac{1}{s^{\alpha +\beta }}+\frac{a_{1}}{a_{3}}\frac{1}{s^{2\left( \alpha
+\beta \right) }}\right) \left( 1+\frac{b_{1}}{b_{2}}\frac{1}{s^{\alpha
+\beta }}\right) ^{-1} \\
& =\frac{a_{3}}{b_{2}}\frac{1}{s^{1+\nu -\alpha }}\left( 1+\frac{a_{2}}{a_{3}%
}\frac{1}{s^{\alpha +\beta }}+\frac{a_{1}}{a_{3}}\frac{1}{s^{2\left( \alpha
+\beta \right) }}\right) \left( 1-\frac{b_{1}}{b_{2}}\frac{1}{s^{\alpha
+\beta }}+\left( \frac{b_{1}}{b_{2}}\right) ^{2}\frac{1}{s^{2\left( \alpha
+\beta \right) }}+O\left( s^{-3\left( \alpha +\beta \right) }\right) \right)
\\
& =\frac{a_{3}}{b_{2}}\frac{1}{s^{1+\nu -\alpha }}\left( 1+\left( \frac{a_{2}%
}{a_{3}}-\frac{b_{1}}{b_{2}}\right) \frac{1}{s^{\alpha +\beta }}+\left( 
\frac{a_{1}}{a_{3}}-\frac{a_{2}}{a_{3}}\frac{b_{1}}{b_{2}}+\left( \frac{b_{1}%
}{b_{2}}\right) ^{2}\right) \frac{1}{s^{2\left( \alpha +\beta \right) }}%
+O\left( s^{-3\left( \alpha +\beta \right) }\right) \right) \\
& =\frac{a_{3}}{b_{2}}\frac{1}{s^{1+\nu -\alpha }}+\frac{a_{3}}{b_{2}}\left( 
\frac{a_{2}}{a_{3}}-\frac{b_{1}}{b_{2}}\right) \frac{1}{s^{1+\beta +\nu }} \\
& \quad +\frac{a_{3}}{b_{2}}\left( \frac{a_{1}}{a_{3}}-\frac{a_{2}}{a_{3}}%
\frac{b_{1}}{b_{2}}+\left( \frac{b_{1}}{b_{2}}\right) ^{2}\right) \frac{1}{%
s^{1+\alpha +2\beta +\nu }}+O\left( s^{^{-\left( 1+2\alpha +3\beta +\nu
\right) }}\right) ,\quad \text{when}\quad s\rightarrow \infty \text{,}
\end{align*}%
that implies the short time asymptotics of the creep compliance as%
\begin{align}
\varepsilon _{cr}\left( t\right) & =\frac{a_{3}}{b_{2}}\frac{t^{\nu -\alpha }%
}{\Gamma \left( 1+\nu -\alpha \right) }+\frac{a_{3}}{b_{2}}\left( \frac{a_{2}%
}{a_{3}}-\frac{b_{1}}{b_{2}}\right) \frac{t^{\beta +\nu }}{\Gamma \left(
1+\beta +\nu \right) }  \notag \\
& \quad +\frac{a_{3}}{b_{2}}\left( \frac{a_{1}}{a_{3}}-\frac{a_{2}}{a_{3}}%
\frac{b_{1}}{b_{2}}+\left( \frac{b_{1}}{b_{2}}\right) ^{2}\right) \frac{%
t^{\alpha +2\beta +\nu }}{\Gamma \left( 1+\alpha +2\beta +\nu \right) }%
+O\left( t^{2\alpha +3\beta +\nu }\right) ,\quad \text{when}\quad
t\rightarrow 0,  \label{cr-sha}
\end{align}%
while the large time asymptotics of the creep compliance%
\begin{align}
\varepsilon _{cr}\left( t\right) & =\frac{a_{1}}{b_{1}}\frac{t^{\beta +\nu }%
}{\Gamma \left( 1+\beta +\nu \right) }+\frac{a_{1}}{b_{1}}\left( \frac{a_{2}%
}{a_{1}}-\frac{b_{2}}{b_{1}}\right) \frac{t^{\nu -\alpha }}{\Gamma \left(
1+\nu -\alpha \right) }  \notag \\
& \quad +\frac{a_{1}}{b_{1}}\left( \frac{a_{3}}{a_{1}}-\frac{a_{2}}{a_{1}}%
\frac{b_{2}}{b_{1}}+\left( \frac{b_{2}}{b_{1}}\right) ^{2}\right) \frac{%
t^{-2\alpha -\beta +\nu }}{\Gamma \left( 1-2\alpha -\beta +\nu \right) }%
+O\left( t^{-\xi }\right) ,\quad \text{when}\quad t\rightarrow \infty ,
\label{cr-la}
\end{align}%
follows from the asymptotics of creep compliance in Laplace domain, obtained
as%
\begin{align*}
\tilde{\varepsilon}_{cr}\left( s\right) & =\frac{a_{1}}{b_{1}}\frac{1}{%
s^{1+\beta +\nu }}\frac{1+\frac{a_{2}}{a_{1}}s^{\alpha +\beta }+\frac{a_{3}}{%
a_{2}}s^{2\left( \alpha +\beta \right) }}{1+\frac{b_{2}}{b_{1}}s^{\alpha
+\beta }} \\
& =\frac{a_{1}}{b_{1}}\frac{1}{s^{1+\beta +\nu }}\left( 1+\frac{a_{2}}{a_{1}}%
s^{\alpha +\beta }+\frac{a_{3}}{a_{1}}s^{2\left( \alpha +\beta \right)
}\right) \left( 1+\frac{b_{2}}{b_{1}}s^{\alpha +\beta }\right) ^{-1} \\
& =\frac{a_{1}}{b_{1}}\frac{1}{s^{1+\beta +\nu }}\left( 1+\frac{a_{2}}{a_{1}}%
s^{\alpha +\beta }+\frac{a_{3}}{a_{1}}s^{2\left( \alpha +\beta \right)
}\right) \left( 1-\frac{b_{2}}{b_{1}}s^{\alpha +\beta }+\left( \frac{b_{2}}{%
b_{1}}\right) ^{2}s^{2\left( \alpha +\beta \right) }+O\left( s^{3\left(
\alpha +\beta \right) }\right) \right) \\
& =\frac{a_{1}}{b_{1}}\frac{1}{s^{1+\beta +\nu }}\left( 1+\left( \frac{a_{2}%
}{a_{1}}-\frac{b_{2}}{b_{1}}\right) s^{\alpha +\beta }+\left( \frac{a_{3}}{%
a_{1}}-\frac{a_{2}}{a_{1}}\frac{b_{2}}{b_{1}}+\left( \frac{b_{2}}{b_{1}}%
\right) ^{2}\right) s^{2\left( \alpha +\beta \right) }+O\left( s^{3\left(
\alpha +\beta \right) }\right) \right) \\
& =\frac{a_{1}}{b_{1}}\frac{1}{s^{1+\beta +\nu }}+\frac{a_{1}}{b_{1}}\left( 
\frac{a_{2}}{a_{1}}-\frac{b_{2}}{b_{1}}\right) \frac{1}{s^{1+\nu -\alpha }}
\\
& \quad +\frac{a_{1}}{b_{1}}\left( \frac{a_{3}}{a_{1}}-\frac{a_{2}}{a_{1}}%
\frac{b_{2}}{b_{1}}+\left( \frac{b_{2}}{b_{1}}\right) ^{2}\right) \frac{1}{%
s^{1-2\alpha -\beta +\nu }}+O\left( s^{-\left( 1-\xi \right) }\right) ,\quad 
\text{when}\quad s\rightarrow 0,
\end{align*}%
with $\xi \in \left( 0,1\right) .$

\begin{figure}[h]
\begin{center}
\begin{minipage}{0.46\columnwidth}
  \subfloat[Creep compliance is a Bernstein function.]{
   \includegraphics[width=\columnwidth]{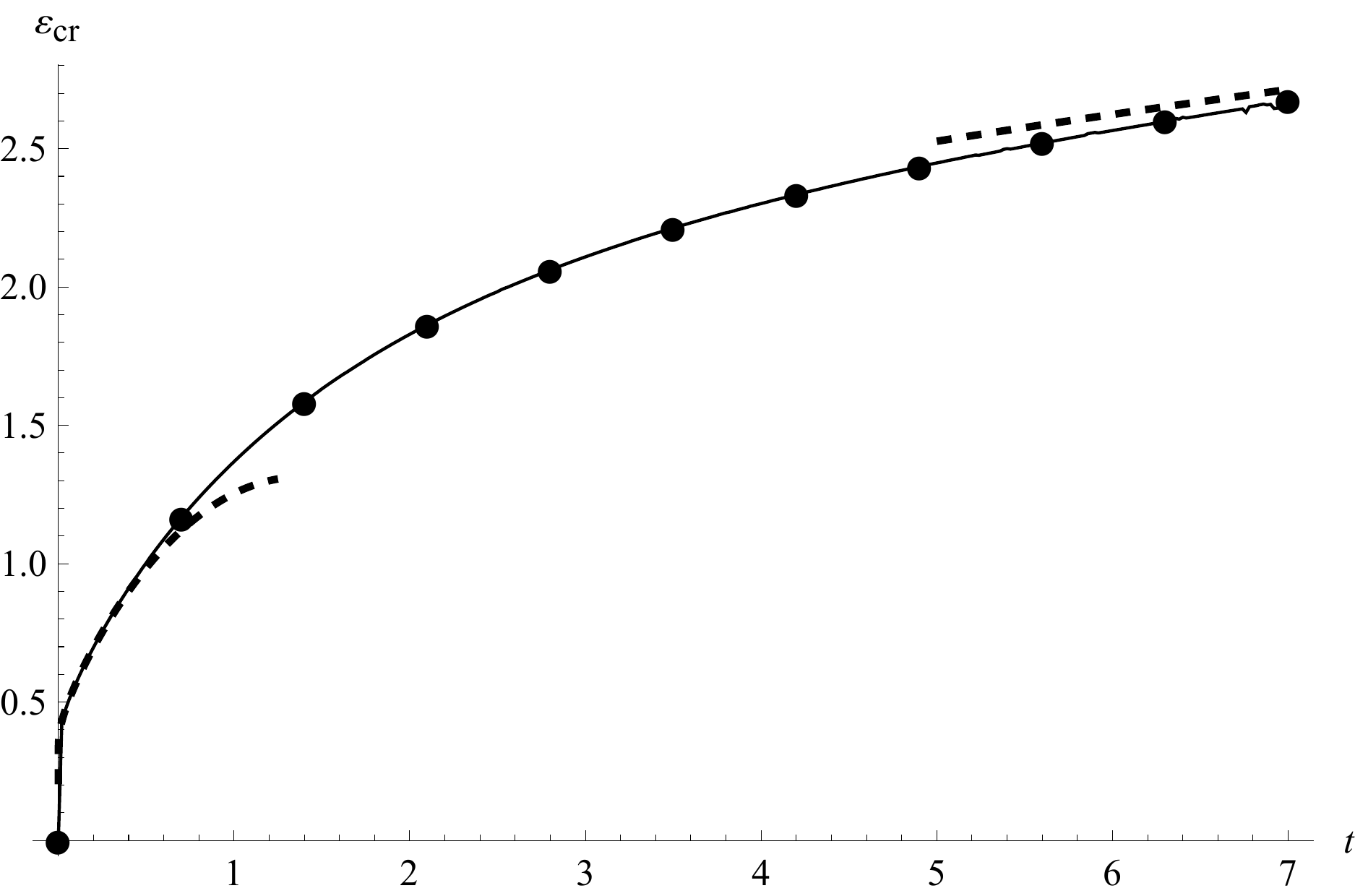}
   \label{cr-np}}
  \end{minipage}
\hfil
\begin{minipage}{0.46\columnwidth}
  \subfloat[Creep compliance is not a Bernstein function..]{
   \includegraphics[width=\columnwidth]{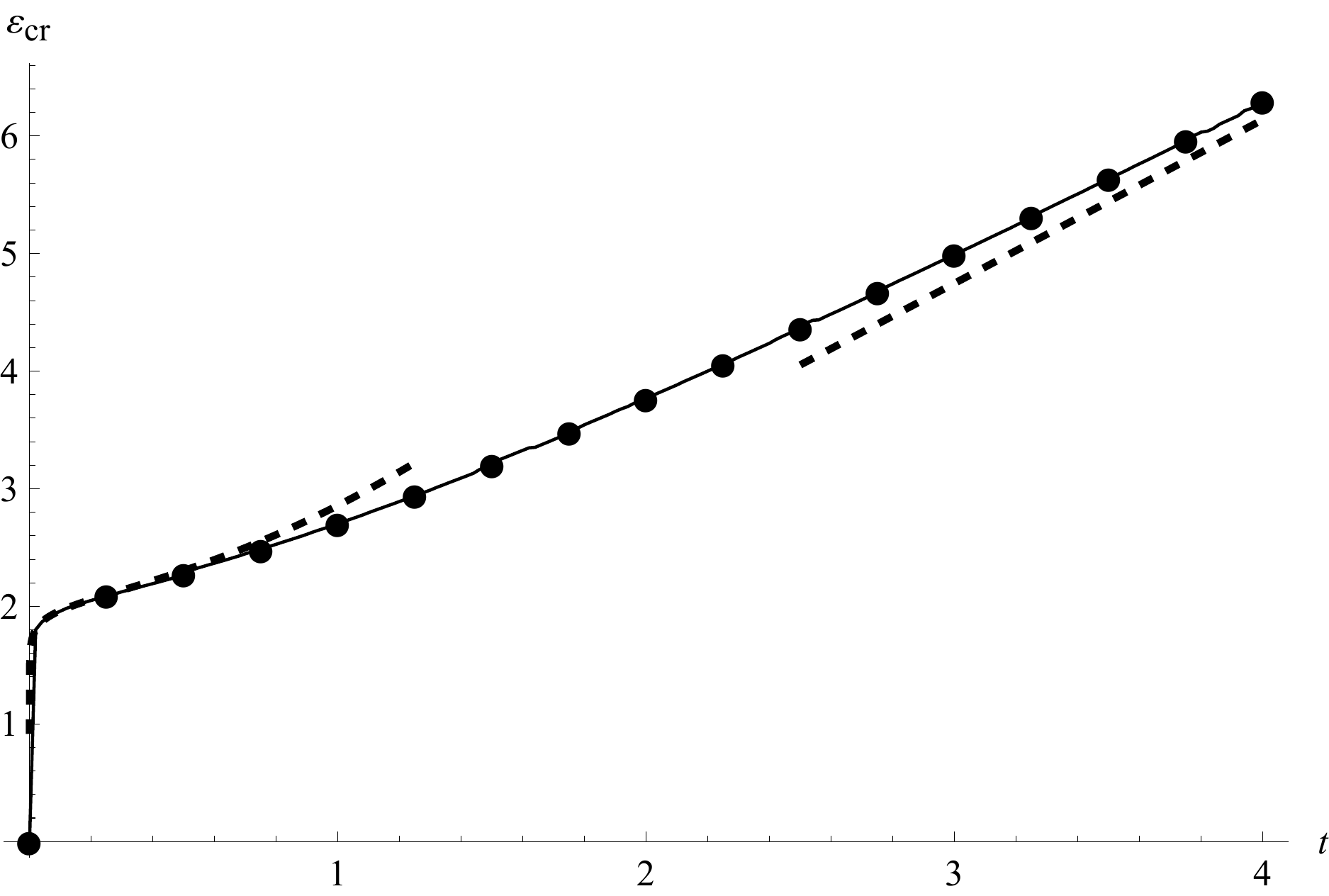}
   \label{cr-ccp}}
  \end{minipage}
\end{center}
\caption{Time evolution of the creep compliance. }
\label{cr}
\end{figure}

Figure \ref{cr} presents time evolution of the creep compliance for
different values of model parameters, so that the creep compliance from
Figure \ref{cr-np}, obtained according to (\ref{cr-preko-jedinice}) for the
same set of parameters as the relaxation modulus from Figure \ref{sr-np}, is
a Bernstein function, while the creep compliance from Figure \ref{cr-ccp},
obtained for the set of parameters in Table \ref{parametri} guaranteeing
that function $\tilde{\sigma}_{sr}$ has a pair of complex conjugated poles,
is obviously not a Bernstein function, since it is a convex, rather than a
concave function. Creep compliance curves from Figure \ref{cr}, obtained
according to the integral representation (\ref{cr-preko-jedinice}),
perfectly overlap with the creep compliance curves represented by dots, that
are calculated according to the Mittag-Leffler representation (\ref%
{cr-preko-ML}). Finally, short and long time asymptotic expressions (\ref%
{cr-sha}) and (\ref{cr-la}) generate curves that show satisfactory
overlapping with the creep compliance curves. 

\section{Position and number of zeros of functions $\protect\phi _{\protect%
\sigma }$ and $\protect\phi _{\protect\varepsilon }$ \label{DPNZ}}

Functions $\phi _{\sigma }$ and $\phi _{\varepsilon }$ in the case of
fractional anti-Zener and Zener models are of three different forms%
\begin{equation*}
\phi \left( s\right) =\left\{ \!\!\!%
\begin{tabular}{l}
$as^{\xi }+b,$ \smallskip \\ 
$as^{1+\xi }+b,$%
\end{tabular}%
\ \right. \quad \phi \left( s\right) =\left\{ \!\!\!%
\begin{tabular}{l}
$as^{\xi }+bs^{\zeta }+c,$ \smallskip \\ 
$as^{1+\xi }+bs^{\zeta }+c,$%
\end{tabular}%
\ \right. \quad \text{and}\quad \phi \left( s\right) =as^{2\xi }+bs^{\xi }+c,
\end{equation*}%
see Table \ref{skupina}, where $a,b,c>0$ and $\xi ,\zeta \in \left(
0,1\right) .$ Functions containing powers of complex variable $s$ that are
in the interval $\left( 0,1\right) ,$ as it is well-known, do not admit
zeros in the main Riemann branch $\arg s\in \left( -\pi ,\pi \right) ,$
while functions having the power of leading term in the interval $\left(
1,2\right) $ may admit no zeros, a negative real zero, and a pair of complex
conjugated zeros with negative real part, as it is proved in the sequel.

\subsection{Two-term function}

Functions $\phi _{\sigma }$ and $\phi _{\varepsilon }$ in the case of model
ID.ID takes the form%
\begin{equation*}
\phi \left( s\right) =as^{1+\xi }+b,\quad s\in 
%TCIMACRO{\U{2102} }%
%BeginExpansion
\mathbb{C}
%EndExpansion
,\quad \arg s\in \left( -\pi ,\pi \right) ,
\end{equation*}%
see Table \ref{skupina}, where $a,b>0$ and $\xi \in \left( 0,1\right) .$ The
function $\phi $ proves to admit a pair of complex conjugated zeros with
negative real part determined by 
\begin{equation*}
s_{\scriptscriptstyle{\mathrm{CCP}}}=\left( \frac{b}{a}\right) ^{\frac{1}{%
1+\xi }}\mathrm{e}^{\mathrm{i}\frac{\pi }{1+\xi }}\quad \text{and}\quad \bar{%
s}_{\scriptscriptstyle{\mathrm{CCP}}}=\left( \frac{b}{a}\right) ^{\frac{1}{%
1+\xi }}\mathrm{e}^{-\mathrm{i}\frac{\pi }{1+\xi }}.
\end{equation*}

The real and imaginary part of function $\phi $, by substituting $s=\rho 
\mathrm{e}^{\mathrm{i}\varphi }$, become%
\begin{equation*}
\func{Re}\phi \left( \rho ,\varphi \right) =a\rho ^{1+\xi }\cos \left(
\left( 1+\xi \right) \varphi \right) +b\quad \text{and}\quad \func{Im}\phi
\left( \rho ,\varphi \right) =a\rho ^{1+\xi }\sin \left( \left( 1+\xi
\right) \varphi \right) ,
\end{equation*}%
so that $\func{Im}\phi \left( \rho ,\varphi \right) =0$ implies $\varphi
_{k}=\frac{k\pi }{1+\xi }.$ By requesting $\varphi _{k}\in \left( -\pi ,\pi
\right) ,$ one finds $-\left( 1+\xi \right) \leqslant k\leqslant \left(
1+\xi \right) $, implying $k\in \left\{ -1,0,1\right\} ,$ transforming the
real part of function $\phi $ into%
\begin{equation*}
\func{Re}\phi \left( \rho ,\varphi _{k}\right) =a\rho ^{1+\xi }\cos \left(
k\pi \right) +b,
\end{equation*}%
that is zero for $k\in \left\{ -1,1\right\} ,$ implying $\rho =\left( \frac{b%
}{a}\right) ^{\frac{1}{1+\xi }},$ while for $k=0$ one has $\func{Re}\phi
\left( \rho ,\varphi _{0}\right) =a\rho ^{1+\xi }+b>0.$ Therefore, the zeros
of function $\phi $ are $s_{1}=\left( \frac{b}{a}\right) ^{\frac{1}{1+\xi }}%
\mathrm{e}^{\mathrm{i}\frac{\pi }{1+\xi }}$ and $s_{-1}=\left( \frac{b}{a}%
\right) ^{\frac{1}{1+\xi }}\mathrm{e}^{-\mathrm{i}\frac{\pi }{1+\xi }}.$

\subsection{Three-term function}

Functions $\phi _{\sigma }$ and $\phi _{\varepsilon }$ in the case of models
IID.IID and IDD.IDD may take the form%
\begin{equation*}
\phi \left( s\right) =as^{1+\xi }+bs^{\zeta }+c,\quad s\in 
%TCIMACRO{\U{2102} }%
%BeginExpansion
\mathbb{C}
%EndExpansion
,\quad \arg s\in \left( -\pi ,\pi \right) ,
\end{equation*}%
see Table \ref{skupina}, where $a,b,c>0$ and $\xi ,\zeta \in \left(
0,1\right) .\ $The function $\phi $ proves to admit: a pair of complex
conjugated zeros with negative real part if $\func{Re}\phi \left( \rho
^{\ast },\pi \right) >0,$ a negative real zero $s=-\rho ^{\ast }$ if $\func{%
Re}\phi \left( \rho ^{\ast },\pi \right) =0,$ no zeros if $\func{Re}\phi
\left( \rho ^{\ast },\pi \right) <0,$ where $\func{Re}\phi $ is given by (%
\ref{re-i-im-psi})$_{1}$ and where $\rho ^{\ast }=\left( \frac{b}{a}\frac{%
\sin \left( \zeta \pi \right) }{\sin \left( \xi \pi \right) }\right) ^{\frac{%
1}{1+\xi -\zeta }}.$

The real and imaginary part of function $\phi $, by substituting $s=\rho 
\mathrm{e}^{\mathrm{i}\varphi }$, become%
\begin{equation}
\func{Re}\phi \left( \rho ,\varphi \right) =a\rho ^{1+\xi }\cos \left(
\left( 1+\xi \right) \varphi \right) +b\rho ^{\zeta }\cos \left( \zeta
\varphi \right) +c\quad \text{and}\quad \func{Im}\phi \left( \rho ,\varphi
\right) =a\rho ^{1+\xi }\sin \left( \left( 1+\xi \right) \varphi \right)
+b\rho ^{\zeta }\sin \left( \zeta \varphi \right) .  \label{re-i-im-psi}
\end{equation}%
According to (\ref{re-i-im-psi})$_{2}$, one has $\func{Im}\phi \left( \rho
,-\varphi \right) =-\func{Im}\phi \left( \rho ,\varphi \right) $ implying
that $s_{0}=\rho _{0}\mathrm{e}^{\mathrm{i}\varphi _{0}}$ and its complex
conjugate $\bar{s}_{0}=\rho _{0}\mathrm{e}^{-\mathrm{i}\varphi _{0}}$ are
both zeros of function $\phi ,$ so that it is sufficient to to consider the
upper complex half-plane only. Clearly, there are no zeros of function $\phi 
$ lying in the right complex half-plane, since $\func{Im}\phi \left( \rho
,\varphi \right) >0$ for $\varphi \in \left( 0,\frac{\pi }{2}\right) $ and
moreover there are no positive zeros of function $\phi $ as well, since $%
\func{Re}\phi \left( \rho ,0\right) >0.$ Therefore, zeros of function $\phi $
may only lie in the left complex half-plane.

In order to show that function $\phi $ admits a single zero in the left
complex quarter-plane, the argument principle and contour $\gamma ,$ shown
in Figure \ref{fig-gama}, is used, stating that if the variable $s$ changes
along the contour $\gamma $ closed in the complex plane, then the number of
zeros $N$ of function $\phi $ in the domain encircled by contour $\gamma $
is given by $\Delta \arg \phi \left( s\right) =2\pi N,$ provided that
function $\phi $ does not have poles in the mentioned domain.

\noindent 
\begin{minipage}{\columnwidth}
\begin{minipage}[c]{0.4\columnwidth}
\centering
\includegraphics[width=0.7\columnwidth]{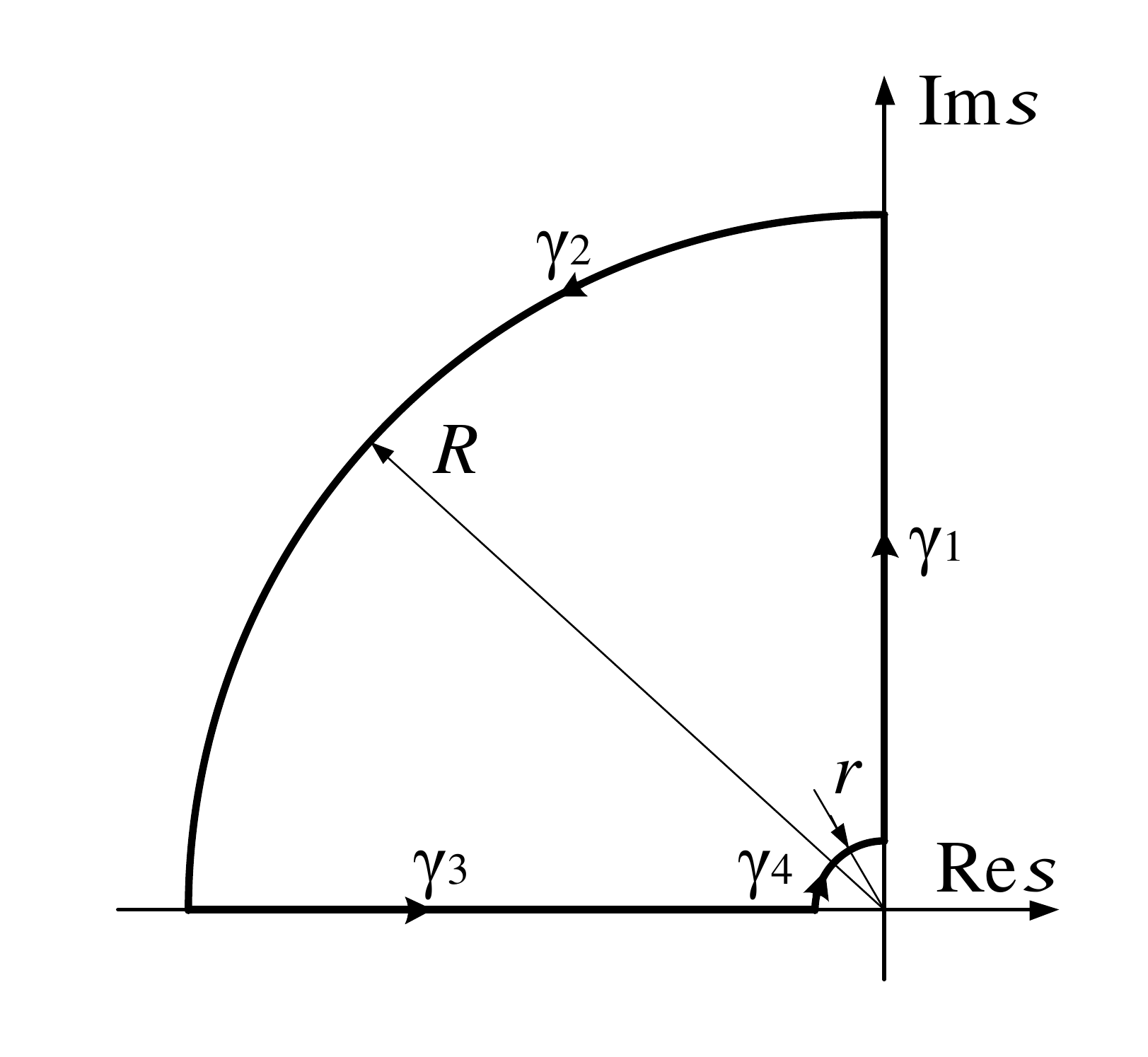}
\captionof{figure}{Contour $\gamma$.}
\label{fig-gama}
\end{minipage}
\hfil
\begin{minipage}[c]{0.55\columnwidth}
\centering
\begin{tabular}{rll}
$\gamma _{1}:$ & $s=\rho \mathrm{e}^{\mathrm{i}\frac{\pi}{2}},$ & $\rho \in \left[ r,R%
\right] ,$ \\
$\gamma _{2}:$ & $s=R\mathrm{e}^{\mathrm{i}\varphi },$ & $\varphi \in \left[ 
\frac{\pi }{2},\pi \right] ,$ \\ 
$\gamma _{3}:$ & $s=\rho \mathrm{e}^{\mathrm{i}\pi },$ & $\rho \in \left[ r,R%
\right] ,$ \\ 
$\gamma _{4}:$ & $s=r\mathrm{e}^{\mathrm{i}\varphi },$ & $\varphi \in \left[ 
\frac{\pi }{2},\pi \right] .$ \\ 
\end{tabular}
\captionof{table}{Parametrization of contour $\gamma$.}
\label{fig-gama-param}
\end{minipage}
\end{minipage}\smallskip

The argument of complex numbers $s$ belonging to the contour $\gamma _{1}$
has a fixed value $\varphi =\frac{\pi }{2},$ while their modulus changes so
that $\rho \in \left( r,R\right) ,$ see the parameterization given in Table %
\ref{fig-gama-param}, so that one has 
\begin{equation*}
\func{Im}\phi \left( \rho ,\frac{\pi }{2}\right) =a\rho ^{1+\xi }\cos \frac{%
\xi \pi }{2}+b\rho ^{\zeta }\sin \frac{\zeta \pi }{2}>0,
\end{equation*}%
by (\ref{re-i-im-psi}), which also implies that function $\phi $ does not
have purely imaginary zeros. The asymptotics of real and imaginary parts of
function $\phi ,$ given by (\ref{re-i-im-psi}), taking the form%
\begin{eqnarray*}
&&\func{Re}\phi \left( r,\frac{\pi }{2}\right) \sim c\quad \text{and}\quad 
\func{Im}\phi \left( r,\frac{\pi }{2}\right) \sim br^{\zeta }\sin \frac{%
\zeta \pi }{2}\rightarrow 0^{+},\quad \text{for}\quad r\rightarrow 0,\quad 
\text{as well as} \\
&&\func{Re}\phi \left( R,\frac{\pi }{2}\right) \sim -aR^{1+\xi }\sin \frac{%
\xi \pi }{2}\rightarrow -\infty \quad \text{and}\quad \func{Im}\phi \left( R,%
\frac{\pi }{2}\right) \sim aR^{1+\xi }\cos \frac{\xi \pi }{2}\rightarrow
\infty ,\quad \text{for}\quad R\rightarrow \infty ,
\end{eqnarray*}%
so that $\left\vert \phi \left( R,\frac{\pi }{2}\right) \right\vert \sim
aR^{1+\xi }\rightarrow \infty $ and $\tan \arg \phi \left( R,\frac{\pi }{2}%
\right) \sim -\tan \frac{\xi \pi }{2},$ implying $\arg \phi \left( R,\frac{%
\pi }{2}\right) \sim \pi -\frac{\xi \pi }{2}\in \left( \frac{\pi }{2},\pi
\right) .$ Therefore, as $s$ changes along the contour $\gamma _{1},$ the
argument of function $\phi $ changes from zero to $\pi -\frac{\xi \pi }{2}.$

Complex number $s$ belonging to the contour $\gamma _{2}$ has large but
fixed modulus $\rho =R,$ while its argument changes along the contour $%
\gamma _{2}$ taking the values $\varphi \in \left[ \frac{\pi }{2},\pi \right]
,$ see Table \ref{fig-gama-param}, so that the asymptotics of (\ref%
{re-i-im-psi}) as $R\rightarrow \infty $ gives%
\begin{equation*}
\func{Re}\phi \left( R,\varphi \right) \sim aR^{1+\xi }\cos \left( \left(
1+\xi \right) \varphi \right) \quad \text{and}\quad \func{Im}\phi \left(
R,\varphi \right) \sim aR^{1+\xi }\sin \left( \left( 1+\xi \right) \varphi
\right) ,
\end{equation*}%
so that $\left\vert \phi \left( R,\varphi \right) \right\vert \sim aR^{1+\xi
}\rightarrow \infty $ and $\arg \phi \left( R,\varphi \right) \sim \left(
1+\xi \right) \varphi ,$ implying $\arg \phi \left( R,\varphi \right) \in %
\left[ \frac{\left( 1+\xi \right) \pi }{2},\left( 1+\xi \right) \pi \right]
. $ Therefore, as $s$ changes along contour $\gamma _{2}$, the argument of
function $\phi $ starts from the second quadrant, since $\arg \phi \left( R,%
\frac{\pi }{2}\right) =\frac{\left( 1+\xi \right) \pi }{2}$ and ends either
in the third, or possibly the fourth quadrant, since $\arg \phi \left( R,\pi
\right) =\left( 1+\xi \right) \pi ,$ however not crossing through the first
quadrant, since the real and imaginary parts of function $\phi $\ cannot be
simultaneously positive for $\varphi \in \left[ \frac{\pi }{2},\pi \right] ,$%
\ since $\frac{\pi }{2}\leqslant \left( 1+\xi \right) \varphi \leqslant 2\pi 
$.

The argument of complex numbers $s$ lying on the contour $\gamma _{3}$ has a
fixed value $\varphi =\pi ,$ while their modulus changes in the interval $%
\rho \in \left( r,R\right) ,$ yielding%
\begin{equation*}
\func{Re}\phi \left( \rho ,\pi \right) =-a\rho ^{1+\xi }\cos \left( \xi \pi
\right) +b\rho ^{\zeta }\cos \left( \zeta \pi \right) +c\quad \text{and}%
\quad \func{Im}\phi \left( \rho ,\pi \right) =-a\rho ^{1+\xi }\sin \left(
\xi \pi \right) +b\rho ^{\zeta }\sin \left( \zeta \pi \right) ,
\end{equation*}%
according to (\ref{re-i-im-psi}), so that the asymptotics of real and
imaginary parts of function $\phi $ is obtained in the form%
\begin{equation*}
\func{Re}\phi \left( R,\pi \right) \sim -aR^{1+\xi }\cos \left( \xi \pi
\right) \rightarrow -\infty \quad \text{and}\quad \func{Im}\phi \left( R,\pi
\right) \sim -aR^{1+\xi }\sin \left( \xi \pi \right) \rightarrow -\infty
,\quad \text{for}\quad R\rightarrow \infty ,
\end{equation*}
as well as 
\begin{equation*}
\func{Re}\phi \left( r,\pi \right) \sim c\quad \text{and}\quad \func{Im}\phi
\left( r,\pi \right) \sim br^{\zeta }\sin \left( \zeta \pi \right)
\rightarrow 0^{+}\quad \text{for}\quad r\rightarrow 0,
\end{equation*}%
implying that $\arg \phi \left( R,\pi \right) \sim \left( 1+\xi \right) \pi $
as $R\rightarrow \infty $, as well as $\tan \arg \phi \left( r,\pi \right)
\sim 0$ as $r\rightarrow 0.$ Obviously, $\func{Im}\phi \left( \rho ,\pi
\right) $ changes sign as $\rho $ changes from $R\rightarrow \infty $ to $%
r\rightarrow 0$ and it has a single zero at $\rho ^{\ast }=\left( \frac{b}{a}%
\frac{\sin \left( \zeta \pi \right) }{\sin \left( \xi \pi \right) }\right) ^{%
\frac{1}{1+\xi -\zeta }},$ since $\rho ^{1+\xi }$ is a convex and $\rho
^{\zeta }$ is a concave function, both monotonically increasing from zero to
infinity. Therefore, the sign of $\func{Re}\phi \left( \rho ^{\ast },\pi
\right) $ determines whether the argument of function $\phi ,$ starting from 
$\arg \phi \left( R,\pi \right) =\left( 1+\xi \right) \pi ,$ as $s$ changes
along contour $\gamma _{3},$ ends at $\arg \phi \left( r,\pi \right) =0$ or
at $\arg \phi \left( r,\pi \right) =2\pi ,$ so that the former is achieved
if $\func{Re}\phi \left( \rho ^{\ast },\pi \right) <0$ and the latter is
achieved if $\func{Re}\phi \left( \rho ^{\ast },\pi \right) >0.$ Moreover,
if $\func{Re}\phi \left( \rho ^{\ast },\pi \right) =0,$ then function $\phi $
admits a negative real zero $s=-\rho ^{\ast }=-\left( \frac{b}{a}\frac{\sin
\left( \zeta \pi \right) }{\sin \left( \xi \pi \right) }\right) ^{\frac{1}{%
1+\xi -\zeta }},$ since $\rho ^{\ast }$ solves $\func{Im}\phi \left( \rho
,\pi \right) =0.$

The modulus of complex number $s$ belonging to the contour $\gamma _{4}$ has
a fixed but small value $\rho =r,$ while its argument changes in the
interval $\varphi \in \left[ \frac{\pi }{2},\pi \right] ,$ so that the
asymptotics of real and imaginary parts of function $\phi $, see (\ref%
{re-i-im-psi}), become%
\begin{equation*}
\func{Re}\phi \left( r,\varphi \right) \sim c\quad \text{and}\quad \func{Im}%
\phi \left( r,\varphi \right) \sim br^{\zeta }\sin \left( \zeta \varphi
\right) \rightarrow 0^{+}\quad \text{for}\quad r\rightarrow 0,
\end{equation*}%
implying that $\arg \phi $ does not really change.

Summing up, the change of argument of function $\phi ,$ as $s$ changes along
the contour $\gamma ,$ is either $\Delta \arg \phi \left( s\right) =2\pi $
if $\func{Re}\phi \left( \rho ^{\ast },\pi \right) >0,$ or $\Delta \arg \phi
\left( s\right) =0$ if $\func{Re}\phi \left( \rho ^{\ast },\pi \right) <0,$
implying that function $\phi $ either has one zero, or has no zeros in the
upper left complex quarter-plane. If $\func{Re}\phi \left( \rho ^{\ast },\pi
\right) =0,$ then function $\phi $ admits a negative real zero.

\subsection{Function quadratic in $s^{\protect\xi }$\label{nule-trinom}}

Function $\phi _{\sigma }$ in the case of models I$^{+}$ID.ID and IDD$^{+}$%
.DD$^{+},$ as well as the function $\phi _{\varepsilon }$ in the case of
model ID.IDD$^{+},$ so as both functions $\phi _{\sigma }$ and $\phi
_{\varepsilon }$ in the case of models I$^{+}$ID.I$^{+}$ID, IDD$^{+}$.IDD$%
^{+},$ and I$^{+}$ID.IDD$^{+},$ take the same form%
\begin{equation}
\phi \left( s\right) =as^{2\xi }+bs^{\xi }+c,\quad s\in 
%TCIMACRO{\U{2102} }%
%BeginExpansion
\mathbb{C}
%EndExpansion
,\quad \arg s\in \left( -\pi ,\pi \right) ,  \label{kvadratna}
\end{equation}%
see Table \ref{skupina}, where $a,b,c>0$ and $\xi \in \left( 0,1\right) .$
If $\xi \in \left( 0,\frac{1}{2}\right] ,$ then all powers of variable $s$
belong to interval $\left( 0,1\right) $ and, as it is well-known, the
function $\phi $ does not have zeros, while if $\xi \in \left( \frac{1}{2}%
,1\right) ,$ then the zeros of function $\phi $ are%
\begin{equation*}
s^{\xi }=-\frac{b}{2a}\left( 1\pm \sqrt{1-\frac{4ac}{b^{2}}}\right) .
\end{equation*}

Assuming that the discriminant of equation quadratic in $s^{\xi }$ is
non-negative, i.e., that $\frac{2\sqrt{ac}}{b}\leqslant 1,$ one has%
\begin{equation*}
s^{\xi }=-\frac{b}{2a}\left( 1\pm \sqrt{1-\frac{4ac}{b^{2}}}\right) <0\quad 
\text{i.e.,}\quad \rho ^{\xi }\,\mathrm{e}^{\mathrm{i}\xi \varphi }=\frac{b}{%
2a}\left( 1\pm \sqrt{1-\frac{4ac}{b^{2}}}\right) \,\mathrm{e}^{\mathrm{i}%
\left( 2k+1\right) \pi },\quad k\in 
%TCIMACRO{\U{2124} }%
%BeginExpansion
\mathbb{Z}
%EndExpansion
,
\end{equation*}%
implying 
\begin{equation}
\rho =\sqrt[\xi ]{\frac{b}{2a}\left( 1\pm \sqrt{1-\frac{4ac}{b^{2}}}\right) }%
\quad \text{and}\quad \varphi =\frac{\left( 2k+1\right) \pi }{\xi }.
\label{de-manje-od-nule}
\end{equation}%
The complex number $s$ with modulus $\rho $ and argument $\varphi ,$ given
by (\ref{de-manje-od-nule}), is not a zero of function $\phi ,$ given by (%
\ref{kvadratna}), since the requirement $\varphi \in \left( -\pi ,\pi
\right) $ cannot be satisfied for any $k\in 
%TCIMACRO{\U{2124} }%
%BeginExpansion
\mathbb{Z}
%EndExpansion
,$ due to $-\pi <\frac{\left( 2k+1\right) \pi }{\xi }<\pi $ implying $%
-1<-\xi <2k+1<\xi <1,$ i.e., $-1<k<0.$

On the other hand, if the discriminant of equation quadratic in $s^{\xi }$
is negative, i.e., if $\frac{2\sqrt{ac}}{b}>1,$ then one has%
\begin{equation}
s^{\xi }=-\frac{b}{2a}\left( 1\pm \mathrm{i}\sqrt{\frac{4ac}{b^{2}}-1}%
\right) =-\lambda \pm \mathrm{i}\eta =z_{\pm }\quad \text{with}\quad \lambda
=\frac{b}{2a}>0\quad \text{and}\quad \eta =\frac{b}{2a}\sqrt{\frac{4ac}{b^{2}%
}-1}>0,  \label{lambda-eta}
\end{equation}%
so that%
\begin{equation}
s^{\xi }=\rho ^{\xi }\,\mathrm{e}^{\mathrm{i}\xi \varphi _{\pm }}=\left\vert
z_{\pm }\right\vert \,\mathrm{e}^{\mathrm{i}\left( \arg z_{\pm }+2k\pi
\right) }\quad \text{with}\quad \left\vert z_{\pm }\right\vert =\sqrt{%
\lambda ^{2}+\eta ^{2}}\quad \text{and}\quad \arg z_{\pm }=\pm \left( 1-%
\frac{1}{\pi }\arctan \frac{\eta }{\lambda }\right) \pi ,  \label{s-na-ksi}
\end{equation}%
implying 
\begin{equation}
s_{\pm }=\rho \,\mathrm{e}^{\mathrm{i}\varphi _{\pm }}\quad \text{with}\quad
\rho =\sqrt[\xi ]{\left\vert z_{\pm }\right\vert }\quad \text{and}\quad
\varphi _{\pm }=\frac{\arg z_{\pm }+2k\pi }{\xi }.  \label{ro-i-fi-pm}
\end{equation}%
In order for $\varphi _{+}$ to be an argument of zero of the function $\phi
, $ given by (\ref{kvadratna}), one requires $\varphi _{+}=\frac{\arg
z_{+}+2k\pi }{\xi }=\frac{2k+1-\delta }{\xi }\pi \in \left( 0,\pi \right) ,$
with $\delta =\frac{1}{\pi }\arctan \frac{\eta }{\lambda }\in \left( 0,\frac{%
1}{2}\right) ,$ so that one obtains $0<\left( 2k+1-\delta \right) \pi <\xi
\pi ,$ transforming into $0<\delta <2k+1<\delta +\xi <\frac{3}{2},$ i.e., $-%
\frac{1}{2}<k<\frac{1}{4},$ yielding $k=0.$ Similar argumentation in the
case of $\varphi _{-}$ also implies that $k=0.$ Therefore, the zeros $s_{\pm
}=\rho \,\mathrm{e}^{\mathrm{i}\varphi _{\pm }}$ of function $\phi ,$ given
by (\ref{kvadratna}), are such that%
\begin{equation*}
\rho =\sqrt[2\xi ]{\lambda ^{2}+\eta ^{2}}\quad \text{and}\quad \varphi
_{\pm }=\pm \frac{1-\frac{1}{\pi }\arctan \frac{\eta }{\lambda }}{\xi }\pi ,
\end{equation*}%
with $\lambda $ and $\eta $ defined by (\ref{lambda-eta}), according to (\ref%
{s-na-ksi}) and (\ref{ro-i-fi-pm}). Actually, the function $\phi $ admits a
pair of complex conjugated zeros%
\begin{equation}
s_{\scriptscriptstyle{\mathrm{CCP}}}=\sqrt[2\xi ]{\lambda ^{2}+\eta ^{2}}\,%
\mathrm{e}^{\mathrm{i}\frac{1-\frac{1}{\pi }\arctan \frac{\eta }{\lambda }}{%
\xi }\pi }\quad \text{and}\quad \bar{s}_{\scriptscriptstyle{\mathrm{CCP}}}=%
\sqrt[2\xi ]{\lambda ^{2}+\eta ^{2}}\,\mathrm{e}^{-\mathrm{i}\frac{1-\frac{1%
}{\pi }\arctan \frac{\eta }{\lambda }}{\xi }\pi },  \label{sccp}
\end{equation}%
having negative real part, if $\tan \left( \xi \pi \right) >-\frac{\eta }{%
\lambda },$ since $\varphi _{+}=\frac{1-\frac{1}{\pi }\arctan \frac{\eta }{%
\lambda }}{\xi }\pi <\pi $ implies $1-\frac{1}{\pi }\arctan \frac{\eta }{%
\lambda }<\xi ,$ reducing to $\tan \left( \xi \pi \right) >-\frac{\eta }{%
\lambda },$ while one already knows that $\arg z_{+}=\pi -\arctan \frac{\eta 
}{\lambda }>\frac{\pi }{2},$ see (\ref{s-na-ksi})$_{3}$, and therefore $%
\varphi _{+}=\frac{\arg z_{+}}{\xi }>\frac{\pi }{2},$ see (\ref{ro-i-fi-pm}%
), due to $\xi <1.$ Note, if $\tan \left( \xi \pi \right) <-\frac{\eta }{%
\lambda },$ then $\varphi _{+}>\pi $ and function $\phi ,$ given by (\ref%
{kvadratna}) does not have zeros for $\arg s\in \left( -\pi ,\pi \right) .$

In the special case when $\varphi _{\pm }=\pm \frac{1-\frac{1}{\pi }\arctan 
\frac{\eta }{\lambda }}{\xi }\pi =\pm \pi ,$ i.e., when $\frac{1-\frac{1}{%
\pi }\arctan \frac{\eta }{\lambda }}{\xi }=1$ reducing to $\tan \left( \xi
\pi \right) =-\frac{\eta }{\lambda },$ one has that the function $\phi ,$
given by (\ref{kvadratna}), admits a negative real zero 
\begin{equation}
s_{\scriptscriptstyle{\mathrm{RP}}}=-\sqrt[2\xi ]{\lambda ^{2}+\eta ^{2}},
\label{srp}
\end{equation}%
according to (\ref{sccp}), with the order $\xi $ obtained as a solution to $%
\tan \left( \xi \pi \right) =-\frac{\eta }{\lambda }.$

In conclusion, the function $\phi ,$ given by (\ref{kvadratna}), does not
have zeros if $\xi \in \left( 0,\frac{1}{2}\right] ,$ or if $\xi \in \left( 
\frac{1}{2},1\right) $ and $\frac{2\sqrt{ac}}{b}\leqslant 1,$ or if $\xi \in
\left( \frac{1}{2},1\right) ,$ $\frac{2\sqrt{ac}}{b}>1,$ and $\tan \left(
\xi \pi \right) <-\sqrt{\frac{4ac}{b^{2}}-1},$ while if $\xi \in \left( 
\frac{1}{2},1\right) $ and $\frac{2\sqrt{ac}}{b}>1,$ then it has either a
negative real zero%
\begin{equation}
s_{\scriptscriptstyle{\mathrm{RP}}}=-\sqrt[2\xi ]{\frac{c}{a}},
\label{nule-rp}
\end{equation}%
according to (\ref{srp}) and (\ref{lambda-eta}), if $\tan \left( \xi \pi
\right) =-\sqrt{\frac{4ac}{b^{2}}-1}$, that actually determines the order $%
\xi ,$ or it has a pair of complex conjugated zeros 
\begin{equation}
s_{\scriptscriptstyle{\mathrm{CCP}}}=\sqrt[2\xi ]{\frac{c}{a}}\,\mathrm{e}^{%
\mathrm{i}\frac{1-\frac{1}{\pi }\arctan \sqrt{\frac{4ac}{b^{2}}-1}}{\xi }\pi
}\quad \text{and}\quad \bar{s}_{\scriptscriptstyle{\mathrm{CCP}}}=\sqrt[2\xi 
]{\frac{c}{a}}\,\mathrm{e}^{-\mathrm{i}\frac{1-\frac{1}{\pi }\arctan \sqrt{%
\frac{4ac}{b^{2}}-1}}{\xi }\pi },  \label{nule-ccp}
\end{equation}%
according to (\ref{sccp}) and (\ref{lambda-eta}), if $\tan \left( \xi \pi
\right) >-\sqrt{\frac{4ac}{b^{2}}-1}.$

\section{Conclusion}

The analysis of the energy balance of one-dimensional viscoelastic body
modeled by the fractional anti-Zener and Zener models, conducted in Section %
\ref{power}, yields two expressions for the power per unit volume of the
viscoelastic body, where the first one is expressed through the strain (\ref%
{snaga-epsilon}) and the second one is expressed through the stress (\ref%
{snaga-sigma}), both of them containing information about the elastic
properties of material, represented by the energy per unit volume stored in
the viscoelastic body (\ref{pot-en-epsilon}) and (\ref{pot-en-sigma}), as
well as about the viscous properties of the material, represented by the
dissipated power per unit volume (\ref{P-epsilon}) and (\ref{P-sigma}).

The positivity of stored energy and dissipated power per unit volume, both
of them consisting of the instantaneous and hereditary terms and depending
either on the relaxation modulus and its derivatives, or on the creep
compliance and its derivatives, is guaranteed if the relaxation modulus,
respectively creep compliance,\ is a completely monotone function,
respectively a Bernstein function, and therefore the Laplace transform
method is used in Section \ref{sr-and-cr} in order to obtain the explicit
forms of the relaxation modulus $\sigma _{sr},$ expressed by (\ref{sr-opste}%
) in terms of functions $\sigma _{sr}^{\left( \mathrm{NP}\right) },$ $\sigma
_{sr}^{\left( \mathrm{RP}\right) },$ and $\sigma _{sr}^{\left( \mathrm{CCP}%
\right) },$ given by (\ref{sigma-NP}), (\ref{sigma-RP}), and (\ref{sigma-CCP}%
), as well as of the creep compliance $\varepsilon _{cr}$, expressed by (\ref%
{cr-opste}) in terms of functions $\varepsilon _{cr}^{\left( \mathrm{NP}%
\right) },$ $\varepsilon _{cr}^{\left( \mathrm{RP}\right) },$ and $%
\varepsilon _{cr}^{\left( \mathrm{CCP}\right) },$ given by (\ref{epsilon-NP}%
), (\ref{epsilon-RP}), and (\ref{epsilon-CCP}), allowing for the relaxation
modulus and creep compliance to be non-monotonic and even oscillatory
functions with exponentially decreasing amplitude. These qualitative
properties are due to the existence of poles of the relaxation modulus and
creep compliance in the Laplace domain (\ref{sigma-sr-ld}) and (\ref%
{epsilon-cr-ld}), examined in Section \ref{DPNZ}. The relaxation modulus is
proved to be a completely monotone function, as well as the creep compliance
proved to be a Bernstein function, if additional restrictions on model
parameters, derived in Section \ref{narrovlje} and narrowing the
thermodynamical restrictions on the parameters of fractional anti-Zener and
Zener models, listed in Appendix \ref{FAZ-ZM}, are posed.

In Section \ref{nmericali}, the time evolution of relaxation modulus and
creep compliance for model I$^{+}$ID.ID is examined numerically for the set
of parameters from Table \ref{parametri} and presented in Figures \ref{sr}
and \ref{cr} along with the short and long time asymptotics, displaying a
good agreement. In the case of model I$^{+}$ID.ID, the relaxation modulus
can be a completely monotone function, as obvious Figure \ref{sr-np}, as
well as a function resembling to a completely monotone function, because of
the choice of model parameters not guaranteeing the mentioned property, see
Figure \ref{sr-rp}, and moreover an oscillatory function having
exponentially decreasing amplitude, see Figures \ref{sr-ccp-srsha} and \ref%
{sr-ccp-srla}. On the other hand, the creep compliance can be at most a
non-monotonic function, with the possibility to be a Bernstein function, see
Figure \ref{cr-np}, and a function that is positive, monotonically
increasing and a convex rather than a concave function, see Figure \ref%
{cr-ccp}. The creep compliance curves, obtained according to the integral
representation (\ref{cr-preko-jedinice}) because of convergence reasons
rather than (\ref{cr-preko-fiova}), are compared with the ones obtained by
the Mittag-Leffler representation (\ref{cr-preko-ML}) and showed a good
agreement. 

\pagebreak

\noindent \textbf{Declaration of Competing Interest} \smallskip

\noindent The authors declare that they have no conflict of interest.\smallskip

\noindent \textbf{Data availability }\smallskip

\noindent No data was used for the research described in the article. \smallskip

\noindent \textbf {Acknowledgements} \smallskip

\noindent The work is supported by the Ministry of Science, Technological Development
and Innovation of the Republic of Serbia under grant
451-03-47/2023-01/200125 (SJ and DZ). 

\appendix

\section{Fractional anti-Zener and Zener models\label{FAZ-ZM}}

\subsection{Symmetric models}

\textbf{Model ID.ID} takes the form%
\begin{equation*}
\left( a_{1}\,_{0}\mathrm{I}_{t}^{\alpha }+a_{2}\,_{0}\mathrm{D}_{t}^{\beta
}\right) \sigma \left( t\right) =\left( b_{1}\,_{0}\mathrm{I}_{t}^{\mu
}+a_{2}\,_{0}\mathrm{D}_{t}^{\alpha +\beta -\mu }\right) \varepsilon \left(
t\right) ,  \label{AZ-ID-ID-less-konacno}
\end{equation*}%
with thermodynamical restrictions%
\begin{gather}
0\leqslant \alpha +\beta -\mu \leqslant 1,\quad \mu \leqslant \alpha ,\quad
\beta +\mu \leqslant 1,  \label{TD-AZ-ID-ID-less-1} \\
-\frac{a_{1}}{a_{2}}\frac{\cos \frac{\left( 2\alpha +\beta -\mu \right) \pi 
}{2}}{\cos \frac{\left( \beta +\mu \right) \pi }{2}}\leqslant \frac{b_{1}}{%
b_{2}}\leqslant \frac{a_{1}}{a_{2}}\frac{\sin \frac{\left( 2\alpha +\beta
-\mu \right) \pi }{2}}{\sin \frac{\left( \beta +\mu \right) \pi }{2}}.
\label{TD-AZ-ID-ID-less-2}
\end{gather}

\noindent \textbf{Model ID.DD}$^{{}^{+}}$\noindent \textbf{\ }takes the form%
\begin{equation*}
\left( a_{1}\,_{0}\mathrm{I}_{t}^{\alpha }+a_{2}\,_{0}\mathrm{D}_{t}^{\beta
}\right) \sigma \left( t\right) =\left( b_{1}\,_{0}\mathrm{D}_{t}^{\mu
}+b_{2}\,_{0}\mathrm{D}_{t}^{\alpha +\beta +\mu }\right) \varepsilon \left(
t\right) ,  \label{AZ-ID-DD-over}
\end{equation*}%
with thermodynamical restrictions 
\begin{gather}
1\leqslant \alpha +\beta +\mu \leqslant 2,\quad \beta \leqslant \mu
\leqslant 1-\alpha ,  \label{TD-AZ-ID-DD-over-1} \\
\frac{a_{1}}{a_{2}}\frac{\left\vert \cos \frac{\left( 2\alpha +\beta +\mu
\right) \pi }{2}\right\vert }{\cos \frac{\left( \mu -\beta \right) \pi }{2}}%
\leqslant \frac{b_{1}}{b_{2}}.  \label{TD-AZ-ID-DD-over-2}
\end{gather}

\noindent \textbf{Model IID.IID} takes the form%
\begin{equation*}
\left( a_{1}\,_{0}\mathrm{I}_{t}^{\alpha }+a_{2}\,_{0}\mathrm{I}_{t}^{\beta
}+a_{3}\,_{0}\mathrm{D}_{t}^{\gamma }\right) \sigma \left( t\right) =\left(
b_{1}\,_{0}\mathrm{I}_{t}^{\alpha +\gamma -\eta }+b_{2}\,_{0}\mathrm{I}%
_{t}^{\beta +\gamma -\eta }+b_{3}\,_{0}\mathrm{D}_{t}^{\eta }\right)
\varepsilon \left( t\right) ,  \label{franjo-less-less-konacno}
\end{equation*}%
with thermodynamical restrictions%
\begin{gather}
\beta <\alpha ,\quad \gamma \leqslant \eta ,\quad 0\leqslant \beta +\gamma
-\eta \leqslant \alpha +2\gamma -\eta \leqslant 1,\quad \alpha +\gamma
\leqslant \beta +\eta ,  \label{franjo-less-less-TD-1} \\
-\frac{b_{3}}{b_{1}}\frac{\cos \frac{\left( \alpha +\eta \right) \pi }{2}}{%
\cos \frac{\left( \alpha +2\gamma -\eta \right) \pi }{2}}\leqslant \frac{%
a_{3}}{a_{1}}\leqslant \frac{b_{3}}{b_{1}}\frac{\sin \frac{\left( \alpha
+\eta \right) \pi }{2}}{\sin \frac{\left( \alpha +2\gamma -\eta \right) \pi 
}{2}},  \label{franjo-less-less-TD-2} \\
-\frac{b_{3}}{b_{2}}\frac{\cos \frac{\left( \beta +\eta \right) \pi }{2}}{%
\cos \frac{\left( \beta +2\gamma -\eta \right) \pi }{2}}\leqslant \frac{a_{3}%
}{a_{2}}\leqslant \frac{b_{3}}{b_{2}}\frac{\sin \frac{\left( \beta +\eta
\right) \pi }{2}}{\sin \frac{\left( \beta +2\gamma -\eta \right) \pi }{2}}.
\label{franjo-less-less-TD-3}
\end{gather}

\noindent \textbf{Model IDD.IDD} takes the form%
\begin{equation*}
\left( a_{1}\,_{0}\mathrm{I}_{t}^{\alpha }+a_{2}\,_{0}\mathrm{D}_{t}^{\beta
}+a_{3}\,_{0}\mathrm{D}_{t}^{\gamma }\right) \sigma \left( t\right) =\left(
b_{1}\,_{0}\mathrm{I}_{t}^{\mu }+b_{2}\,_{0}\mathrm{D}_{t}^{\alpha +\beta
-\mu }+b_{3}\,_{0}\mathrm{D}_{t}^{\alpha +\gamma -\mu }\right) \varepsilon
\left( t\right) ,  \label{IDD-IDD-less-less-konacno}
\end{equation*}%
with thermodynamical restrictions%
\begin{gather}
0\leqslant \alpha +\gamma -\mu \leqslant 1,\quad \beta <\gamma ,\quad \mu
\leqslant \alpha ,\quad \gamma +\mu \leqslant \alpha +\beta ,\quad \gamma
+\mu \leqslant 1,  \label{IDD-IDD-less-less-TD-1} \\
-\frac{b_{2}}{b_{1}}\frac{\cos \frac{\left( 2\alpha +\beta -\mu \right) \pi 
}{2}}{\cos \frac{\left( \beta +\mu \right) \pi }{2}}\leqslant \frac{a_{2}}{%
a_{1}}\leqslant \frac{b_{2}}{b_{1}}\frac{\sin \frac{\left( 2\alpha +\beta
-\mu \right) \pi }{2}}{\sin \frac{\left( \beta +\mu \right) \pi }{2}},
\label{IDD-IDD-less-less-TD-2} \\
-\frac{b_{3}}{b_{1}}\frac{\cos \frac{\left( 2\alpha +\gamma -\mu \right) \pi 
}{2}}{\cos \frac{\left( \gamma +\mu \right) \pi }{2}}\leqslant \frac{a_{3}}{%
a_{1}}\leqslant \frac{b_{3}}{b_{1}}\frac{\sin \frac{\left( 2\alpha +\gamma
-\mu \right) \pi }{2}}{\sin \frac{\left( \gamma +\mu \right) \pi }{2}}.
\label{IDD-IDD-less-less-TD-3}
\end{gather}

\noindent \textbf{Model IID.IDD} takes the form%
\begin{equation*}
\left( a_{1}\,_{0}\mathrm{I}_{t}^{\alpha }+a_{2}\,_{0}\mathrm{I}_{t}^{\beta
}+a_{3}\,_{0}\mathrm{D}_{t}^{\gamma }\right) \sigma \left( t\right) =\left(
b_{1}\,_{0}\mathrm{I}_{t}^{\mu }+b_{2}\,_{0}\mathrm{D}_{t}^{\nu }+b_{3}\,_{0}%
\mathrm{D}_{t}^{\alpha +\gamma -\mu }\right) \varepsilon \left( t\right) ,
\label{frale-less-less-less-less-konacno}
\end{equation*}%
with thermodynamical restrictions%
\begin{gather}
\mu \leqslant \beta <\alpha ,\quad \gamma \leqslant \nu ,\quad \alpha +\beta
+\gamma \leqslant 1+\mu ,\quad \mu +\nu -\gamma <\alpha \leqslant 1-\nu ,
\label{TD-frale-less-less-less-less-1} \\
-\frac{b_{3}}{b_{1}}\frac{\cos \frac{\left( 2\alpha +\gamma -\mu \right) \pi 
}{2}}{\cos \frac{\left( \gamma +\mu \right) \pi }{2}}\leqslant \frac{a_{3}}{%
a_{1}}\leqslant \frac{b_{3}}{b_{1}}\frac{\sin \frac{\left( 2\alpha +\gamma
-\mu \right) \pi }{2}}{\sin \frac{\left( \gamma +\mu \right) \pi }{2}}.
\label{TD-frale-less-less-less-less-2}
\end{gather}

\noindent \textbf{Model I}$^{{}^{+}}$\textbf{ID.I}$^{{}^{+}}$\textbf{ID}
takes the form%
\begin{equation*}
\left( a_{1}\,_{0}\mathrm{I}_{t}^{1+\alpha }+a_{2}\,_{0}\mathrm{I}_{t}^{%
\frac{1+\alpha -\gamma }{2}}+a_{3}\,_{0}\mathrm{D}_{t}^{\gamma }\right)
\sigma \left( t\right) =\left( b_{1}\,_{0}\mathrm{I}_{t}^{1+\mu }+b_{2}\,_{0}%
\mathrm{I}_{t}^{\frac{1+\mu -\left( \alpha +\gamma -\mu \right) }{2}%
}+b_{3}\,_{0}\mathrm{D}_{t}^{\alpha +\gamma -\mu }\right) \varepsilon \left(
t\right) ,  \label{franjo-over-over-konacno}
\end{equation*}%
with thermodynamical restrictions%
\begin{gather}
\mu \leqslant \alpha ,\quad \alpha +\gamma +2\left( \alpha -\mu \right)
=3\alpha +\gamma -2\mu \leqslant 1,  \label{franjo-over-over-TD-1} \\
\frac{a_{2}}{a_{1}}\leqslant \frac{b_{2}}{b_{1}}\frac{\cos \frac{\left(
1-3\alpha -\gamma +2\mu \right) \pi }{4}}{\cos \frac{\left( 1+\alpha -\gamma
-2\mu \right) \pi }{4}},\quad \frac{a_{3}}{a_{2}}\leqslant \frac{b_{3}}{b_{2}%
}\frac{\sin \frac{\left( 1+3\alpha +\gamma -2\mu \right) \pi }{4}}{\sin 
\frac{\left( 1-\alpha +\gamma +2\mu \right) \pi }{4}},
\label{franjo-over-over-TD-2} \\
a_{3}b_{1}\cos \frac{\left( \gamma +\mu \right) \pi }{2}-a_{2}b_{2}\sin 
\frac{\left( \alpha -\mu \right) \pi }{2}\leqslant a_{1}b_{3}\cos \frac{%
\left( 2\alpha +\gamma -\mu \right) \pi }{2},  \label{franjo-over-over-TD-3}
\\
a_{1}b_{3}\sin \frac{\left( 2\alpha +\gamma -\mu \right) \pi }{2}\leqslant
a_{2}b_{2}\cos \frac{\left( \alpha -\mu \right) \pi }{2}-a_{3}b_{1}\sin 
\frac{\left( \gamma +\mu \right) \pi }{2}.  \label{franjo-over-over-TD-4}
\end{gather}

\noindent \textbf{Model IDD}$^{{}^{+}}$\textbf{.IDD}$^{{}^{+}}$ takes the
form%
\begin{equation*}
\left( a_{1}\,_{0}\mathrm{I}_{t}^{\alpha }+a_{2}\,_{0}\mathrm{D}_{t}^{\frac{%
1+\gamma -\alpha }{2}}+a_{3}\,_{0}\mathrm{D}_{t}^{1+\gamma }\right) \sigma
\left( t\right) =\left( b_{1}\,_{0}\mathrm{I}_{t}^{\alpha +\gamma -\eta
}+b_{2}\,_{0}\mathrm{D}_{t}^{\frac{1+\eta -\left( \alpha +\gamma -\eta
\right) }{2}}+b_{3}\,_{0}\mathrm{D}_{t}^{1+\eta }\right) \varepsilon \left(
t\right) ,  \label{IDD-IDD-over-over-konacno}
\end{equation*}%
with thermodynamical restrictions%
\begin{gather}
0\leqslant \alpha +\gamma -\eta \leqslant 1,\quad \gamma \leqslant \eta
,\quad \alpha +\eta +\left( \eta -\gamma \right) =\alpha -\gamma +2\eta
\leqslant 1,  \label{TD-IDD-IDD-over-over-1} \\
\frac{a_{2}}{a_{1}}\leqslant \frac{b_{2}}{b_{1}}\frac{\sin \frac{\left(
1+\alpha -\gamma +2\eta \right) \pi }{4}}{\sin \frac{\left( 1+\alpha
+3\gamma -2\eta \right) \pi }{4}},\quad \frac{a_{3}}{a_{2}}\leqslant \frac{%
b_{3}}{b_{2}}\frac{\cos \frac{\left( 1-\alpha +\gamma -2\eta \right) \pi }{4}%
}{\cos \frac{\left( 1-\alpha -3\gamma +2\eta \right) \pi }{4}},
\label{TD-IDD-IDD-over-over-3} \\
a_{3}b_{1}\cos \frac{\left( \alpha +2\gamma -\eta \right) \pi }{2}%
-a_{2}b_{2}\sin \frac{\left( \eta -\gamma \right) \pi }{2}\leqslant
a_{1}b_{3}\cos \frac{\left( \alpha +\eta \right) \pi }{2},
\label{TD-IDD-IDD-over-over-4} \\
a_{1}b_{3}\sin \frac{\left( \alpha +\eta \right) \pi }{2}\leqslant
a_{2}b_{2}\cos \frac{\left( \eta -\gamma \right) \pi }{2}-a_{3}b_{1}\sin 
\frac{\left( \alpha +2\gamma -\eta \right) \pi }{2}.
\label{TD-IDD-IDD-over-over-5}
\end{gather}

\noindent \textbf{Model I}$^{{}^{+}}$\textbf{ID.IDD}$^{{}^{+}}$ takes the
form%
\begin{equation*}
\left( a_{1}\,_{0}\mathrm{I}_{t}^{1+\alpha }+a_{2}\,_{0}\mathrm{I}_{t}^{%
\frac{1+\alpha -\gamma }{2}}+a_{3}\,_{0}\mathrm{D}_{t}^{\gamma }\right)
\sigma \left( t\right) =\left( b_{1}\,_{0}\mathrm{I}_{t}^{\alpha +\gamma
-\eta }+b_{2}\,_{0}\mathrm{D}_{t}^{\frac{1+\eta -\left( \alpha +\gamma -\eta
\right) }{2}}+b_{3}\,_{0}\mathrm{D}_{t}^{1+\eta }\right) \varepsilon \left(
t\right) ,  \label{frale-over-less-less-over-konacno}
\end{equation*}%
with thermodynamical restrictions%
\begin{gather}
\eta \leqslant \gamma ,\quad \alpha +\gamma +2\left( \gamma -\eta \right)
=\alpha +3\gamma -2\eta \leqslant 1,  \label{TD-frale-over-less-less-over-1}
\\
\frac{a_{1}}{b_{1}}\frac{\sin \frac{\left( 1+\alpha -\gamma +2\eta \right)
\pi }{4}}{\cos \frac{\left( 1-\alpha -3\gamma +2\eta \right) \pi }{4}}%
\leqslant \frac{a_{2}}{b_{2}}\leqslant \frac{a_{3}}{b_{3}}\frac{\cos \frac{%
\left( 1-\alpha -3\gamma +2\eta \right) \pi }{4}}{\sin \frac{\left( 1+\alpha
-\gamma +2\eta \right) \pi }{4}},  \label{TD-frale-over-less-less-over-2} \\
a_{1}b_{3}\cos \frac{\left( \alpha +\eta \right) \pi }{2}-a_{2}b_{2}\sin 
\frac{\left( \gamma -\eta \right) \pi }{2}\leqslant a_{3}b_{1}\cos \frac{%
\left( \alpha +2\gamma -\eta \right) \pi }{2},
\label{TD-frale-over-less-less-over-3} \\
a_{3}b_{1}\sin \frac{\left( \alpha +2\gamma -\eta \right) \pi }{2}\leqslant
a_{2}b_{2}\cos \frac{\left( \gamma -\eta \right) \pi }{2}-a_{1}b_{3}\sin 
\frac{\left( \alpha +\eta \right) \pi }{2}.
\label{TD-frale-over-less-less-over-4}
\end{gather}

\subsection{Asymmetric models\label{asimetricni modeli}}

\noindent \textbf{Model IID.ID} takes the form%
\begin{equation*}
\left( a_{1}\,_{0}\mathrm{I}_{t}^{\alpha +\beta -\gamma }+a_{2}\,_{0}\mathrm{%
I}_{t}^{\nu }+a_{3}\,_{0}\mathrm{D}_{t}^{\gamma }\right) \sigma \left(
t\right) =\left( b_{1}\,_{0}\mathrm{I}_{t}^{\alpha }+b_{2}\,_{0}\mathrm{D}%
_{t}^{\beta }\right) \varepsilon \left( t\right) ,
\label{(A)Z-IID-ID-less-konacno}
\end{equation*}%
with thermodynamical restrictions%
\begin{gather}
0\leqslant \alpha \leqslant \nu <\alpha +\beta -\gamma \leqslant 1,\quad
\beta +\nu \leqslant 1,  \label{TD-(A)Z-IID-ID-less-1} \\
-\frac{a_{1}}{a_{3}}\frac{\cos \frac{\left( \alpha +2\beta -\gamma \right)
\pi }{2}}{\cos \frac{\left( \alpha +\gamma \right) \pi }{2}}\leqslant \frac{%
b_{1}}{b_{2}}\leqslant \frac{a_{1}}{a_{3}}\frac{\sin \frac{\left( \alpha
+2\beta -\gamma \right) \pi }{2}}{\sin \frac{\left( \alpha +\gamma \right)
\pi }{2}}.  \label{TD-(A)Z-IID-ID-less-2}
\end{gather}

\noindent \textbf{Model IDD.DD}$^{{}^{+}}$ takes the form%
\begin{equation*}
\left( a_{1}\,_{0}\mathrm{I}_{t}^{\alpha }+a_{2}\,_{0}\mathrm{D}_{t}^{\beta
}+a_{3}\,_{0}\mathrm{D}_{t}^{\gamma }\right) \sigma \left( t\right) =\left(
b_{1}\,_{0}\mathrm{D}_{t}^{\mu }+b_{2}\,_{0}\mathrm{D}_{t}^{\alpha +\beta
+\mu }\right) \varepsilon \left( t\right) ,
\end{equation*}%
with thermodynamical restrictions%
\begin{gather}
1\leqslant \alpha +\beta +\mu \leqslant 2,\quad \beta <\gamma \leqslant \mu
\leqslant 1-\alpha ,  \label{TD-AZ-IDD-DD-less-over-1} \\
\frac{a_{1}}{a_{2}}\frac{\left\vert \cos \frac{\left( 2\alpha +\beta +\mu
\right) \pi }{2}\right\vert }{\cos \frac{\left( \mu -\beta \right) \pi }{2}}%
\leqslant \frac{b_{1}}{b_{2}}.  \label{TD-AZ-IDD-DD-less-over-2}
\end{gather}

\noindent \textbf{Model I}$^{{}^{+}}$\textbf{ID.ID} takes the form%
\begin{equation*}
\left( a_{1}\,_{0}\mathrm{I}_{t}^{\alpha +\beta +\nu }+a_{2}\,_{0}\mathrm{I}%
_{t}^{\nu }+a_{3}\,_{0}\mathrm{D}_{t}^{\alpha +\beta -\nu }\right) \sigma
\left( t\right) =\left( b_{1}\,_{0}\mathrm{I}_{t}^{\alpha }+b_{2}\,_{0}%
\mathrm{D}_{t}^{\beta }\right) \varepsilon \left( t\right) ,
\label{(A)Z-IID-ID-over-konacno}
\end{equation*}%
with thermodynamical restrictions%
\begin{gather}
0\leqslant \alpha +\beta -\nu \leqslant 1,\quad 1\leqslant \alpha +\beta
+\nu \leqslant 2,\quad \alpha \leqslant \nu \leqslant 1-\beta ,
\label{TD-(A)Z-IID-ID-over-1} \\
\frac{a_{1}}{a_{2}}\frac{\left\vert \cos \frac{\left( \alpha +2\beta +\nu
\right) \pi }{2}\right\vert }{\cos \frac{\left( \nu -\alpha \right) \pi }{2}}%
\leqslant \frac{b_{1}}{b_{2}}\leqslant \frac{a_{2}}{a_{3}}\frac{\sin \frac{%
\left( \beta +\nu \right) \pi }{2}}{\sin \frac{\left( 2\alpha +\beta -\nu
\right) \pi }{2}}.  \label{TD-(A)Z-IID-ID-over-2}
\end{gather}

\noindent \textbf{Model IDD}$^{{}^{+}}$\textbf{.DD}$^{{}^{+}}$ takes the form%
\begin{equation*}
\left( a_{1}\,_{0}\mathrm{I}_{t}^{\alpha }+a_{2}\,_{0}\mathrm{D}_{t}^{\beta
}+a_{3}\,_{0}\mathrm{D}_{t}^{\alpha +2\beta }\right) \sigma \left( t\right)
=\left( b_{1}\,_{0}\mathrm{D}_{t}^{\mu }+b_{2}\,_{0}\mathrm{D}_{t}^{\alpha
+\beta +\mu }\right) \varepsilon \left( t\right) ,
\end{equation*}%
with thermodynamical restrictions%
\begin{gather}
1\leqslant \alpha +2\beta \leqslant 2,\quad 1\leqslant \alpha +\beta +\mu
\leqslant 2,\quad \beta \leqslant \mu \leqslant 1-\alpha ,
\label{TD-AZ-IDD-DD-over-over-1} \\
\frac{a_{1}}{a_{2}}\frac{\left\vert \cos \frac{\left( 2\alpha +\beta +\mu
\right) \pi }{2}\right\vert }{\cos \frac{\left( \mu -\beta \right) \pi }{2}}%
\leqslant \frac{b_{1}}{b_{2}}\leqslant \frac{a_{2}}{a_{3}}\frac{\sin \frac{%
\left( \alpha +\mu \right) \pi }{2}}{\sin \frac{\left( \alpha +2\beta -\mu
\right) \pi }{2}}.  \label{TD-AZ-IDD-DD-over-over-2}
\end{gather}

\noindent \textbf{Model ID.IDD} takes the form%
\begin{equation*}
\left( a_{1}\,_{0}\mathrm{I}_{t}^{\mu }+a_{2}\,_{0}\mathrm{D}_{t}^{\nu
}\right) \sigma \left( t\right) =\left( b_{1}\,_{0}\mathrm{I}_{t}^{\alpha
}+b_{2}\,_{0}\mathrm{D}_{t}^{\beta }+b_{3}\,_{0}\mathrm{D}_{t}^{\mu +\nu
-\alpha }\right) \varepsilon \left( t\right) ,
\label{(A)Z-ID-IDD-less-konacno}
\end{equation*}%
with thermodynamical restrictions%
\begin{gather}
0\leqslant \nu \leqslant \beta <\mu +\nu -\alpha \leqslant 1,\quad \alpha
\leqslant \mu \leqslant 1-\beta ,  \label{TD-(A)Z-ID-IDD-less-1} \\
-\frac{a_{1}}{a_{2}}\frac{\cos \frac{\left( 2\mu +\nu -\alpha \right) \pi }{2%
}}{\cos \frac{\left( \alpha +\nu \right) \pi }{2}}\leqslant \frac{b_{1}}{%
b_{3}}\leqslant \frac{a_{1}}{a_{2}}\frac{\sin \frac{\left( 2\mu +\nu -\alpha
\right) \pi }{2}}{\sin \frac{\left( \alpha +\nu \right) \pi }{2}}.
\label{TD-(A)Z-ID-IDD-less-2}
\end{gather}

\noindent \textbf{Model ID.DDD}$^{{}^{+}}$ takes the form%
\begin{equation*}
\left( a_{1}\,_{0}\mathrm{I}_{t}^{\alpha }+a_{2}\,_{0}\mathrm{D}_{t}^{\beta
}\right) \sigma \left( t\right) =\left( b_{1}\,_{0}\mathrm{D}_{t}^{\mu
}+b_{2}\,_{0}\mathrm{D}_{t}^{\nu }+b_{3}\,_{0}\mathrm{D}_{t}^{\alpha +\beta
+\nu }\right) \varepsilon \left( t\right) ,
\end{equation*}%
with thermodynamical restrictions%
\begin{gather}
1\leqslant \alpha +\beta +\nu \leqslant 2,\quad \beta \leqslant \mu <\nu
\leqslant 1-\alpha ,  \label{TD-(A)Z-ID-DDD-over-1} \\
\frac{a_{1}}{a_{2}}\frac{\left\vert \cos \frac{\left( 2\alpha +\beta +\nu
\right) \pi }{2}\right\vert }{\cos \frac{\left( \nu -\beta \right) \pi }{2}}%
\leqslant \frac{b_{2}}{b_{3}}.  \label{TD-(A)Z-ID-DDD-over-2}
\end{gather}

\noindent \textbf{Model ID.IDD}$^{{}^{+}}$ takes the form%
\begin{equation*}
\left( a_{1}\,_{0}\mathrm{I}_{t}^{\alpha }+a_{2}\,_{0}\mathrm{D}_{t}^{\beta
}\right) \sigma \left( t\right) =\left( b_{1}\,_{0}\mathrm{I}_{t}^{\alpha
+\beta -\nu }+b_{2}\,_{0}\mathrm{D}_{t}^{\nu }+b_{3}\,_{0}\mathrm{D}%
_{t}^{\alpha +\beta +\nu }\right) \varepsilon \left( t\right) ,
\label{(A)Z-ID-IDD-over-konacno}
\end{equation*}%
with thermodynamical restrictions%
\begin{gather}
0\leqslant \alpha +\beta -\nu \leqslant 1,\quad 1\leqslant \alpha +\beta
+\nu \leqslant 2,\quad \beta \leqslant \nu \leqslant 1-\alpha ,
\label{TD-(A)Z-ID-IDD-over-1} \\
\frac{b_{1}}{b_{2}}\frac{\sin \frac{\left( \alpha +2\beta -\nu \right) \pi }{%
2}}{\sin \frac{\left( \alpha +\nu \right) \pi }{2}}\leqslant \frac{a_{1}}{%
a_{2}}\leqslant \frac{b_{2}}{b_{3}}\frac{\cos \frac{\left( \nu -\beta
\right) \pi }{2}}{\left\vert \cos \frac{\left( 2\alpha +\beta +\nu \right)
\pi }{2}\right\vert }.  \label{TD-(A)Z-ID-IDD-over-2}
\end{gather}

%\bibliographystyle{plain}
%\bibliography{dz}

\end{document}

%% file: tabela-svih-modela.tex
 \begin{table}[h]
 \begin{center}
 \begin{tabular}{|c|c|c|c|c|c|c|c|}

 \hline \xrowht{14pt}

 & Function  $\phi _{\sigma }$ & Function $\phi _{\varepsilon }$ & Model &  & Order $\xi $ & Order $\lambda $ & Order $\kappa $ \\ 

 \hhline{|=|=|=|=|=|=|=|=|} \xrowht{14pt}

\multirow{8}{*}{\multirowcell{8}{\rotatebox{90}{Symmetric models}}} & \multirow{2}{*}{\multirowcell{2}{$ a_{1}+a_{2}s^{\alpha +\beta}$ }}  & \multirow{2}{*}{\multirowcell{2}{$ b_{1}+b_{2}s^{\alpha +\beta}$} } 
& ID.ID & $*$ & $\alpha -\mu $ & $-$ & $-$ \\

\cline{4-8} \xrowht{14pt}

 & & & ID.DD$^{+}$ &  & $\alpha +\mu <1$ & $-$ & $-$ \\

\Xcline{2-8}{2\arrayrulewidth} \xrowht{14pt}

& \multirow{3}{*}{\multirowcell{3}{$a_{1}+a_{2}s^{\lambda }+a_{3}s^{\alpha +\gamma }$}} &  \multirow{2}{*}{\multirowcell{2}{$b_{1}+b_{2}s^{\lambda }+b_{3}s^{\alpha +\gamma }$}} 
 & IID.IID & $*$ & $\eta-\gamma $ & $\alpha -\beta $ & $-$ \\

\cline{4-8} \xrowht{14pt}

 & & & IDD.IDD & $*$ & $\alpha -\mu $ & $\alpha +\beta <1$ & $-$ \\

\cline{3-8} \xrowht{14pt}

 & & $b_{1}+b_{2}s^{\kappa }+b_{3}s^{\alpha +\gamma }$ & IID.IDD &  & $\alpha -\mu $ & $\mu +\nu <1$ & $\alpha -\beta $ \\ 

\cline{2-8} \xrowht{14pt}

 & \multirow{3}{*}{\multirowcell{3}{$a_{1}+a_{2}s^{\frac{1+\alpha +\gamma }{2}}+a_{3}s^{1+\alpha +\gamma }$}} 
 &\multirow{3}{*}{\multirowcell{3}{$b_{1}+b_{2}s^{\frac{1+\alpha +\gamma }{2}}+b_{3}s^{1+\alpha +\gamma }$}} & I$^{+}$ID.I$^{+}$ID &  & $\alpha -\mu $ & $-$ & $-$ \\

\cline{4-8} \xrowht{14pt} 

 &  &  & IDD$^{+}$.IDD$^{+}$ &  & $\eta -\gamma $ & $-$ & $-$ \\ 

\cline{4-8} \xrowht{14pt}

& & & I$^{+}$ID.IDD$^{+}$ &  & $1-\left( \gamma -\eta \right) $ & $-$ & $-$ \\ 

\Xhline{4\arrayrulewidth} \xrowht{14pt}

 \multirow{7}{*}{\multirowcell{7}{\rotatebox{90}{Asymmetric models}}}  
 & \multirow{2}{*}{\multirowcell{2}{$a_{1}+a_{2}s^{\lambda }+a_{3}s^{\kappa }$}} 
 & \multirow{4}{*}{\multirowcell{4}{$b_{1}+b_{2}s^{\alpha +\beta }$ }}
 & IID.ID &  & $\beta -\gamma $ & $\left( \alpha +\beta \right)-\left( \nu +\gamma \right) $ & $\alpha +\beta <1$ \\ 

\cline{4-8} \xrowht{14pt}

 & & & IDD.DD$^{+}$ &  & $\alpha +\mu <1$ & $\alpha +\beta <1$ & $\alpha +\mu <1$ \\ 

\cline{2-2} \cline{4-8} \xrowht{14pt}

 & \multirow{2}{*}{\multirowcell{2}{$a_{1}+a_{2}s^{\alpha +\beta }+a_{3}s^{2\left( \alpha +\beta \right) }$}}
 & & I$^{+}$ID.ID &  & $\beta +\nu <1$ & $-$ & $-$ \\

\cline{4-8} \xrowht{14pt}

 & & & IDD$^{+}$.DD$^{+}$ &  & $\alpha +\mu <1$ & $-$ & $-$ \\

\Xcline{2-8}{2\arrayrulewidth} \xrowht{14pt}

 & \multirow{3}{*}{\multirowcell{3}{$a_{1}+a_{2}s^{\alpha +\beta }$}} & \multirow{2}{*}{\multirowcell{2}{$b_{1}+b_{2}s^{\lambda }+b_{3}s^{\kappa }$}} & ID.IDD &  & $\alpha-\mu $ & $\mu +\nu <1$ & $\alpha +\beta <1$ \\ 

\cline{4-8} \xrowht{14pt}

 & & & ID.DDD$^{+}$ &  & $\alpha +\mu <1$ & $\nu -\mu $ & $\alpha +\nu-\left( \mu -\beta \right) $ \\

\cline{3-8} \xrowht{14pt} 

  & & $b_{1}+b_{2}s^{\alpha +\beta}+b_{3}s^{2\left( \alpha +\beta \right) }$ & ID.IDD$^{+}$ &  & $\nu -\beta $ & $-$ & $-$ \\

\hline

 \end{tabular}
 \end{center}
 \caption{Summary of constitutive functions $\phi _{\sigma }$ and $\phi _{\varepsilon }$, along with the order $\xi$, corresponding to the thermodynamically consistent fractional anti-Zener and Zener models. The notation $*$ means that the orders $\alpha +\beta$ and $\alpha +\gamma$ belong either to the interval $(0,1)$ or interval $(1,2)$.}
 \label{skupina}
 \end{table}

%% file: parametri.tex
 \begin{table}[h]
 \begin{center}

\begin{tabular}{|c|c|c|c|c|c|c|c|c|}

 \hline \xrowht{14pt}

Case when $\tilde{\sigma}_{sr}$ has& $\alpha $ & $\beta $ & $\nu $ & $a_{1}$ & $a_{2}$ & $a_{3}$ & $b_{1}$ & $b_{2}$ \\

 \hhline{|=|=|=|=|=|=|=|=|=|} \xrowht{14pt}

no poles & \multirow{3}{*}{\multirowcell{3}{$0.35$}} & \multirow{3}{*}{\multirowcell{3}{$0.55$}} & \multirow{3}{*}{\multirowcell{3}{$0.4$}} & $0.05$ & $1.5$ & $0.45$ & $0.7$ & $0.95$ \\  \xrowht{14pt}
a negative real pole &  &  &  & $11$ & $28.4029\dots $ & $20.27$ & $7$ & $9.5$ \\  \xrowht{14pt}
\makecell{a pair of complex\\ conjugated poles} & &  &  & $11$ & $15$ & $20.27$ & $7$ & $9.5$\\ 

\hline

\end{tabular}

 \end{center}
 \caption{Model parameters used for numerical examples.}
 \label{parametri}
 \end{table}